\documentclass[11pt,reqno]{amsart}
\usepackage[french,english]{babel}
\usepackage[utf8]{inputenc}
\usepackage[T1]{fontenc}            
\usepackage{amsmath}
\usepackage{graphicx}
\usepackage{amssymb}
\usepackage{amsthm}
\usepackage{comment}
\usepackage{pgf, tikz}
\usepackage{amsfonts,mathrsfs}
\usepackage{thmtools}
\usepackage{geometry}
\usepackage{stmaryrd}
\usepackage{color}
\usepackage{hyperref}
\usepackage{bm}
\hypersetup{colorlinks = true , linkcolor=blue, citecolor=blue}
\numberwithin{equation}{section}
\usepackage{comment}
\theoremstyle{definition}

\declaretheorem[numberwithin=,title=Definition]{Def}
\declaretheorem[numberwithin=,title=Proposition]{Prop}

\declaretheorem[numberwithin=,title=Theorem]{Thm}
\declaretheorem[numberwithin=,title=Lemma]{Lemme}

\declaretheorem[numberwithin=,title=Remark]{Rq}

\newtheorem{Assumption}{Assumption}

\renewcommand{\P}{\mathbb{P}}
\renewcommand{\H}{\ensuremath{\mathcal{H}}\xspace}
\newcommand{\p}{\ensuremath{\widehat{p}}\xspace}
\newcommand{\pa}{\ensuremath{\partial}\xspace}
\newcommand{\PP}{\mathcal{P}_{\beta}(\R^d)}
\newcommand{\PPP}{\mathcal{P}(\R^d)}
\newcommand{\FF}{\mathcal{F}}
\newcommand{\E}{\mathbb{E}}

\newcommand{\LL}{\mathcal{L}}

\newcommand{\KK}{\mathcal{K}}
\newcommand{\R}{\ensuremath{\mathbb{R}}\xspace}

\newcommand{\N}{\ensuremath{\mathbb{N}}\xspace}
\newcommand{\C}{\ensuremath{\mathbb{C}}\xspace}

\newcommand{\eps}{\varepsilon}
\newcommand{\CC}{\mathcal{C}}
\newcommand{\BB}{\mathcal{B}}

\newcommand{\1}{\mathbf{1}}

\newcommand{\NN}{\ensuremath{\mathcal{N}}\xspace}
\newcommand{\NNN}{\ensuremath{\widetilde{\mathcal{N}}}\xspace}

\newcommand{\delb}{\frac{\delta }{\delta  m}b}
\newcommand{\delu}{\frac{\delta }{\delta  m}u}
\newcommand{\del}{\frac{\delta }{\delta  m}}
\newcommand{\dell}{\frac{\delta^2 }{\delta  m^2}}
\newcommand{\muu}{\overline{\mu}}
\newcommand{\x}{\bm{x}}

\begin{document}

	\title[Propagation of chaos for stable-driven McKean-Vlasov SDEs]{Quantitative weak propagation of chaos for McKean-Vlasov SDEs driven by stable processes}

	\author{Thomas Cavallazzi}
	\address{Université Paris-Saclay, CentraleSupélec, MICS and CNRS FR-3487, France}
	\email{thomas.cavallazzi@centralesupelec.fr}
	\date{January 24, 2024}
	%\urladdr{dfdf} 
	
	\begin{abstract}
		We consider a general McKean-Vlasov stochastic differential equation driven by a rotationally invariant $\alpha$-stable process on $\R^d$, with $\alpha \in (1,2)$. We assume that the diffusion coefficient is the identity matrix and that the drift is bounded and Hölder continuous with respect to both space and measure variables, in some precise sense. The main goal of this work is to prove new weak propagation of chaos estimates for the associated mean-field interacting particle system. We also establish a pointwise control on the difference between the density of one particle and the density of the limiting McKean-Vlasov SDE. Our study relies on the regularizing properties and the dynamics of the semigroup associated with the McKean-Vlasov stochastic differential equation, which acts on functions defined on $\PP$, the space of probability measures on $\R^d$ having a finite moment of order $\beta \in (1,\alpha)$. More precisely, the dynamics of the semigroup is described by a backward Kolmogorov partial differential equation defined on the strip $[0,T]\times \PP$.
	\end{abstract}

 \maketitle
 
 \tableofcontents
 
 \section{Introduction}
 
 In this work, we are interested in the following McKean-Vlasov Stochastic Differential Equation (SDE) driven by $Z= (Z_t)_t$ a rotationally invariant $\alpha$-stable process on $\R^d$, with $\alpha \in (1,2)$,
 
 \begin{equation}\label{McKVSDE_intro}
 	\begin{cases}
 		dX_t= b(t,X_t,[X_t])\, dt +dZ_t, \quad t \in [0,T],\\ \mu_t := [X_t], \\ X_0=\xi, \quad [\xi] = \mu \in \PPP,
 	\end{cases} 
 \end{equation}
 where  $T>0$ is fixed, $[\xi]$ denotes the distribution of the random variable $\xi$ which is independent of $Z$, $\PPP$ is the space of probability measures on $\R^d$ and where $b : [0,T]\times \R^d \times \PPP \to \R^d$ is the drift coefficient. This kind of SDEs are also called nonlinear SDEs since the associated Fokker-Planck equation solved by the flow of marginal distributions $(\mu_t)_t$ is nonlinear. We also focus on the associated mean-field interacting particle system defined, for all $N \geq 1$, by \begin{equation}\label{edsparticles_intro}
 	\left\{  \begin{array}{lll}
 		&dX^{i,N}_t = b(t,X^{i,N}_t,\muu^N_t) \,dt +  dZ^i_t, \quad t \in [0,T],\quad i \in \{1,\dots,N\}, \\ &\muu^N_t := \frac{1}{N} \displaystyle\sum_{j=1}^N \delta_{X^{j,N}_t},\\ &X^{i,N}_0 = \xi^i,
 	\end{array}\right.\end{equation}
 where $(Z^i,\xi^i)_{ i \geq 1}$ are i.i.d.\ copies of $(Z,\xi)$. The connection between \eqref{McKVSDE_intro} and \eqref{edsparticles_intro} is that the McKean-Vlasov SDE \eqref{McKVSDE_intro} describes the dynamics of one particle of the system \eqref{edsparticles_intro}, when the number of particles $N$ tends to infinity. A stronger property that is also expected to hold true in general is the so-called propagation of chaos. It states that for all $k \geq 1$, the dynamics of $k$ particles is described, when $N$ tends to infinity, by $k$ independent copies of \eqref{McKVSDE_intro}. This was originally studied by McKean \cite{McKeanPropchaos67} and then investigated by Sznitman \cite{SznitmanTopicspropagationchaos1991}, when $Z$ is a Brownian motion. These mean-field
 systems have many applications for example in physics (kinetic theory), in biology to study the motion of a
 cell population for example, in neuroscience to model the interactions between neurons, in social sciences to
 describe self-organization behaviors and also in the Mean-Field Games theory.\\
 
  Propagation of chaos can be considered in the weak sense, i.e.\ in distribution, in particular through the convergence of the empirical measure $\muu^N$ (see \cite{SznitmanTopicspropagationchaos1991}), or in the strong sense, i.e.\ at the level of paths by coupling. The terminology of weak and strong propagation of chaos used here is based of the corresponding properties for numerical schemes for SDEs. It can be qualitative or quantitative when a rate of convergence is shown. By quantitative weak propagation of chaos, we precisely mean to find a rate of convergence of $\E |\phi(\muu^N_t) - \phi(\mu_t)|$ and $|\E(\phi(\muu^N_t) - \phi(\mu_t))|$, for $\phi$ in a large class of functions defined on the space of probability measures.\\

    When $Z$ is a Brownian motion, McKean-Vlasov SDEs and mean-field particle systems have been widely studied. In addition to  McKean \cite{McKeanPropchaos67} and Sznitman \cite{SznitmanTopicspropagationchaos1991}, one can mention for example Gärtner \cite{Gartner}, Méléard \cite{Meleard1996POC}, Malrieu \cite{malrieu:hal-01282602}, Mischler et al.\ \cite{Mischlernewapproachquantitative2015}, Jabin and Wang \cite{Jabin_2018}, Lacker \cite{lacker2018strong, Lacker_entropy}, Toma\v{s}evi\'c \cite{tomasevic:hal-03086253} and Jabir \cite{Jabir_girsanov}.  For a detailed review on the topic of propagation of chaos, we refer to \cite{chaintron:hal-03585065,chaintron:hal-03585067}. The recent development of Mean-Field Games has given a new impulsion to study these systems and the associated propagation of chaos. In particular, it provides a new formalism and new tools, such as the notion of master equations, which are Partial Differential Equations (PDEs) on the space of measures, that allow to revisit or generalize some of the preceding works. In this direction, one can mention for example the book of Carmona and Delarue \cite{CarmonaProbabilisticTheoryMean2018}, Chaudru de Raynal and Frikha \cite{deraynal2021wellposedness, frikha2021backward}, Chassagneux et al.\ \cite{chassagneux2019weak}, Delarue and Tse \cite{delarue2021uniform} and Jourdain and Tsé \cite{JourdainCentrallimittheorem2021}. Such new tools are at the core of this work.\\
 
  It is natural to consider other types of noise such as Lévy processes, which are also largely used to model physical systems (Lévy flights and anomalous diffusion), see e.g. \cite{Mann_fractal_fractional_diffusions} for the physical point of view and \cite{Jourdain_fractional_diffusion} for
  the mathematical point of view. The propagation of chaos phenomenon has also been studied for McKean-Vlasov SDEs driven by Lévy processes. In \cite{GRAHAM199269},
 following the approach of Sznitman \cite{SznitmanTopicspropagationchaos1991}, Graham proves qualitative weak propagation of chaos
 under Lipschitz assumptions (with respect to the $1$- Wasserstein metric $W_1$ for the measure variable) for a mean-field system driven by a Poisson random measure and its compensated
 measure. He works in the $L^1$ framework, i.e.\ the Poisson random measure is associated with a Poisson process
 having a finite moment of order $1$. When the driving noise is a Lévy process having a finite moment of order $2$, we refer to Jourdain et al. \cite{JourdainNonlinearSDEsdriven2007}. The authors prove quantitative strong propagation of chaos in $L^2$ under standard Lipschitz assumptions on the drift and diffusion coefficients (with respect to the $2$-Wasserstein metric $W_2$). We also refer to Neelima et al. \cite{NeelimaWellposednesstamedEuler2020} where quantitative strong propagation of chaos is proved in $L^2$, relaxing the assumptions of \cite{JourdainNonlinearSDEsdriven2007}. We
 also mention \cite{Mischlernewapproachquantitative2015}, where Mischler, Mouhot and Wennberg exhibit conditions leading to weak propagation of chaos estimates. As an application, they study an inelastic Boltzmann collision jump process. In the one-dimensional case, Frikha and Li \cite{FrikhaWellposednessapproximationonedimensional2020} study a McKean-Vlasov
 SDE driven by a compensated Poisson random measure with positive jumps. They prove quantitative strong propagation of chaos in $L^1$ under one-sided Lipschitz assumptions on the coefficients (with respect to $W_1$ for the measure variable).\\
 
The aim of this work is to prove quantitative weak propagation of chaos for \eqref{McKVSDE_intro}, assuming that the drift $b$ is bounded and Hölder continuous in some precise sense with respect
 to both space and measure variables (see Assumption \ref{Assumption2} for the precise assumptions). The method that we use relies on the semigroup acting on functions defined on $\PPP$, associated with the McKean-Vlasov SDE. More details are given below. The study in the Brownian case has been made by Chaudru de Raynal and Frikha \cite{deraynal2021wellposedness,frikha2021backward}, under similar assumptions. The main difficulties here are the following. Firstly, the dependence of the drift $b$ with respect to the measure variable is general, so that differential calculus for functions defined on $\PP$ is needed. We use the notion of linear (functional) derivative (see Definition \ref{deflinearderivative}). Secondly, since the coefficients are not Lipschitz continuous, one needs to benefit from a regularization by noise phenomenon, in particular to prove the well-posedness of \eqref{McKVSDE_intro}. Thirdly, since $Z$ is a Lévy process having an infinite moment of order $2$, another difficulty here is the impossibility of working in $L^2.$ It is thus impossible to rely on the tools, mentioned above, developed for Mean-Field Games in the $L^2$ framework. Moreover, the presence of jumps induces supplementary technical difficulties to develop such tools.\\

 In the first part of this work, the weak well-posedness of \eqref{McKVSDE_intro} is established through the related nonlinear martingale problem by using the Banach fixed point theorem applied on a suitable complete metric space (see Assumption \ref{Assumption1} and Theorem \ref{Thm_wellposedness}). Notice that the well-posedness of \eqref{McKVSDE_intro} was already shown in \cite{frikha_menozzi_konakov} by using the same fixed point argument. The fixed point theorem is only proved to hold in small time in \cite{frikha_menozzi_konakov}, which is enough for the well-posedness. However, in the following, we need to apply the Banach fixed point theorem on the whole time interval $[0,T]$ in order to use Picard iterations to approximate the solution to \eqref{McKVSDE_intro} on $[0,T]$. That is why we do the proof a bit differently to obtain this result.\\
 
We then study the regularity of the transition density associated with \eqref{McKVSDE_intro}, in particular the differentiability with respect to the initial distribution $\mu \in \PP$, with $\beta \in (1,\alpha)$, where $\PP$ is the space of probability measures $\mu$ on $\R^d$ having a finite moment of order $\beta$. Let us introduce, for $x \in \R^d$ and $s \in [0,T)$, the following decoupled
 stochastic flow associated to SDE \eqref{McKVSDE_intro}

 \begin{equation}\label{decoupled_SDE_intro}
 	\begin{cases}
 		dX^{s,x,\mu}_t= b(t,X_t^{s,x,\mu},[X^{s,\mu}_t])\, dt +dZ_t,\quad t \in [s,T], \\ X_s^{s,x,\mu}=x,
 	\end{cases}
 \end{equation}
where $[X^{s,\mu}_t]$ denotes the distribution of the solution to \eqref{McKVSDE_intro} at time $t$, starting at time $s$ by any random variable $\xi$ with distribution $\mu \in \PP$. The notation $[X^{s,\mu}_t]$ makes sense since by weak well-posedness, the distribution of the solution to \eqref{McKVSDE_intro} depends on the initial condition only through its distribution $\mu$. The distributions of $X^{s,x,\mu}_t$ and $X^{s,\mu}_t$ admit densities respectively denoted by $p(\mu,s,t,x,\cdot)$ and $p(\mu,s,t,\cdot)$. Moreover, they satisfy the following relation,  stemming from the well-posedness of the martingale problem associated with \eqref{McKVSDE_intro},
 \begin{equation}\label{decoupleddensities_link_intro}
	p(\mu,s,t,y) = \int_{\R^d} p(\mu,s,t,x,y) \, d\mu(x).
\end{equation}
 Let us fix $t \in (0,T]$ and $y \in \R^d$. We study the regularity of the map $(s,x,\mu) \in [0,t) \times \R^d \times \PP \mapsto p(\mu,s,t,x,y)$ and we prove in Theorem \ref{Thm_density_McKean_light} that $p$ is solution to the following backward Kolmogorov PDE \begin{equation}\label{EDPKolmogorovfundamental_intro}\begin{cases}\pa_s p(\mu,s,t,x,y) + \LL_sp(\cdot,s,t,\cdot,y)(\mu,x) =0, \quad \forall(\mu,s,x) \in \PP \times [0,t) \times \R^d, \\ p(\mu,s,t,x,\cdot) \underset{s \rightarrow t-}{\longrightarrow} \delta_x, \quad \text{in the weak sense,}\end{cases}\end{equation}
 where $\LL_s$ is defined, for smooth enough functions $h$ on $ \PP \times \R^d$, by \begin{align}\label{def_operator_McKV_intro}
 	\LL_sh(\mu,x)\notag&:= b(s,x,\mu)\cdot \pa_x h(\mu,x) +  \int_{\R^d} \left[ h(\mu,x+z) - h(\mu,x) - z\cdot \pa_x h(\mu,x) \right] \, \frac{dz}{|z|^{d+\alpha}} \\ &\quad + \int_{\R^d}  b(s,v,\mu) \cdot \pa_v\del h(\mu,x)(v) \, d\mu(v) \\ \notag&\quad + \int_{\R^d} \int_{\R^d} \left[\del h(\mu,x)(v+z) - \del h(\mu,x)(v) - z\cdot \pa_v \del h(\mu,x)(v) \right] \, \frac{dz}{|z|^{d+\alpha}} \, d\mu(v),
 \end{align}
where $\del$ denotes the linear derivative (see Definition \ref{deflinearderivative}). Estimates and Hölder controls on the derivatives of $p$ with respect to $s$, $x$ and $\mu$ are proved in Theorem \ref{Thm_density_McKean} under Assumption \ref{Assumption2} by using Picard iterations and the parametrix method.\\
 
 Once that the regularity of the transition density has been studied, we focus on the regularizing properties of the semigroup, acting on functions defined on $\PP$, associated with \eqref{McKVSDE_intro}. It is at the core of the method to prove the propagation of chaos. For a function $\phi : \PP \to \R$, the action of the McKean-Vlasov semigroup on $\phi$ is given by the function $U$ defined by \begin{equation}\label{defsg_intro}U:(t,\mu) \in [0,T]\times \PP \mapsto \phi([X^{t,\mu}_T]),\end{equation}
 where $[X^{t,\mu}_T]$ denotes the distribution of the solution to \eqref{McKVSDE_intro} at time $T$ and starting at time $t$ with a random variable $\xi$ with distribution $\mu \in \PP$. More precisely, we prove that if $\phi : \PP \to \R$ has a linear derivative that is uniformly Hölder continuous, then the map  $U$ defined by \eqref{defsg_intro} is more regular. Moreover, estimates and Hölder controls on the derivative of $U$ with respect to the measure variable are proved. Inspired by \cite{CarmonaProbabilisticTheoryMean2018} (see $(5.128)$), \cite{Cardaliaguetmasterequationconvergence2015} and \cite{deraynal2021wellposedness}, we describe the dynamics of the semigroup in Theorem \ref{Thm_EDP} by showing that $U$ is the unique classical solution to the backward Kolmogorov PDE \begin{equation}\label{EDP_KOLMO_intro}
 	\begin{cases}
 		\pa_t U(t,\mu) + \mathscr{L}_tU(t,\cdot)(\mu) =0, \quad \forall (t,\mu) \in [0,T) \times \PP, \\ 
 		U(T,\mu) = \phi(\mu), \quad \forall \mu \in \PP,
 	\end{cases}
 \end{equation} 
where $\mathscr{L}_t$ is defined, for smooth enough functions $h$ on $\PP$, by \begin{align}\label{def_operator_McKV_measure_intro}
	\mathscr{L}_th(\mu)&:= \int_{\R^d}  b(t,v,\mu) \cdot \pa_v\del h(\mu)(v) \, d\mu(v) \\ \notag&\quad + \int_{\R^d} \int_{\R^d} \left[\del h(\mu)(v+z) - \del h(\mu)(v) - z\cdot \pa_v \del h(\mu)(v) \right] \, \frac{dz}{|z|^{d+\alpha}} \, d\mu(v).
\end{align}\\

Then, we use the preceding results to prove quantitative propagation of chaos for the particle system \eqref{edsparticles_intro} towards the limiting McKean-Vlasov SDE \eqref{McKVSDE_intro}.\
 
Firstly, we prove in Theorem \ref{Thm_POC} weak propagation of chaos estimates for the particle system \eqref{edsparticles_intro}, as defined above. We refer the reader to Remark \ref{RQ_comparison_POC} for a comparison with the existing literature. The method that we use relies on the solution to the backward Kolmogorov PDE \eqref{EDP_KOLMO_intro} $U$ defined in \eqref{defsg_intro}, and thus on the McKean-Vlasov semigroup. This strategy was originally described in \cite[pages $506-508$]{CarmonaProbabilisticTheoryMean2018}, inspired by \cite{Cardaliaguetmasterequationconvergence2015} and \cite{Mischlernewapproachquantitative2015}, and was employed for example in \cite{chassagneux2019weak,delarue2021uniform,frikha2021backward}. Let us describe the main ideas. We begin by computing the time derivative of the map $ t \in [0,T) \mapsto U(t,\muu^N_t)$ by applying the standard Itô formula for the empirical projection $(t,x_1,\dots,x_N)~\in ~[0,T] \times (\R^d)^N \mapsto U\left(t, \frac1N \sum_{k=1}^N \delta_{x_k}\right)$ and for the particle system. Noting that $t \in[0,T] \mapsto U(t,\mu_t)$ is constant, we naturally expect that the time derivative previously computed tends to $0$ as $N$ converges to infinity. This convergence has to be shown with an explicit rate of convergence using the PDE satisfied by $U$ and some estimates on $U$. Finally, we express $\phi(\muu^N_T) - \phi(\mu_T) = U(T,\muu^N_T) - U(T,\mu_T)$ as the sum of $ U(0,\muu^N_0) - U(0,\mu_0)$ plus a remainder term related to the time derivative which was previously estimated. Since the initial data are i.i.d., the first term is controlled by standard estimates, for example those in \cite{fournier:hal-00915365}.\

Secondly, we focus on the approximation of the distribution of one particle of the system \eqref{edsparticles_intro} by the limiting McKean-Vlasov process at the level of densities in Theorem \ref{Thm_poc_density}. We prove a pointwise error bound between the density of one particle and the density of the limiting McKean-Vlasov SDE. This allows to quantify the decrease with respect to $N$ of the total variation distance between the distributions of one particle and the solution to the McKean-Vlasov SDE. The method used to prove this result relies on similar ideas than presented in the preceding paragraph. It was used in the Brownian case in \cite{frikha2021backward} (see also references therein). Instead of considering the semigroup associated with the McKean-Vlasov SDE, we work with the not decoupled density of the solution to the McKean-Vlasov SDE \eqref{McKVSDE_intro}, which is somehow a fundamental solution to the Backward Kolmogorov PDE \eqref{EDP_KOLMO_intro} (see Theorem \ref{Thm_EDP_densité_McKV}). Denoting by $p(\mu,s,t,\cdot)$ the density of the McKean-Vlasov SDE \eqref{McKVSDE_intro} at time $t$ and starting at time $s$ from a random variable with distribution $\mu \in \PP$, the idea is to study the dynamics of the map $s \in [0,t) \mapsto p(\muu^N_s,s,t,y)$, for any $y \in \R^d$. The two ingredients are the following. On the one hand, we quantify the error with respect to $N$ between $\E p(\muu^N_s,s,t,y)$ and $p(\mu_s,s,t,y)=p(\mu_0,0,t,y)$. On the other hand, we prove that when $s$ tends to $t$, then $ \E p(\muu^N_s,s,t,\cdot)$ converges pointwise to the density of one particle at time $t$ and starting at time $0$ from $\mu_0$.\\ 

The paper is organized as follows. In Section \ref{section_overview}, we present our results and we comment them. The proofs are given in the next sections. Section \ref{section_well_posedness} is dedicated to prove the weak well-posedness of the McKean-Vlasov SDE \eqref{McKVSDE_intro} (Theorem \ref{Thm_wellposedness}). Then, we study the regularity of its transition density in Section \ref{section_properties_transition_densit} (Theorem \ref{Thm_density_McKean}). Section \ref{section_PDE} is devoted to the proof of the regularizing properties of the semigroup and the backward Kolmogorov PDE that describes its dynamics (Theorem \ref{Thm_EDP}). In Section \ref{section_prop_chaos}, we prove the quantitative weak propagation of chaos result (Theorem \ref{Thm_POC}) and in Section \ref{section_poc_density} the error bound for the approximation of the distribution of one particle by the limiting McKean-Vlasov SDE at the level of densities (Theorem \ref{Thm_poc_density}). We prove in Section \ref{section_proof_prop} the technical proposition leading to the regularity of the transition density and the related estimates stated in Theorem \ref{Thm_density_McKean}. In Appendix \ref{section_appendix_diff_calculus}, we gather some definitions and propositions related to differential calculus for functions defined on $\PPP$ or $\PP$. Finally, Appendix \ref{section_appendix_parametrix} aims at presenting the parametrix method and estimates on the density of the solution to a linear stable-driven SDE, which are the core of the proof of Theorem \ref{Thm_density_McKean}.\\

Let us finally introduce some notations used several times in the article.\\

\textbf{Notations}

\begin{enumerate}
	
	\item[-] $\PPP$ denotes the set of probability measures on $\R^d$.
	\item[-] $d_{TV}$ is the total variation metric on $\PPP$.
	
	\item[-] $\PP$ denotes the set of probability measures $\mu$ on $\R^d$ such that $\int_{\R^d} |x|^\beta \, d\mu(x) < + \infty,$ for $\beta \geq 1.$ It is equipped with the Wasserstein metric of order $\beta$ denoted by $W_\beta,$ which makes it complete. Denoting by $\Pi(\mu,\nu)$ the set of couplings between two probability measures $\mu, \nu \in \PP$, the metric $W_\beta$ is defined by $$W_{\beta}(\mu,\nu) = \inf_{\pi \in\Pi(\mu,\nu) } \left(\int_{\R^d\times \R^d} |x-y|^\beta \, d\pi(x,y)\right)^{\frac{1}{\beta}}.$$ 
	\item[-] $M_\beta(\mu)$ denotes the moment of order $\beta$ of $\mu \in \PP$ defined by $$M_\beta(\mu) := \left( \int_{\R^d} |x|^\beta \, d\mu(x) \right)^{\frac1\beta }.$$
	\item[-] $[\xi]$ denotes the distribution of the random variable $\xi$.
	\item[-]  $\muu^N_{\bm{x}}:= \frac1N \sum_{k=1}^N \delta_{x_k}$ denotes the empirical measure, for $\bm{x} = (x_1,\dots,x_N) \in (\R^d)^N$.
	\item[-]  $\bm{\tilde{z}_k} := (0, \dots,z,\dots,0) \in (\R^d)^N$ for $z \in \R^d,$ where $z$ appears in the $k$-th position.

	\item[-] $B_r$ is the open ball in $\R^d$ centered at $0$ and of radius $r$ for the euclidean norm.
	
	\item[-] $B_r^c$ denotes the complementary of $B_r$.

	\item[-] $a\wedge b$ denotes the minimum between $a$ and $b$.
	\item[-] $a \vee b $ denotes the maximum between $a$ and $b$.
	\item[-] $\BB$ denotes the Beta function defined, for all $x,y >0$, by \[ \BB(x,y) := \int_0^1 (1-t)^{-1+x} t^{-1+y}  \, dt. \]
	
	\item[-] $C$ is a generic constant that may depend only on the fixed parameters of the problem and which may change from line to line.
	
\end{enumerate}

 \section{Overview on the main results and comments}\label{section_overview}
 
 Let us fix $Z=(Z_t)_t$ a rotationally invariant $\alpha$-stable process on $\R^d$ with $\alpha \in (1,2).$ Its associated Poisson random measure is denoted by $\NN$ and the compensated Poisson random measure by $\NNN$. Since $\alpha \in (1,2)$, by \cite[Remark $14.6$ and Theorem $14.7$]{SatoLevyprocessesinfinitely1999}, we can write for all $t \geq 0$ \[ Z_t = \int_0^t\int_{\R^d} z \, \NNN(ds,dz). \] The density of $Z_t$ is denoted by $q(t,\cdot)$ and the Lévy measure $\nu$ of $Z$ is given by \[ d\nu(z) := \frac{d z}{|z|^{d+\alpha}}.\]
The generator of $Z$ is the fractional Laplacian $\Delta^{\frac{\alpha}{2}}$ which is defined for all $f \in \CC^{1+\gamma}_b(\R^d;\R)$, with $\gamma > \alpha -1$ (i.e.\ $f$ belongs to $\CC^1_b(\R^d;\R)$ and $\nabla f$ is $\gamma$-Hölder continuous) and for all $x \in \R^d$, by \begin{equation}\label{stableoperatorMK}
 	\Delta^{\frac{\alpha}{2}} f(x):= \int_{\R^d} (f(x+z) - f(x) - \nabla f (x) \cdot z )\, d\nu(z).
 \end{equation}
 
 We are interested in the following stable-driven McKean-Vlasov SDE \begin{equation}\label{McKVSDE_existence}
 	\begin{cases}
 		dX_t^{s,\xi}= b(t,X_t^{s,\xi},[X_t^{s,\xi}])\, dt +dZ_t, \quad  t \in [s,T], \\ X_s^{s,\xi}=\xi, \quad [\xi] = \mu \in \PPP,
 	\end{cases} 
 \end{equation}
 where $[\xi]$ denotes the distribution of the random variable $\xi$ and $s \in [0,T)$.\\
 
  Let us define, for $k > - \alpha $ the function $\rho^k$ on $(0,+\infty)\times \R^d$ by \begin{equation}\label{defdensityrefMK}
 	\forall t>0, \, x \in \R^d,\, \rho^k(t,x):= t^{-\frac{d}{\alpha}}(1+t^{-\frac{1}{\alpha}}|x|)^{-d-\alpha -k}.
 \end{equation}
 These functions are related to gradient estimates on the transition density $q$ of the stable process $Z$ (see Lemma \ref{Lemmegradientestimatestable} in the Appendix). Some useful properties of these functions are given in Lemma \ref{Lemmedensityreferencecontrol}.\\
 
 \subsection{Well-posedness of the nonlinear martingale problem and Picard iterations}
 
The first point is to prove the existence and uniqueness, in the weak sense, of the solution to \eqref{McKVSDE}. Let us first recall the definition of the nonlinear martingale problem associated with \eqref{McKVSDE}.
 
 \begin{Def}\label{def_pb_ptgle} Let us fix $(s,\mu) \in [0,T) \times \PPP$. We say that a probability measure $\P$ on the Skorokhod space
 	$\mathcal{D}([s,T];\R^d)$, endowed with its canonical filtration $(\FF_t)_t$, with time marginal distributions $(\P_t)_{t\in[s,T]} \in \CC^0([s,T];\PPP)$ solves the
 	nonlinear martingale problem associated to the SDE \eqref{McKVSDE_existence} with initial distribution $\mu$ at time $s$ if the
 	canonical process $(y_t)_{t\in[s,T]}$ satisfies the two following conditions.
 	
 	\begin{enumerate}
 		\item We have $\P_s = \mu$.
 		\item For any $\phi \in \CC^{1,2}_b([s,T]\times \R^d)$, the process defined for $t \in [s,T]$ by 
 		\begin{equation*}
 			\phi(t,y_t) - \phi(s,y_s) - \int_s^t \left(\pa_r + L^{\P}_r\right) \phi(r,y_r) \, dr
 		\end{equation*}
 		is a $(\mathcal{D}([s,T];\R^d),(\FF_t)_t,\P)$-martingale starting from $0$ at time $t=s$ and where 
 		
 		\begin{equation}\label{def_operator_mtgle_pb}
 			L^{\P}_rf(t,x) := b(r,x,\P_r)\cdot \pa_x f(t,x) + \Delta^{\frac{\alpha}{2}} f(t,\cdot)(x).
 		\end{equation}
 	\end{enumerate}
 \end{Def}
 
 We assume that the drift $b:[0,T]\times \R^d \times \PPP  \to \R^d$ satisfies the following properties.
 
 \begin{Assumption}\label{Assumption1}
 	\
 	
 	\begin{enumerate}
 		\item The drift $b$ is measurable and globally bounded on $[0,T]\times \R^d \times \PPP$.
 		\item For any $(t,\mu)\in [0,T]\times \PPP$, the map $b(t,\cdot,\mu)$ is $\eta$-Hölder continuous on $\R^d$ uniformly with respect to $(t,\mu) \in [0,T]\times \PPP$, with $\eta \in ( 0, 1]$, i.e.\ there exists $C>0$ such that for all $t\in [0,T]$, $\mu \in\PPP$ and $x_1,x_2 \in \R^d$ \[|b(t,x_1,\mu)-b(t,x_2,\mu)| \leq C |x_1 - x_2|^\eta.\]
 		\item For any $(t,x) \in [0,T] \times \R^d$, the map $b(t,x,\cdot)$ is Lipschitz continuous with respect to the total variation metric $d_{TV}$ uniformly with respect to $(t,x)\in[0,T]\times \R^d$, i.e.\ there exists $C>0$ such that for all $t\in [0,T]$, $x \in\R^d$ and $\mu_1,\mu_2\in \PPP$ \[|b(t,x,\mu_1)-b(t,x,\mu_2)| \leq C d_{TV}(\mu_1,\mu_2).\]
 	\end{enumerate}
 	
 \end{Assumption}

We can now state the weak well-posedness result.

 \begin{Thm}[Weak well-posedness]\label{Thm_wellposedness}
 	Under Assumption \ref{Assumption1}, the martingale problem associated with the McKean-Vlasov SDE \eqref{McKVSDE_existence} is well-posed for all initial distribution $\mu \in \PPP$. In particular, the SDE \eqref{McKVSDE_existence} is well-posed in the weak sense. Moreover, for any $P \in \CC^0([s,T];\PPP)$ with $P_s = \mu \in \PPP$, we can define recursively $\overline{X}^{(m)}$ as the unique weak solution to \begin{equation}\label{SDE_iteree_thm_existence}
 		\begin{cases}
 			d\overline{X}^{(m)}_t = b(t,\overline{X}^{m}_t,[\overline{X}^{(m-1)}_t]))\,dt + dZ_t, \quad t \in [s,T], \\ \overline{X}^{(m)}_s = \xi,
 		\end{cases}
 	\end{equation}
 	with $[\xi] =\mu$ and $([\overline{X}^{(0)}_t])_{t \in [s,T]} = P$. Then, denoting by $P^*$ the unique solution to the martingale problem associated with \eqref{McKVSDE_existence}, with initial condition $\mu\in \PPP$ at time $s$, we have \begin{equation}\label{convergence_Thm_existence}
 		\sup_{t \in [s,T]} d_{TV} (P^*_t,[\overline{X}^{(m)}_t]) \underset{m \to + \infty}{\longrightarrow}0,
 	\end{equation}
 	where $d_{TV}$ is the total variation metric.
 \end{Thm}

This result is proved in Section \ref{section_well_posedness}.

\begin{Rq}The well-posedness of the nonlinear martingale problem was already proved in \cite{frikha_menozzi_konakov}. Therein, the fixed point theorem is only shown to hold in small time since it is enough to prove the well-posedness. However, we will need in the following the convergence of the sequence of Picard iterations on the whole interval $[s,T]$. That is why we do the proof a bit differently. Instead of using the implicit parametrix expansion, we write the complete parametrix expansion, which yields this global result.
\end{Rq}

\subsection{Cauchy problem for the transition density associated with the McKean-Vlasov SDE}

 Now that the well-posedness has been established, we focus on the regularity of the transition density associated with \eqref{McKVSDE}. In particular, our aim is to study its regularity with respect to the initial distribution $\mu$. To do this, we fix $s\in [0,T)$, $\beta \in (1,\alpha)$ and we consider the following stable-driven McKean-Vlasov SDE \begin{equation}\label{McKVSDE}
 	\begin{cases}
 		dX^{s,\xi}_t= b(t,X_t^{s,\xi},[X_t^{s,\xi}])\, dt +dZ_t, \quad t \in [s,T], \\ X_s^{s,\xi}=\xi, \quad [\xi] = \mu \in \PP.
 	\end{cases} 
 \end{equation}
 The associated martingale problem is well-posed using Theorem \ref{Thm_wellposedness} and there is weak existence and uniqueness for the SDE \eqref{McKVSDE}. Moreover the distribution of $X^{s,\xi}_t$ depends only on $\mu$ and not on the choice of the random variable $\xi$ such that $[\xi]=\mu$ and we denote it by \begin{equation}\label{notationX^mu}[X^{s,\mu}_t] := [X^{s,\xi}_t].  \end{equation} 
 We introduce, for $x \in \R^d$, the following decoupled
 stochastic flow associated to SDE \eqref{McKVSDE}

 \begin{equation}\label{decoupled_SDE}
 	\begin{cases}
 		dX^{s,x,\mu}_t= b(t,X_t^{s,x,\mu},[X^{s,\mu}_t])\, dt +dZ_t,\quad  t \in [s,T], \\ X_s^{s,x,\mu}=x \in \R^d.
 	\end{cases}
 \end{equation}
 
 The distribution of $X^{s,x,\mu}_t$ admits a density denoted by $p(\mu,s,t,x,\cdot)$ by using the parametrix expansion given in Theorem \ref{Thmdensityparametrix}. The distribution of $X^{s,\mu}_t$ has also a density denoted by $p(\mu,s,t,\cdot)$. Moreover, it satisfies the following relation stemming from the well-posedness of the martingale problem \begin{equation}\label{decoupleddensities_link}
 	p(\mu,s,t,y) = \int_{\R^d} p(\mu,s,t,x,y) \, d\mu(x).
 \end{equation}

 We introduce the following assumption, which is stronger than Assumption \ref{Assumption1}. 

\begin{Assumption}\label{Assumption2}
	\
	
	\begin{enumerate} 
		\item The drift $b$ is jointly continuous and globally bounded on $[0,T]\times \R^d \times \PPP$.
		\item For any $(t,\mu)\in [0,T]\times \PPP$, the map $b(t,\cdot,\mu)$ is $\eta$-Hölder continuous on $\R^d$, for some fixed $\eta \in (0,1]$, uniformly with respect to $(t,\mu) \in [0,T]\times \PPP$, i.e.\ there exists $C>0$ such that for all $t\in [0,T]$, $\mu \in\PPP$ and $x_1,x_2 \in \R^d$ \[|b(t,x_1,\mu)-b(t,x_2,\mu)| \leq C |x_1 - x_2|^\eta.\]
		
		\item For any $(t,x) \in [0,T] \times \R^d$, the map $\mu \in \PPP \mapsto b(t,x,\mu)$ has a linear derivative such that $\del b(t,x,\mu)(\cdot)$ is $\eta$-Hölder continuous on $\R^d$ uniformly with respect to $(t,x,\mu) \in [0,T] \times \R^d\times \PPP$ and $\del b$ is bounded on $[0,T]\times \R^d \times \PPP \times \R^d$.
		
		\item For any $(t,x,v) \in [0,T] \times (\R^d)^2$, the map $\mu \in \PPP \mapsto \del b(t,x,\mu)(v)$ has a linear derivative such that $\dell b(t,x,\mu)(v,\cdot)$ is $\eta$-Hölder continuous uniformly with respect to $(t,x,\mu,v) \in [0,T]\times \R^d \times \PPP\times \R^d$ and $\dell b$ is bounded on $[0,T] \times \R^d \times \PPP \times (\R^d)^2$.
	\end{enumerate}
	
\end{Assumption}

We need this stronger assumptions in order to study the regularity with respect to the initial distribution of the transition density $p$.\\

Let us now state in the next theorem the backward Kolmogorov PDE satisfied by the decoupled transition density $p$ as well as some important gradient estimates on $p$.

\begin{Thm}\label{Thm_density_McKean_light}
	Let us fix $(t,y) \in (0,T]\times \R^d$. Under Assumption \ref{Assumption2}, the mapping $(\mu,s,x) \in \PP \times [0,t) \times \R^d \mapsto p(\mu,s,t,x,y)$ belongs to $\CC^{1}(\PP \times [0,t)\times \R^d )$ (see Definition \ref{def_C1_space}), and is solution to the following backward Kolmogorov PDE \begin{equation}\label{EDPKolmogorovfundamental}\begin{cases}\pa_s p(\mu,s,t,x,y) + \LL_sp(\cdot,s,t,\cdot,y)(\mu,x) =0, \quad \forall(\mu,s,x) \in \PP \times [0,t) \times \R^d, \\ p(\mu,s,t,x,\cdot) \underset{s \rightarrow t-}{\longrightarrow} \delta_x, \quad \text{in the weak sense,}\end{cases}\end{equation}
	where $\LL_s$ is defined, for smooth enough function $h$ on $\PP \times \R^d$, by \begin{align}\label{def_operator_McKV}
		\LL_sh(\mu,x)\notag&:= b(s,x,\mu)\cdot \pa_x h(\mu,x) +  \int_{\R^d} \left[ h(\mu,x+z) - h(\mu,x) - z\cdot \pa_x h(\mu,x) \right] \, \frac{dz}{|z|^{d+\alpha}} \\ &\quad + \int_{\R^d}  b(s,v,\mu) \cdot \pa_v\del h(\mu,x)(v) \, d\mu(v) \\ \notag&\quad + \int_{\R^d} \int_{\R^d} \left[\del h(\mu,x)(v+z) - \del h(\mu,x)(v) - z\cdot \pa_v \del h(\mu,x)(v) \right] \, \frac{dz}{|z|^{d+\alpha}} \, d\mu(v).
	\end{align}
	
	Moreover, there exists a positive constant $C>0$ such that for all $j \in \{0,1\}$, $\mu \in \PP$, $0 \leq s < t \leq T$ and $x,y,v \in \R^d$
			
			\begin{equation}\label{density_bound_thm_light}
				|\pa_x^j p(\mu,s,t,x,y) |
				\leq C (t-s)^{-\frac{j}{\alpha}} \rho^{j}(t-s,y-x),
			\end{equation}
		
			\begin{equation}\label{density_bound_laplacian_thm_light}
			|\Delta^{\frac{\alpha}{2}} p(\mu,s,t,\cdot,y)(x) |
			\leq C (t-s)^{-1} \rho^{0}(t-s,y-x),
		\end{equation}

			 \begin{equation}\label{gradient_linear_der_bound_thm_light}
				\left\vert \pa_v^j\del p(\mu,s,t,x,y)(v)\right\vert \leq C (t-s)^{j\frac{\eta -1}{\alpha} + 1 - \frac{1}{\alpha}}\rho^0(t-s,y-x),
			\end{equation}
		and

		\begin{equation}\label{density_op_bound_thm}
			\left\vert\Delta^{\frac{\alpha}{2}} \left[\del p(\mu,s,t,x,y) \right](v)\right\vert \leq C(t-s)^{- \frac{1}{\alpha}} \rho^{0}(t-s,y-x),
		\end{equation}
	where $\rho^j$ was defined by \eqref{defdensityrefMK}.

\end{Thm}

The proof of this theorem is postponed to Section \ref{section_properties_transition_densit}. More precisely, certain other estimates (see Theorem \ref{Thm_density_McKean}), in particular Hölder controls, are needed not only to prove the preceding theorem, but also to exhibit the regularizing properties of the semigroup associated with \eqref{McKVSDE} which are presented in the next subsection.

\begin{Rq}Contrary to the Brownian case ($\alpha =2$) studied in \cite{deraynal2021wellposedness,frikha2021backward}, we do not need to prove $\CC^2$ regularity with respect to $x \in \R^d$ and $\mu \in \PP$. This would impose stronger assumptions.
\end{Rq}

\begin{Rq}
	Let us introduce, for $(s,\mu) \in[0,T]\times\PP$, the operator $L_s^{\mu}$ associated with the linear SDE related to the McKean-Vlasov SDE \eqref{McKVSDE}, i.e.\ where the measure argument is frozen equal to $\mu$. It is defined, for smooth enough function $f$ on $\R^d$ by $$ \forall x \in \R^d, \quad L_s^{\mu} f(x) := b(s,x,\mu)\cdot \pa_x f(x) +  \Delta^{\frac{\alpha}{2}} f(x).$$
	Then, we have for all smooth enough function $h$ on $\R^d \times \PP$ $$ \LL_s h(\mu,x) = L_s^\mu h(\mu,\cdot)(x) + \int_{\R^d} L^{\mu}_s \left[\del h(\mu,x)\right](v) \, d\mu(v).$$
\end{Rq}

Let us now deal with the not decoupled density of the McKean-Vlasov SDE \eqref{McKVSDE}.

\begin{Thm}\label{Thm_EDP_densité_McKV}
	Let us fix $(t,y) \in (0,T]\times \R^d$. Under Assumption \ref{Assumption2}, the mapping $(\mu,s) \in \PP \times [0,t)  \mapsto p(\mu,s,t,y)$ belongs to $\CC^{1}(\PP \times [0,t))$ (see Definition \ref{def_C1_space}), and is solution to the following backward Kolmogorov PDE \begin{equation}\label{EDPKolmogorovfundamental_not_dec}\begin{cases}\pa_s p(\mu,s,t,y) + \mathscr{L}_sp(\cdot,s,t,y)(\mu) =0, \quad \forall(\mu,s) \in \PP \times [0,t), \\ p(\mu,s,t,\cdot) \underset{s \rightarrow t-}{\longrightarrow} \mu, \quad \text{with respect to $W_1$,}\end{cases}\end{equation}
	where $\mathscr{L}_s$ is defined, for smooth enough functions $h$ on $\PP$, by \begin{align}\label{def_operator_McKV_measure}
		\mathscr{L}_sh(\mu)&:= \int_{\R^d}  b(s,v,\mu) \cdot \pa_v\del h(\mu)(v) \, d\mu(v) \\ \notag&\quad + \int_{\R^d} \int_{\R^d} \left[\del h(\mu)(v+z) - \del h(\mu)(v) - z\cdot \pa_v \del h(\mu)(v) \right] \, \frac{dz}{|z|^{d+\alpha}} \, d\mu(v).
	\end{align}
\end{Thm}

We do not give the proof of this result since it follows from the same reasoning as in the proof of $\eqref{EDPKolmogorovfundamental}$ in Theorem \ref{Thm_density_McKean_light} presented in Section \ref{section_properties_transition_densit} (see also \cite[Theorem $3.3$]{frikha2021backward} for the proof in the Brownian case).

\subsection{Backward Kolmogorov PDE on the space of measures}

We can now focus on the study of the semigroup associated with \eqref{McKVSDE}, acting on functions defined on $\PP$. Let us recall that $\beta \in (1,\alpha)$ is fixed. For a fixed function $\phi : \PP \to \R$, the action of the semigroup on $\phi$ is given by the map $U$ defined by \begin{equation}\label{def_sol_PDE}
	U(t,\mu) := \phi([X^{t,\mu}_T]),\quad \forall (t,\mu) \in [0,T]\times \PP,
\end{equation}
where $[X^{t,\mu}_T]$ is the flow of marginal distributions of \eqref{McKVSDE}, where the initial distribution is equal to $\mu$ at time $t$ and was defined in \eqref{notationX^mu}. We aim at studying the regularizing properties of the semigroup, i.e.\ the gain of regularity between $\phi$ and $U$ with respect to the measure variable. This will be crucial to prove the propagation of chaos. The regularization of $\phi$ by a smooth flow of probability measures is presented in Proposition \ref{Prop_lin_der_along_flow_density}.\\

 We define the space of functions on which the semigroup acts.

\begin{Def}\label{def_holder_test_functions}
	Let us fix $\delta \in (0,1]$. The space $\CC^{(1,\delta)}(\PP)$ is defined as the set of continuous functions $\phi : \PP \rightarrow \R$ admitting a linear derivative such that there exists a positive constant $C$ such that for all $\mu \in \PP$, $v_1,v_2 \in \R^d$ \[ \left\vert \del \phi (\mu)(v_1) -  \del \phi (\mu)(v_2) \right\vert \leq C |v_1 - v_2|^\delta.\]

\end{Def}

We state in the next theorem the regularizing properties of the semigroup acting on $\CC^{(1,\delta)}(\PP)$ and we describe its dynamics through the backward Kolmogorov PDE that it satisfies.  

\begin{Thm}[Backward Kolmogorov PDE]\label{Thm_EDP}
	
	Let us fix $\phi\in \CC^{(1,\delta)}(\PP)$. Then, under Assumption \ref{Assumption2}, the function $U$ defined in \eqref{def_sol_PDE} belongs to $\CC^0([0,T] \times \PP) \cap \CC^1([0,T)\times \PP)$ (see Definition \ref{def_C1_space}) and satisfies the following properties.
	
	\begin{itemize}
%		\item There exists a positive constant $C$ such that for all $t \in [0,T)$, $\mu \in \PP$, $v \in \R^d$ 
%		
%		\begin{equation}\label{lin_der_bound_sol_PDE}
%			\left\vert \del U (t,\mu)(v)\right\vert \leq C(1+|v|^\delta).
%		\end{equation} 
		
		\item There exists a positive constant $C$ such that for all $t \in [0,T)$, $\mu \in \PP$, $v \in \R^d$ 
		
		\begin{equation}\label{gradient_lin_der_bound_sol_PDE}
			\left\vert \pa_v\del U (t,\mu)(v)\right\vert \leq C(T-t)^{\frac{\delta -1}{\alpha}}.
		\end{equation}
		
		\item For all $\gamma \in (0,1]\cap(0,(2\alpha -2) \wedge (\eta + \alpha -1))$, there exists a positive constant $C$ such that for all $t \in [0,T)$, $\mu \in \PP$, $v_1,v_2 \in \R^d$ 
		
		\begin{equation}\label{gradient_lin_der_holder_sol_PDE}
			\left\vert \pa_v\del U (t,\mu)(v_1) - \pa_v\del U (t,\mu)(v_2)\right\vert \leq C(T-t)^{\frac{\delta -1 - \gamma}{\alpha}}|v_1-v_2|^\gamma.
		\end{equation}
	\end{itemize}
	
	Moreover, $U$ is solution to the following backward Kolmogorov PDE
	
	\begin{equation}\label{EDP_Kolmo_backward}
		\begin{cases}
			\pa_t U(t,\mu) + \mathscr{L}_tU(t,\cdot)(\mu) =0, \quad \forall (t,\mu) \in [0,T) \times \PP, \\ 
			U(T,\mu) = \phi(\mu), \quad \forall \mu \in \PP,
		\end{cases}
	\end{equation}
where $\mathscr{L}_t$ was defined in \eqref{def_operator_McKV_measure}. It is the unique solution to \eqref{EDP_Kolmo_backward} among all functions in $\CC^0([0,T] \times \PP) \cap \CC^1([0,T)\times \PP)$ satisfying \eqref{gradient_lin_der_bound_sol_PDE} and \eqref{gradient_lin_der_holder_sol_PDE}.
\end{Thm}

We prove this result in Section \ref{section_PDE}.

\subsection{Quantitative weak propagation of chaos}
We are now going to use the regularizing properties and the dynamics of the semigroup given in Theorem \ref{Thm_EDP} to prove quantitative weak propagation of chaos for the mean-field interacting particle system associated with \eqref{McKVSDE}.
 Let us introduce $(Z^n)_n$ an i.i.d.\ sequence of $\alpha$-stable processes having the same distribution as $Z$ and $(X^n_0)_n$ an i.i.d.\ sequence of random variables with common distribution $\mu_0 \in \PP$, where $\beta \in (1,\alpha)$ is still fixed. For any $N \geq 1,$ the system of $N$ particles associated with \eqref{McKVSDE} is defined as the unique solution to the following classical SDE \begin{equation}\label{edsparticles}
	\left\{  \begin{array}{lll}
		&dX^{i,N}_t = b(t,X^{i,N}_t,\muu^N_t) \,dt +  dZ^i_t, \quad t \in [0,T],\quad i \in \{1,\dots,N\}, \\ &\muu^N_t := \frac{1}{N} \displaystyle\sum_{j=1}^N \delta_{X^{j,N}_t},\\ &X^{i,N}_0 = X^i_0.
	\end{array}\right.\end{equation}

This linear SDE on $(\R^d)^N$ is well-posed in the weak sense using \cite[Corollary $1.4$ $(iii)$]{Cheng_stable_homeo} since its drift coefficient is Hölder continuous in space uniformly in time. The limiting McKean-Vlasov SDE is \eqref{McKVSDE} starting at time $s=0$ from any random variable $\xi$ with distribution $\mu_0 \in \PP$. We denote by $(\mu_t)_{t \in [0,T]}$ the flow of marginal distributions of its solution. Let us define the space of test functions that we use to quantify the weak propagation of chaos of the particle system \eqref{edsparticles} towards the McKean-Vlasov SDE \eqref{McKVSDE}.
 
 \begin{Def}\label{def_space_prop_chaos}
 	For $\delta \in (0,1]$ and $L>0$, we define the space $\CC^{(2,\delta)}_L(\PP)$ as the set of continuous functions  $\phi : \PP \rightarrow \R$ admitting two linear derivatives $\del \phi$ and $\dell \phi$ such that for all $\mu \in \PP$, $v_1,v_2,v_1',v_2' \in \R^d$ \[ \left\vert \del \phi (\mu)(v_1) -  \del \phi (\mu)(v_2) \right\vert \leq L |v_1 - v_2|^\delta,\] and   \[ \left\vert \dell \phi (\mu)(v_1,v_1') -  \dell \phi (\mu)(v_2,v_2') \right\vert \leq L (|v_1 - v_2|^\delta + |v_1' - v_2'|^\delta).\]
 \end{Def}

\begin{Rq}
	Note that $\CC^{(2,\delta)}_L(\PP)$ is a subspace of $\CC^{(1,\delta)}(\PP)$ defined in Definition \ref{def_holder_test_functions}.
\end{Rq}

We now state our quantitative propagation of chaos result. 

\begin{Thm}[Quantitative weak propagation of chaos]\label{Thm_POC} \
	
	\noindent Let us fix $\delta \in (0,1]$, $L>0$ and $\gamma \in (0,1] \cap (0,(\delta + \alpha -1)\wedge(2\alpha -2)\wedge (\eta +\alpha -1))$. Then, under Assumption \ref{Assumption2}, there exists a positive constant $C = C(d,T,\alpha,\beta,\ref{Assumption2},\gamma,\delta,L)$, non-decreasing with respect to $T$, such that for all $\phi \in \CC^{(2,\delta)}_L(\PP)$ and $N \geq 1$, it holds 
	
	\begin{equation}\label{eq1_thm_poc}
		\E|\phi(\muu^N_T) - \phi(\mu_T)| \leq CT^{\frac{\delta-1}{\alpha}}\E W_1(\muu^N_0,\mu_0) + \frac{C}{N^{1-\frac{1}{\beta}}}.
	\end{equation}
Moreover, there exists a positive constant $C = C(d,T,\alpha,\beta,\ref{Assumption2},\gamma,\delta,L,M_1(\mu_0))$, non-decreasing with respect to $T$, such that for all $\phi \in \CC^{(2,\delta)}_L(\PP)$ and $N \geq 1$, it holds 
	
	\begin{equation}\label{eq3_thm_poc}
		|\E(\phi(\muu^N_T) - \phi(\mu_T))| \leq  \frac{C}{N^\gamma}T^{\frac{\delta-1}{\alpha}}.
	\end{equation}
\end{Thm}

The proof is given in Section \ref{section_prop_chaos}.

\begin{Rq} \begin{enumerate}
		
		\item The initial data term in \eqref{eq1_thm_poc} can be handled using Fournier and Guillin \cite{fournier:hal-00915365}, in particular in the case where $\mu_0$ has more moments than $\beta$. Indeed, one has if $\mu_0 \in \mathcal{P}_q(\R^d)$ with $q\geq 1$ \begin{equation}\label{FG}
			\E W_1(\muu^N_0,\mu_0) \leq C \left\{\begin{array}{lll}
				&N^{-\frac12} + N^{-\left(1-\frac{1}{q}\right)}, \quad &\text{if}\quad d=1\quad \text{and}\quad q \neq 2,\\ 
				&N^{-\frac12}\ln(1+N) +N^{-\left(1-\frac{1}{q}\right)}, \quad &\text{if}\quad d=2 \quad \text{and}\quad q \neq 2, \\  &N^{-\frac1d} + N^{-\left(1-\frac{1}{q}\right)}, \quad &\text{if}\quad d\geq 3\quad \text{and} \quad q \neq \frac{d}{d -1}.
			\end{array}\right.
		\end{equation}
		\item Let us justify why the estimate \eqref{eq3_thm_poc} is interesting. Firstly, this result quantifies the approximation of semigroup associated with the McKean-Vlasov SDE by its empirical projection defined in Definition \ref{def_empirical_proj} applied to the particle system. Secondly, it allows to quantify the approximation of the distribution of one particle by the distribution of the limiting McKean-Vlasov process with respect to $W_1$. Indeed, denoting by $\Vert \varphi \Vert_{\text{Lip}} := \sup_{x\neq y} \frac{|\varphi(x) -\varphi(y)|}{|x-y|}$ for $\varphi : \R^d \to \R$, the set $$\mathscr{C}:= \left\{ \phi : \PP \to \R, \, \exists \varphi : \R^d \to \R, \, \text{with}\, \Vert\varphi \Vert_{\text{Lip}}  \leq 1,\, \text{and}\, \phi(\mu) = \int_{\R^d} \varphi \,d\mu, \, \forall \mu \in \PP \right\}$$ is contained in $\CC^{(2,1)}_1(\PP)$. Thus, the estimate \eqref{eq3_thm_poc} and Kantorovich-Rubinstein's theorem \cite[Remark $6.5$]{Villaniold} ensure that \begin{align}\label{mean_field_limit_eq} 
		\notag\sup_{t \in [0,T]}W_1([X^{1,N}_t],\mu_t) & =\sup_{t \in [0,T]}\sup_{\varphi,\, \|\varphi\|_{\text{Lip}}\leq1} \left\vert 
		\E \varphi(X^{1,N}_t) -  \int_{\R^d} \varphi \, d\mu_t \right\vert \\	\notag &= \sup_{t \in [0,T]}\sup_{\varphi,\, \|\varphi\|_{\text{Lip}}\leq1} \left\vert 
		\E \left(\frac1N \sum_{k=1}^N \varphi(X^{k,N}_t)\right) -  \int_{\R^d} \varphi \, d\mu_t \right\vert\\	\notag &=
		\sup_{t \in [0,T]}\sup_{\varphi,\, \|\varphi\|_{\text{Lip}}\leq 1} \left\vert 
		\E \int_{\R^d} \varphi \, d\muu^N_t - \E \int_{\R^d} \varphi \, d\mu_t \right\vert \\ 	\notag&=
		\sup_{t \in [0,T]}\sup_{\phi \in \mathscr{C}} \left\vert 
		\E \phi(\muu^N_t) - \E \phi (\mu_t) \right\vert\\&\leq \frac{C_T}{N^\gamma},
	\end{align}
where $\gamma \in (0,1] \cap (0,(2\alpha -2)\wedge (\eta +\alpha -1))$ and since the constant $C$ in Theorem \ref{Thm_POC} is non-decreasing with respect to $T$.

\end{enumerate}
\end{Rq}
Let us now compare our result with the existing literature.
\begin{Rq}\label{RQ_comparison_POC}
	\begin{enumerate}
		\item Let us formally take $\alpha =2$, which corresponds to the Brownian case treated in \cite{frikha2021backward}. Then, we can take $\gamma =1$ and $\beta =2$ in Theorem \ref{Thm_POC}. The rates of convergence proved in our theorem are precisely those proved in \cite[Theorem $3.6$]{frikha2021backward}, i.e.\ $N^{-\frac12}$ for \eqref{eq1_thm_poc} and $N^{-1}$ for \eqref{eq3_thm_poc}.
		\item In dimension $d=1$, we recover with \eqref{eq1_thm_poc} the same rate of convergence obtained in \cite{FrikhaWellposednessapproximationonedimensional2020}, for the strong propagation of chaos in $L^1$, since $\E W_1(\muu^N_0,\mu_0) \leq CN^{\frac{1}{\beta}-1}$ by \eqref{FG}. 
		\item In \cite{Cavallazzi_Ito_jump}, the example of a nonlinear Ornstein-Uhlenbeck process is treated using the same method. It corresponds to take $$b(t,x,\mu):= x + \int_{\R^d} y \, d\mu(y).$$  However, there is preliminary step in the proof, which consists in removing the jumps larger than the number of particles $N$ from all the noises. This is due to the unboundedness of $b$ with respect to both space and measures variables in this case, which yields weaker estimates on the semigroup. It is proved in \cite{Cavallazzi_Ito_jump} that there exists a positive constant $C$ such that for any $\phi \in \CC^{(2,1)}_1(\PP)$ \begin{equation}\label{eq1_thm_poc_OU}
			\E \left\vert \phi(\muu^N_T) - \phi(\mu_T) \right\vert \leq 	C\,\E W_1(\muu^N_0,\mu_0) +  C\frac{\ln(N)^{\frac{1}{\alpha}}}{N^{1-\frac{1}{\alpha}}},\end{equation} which is better than \eqref{eq1_thm_poc} in spite of our stronger assumption of boundedness on $b$. Indeed, by removing the large jumps in a first step, we can take $\beta = \alpha$ in \eqref{eq1_thm_poc}, up to the logarithmic factor present in \eqref{eq1_thm_poc_OU}. This factor precisely comes from the fact that $\int_{1 \leq |z| \leq N } |z|^\alpha \, d\nu(z) \underset{N \to +\infty}{\sim} \ln(N),$ which is the price to pay to take $\beta = \alpha$ in \eqref{eq1_thm_poc}. The estimate \eqref{eq3_thm_poc} is better in our framework since the rate of the corresponding estimate in \cite{Cavallazzi_Ito_jump} is $N^{1- \alpha}$ and we can take $\gamma >\alpha -1$ in Theorem \ref{Thm_POC}. This is natural since the drift is unbounded in \cite{Cavallazzi_Ito_jump}.
	\end{enumerate}
\end{Rq}

\subsection{Quantitative approximation of the distribution of one particle by the limiting McKean-Vlasov process at the level of densities}
We keep the same notations as in the preceding subsection. We present here the quantitative propagation of chaos result at the level of densities for the particle system \eqref{edsparticles}. Let us first introduce the following assumption which deals with the existence of a density for \eqref{edsparticles}.

\begin{Assumption}\label{Assumption_ex_density}
We assume that for any $t>0$, the particles $(X^{1,N}_t, \dots, X^{N,N}_t)$ defined by \eqref{edsparticles} with initial distribution $\mu_0 \in \PP$ at time $0$ have a density on $(\R^d)^N$ denoted by $\bm{p}^N(\mu_0,0,t,\cdot)$. 
\end{Assumption}

This assumption is implied by Assumption \ref{Assumption2} in the two following cases.

\begin{Rq}
\begin{itemize} 
 \item When $d=1$, Assumption \ref{Assumption_ex_density} is implied by Assumption \ref{Assumption2} using \cite[Theorem $1.1$]{Chen_existence_density}. Indeed, the map $\bm{x} = (x_1,\dots,x_N) \in (\R^d)^N \mapsto (b(t,x_i,\muu^N_{\bm{x}}))_{i \in \{1,\dots,N\}}$ is Hölder continuous uniformly with respect to $t \in [0,T]$. 
 \item	When $b$ does not depend on time and $d >1$, Assumption \ref{Assumption_ex_density} is implied by Assumption \ref{Assumption2} by \cite[Theorem $1.1$]{DEBUSSCHE20131757}. If we want to avoid completely Assumption \ref{Assumption_ex_density}, we need to study if we can extend \cite[Theorem $1.1$]{Chen_existence_density} when $d >1$.
 \end{itemize}

\end{Rq}

 Under Assumption \ref{Assumption_ex_density}, the density of the first particle exists and is denoted by $p^{1,N}(\mu_0,s,t,\cdot).$ It is given, for all $y_1 \in \R^d$, by $$ p^{1,N}(\mu_0,0,t,y_1) = \int_{(\R^d)^{N-1}} \bm{p}^N(\mu_0,0,t,y_1,y_2,\dots,y_N) \, dy_2 \dots dy_N.$$ Note that by exchangeability, all the particles have the same distribution. Let us recall that $p(\mu_0,0,t,\cdot)$, defined in \eqref{decoupleddensities_link}, is the density on $\R^d$ of $\mu_t$, which is the distribution of the solution to the McKean-Vlasov SDE \eqref{McKVSDE} at time $t$ with initial distribution $\mu_0 \in \PP$ at time $0$. 
 
 \begin{Thm}\label{Thm_poc_density}
 	
 	Let us fix $\gamma \in (\alpha -1,1] \cap (\alpha -1,(2\alpha -2) \wedge (\eta +\alpha -1))$ and $\gamma' \in [\alpha -1,1]$. We define for $\mu \in \PP$, $t \in (0,T]$ and $y,z \in \R^d$ $$q_0(\mu,0,t,y) := \int_{\R^d} \rho^0(t,y-x) \, d\mu(x), \quad \text{and}\quad f(z):= |z|^{1+\gamma} \1_{B_1}(z) + |z| \1_{B_1^c}(z),$$ where $\rho^0$ was defined by \eqref{defdensityrefMK}. We also set $\zeta := - \left(1 - \frac{2+\gamma}{\alpha}\right) \in (0,1)$ since $\gamma \in (0, 2\alpha - 2)$, and we denote by $\BB$ the Beta function defined, for all $x,y >0$, by \[ \BB(x,y) := \int_0^1 (1-t)^{-1+x} t^{-1+y}  \, dt.\] Then, under Assumptions \ref{Assumption2} and \ref{Assumption_ex_density}, the following properties are satisfied. \begin{itemize}
 		 \item \textbf{(Upper-bound for the density of one particle).}  There exists $C= C(d,T,\alpha,\beta,\ref{Assumption2},\gamma)>0$ such that for all $\mu_0 \in \PP$, $t \in (0,T]$, $y \in \R^d$ and $N\geq 1$
 	
 	\begin{align}\label{density_estimate_one_particle}
 		 p^{1,N}(\mu_0,0,t,y) &\leq C q_0(\mu_0,0,t,y) + \sum_{k=1}^{\infty} \frac{C^{k+1}}{N^{k\gamma}} \frac{t^{k(1-\zeta)}}{k(1 - \zeta)} \left(\prod_{j=1}^{k-1} \BB(j(1-\zeta), 1-\zeta)\right)\\ &\notag \hspace{5cm} \sum_{I \in P_{k}} \int_{(\R^d)^k} q_0\left(\mu_0,0,t,y- \sum_{i \in I} z_i\right) \prod_{j=1}^{k} f(z_j) \, d\nu(z_j),
 	\end{align}
 where $P_k$ denotes the set of all subsets of $\{1,\dots,k\}$ and by convention $ q_0\left(\mu_0,0,t,y- \sum_{i \in \emptyset} z_i\right) := q_0\left(\mu_0,0,t,y\right)$.
 
 \item \textbf{(Pointwise estimate for the approximation of the density of one particle).}  There exists $C= C(d,T,\alpha,\beta,\ref{Assumption2},\gamma,\gamma')>0$ such that for all $\mu_0 \in \PP$, $t \in (0,T]$, $y \in \R^d$ and $N \geq 1$
 	
 	\begin{align}\label{eq_thm_poc_density}
 		\notag\left\vert p^{1,N}(\mu_0,0,t,y) - p(\mu_0,0,t,y) \right\vert &\leq \frac{C}{N^{\gamma'}}t^{1-\frac{1+\gamma'}{\alpha}} (1+ M_{\gamma'}(\mu_0)) \int_{\R^d} (1 + |x|^{\gamma'}) \rho^0(t,y-x) \, d\mu_0(x) \\ &\quad+  \sum_{k=1}^{\infty} \frac{C^{k+1}}{N^{k\gamma}} t^{k(1-\zeta)} \left(\prod_{j=1}^{k-1} \BB(j(1-\zeta), 1-\zeta) \right)\BB(1 + k(1-\zeta),1-\zeta) \\ & \notag\hspace{3cm} \sum_{I \in P_{k}} \int_{(\R^d)^{k}} q_0\left(\mu_0,0,t,y- \sum_{i \in I} z_i\right) \prod_{j=1}^{k} f(z_j) \, d\nu(z_j).
 	\end{align}
 \item \textbf{(Estimate for the approximation of the distribution of one particle in total variation).} There exists $C= C(d,T,\alpha,\beta,\ref{Assumption2},\gamma)>0$ such that for all initial distribution $\mu_0 \in \PP$, $t \in (0,T]$ and $N \geq 1$  \begin{equation}\label{eq_thm_poc_density_TV1}
 	d_{TV}([X^{1,N}_t],\mu_t) \leq \frac{C}{N^\gamma}t^{1-\frac{1+\gamma}{\alpha}}(1+M_\gamma(\mu_0)),
 \end{equation}
and \begin{equation}\label{eq_thm_poc_density_TV2}
	\sup_{t\in [0,T]} d_{TV}([X^{1,N}_t],\mu_t) \leq \frac{C}{N^{\alpha -1}}(1 + M_{\alpha-1}(\mu_0)).
\end{equation}

	\end{itemize} 
 \end{Thm}

The proof is given in Section \ref{section_poc_density}.

\begin{Rq}
Let us compare this result with the estimate \eqref{mean_field_limit_eq} quantified with the Wasserstein metric $W_1$. We find the same rate of convergence $N^{-\gamma}$ with $\gamma \in (0,1] \cap (0,(2\alpha -2)\wedge (\eta +\alpha -1)).$ However in \eqref{eq_thm_poc_density_TV1}, there is a time-integrable singularity, similarly to \cite[Theorem $3.5$]{frikha2021backward}, which is not present in \eqref{mean_field_limit_eq}. This time singularity can be removed with the slower rate of convergence $N^{1-\alpha}$. We recover the rate of convergence $N^{-1}$ shown in \cite[Theorem $3.5$]{frikha2021backward} in the Brownian case by taking formally $\alpha =2$ in the preceding result. 
\end{Rq}

 \section{Well-posedness of the nonlinear martingale problem and Picard iterations}\label{section_well_posedness}

This section is dedicated to the proof of Theorem \ref{Thm_wellposedness}. The proof is based on the Banach fixed point theorem on a suitable complete metric space. \\

\noindent \textbf{Introduction of the complete space and parametrix expansion.} Let us consider the space $\CC^0([s,T];\PPP)$ which is complete under the uniform metric $d_{s,T}$ associated to the total variation metric $d_{TV}$ defined, for $P,Q \in \CC^0([s,T];\PPP)$ by \[d_{s,T}(P,Q):= \sup_{r\in[s,T]} d_{TV}(P_r,Q_r).\]   We introduce the space \begin{equation*}
		\mathcal{A}_{s,T,\mu}:= \left\{ P \in \CC^0([s,T];\PPP), \, P_s=\mu\right\}.
	\end{equation*} 
Note that it is a closed subspace of $(\CC^0([s,T];\PPP),d_{s,T})$ and thus $(\mathcal{A}_{s,T,\mu},d_{s,T})$ is complete. For any $P \in \mathcal{A}_{s,T,\mu}$, we consider the following linear time-inhomogeneous SDE

\begin{equation*}
	\begin{cases}
		d\overline{X}^{s,\xi,P}_t = b(t,\overline{X}^{s,\xi,P}_t, P_t)\,dt + dZ_t, \quad  t \in [s,T], \\ \overline{X}^{s,\xi,P}_s = \xi.
	\end{cases}
\end{equation*}
Notice that this SDE is well-posed in the weak sense since it is the case for the related linear martingale problem by \cite{Mikulevicius_MP_stable}. Its flow of marginal distributions $([X_t^{s,\xi,P}])_{t \in [s,T]}$ belongs to $\mathcal{A}_{s,T,\mu}$. We can thus define a map $\mathcal{I} : \mathcal{A}_{s,T,\mu} \rightarrow \mathcal{A}_{s,T,\mu}$ such that for any $P \in \mathcal{A}_{s,T,\mu}$, $\mathcal{I}(P)_t = [X_t^{s,\xi,P}].$ We remark that a probability measure $\P$ on the Skorokhod space $\mathcal{D}([s,T];\R^d)$ solves the martingale problem related to the McKean-Vlasov SDE \eqref{McKVSDE_existence} if and only if its flow of marginal distributions $(\P_t)_{t \in [s,T]}$ is a fixed point of $\mathcal{I}$. Our goal is thus to prove that for some $m \geq 1$, the $m$-th iterate $\mathcal{I}^{(m)}$ is a contraction on $(\mathcal{A}_{s,T,\mu},d_{s,T})$. We fix $P^1,P^2 \in \mathcal{A}_{s,T,\mu}$ and we define recursively for all $m \geq 1$, $\overline{X}^{1,(m)}$ and  $\overline{X}^{2,(m)}$ as the unique weak solutions to \begin{equation}\label{SDE_iteree_existence}
	\begin{cases}
		d\overline{X}^{i,(m)}_t = b(t,\overline{X}^{i,m}_t,[\overline{X}^{i,(m-1)}_t]))\,dt + dZ_t, \quad  t \in [s,T], \quad   i\in \{1,2\},\\ \overline{X}^{i,(m)}_s = \xi,
	\end{cases}
\end{equation}
with $([\overline{X}^{i,(0)}_t])_{t \in [s,T]} = P^i$. We also introduce the associated decoupling fields. Namely they are the weak solutions to 
\begin{equation}\label{SDE_decoupled_iteree_existence}
	\begin{cases}
		d\overline{X}^{x,i,(m)}_t = b(t,\overline{X}^{x,i,m}_t,[\overline{X}^{i,(m-1)}_t]))\,dt + dZ_t, \quad t \in [s,T], \quad i\in \{1,2\},\\ \overline{X}^{x,i,(m)}_s = x.
	\end{cases}
\end{equation}

Thanks to Theorem \ref{Thmdensityparametrix}, the distribution of $\overline{X}^{x,i,(m)}_t$ has a density with respect to the Lebesgue measure denoted by $p_{i,m}(\mu,s,t,x,\cdot)$. Remark that the notation makes sense since, by weak well-posedness, the distribution of $\overline{X}^{x,i,(m)}_t$ depends on the initial condition $\xi$ only through its distribution $\mu$. Moreover, by weak well-posedness of SDE \eqref{SDE_iteree_existence}, $\overline{X}^{i,(m)}_t$ has a density $p_{i,m}(\mu,s,t,\cdot)$ such that for all $y \in \R^d$ \[ p_{i,m}(\mu,s,t,y) = \int_{\R^d} p_{i,m}(\mu,s,t,x,y) \, d\mu(x).\] Let us give the implicit parametrix representation of $p_{i,m}(\mu,s,t,x,y)$, which is given in Appendix \ref{section_appendix_parametrix}. We define for all $0\leq s \leq r< t \leq T$ and $x,y \in \R^d$ \begin{align}\label{def_proxy_kernel_existence}
	& \p(r,t,x,y) := q(t-r,y-x), \\ &\notag \H_{i,m}(\mu,r,t,x,y) := b(r,x,[\overline{X}^{i,(m-1)}_r]) \cdot\pa_x \p(r,t,x,y),
\end{align}
where $q(t,\cdot)$ is the density of $Z_t$. The space-time convolution between to functions $f$ and $g$ is given by \begin{equation}\label{defconvolopMK_existence} f \otimes g (\mu,r,t,x,y) := \int_r^t \int_{\R^d} f(\mu,r,r',x,z) g(\mu,r',t,z,y) \, dz \, dr',
\end{equation}
when it is well-defined. The convolution iterates of $\H_{i,m}$ are defined recursively by $\H_{i,m}^{k+1} = \H_{i,m} \otimes \H_{i,m}^{k}$. By convention $\H_{i,m}^0 = \text{Id}$. By Assumption \ref{Assumption1}, we can apply Theorem \ref{Thmdensityparametrix} which ensures that  \begin{equation}\label{rep_density_existence}
 	p_{i,m}(\mu,s,t,x,y) = \p(s,t,x,y) + p_{i,m}\otimes \H_{i,m}(\mu,s,t,x,y).
 \end{equation}
Using Theorem \ref{Thmdensityparametrix} and Proposition \ref{Prop_H^k}, we deduce that there exists a positive constant $K$ such that for all $i \in \{1,2\}$, $m\geq 1$, $k \geq1$, $0 \leq s \leq r <t \leq T$ and $x,y \in \R^d$ \begin{equation}\label{control_density_kernel_existence}\begin{cases}
\left\vert p_{i,m}(\mu,s,t,x,y)\right\vert \leq K \rho^0(t-s,y-x), \\ \left\vert \H_{i,m}^k(\mu,r,t,x,y)\right\vert \leq  K^k (t-r)^{-\frac{1}{\alpha} + (k-1)\left( 1 - \frac{1}{\alpha} \right)} \displaystyle\prod_{j=1}^{k-1} \BB \left( j\left( 1 - \frac{1}{\alpha} \right),1 - \frac{1}{\alpha} \right) \rho^1(t-r,y-x),\end{cases}
\end{equation}
where the functions $\rho^j$ were defined in \eqref{defdensityrefMK} and $\BB$ is the Beta function defined, for all $x,y >0$ by \[ \BB(x,y) := \int_0^1 (1-t)^{-1+x} t^{-1+y}  \, dt = \frac{\Gamma(x) \Gamma(y)}{\Gamma(x+y)},\] $\Gamma$ being the Gamma function. For any $ 0 \leq s \leq r <t \leq T$, $x,y \in \R^d$, we define \begin{align*}\Delta p_{m}(\mu,s,t,x,y) &:= p_{1,m}(\mu,s,t,x,y) - p_{2,m}(\mu,s,t,x,y),\\ \Delta \H_{m}(\mu,r,t,x,y)&:= \H_{1,m}(\mu,r,t,x,y) - \H_{2,m}(\mu,r,t,x,y).\end{align*}

It follows from the definition of $\H_{i,m}$, the Lipschitz continuity of $b(s,x,\cdot)$ with respect to the total variation metric (and its expression as the $L^1$ norm of the difference between the densities) and \eqref{gradientestimatestable}, that for some constant $C>0$, one has for all $m \geq 1$ \begin{align}\label{proof_existence_eq2}
\left\vert \Delta \H_{m+1} (\mu,r,t,x,y)\right\vert &\leq C d_{TV}([\overline{X}^{1,(m)}_r],[\overline{X}^{2,(m)}_r])(t-r)^{-\frac{1}{\alpha}} \rho^1(t-r,y-x) \\ \notag&\leq C  \int_{\R^{2d}} \left\vert \Delta p_m(\mu,s,r,x',y')\right\vert d\mu(x') \, dy' \,(t-r)^{-\frac{1}{\alpha}} \rho^1(t-r,y-x).
\end{align}
We similarly obtain for $m=0$ that  \begin{align}\label{proof_existence_eq3}
	\left\vert \Delta \H_{1} (\mu,r,t,x,y)\right\vert &\leq C d_{s,r}(P^1,P^2) (t-r)^{-\frac{1}{\alpha}} \rho^1(t-r,y-x) .
\end{align} 

Let us prove that for all $m \geq 1$, the following representation formula holds true \begin{equation}\label{representation_difference_densities_existence}
	\Delta p_m (\mu,s,t,x,y) =\sum_{k=0}^{\infty} p_{2,m}\otimes \Delta \H_{m}\otimes \H_{1,m}^k(\mu,s,t,x,y),
\end{equation}
the series being absolutely convergent. Starting from the implicit parametrix representation formula \eqref{rep_density_existence}, we have \begin{equation*}
	\Delta p_m(\mu,s,t,x,y) = p_{2,m}\otimes \Delta \H_{m} (\mu,s,t,x,y) =  \Delta p_m \otimes \H_{1,m}(\mu,s,t,x,y).
\end{equation*}
Iterating this procedure, we easily prove by induction that for all $N \geq 1$ 

\begin{equation}\label{proof_existence_eq1}
	\Delta p_m(\mu,s,t,x,y) =  \sum_{k=0}^N p_{2,m}\otimes \Delta \H_{m} \otimes \H_{1,m}^k(\mu,s,t,x,y) + \Delta p_m \otimes \H_{1,m}^{N+1}(\mu,s,t,x,y).
\end{equation}
We want to take the limit $N \to +\infty$ in \eqref{proof_existence_eq1}. By using \eqref{control_density_kernel_existence} and the convolution inequality \eqref{convolineqdensityref}, we obtain that for some positive constant $C$

\begin{align*}
	&\left\vert \Delta p_m \otimes \H_{1,m}^{N+1}(\mu,s,t,x,y) \right\vert \\ &\leq \int_s^t \int_{\R^d} C\rho^0(r-s,z-x)  C^{N+1} (t-r)^{-\frac{1}{\alpha} + N\left( 1 - \frac{1}{\alpha} \right)} \prod_{j=1}^{N} \BB \left( j\left( 1 - \frac{1}{\alpha} \right),1 - \frac{1}{\alpha} \right) \rho^1(t-r,y-z) \, dz \,dr \\ &\leq  C^{N+2} (t-s)^{(N+1)\left( 1 - \frac{1}{\alpha} \right)} \prod_{j=1}^{N} \BB \left( j\left( 1 - \frac{1}{\alpha} \right),1 - \frac{1}{\alpha} \right) \rho^0(t-s,y-x).
\end{align*}
The upper-bound converges to $0$ as $N$ tends to infinity thanks to the asymptotic behavior of the Beta function. Following the same lines, we prove using again \eqref{control_density_kernel_existence} that the series appearing in \eqref{proof_existence_eq1} is absolutely convergent. Letting $N$ tend to infinity in \eqref{proof_existence_eq1} yields the representation formula \eqref{representation_difference_densities_existence}. \\

We are now going to prove by induction that there exists a positive constant $C$ such that for all $m \geq 1$, $ t \in (s,T]$, $x,y \in \R^d$ \begin{equation}\label{proof_existence_eq4}
	\left\vert \Delta p_m(\mu,s,t,x,y) \right\vert \leq C^m(t-s)^{m\left(1 - \frac{1}{\alpha}\right)} \prod_{j=1}^{m-1} \BB\left(1 + j\left(1-\frac{1}{\alpha}\right),1 - \frac{1}{\alpha}\right)d_{s,t}(P^1,P^2) \rho^{0}(t-s,y-x).
\end{equation}

\noindent\textbf{Base case $\bm{m=1}$.} It follows from \eqref{control_density_kernel_existence}, \eqref{proof_existence_eq3} and the convolution inequality \eqref{convolineqdensityref} that \begin{align}\label{proof_existence_eq8}
	\left\vert p_{2,1}\otimes \Delta \H_{1}(\mu,s,t,x,y) \right\vert\notag &\leq C \int_s^t \int_{\R^d} \rho^0(r-s,z-x) (t-r)^{-\frac{1}{\alpha}}d_{s,r}(P^1,P^2) \rho^1(t-r,y-z) \, dz\,dr \\ &\leq C (t-s)^{1- \frac{1}{\alpha}} d_{s,t}(P^1,P^2) \rho^0(t-s,y-x),
\end{align}
since $d_{s,r}(P^1,P^2) \leq d_{s,t}(P^1,P^2)$ for $r \in [s,t]$. Following the same lines and using the bound \eqref{control_density_kernel_existence} on $\H_{1,1}^k$, we show that for some constant $C>0$, one has for all $k\geq 1$
\begin{align*}
	&\left\vert p_{2,1} \otimes \Delta \H_{1} \otimes \H_{1,1}^k(\mu,s,t,x,y) \right\vert \\ &\leq  C(t-s)^{1 - \frac{1}{\alpha}}d_{s,t}(P^1,P^2)  \rho^0(t-s,y-x) C^k(t-s)^{k\left(1- \frac{1}{\alpha}\right)} \prod_{l=1}^{k-1} \BB\left(l\left(1- \frac{1}{\alpha}\right),1- \frac{1}{\alpha}\right).
\end{align*} 
By summing over $k \geq 1$, we deduce using the asymptotic behavior of the Beta function that for some constant $C>0$, we have \begin{align}\label{proof_existence_eq7}
	&\sum_{k=1}^{\infty}\left\vert p_{2,1} \otimes \Delta \H_{1} \otimes \H_{1,1}^k(\mu,s,t,x,y) \right\vert \\ \notag&\leq  C(t-s)^{1 - \frac{1}{\alpha}} d_{s,t}(P^1,P^2)  \rho^0(t-s,y-x).
\end{align}

Plugging \eqref{proof_existence_eq8} and \eqref{proof_existence_eq7} into the representation formula \eqref{representation_difference_densities_existence} for $m=1$ concludes the proof of the base case.\\
 
\noindent\textbf{Induction step.} We assume that $\eqref{proof_existence_eq4}$ holds true at step $m$ for a certain constant $C$ that will be chosen at the end of the induction step to ensure that \eqref{proof_existence_eq4} is verified at step $m+1$. We denote by $K$ any constant independent of $m$ and $C$ appearing in the induction step. By using \eqref{proof_existence_eq2}, one has \begin{multline*}
\left\vert \Delta \H_{m+1} (\mu,r,t,x,y) \right\vert \leq K (t-r)^{-\frac{1}{\alpha}}  \rho^{1}(t-r,y-x) \\C^m(r-s)^{m\left(1 - \frac{1}{\alpha}\right)} \prod_{j=1}^{m-1} \BB\left(1 + j\left(1-\frac{1}{\alpha}\right),1 - \frac{1}{\alpha}\right)d_{s,r}(P^1,P^2).
\end{multline*}

This inequality combined with \eqref{control_density_kernel_existence} yields

\begin{align}\label{proof_existence_eq5}
&\left\vert p_{2,m+1} \otimes \Delta \H_{m+1}(\mu,s,t,x,y) \right\vert \\ \notag&\leq K\int_s^t \int_{\R^d} \rho^0(r-s,z-x)  (t-r)^{-\frac{1}{\alpha}}  \rho^{1}(t-r,y-z) C^m(r-s)^{m\left(1 - \frac{1}{\alpha}\right)} \\ \notag&\hspace{6cm}\prod_{j=1}^{m-1} \BB\left(1 + j\left(1-\frac{1}{\alpha}\right),1 - \frac{1}{\alpha}\right)d_{s,r}(P^1,P^2) \, dz \,dr \\ \notag&\leq K C^m(t-s)^{(m+1)\left(1 - \frac{1}{\alpha}\right)} \prod_{j=1}^{m} \BB\left(1 + j\left(1-\frac{1}{\alpha}\right),1 - \frac{1}{\alpha}\right)d_{s,t}(P^1,P^2)  \rho^0(t-s,y-x).
\end{align} 

Following the same lines and using the bound \eqref{control_density_kernel_existence} on $\H_{1,m+1}^k$, we show that for some constant $K>0$, one has  for all $k\geq 1$
\begin{align*}
	&\left\vert p_{2,m+1} \otimes \Delta \H_{m+1} \otimes \H_{1,m+1}^k(\mu,s,t,x,y) \right\vert \\ &\leq K C^m(t-s)^{(m+1)\left(1 - \frac{1}{\alpha}\right)} \prod_{j=1}^{m} \BB\left(1 + j\left(1-\frac{1}{\alpha}\right),1 - \frac{1}{\alpha}\right)d_{s,t}(P^1,P^2)  \rho^0(t-s,y-x) \\ &\hspace{3cm} K^k(t-s)^{k\left(1- \frac{1}{\alpha}\right)} \prod_{l=1}^{k-1} \BB\left(l\left(1- \frac{1}{\alpha}\right),1- \frac{1}{\alpha}\right).
\end{align*} 
By summing over $k \geq 1$, we deduce, thanks to the asymptotic behavior of the Beta function, that for some constant $K>0$, we have \begin{align}\label{proof_existence_eq6}
	&\sum_{k=1}^{\infty}\left\vert p_{2,m+1} \otimes \Delta \H_{m+1} \otimes \H_{1,m+1}^k(\mu,s,t,x,y) \right\vert \\ \notag&\leq K C^m(t-s)^{(m+1)\left(1 - \frac{1}{\alpha}\right)} \prod_{j=1}^{m} \BB\left(1 + j\left(1-\frac{1}{\alpha}\right),1 - \frac{1}{\alpha}\right)d_{s,t}(P^1,P^2)  \rho^0(t-s,y-x).
\end{align} 
Finally, plugging \eqref{proof_existence_eq5} and \eqref{proof_existence_eq6} into the representation formula \eqref{representation_difference_densities_existence}, we obtain 
 \begin{equation*}
	\left\vert \Delta p_{m+1}(\mu,s,t,x,y) \right\vert \leq KC^m(t-s)^{(m+1)\left(1 - \frac{1}{\alpha}\right)} \prod_{j=1}^{m} \BB\left(1 + j\left(1-\frac{1}{\alpha}\right),1 - \frac{1}{\alpha}\right)d_{s,t}(P^1,P^2) \rho^{0}(t-s,y-x).
\end{equation*}
This ends the proof of the induction step provided that we choose $C\geq K$ in \eqref{proof_existence_eq4}, which is possible since $K$ does not depend on $m$.\\

\noindent\textbf{Conclusion of the proof of Theorem \ref{Thm_wellposedness}.} Using the asymptotic behavior of the Beta function, we get that for $m$ large enough, we have for all $t \in [s,T]$, $x,y \in \R^d$ \[ \left\vert \Delta p_m(\mu,s,t,x,y)\right\vert \leq \eps d_{s,t}(P^1,P^2) \rho^0(t-s,y-x),\]
where $\eps >0$ is such that \begin{equation*}
	\eps \int_{\R^d}\rho^0(t-s,y) \, dy = \eps\int_{\R^d} \rho^0(1,y) \, dy = \frac{1}{2}.
\end{equation*}
Finally, we have 

\begin{align*}
 d_{s,T}(\mathcal{I}^{(m)}(P^1),\mathcal{I}^{(m)}(P^2)) &= \sup_{t \in (s,T]} \sup_{h,\, \|h\|_{\infty}\leq1} \left\vert \int_{\R^d} h(y) (p_{1,m}(\mu,s,t,y) - p_{2,m}(\mu,s,t,y)) \, dy \right\vert \\ &\leq  \sup_{t \in (s,T]} \int_{\R^{2d}} \left\vert \Delta p_m(\mu,s,t,x,y) \right\vert d\mu(x)\, dy \\ &\leq \frac12 d_{s,T}(P^1,P^2).
\end{align*}

The Banach fixed point theorem ensures that $\mathcal{I}$ has a unique fixed point in $\mathcal{A}_{s,T,\mu}$. Thus, the martingale problem associated to the McKean-Vlasov SDE \eqref{McKVSDE_existence} is well-posed. Moreover, we know that for any initial data $P \in \mathcal{A}_{s,T,\mu}$, the sequence $(\mathcal{I}^{(m)}(P))_{m \geq 1}$ converges towards the solution to the martingale problem with respect to the metric $d_{s,T}$. This proves \eqref{convergence_Thm_existence}.

 \section{Properties of the transition density}\label{section_properties_transition_densit}
 
  This section is devoted to prove Theorem \ref{Thm_density_McKean_light}. It will be a direct consequence of the following result.
  
  \begin{Thm}[Regularity estimates on the decoupled transition density]\label{Thm_density_McKean}
  	Let us fix $0 \leq s < t \leq T$ and $y\in \R^d$. Under Assumption \ref{Assumption2}, the mapping $(\mu,x) \in \PP \times \R^d \mapsto p(\mu,s,t,x,y)$ belongs to $\CC^{1}(\PP \times  \R^d)$ (see Definition \eqref{def_C1_space}). Moreover, it satisfies the following properties.
  	
  	\begin{itemize}

  		\item There exists $C>0$ such that for all $j \in \{0,1\}$, $\mu \in \PP$, $0 \leq s < t \leq T$ and $x,y \in \R^d$
  		
  		\begin{equation}\label{density_bound_thm}
  			|\pa_x^j p(\mu,s,t,x,y) |
  			\leq C (t-s)^{-\frac{j}{\alpha}} \rho^{j}(t-s,y-x).
  		\end{equation}

  		\item For all $j \in \{ 0,1\}$ and $\gamma \in (0, 1]$ with $\gamma \in (0, (2\alpha -2) \wedge (\eta + \alpha -1))$ if $j=1$, there exists $C>0$ such that for all $\mu \in \PP$, $0 \leq s < t \leq T$ and $x_1,x_2,y \in \R^d$
  		
  		\begin{equation}\label{gradientdensity_holder_thm}
  			|\pa_x^j p(\mu,s,t,x_1,y) -  \pa^j_x p(\mu,s,t,x_2,y)| \leq C (t-s)^{-\frac{ j+\gamma}{\alpha}} |x_1-x_2|^\gamma \left[\rho^{j}(t-s,y-x_1) +\rho^{j}(t-s,y-x_2)\right].
  		\end{equation}

  		\item There exists $C>0$ such that for all $\mu \in \PP$, $ 0 \leq s < t \leq T$, $x,y,v \in \R^d$ \begin{equation}\label{linear_der_bound_thm}
  			\left\vert \del p(\mu,s,t,x,y)(v)\right\vert \leq C (t-s)^{1 - \frac{1}{\alpha}}\rho^0(t-s,y-x).
  		\end{equation}
  		
  		\item There exists $C>0$ such that for all $\mu \in \PP$, $ 0 \leq s < t \leq T$, $x,y,v \in \R^d$ \begin{equation}\label{gradient_linear_der_bound_thm}
  			\left\vert \pa_v\del p(\mu,s,t,x,y)(v)\right\vert \leq C (t-s)^{\frac{\eta -1}{\alpha} + 1 - \frac{1}{\alpha}}\rho^0(t-s,y-x).
  		\end{equation}

  		\item For all $\gamma \in (0,1]\cap(0,(2\alpha -2)\wedge (\eta + \alpha -1))$, there exists $C>0$ such that for all $\mu \in \PP$, $ 0 \leq s < t \leq T$, $x,y,v_1,v_2 \in \R^d$ \begin{equation}\label{gradient_linear_der_holder_v_thm}
  			\left\vert \pa_v\del p(\mu,s,t,x,y)(v_1) - \pa_v \del p(\mu,s,t,x,y)(v_2)\right\vert \leq C (t-s)^{\frac{\eta -1 - \gamma}{\alpha} + 1 - \frac{1}{\alpha}}|v_1-v_2|^\gamma\rho^0(t-s,y-x).
  		\end{equation}
  		
  		\item For all $\gamma \in (0,1]$, there exists $C>0$ such that for all , $\mu \in \PP$, $ 0 \leq s < t \leq T$, $x,y,v_1,v_2 \in \R^d$ \begin{equation}\label{linear_der_holder_v_thm}
  			\left\vert \del p(\mu,s,t,x,y)(v_1) - \del p(\mu,s,t,x,y)(v_2)\right\vert \leq C (t-s)^{ 1 - \frac{1+\gamma }{\alpha}}|v_1-v_2|^\gamma\rho^0(t-s,y-x).
  		\end{equation}
  		
  		\item For all $\gamma \in (0,1]$, there exists $C>0$ such that for all $\mu \in \PP$, $ 0 \leq s < t \leq T$, $x_1,x_2,y,v\in \R^d$ \begin{multline}\label{linear_der_holder_x_thm}
  			\left\vert \del p(\mu,s,t,x_1,y)(v) - \del p(\mu,s,t,x_2,y)(v)\right\vert \leq C (t-s)^{ 1 - \frac{1+\gamma}{\alpha}}|x_1-x_2|^\gamma \\ \left[\rho^0(t-s,y-x_1) + \rho^0(t-s,y-x_2)\right].
  		\end{multline}
  		
  		\item For all $\gamma \in (0,1]\cap(0,\eta + \alpha -1)$, there exists $C>0$ such that for all $\mu \in \PP$, $ 0 \leq s < t \leq T$, $x_1,x_2,y,v\in \R^d$ \begin{multline}\label{gradient_linear_der_holder_x_thm}
  			\left\vert \pa_v\del p(\mu,s,t,x_1,y)(v) - \pa_v \del p(\mu,s,t,x_2,y)(v)\right\vert \leq  C(t-s)^{\frac{\eta -1 - \gamma}{\alpha} + 1 - \frac{1}{\alpha}}|x_1-x_2|^\gamma \\ \left[\rho^0(t-s,y-x_1) + \rho^0(t-s,y-x_2)\right].
  		\end{multline}

  		\item For all $j \in \{0,1\}$ and $\gamma \in (0,1]$ with $\gamma \in (0, \alpha -1 + \eta)$ if $j=1$, there exists $C>0$ such that for all $0\leq s <t \leq T$, $\mu_1,\mu_2 \in \PP$, $x,y \in \R^d$ 
  		
  		\begin{equation}\label{gradient_density_holder_measure_thm}
  			\left\vert \pa_x^j p(\mu_1,s,t,x,y) -  \pa_x^j p(\mu_2,s,t,x,y) \right\vert \leq C(t-s)^{1-\frac{1+\gamma +j}{\alpha}} W_1^\gamma (\mu_1,\mu_2) \rho^j(t-s,y-x).
  		\end{equation}
  		
  		\item For all $\gamma \in (0,1]$, there exists $C>0$ such that for all $0 \leq s < t \leq T$, $\mu_1,\mu_2 \in \PP$, $x,y,v \in \R^d$
  		
  		\begin{equation}\label{linear_derivative_holder_measure_thm}
  			\left\vert \del p (\mu_1,s,t,x,y)(v) -  \del p (\mu_2,s,t,x,y)(v)\right\vert \leq C(t-s)^{-\frac{\gamma}{\alpha} + 1 - \frac{1}{\alpha}}W_1^\gamma(\mu_1,\mu_2)\rho^0(t-s,y-x).
  		\end{equation}
  		
  		\item For all $\gamma \in (0,1]$, there exists $C>0$ such that for all $0 \leq s < t \leq T$, $\mu_1,\mu_2 \in \PP$, $x,y,v \in \R^d$
  		
  		\begin{equation}\label{gradient_linear_derivative_holder_measure_thm}
  			\left\vert \pa_v\del p (\mu_1,s,t,x,y)(v) -  \pa_v\del p (\mu_2,s,t,x,y)(v)\right\vert \leq C (t-s)^{\frac{\eta -\gamma -1}{\alpha} + 1 - \frac{1}{\alpha}}W_1^\gamma(\mu_1,\mu_2)\rho^0(t-s,y-x).
  		\end{equation}

  		\item For all $j \in \{0,1\}$, $ \gamma \in \left(0,1 - \frac{j}{\alpha}\right)$, there exists a constant $C>0$ such that for all $t \in (0,T],$ $s_1,s_2 \in [0,t)$, $\mu \in \PP$, $x,y \in \R^d$ 
  		
  		\begin{multline}\label{gradient_density_time_Holder_thm}
  			|\pa_x^j p(\mu,s_1,t,x,y) - \pa_x^jp(\mu,s_2,t,x,y)| \\ \leq C \left[\frac{|s_1-s_2|^\gamma}{(t - s_1)^{\gamma + \frac{j}{\alpha}}}\rho^{j}(t-s_1,y-x) + \frac{|s_1-s_2|^\gamma}{(t-s_2)^{\gamma+ \frac{j}{\alpha}}}\rho^{j}(t-s_2,y-x)\right].
  		\end{multline}
  		
  		\item For all $\gamma \in \left(0,1\right)$, there exists a constant $C>0$ such that for all $t \in (0,T]$, $s_1,s_2 \in[0,t)$, $\mu \in \PP$, $x,y,v \in \R^d$
  		
  		\begin{align}\label{lin_der_density_holder_time_thm}
  			\notag&\left\vert \del p(\mu,s_1,t,x,y)(v) - \del p(\mu,s_2,t,x,y)(v)\right\vert\\ &\leq C\left[\frac{|s_1-s_2|^\gamma}{(t-s_1\vee s_2)^{\gamma + \frac{1}{\alpha}-1}} \rho^{0}(t-s_1\vee s_2,y-x) + \frac{|s_1-s_2|^\gamma }{(t-s_1\wedge s_2)^{\gamma +\frac{1}{\alpha}-1}} \rho^{0}(t-s_1\wedge s_2,y-x) \right].
  		\end{align}

  		\item For all $\gamma \in \left(0,1+ \frac{\eta -1}{\alpha}\right)$, there exists a constant $C>0$ such that for all $t \in (0,T]$, $s_1,s_2 \in[0,t)$, $\mu \in \PP$, $x,y,v \in \R^d$
  		
  		\begin{align}\label{gradient_lin_der_density_holder_time_thm}
  			&\left\vert \pa_v\del p(\mu,s_1,t,x,y)(v) - \pa_v\del p(\mu,s_2,t,x,y)(v)\right\vert\\\notag &\leq C \left[\frac{|s_1-s_2|^\gamma}{(t-s_1\vee s_2)^{\gamma + \frac{1}{\alpha}-1 + \frac{1-\eta}{\alpha}}} \rho^{0}(t-s_1\vee s_2,y-x) + \frac{|s_1-s_2|^\gamma }{(t-s_1\wedge s_2)^{\gamma +\frac{1}{\alpha}-1+ \frac{1-\eta}{\alpha}}} \rho^{0}(t-s_1\wedge s_2,y-x) \right].
  		\end{align}
  		
  		\item For all $\mu \in \PP$, $0 \leq s < t \leq T$, $x,v \in \R^d$, we have
  		
  		\begin{equation}\label{centering_prop_lin_der_density}
  			\int_{\R^d} \del p(\mu,s,t,x,y)(v) \, dy = \int_{\R^d} \pa_v \del p(\mu,s,t,x,y)(v) \, dy =0.
  		\end{equation}
  	\end{itemize}

  \end{Thm}

   Before proving this result, let us introduce the parametrix method that is at the core of the proof (see Appendix \ref{section_appendix_parametrix} for more details). We denote by $q(t,\cdot)$ the density of $Z_t$. We define for all $0\leq s\leq r < t \leq T$, $\mu \in \PP$ and $x,y \in \R^d$ \begin{align}\label{defproxyMK}
 	& \p(\mu,s,r,t,x,y) = \p(r,t,x,y) := q(t-r,y-x), \\ &\notag \H(\mu,s,r,t,x,y) := b(r,x,[X^{s,\mu}_r]) \cdot\pa_x \p(s,r,t,x,y).
 \end{align}
 Note that the proxy $\p(s,r,t,x,\cdot)$ does not depend on $\mu$ and $s$ and is the density at time $t >r$ of the solution to \begin{equation}\label{SDEproxyMK}
 	\begin{cases}
 		d\widehat{X}^{r,x}_t= dZ_t, \\ \widehat{X}^{r,x}_r=x \in \R^d,
 	\end{cases}
 \end{equation}
 \\
 and $\H$ is the associated parametrix kernel. We also define the space-time convolution operator between to functions $f$ and $g$ by \begin{equation}\label{defconvolopMK} f \otimes g (\mu,s,r,t,x,y) := \int_r^t \int_{\R^d} f(\mu,s,r,r',x,z) g(\mu,s,r',t,z,y) \, dz \, dr',
 \end{equation}
 when it is well-defined. The space-time convolution iterates $\H^k$ of $\H$ are defined recursively by $\H^1 := \H$ and $\H^{k+1} := \H \otimes \H^{k}.$ By convention $f \otimes \H^0$ is equal to $f$. In order to simplify a bit the notations, we will write $f\otimes g (\mu,s,t,x,y):= f \otimes g (\mu,s,s,t,x,y)$, $\H(\mu,s,t,x,y):= \H(\mu,s,s,t,x,y)$, and the same for other maps. Finally, we denote by $\Phi$ the solution to the following Volterra integral equation \[ \Phi(\mu,s,r,t,x,y) = \H(\mu,s,r,t,x,y) + \H \otimes \Phi(\mu,s,r,t,x,y),\] which is given by the uniform convergent series \begin{equation}\label{defvolterraMK}
 	\Phi(\mu,s,r,t,x,y) = \sum_{k=1}^{\infty} \H^k(\mu,s,r,t,x,y).
 \end{equation}
Then, we have using Theorem \ref{Thmdensityparametrix}, that for all $\mu \in \PP$, $0 \leq s < t \leq T$, and $x,y \in \R^d$ \begin{equation}\label{param_rep_density_thm}
	 p(\mu,s,t,x,y) = \p(s,t,x,y) + \sum_{k=0}^{\infty} \p \otimes \H^k(\mu,s,t,x,y).\end{equation}
 
	 Let us now prove Theorem \ref{Thm_density_McKean}.\\

	\textbf{Step 1: Properties of the Picard approximation of the transition density associated to \eqref{decoupled_SDE}.} In order to study the regularity with respect to $\mu$ of $p(\mu,s,t,x,y)$, we consider an approximation sequence based on Picard iteration.
We fix a measure $\nu \in \PP$ and $s\in [0,T)$ and we start by considering the following stable-driven McKean-Vlasov SDE  \begin{equation}\label{McKVSDE_Picard_initialisation}
	\begin{cases}
		dX^{s,\xi,(1)}_t= b(t,X_t^{s,\xi,(1)},\nu)\, dt +dZ_t, \quad  t \in [s,T],\\ X_s^{s,\xi,(1)}=\xi, \quad [\xi] = \mu \in \PP.
	\end{cases} 
\end{equation}
The associated martingale problem is well-posed by \cite{Mikulevicius_MP_stable} and there is weak existence and uniqueness for SDE \eqref{McKVSDE_Picard_initialisation}. As previously, the distribution of $X^{s,\xi,(1)}_t$ is denoted by $[X^{s,\mu,(1)}_t]$. 
We also introduce, for $x \in \R^d$, the following decoupled
stochastic flow associated to SDE \eqref{McKVSDE_Picard_initialisation}

\begin{equation}\label{decoupled_SDE_Picard_initialisation}
	\begin{cases}
		dX^{s,x,\mu,(1)}_t= b(t,X_t^{s,x,\mu,(1)},[X^{s,\mu,(1)}_t])\, dt +dZ_t, \quad  t \in [s,T],\\ X_s^{s,x,\mu,(1)}=x \in \R^d.
	\end{cases}
\end{equation}

Then, for all $m \geq 1$, we define recursively

\begin{equation}\label{McKVSDE_Picard}
	\begin{cases}
		dX^{s,\xi,(m+1)}_t= b(t,X_t^{s,\xi,(m+1)},[X^{s,\mu,(m)}_t])\, dt +dZ_t, \quad  t \in [s,T], \\ X_s^{s,\xi,(m+1)}=\xi, \quad [\xi] = \mu \in \PP.
	\end{cases} 
\end{equation}

and 

\begin{equation}\label{decoupled_SDE_Picard}
	\begin{cases}
		dX^{s,x,\mu,(m+1)}_t= b(t,X_t^{s,x,\mu,(m+1)},[X^{s,\mu,(m)}_t])\, dt +dZ_t,\quad  t \in [s,T], \\ X_s^{s,x,\mu,(m+1)}=x \in \R^d.
	\end{cases}
\end{equation}

Note that these are not McKean-Vlasov SDEs but linear SDEs. The densities of $[X^{s,\mu,(m)}_t]$ and $[X^{s,x,\mu,(m)}_t]$ exist by \ref{Thmdensityparametrix} and are denoted by $p_m(\mu,s,t,\cdot)$ and $p_m(\mu,s,t,x,\cdot)$ and satisfy \begin{equation}\label{decoupleddensities_link_Picard}
	p_m(\mu,s,t,y) = \int_{\R^d} p_m(\mu,s,t,x,y) \, d\mu(x).
\end{equation}

Let us now give their parametrix expansions. We define for all $0\leq s < t \leq T$, $r \in [s,t)$, $\mu \in \PP$ and $x,y \in \R^d$ \begin{align}\label{def_proxy_kernel_Picard}
	& \p_m(s,r,t,x,y) := \p(r,t,x,y) = q(t-r,y-x), \\ &\notag \H_m(\mu,s,r,t,x,y) := b(r,x,[X^{s,\mu,(m-1)}_r]) \cdot\pa_x \p(s,r,t,x,y).
\end{align}
Note that the proxy $\p_m(s,r,t,x,\cdot)$ does not depend on $m$, $\mu$ and $s$. We denote by $\Phi_m$ the solution to the following Volterra integral equation \begin{equation}\label{Volterra_eq_Picard}
	\Phi_m(\mu,s,r,t,x,y) = \H_m(\mu,s,r,t,x,y) + \H_m \otimes \Phi_m(\mu,s,r,t,x,y),\end{equation} which is given by the uniform convergent series \begin{equation}\label{defvolterraMK_Picard}
	\Phi_m(\mu,s,r,t,x,y) = \sum_{k=1}^{\infty} \H^k_m(\mu,s,r,t,x,y).
\end{equation}
Let us recall that the Beta function $\BB$ is defined, for all $x,y >0$ by \[ \BB(x,y) := \int_0^1 (1-t)^{-1+x} t^{-1+y}  \, dt = \frac{\Gamma(x) \Gamma(y)}{\Gamma(x+y)},\] where $\Gamma$ is the Gamma function.\\

Applying Proposition \ref{Prop_H^k} and Theorem \ref{Thmdensityparametrix}, since all the controls are uniform with respect to the measure argument and thus on $m\geq 1$ too, we deduce the following two propositions.

\begin{Prop} \begin{itemize}
		\item There exists $C>0$ such that for all $k \geq 1$, $m \geq 1$, $\mu \in \PP$, $0 \leq s\leq r < t\leq T$ and $x,y \in \R^d$
		
		\begin{equation}\label{H^k_m_bound}
			|\H^k_m(\mu,s,r,t,x,y)| \leq C^k (t-r)^{-\frac{1}{\alpha} + (k-1)\left( 1 - \frac{1}{\alpha} \right)} \prod_{j=1}^{k-1} \BB \left( j\left( 1 - \frac{1}{\alpha} \right),1 - \frac{1}{\alpha} \right) \rho^1(t-r,y-x).
		\end{equation}
		
		\item For $\gamma \in (0, \eta]$ such that $\gamma < \alpha -1$, there exists $C>0$ depending on $\gamma$ such that for all $k\geq 1$, $m\geq 1$, $\mu \in \PP$, $0 \leq s\leq r < t \leq T$ and $x_1,x_2,y \in \R^d$

		\begin{multline}\label{H^k_m_Holder}
			| \H^k_m(\mu,s,r,t,x_1,y) - \H^k_m(\mu,s,r,t,x_2,y) | \leq C^k(t-r)^{-\frac{\gamma +1}{\alpha} + (k-1)\left( 1 - \frac{1}{\alpha} \right)} |x_1-x_2|^\gamma \\ \prod_{j=1}^{k-1} \BB \left( -\frac{\gamma}{\alpha} + j\left( 1 - \frac{1}{\alpha} \right),1 - \frac{1}{\alpha} \right) \left[\rho^{1 }(t-r,y-x_1) + \rho^{1}(t-r,y-x_2) \right].
		\end{multline}

		\item The series \eqref{defvolterraMK_Picard} defining $\Phi_m$ is absolutely convergent and there exists $C>0$ such that for all $m\geq 1$, $\mu \in \PP$, $0 \leq s \leq r<t\leq T$ and $x,y \in \R^d$
		
		\begin{equation}\label{Phi_m_bound}
			|\Phi_m(\mu,s,r,t,x,y)| \leq C (t-r)^{-\frac{1}{\alpha}} \rho^1(t-r,y-x).
		\end{equation}
		
		\item For $\gamma \in (0, \eta]$ such that $\gamma < \alpha -1$, there exists $C>0$ depending on $\gamma$ such that for all $m \geq 1$, $\mu \in \PP$, $0\leq s\leq r < t\leq T$ and $x_1,x_2,y \in \R^d$

		\begin{equation}\label{Phi_m_Holder}
			| \Phi_m(\mu,s,r,t,x_1,y) - \Phi_m(\mu,s,r,t,x_2,y) | \leq C(t-r)^{-\frac{\gamma +1}{\alpha}} |x_1-x_2|^\gamma \left[\rho^{1 }(t-r,y-x_1) + \rho^{1 }(t-r,y-x_2) \right].
		\end{equation}
		
	\end{itemize}
	
\end{Prop}

\begin{Prop}\label{Prop_control_density_Picard_classical}
	
	For any $m\geq 1$, $\mu \in \PP$, $0\leq s<t \leq T$ and $x \in \R^d$, the distribution of $X^{s,x,\mu,(m)}_t$ has a density with respect to the Lebesgue measure denoted by $p_m(\mu,s,t,x,\cdot)$ and given by the absolutely convergent parametrix series \begin{align}\label{representationdensityparametrix_Picard}
		\notag p_m(\mu,s,t,x,y) &= \p(s,t,x,y) + \sum_{k=1}^{\infty} \p \otimes \H^k_m( \mu,s,t,x,y)\\ &= \p(s,t,x,y) + \p \otimes \Phi_m(s,t,x,y).
	\end{align}
	
	For any $m\geq  1$, $\mu \in \PP$, $0\leq s <t \leq T$ and $y \in \R^d$, $p_m(\mu,s,t,\cdot,y)$ is of class $\CC^1$ on $ \R^d$. Moreover, the following properties hold true.\\
	
	\begin{itemize}
		
		\item There exists $C>0$ such that for all $j \in \{0,1\}$, $m\geq 1$, $\mu \in \PP$, $0 \leq s < t \leq T$ and $x,y \in \R^d$
		
		\begin{equation}\label{density_bound_Picard}
			|\pa_x^j p_m(\mu,s,t,x,y) |
			\leq C (t-s)^{-\frac{j}{\alpha}} \rho^{j}(t-s,y-x).
		\end{equation}

		\item For all $j \in \{ 0,1\}$ and $\gamma \in (0, 1]$ with $\gamma \in (0, (2\alpha -2) \wedge (\eta + \alpha -1))$ if $j=1$, there exists $C>0$ such that for all $m\geq 1$, $\mu \in \PP$, $0 \leq s < t \leq T$ and $x_1,x_2,y \in \R^d$
		
		\begin{equation}\label{gradientdensity_holder_Picard}
			|\pa_x^j p_m(\mu,s,t,x_1,y) -  \pa^j_x p_m(\mu,s,t,x_2,y)| \leq C (t-s)^{-\frac{ j+\gamma}{\alpha}} |x_1-x_2|^\gamma \left[\rho^{j}(t-s,y-x_1) +\rho^{j}(t-s,y-x_2)\right].
		\end{equation}

	\end{itemize}
	
\end{Prop}

We state in the following proposition all the properties satisfied by the transition densities $p_m$ which are used to prove Theorem \ref{Thm_density_McKean}. As the proof is rather long and technical, it is postponed to Section \ref{section_proof_prop}.

\begin{Prop}\label{Prop_density_Picard}
	For any $m \geq 1$, $t \in (0,T]$, $y \in \R^d$, the map $(\mu,s,x) \in \PP \times [0,t) \times \R^d \mapsto p_m(\cdot,\cdot,t,\cdot,y)$ belongs to $\CC^1(\PP \times [0,t) \times \R^d)$ and satisfies the following properties. 
	\begin{itemize}
		
		\item There exists $C>0$ such that for all $m\geq 1$, $\mu \in \PP$, $ 0 \leq s < t \leq T$, $x,y,v \in \R^d$ \begin{multline}\label{linear_der_bound_Picard}
			\left\vert \del p_m(\mu,s,t,x,y)(v)\right\vert \leq (t-s)^{1 - \frac{1}{\alpha}}\rho^0(t-s,y-x)\\  \left( \sum_{k=1}^m C^k(t-s)^{(k-1)\left(1 - \frac{1}{\alpha}\right)} \prod_{j=1}^{k-1} \BB\left(1 + j \left(1 - \frac{1}{\alpha}\right), 1 - \frac{1}{\alpha}\right)\right).
		\end{multline}
		
		\item There exists $C>0$ such that for all $m\geq 1$, $\mu \in \PP$, $ 0 \leq s < t \leq T$, $x,y,v \in \R^d$ \begin{multline}\label{gradient_linear_der_bound_Picard}
			\left\vert \pa_v\del p_m(\mu,s,t,x,y)(v)\right\vert \leq  (t-s)^{\frac{\eta -1}{\alpha} + 1 - \frac{1}{\alpha}}\rho^0(t-s,y-x)\\ \left( \sum_{k=1}^m C^k(t-s)^{(k-1)\left(1 + \frac{\eta -1}{\alpha}\right)} \prod_{j=1}^{k-1} \BB\left(1 + \frac{\eta -1}{\alpha} + j \left(1 + \frac{\eta -1}{\alpha}\right), 1 - \frac{1}{\alpha}\right)\right).
		\end{multline}
	
		\item For all $\tilde{\eta} \in (0,\eta \wedge(\alpha-1))$ there exists a constant $C>0$ such that for all $m\geq1$, $0\leq s < t \leq T$, $\mu \in \PP$, $x,y \in \R^d$ 
	
	\begin{multline}\label{density_time_der_bound_Picard}
		|\pa_s p_m(\mu,s,t,x,y)| \leq (t-s)^{-1} \rho^{-\tilde{\eta}}(t-s,y-x) \sum_{k=1}^m C^k (t-s)^{(k-1)\left(1 + \frac{\eta-1}{\alpha}\right)} \\ \prod_{j=1}^{k-1} \BB\left(\frac{\eta}{\alpha} + (j-1)\left(1 + \frac{\eta -1}{\alpha}\right),1 - \frac{1}{\alpha}\right).
	\end{multline}
		
		\item For all $\gamma \in (0,1]\cap(0,(2\alpha -2)\wedge (\eta + \alpha -1))$, there exists $C>0$ such that for all $ m\geq 1$, $\mu \in \PP$, $ 0 \leq s < t \leq T$, $x,y,v_1,v_2 \in \R^d$ \begin{multline}\label{gradient_linear_der_holder_v_Picard}
		\left\vert \pa_v\del p_m(\mu,s,t,x,y)(v_1) - \pa_v \del p_m(\mu,s,t,x,y)(v_2)\right\vert \leq  (t-s)^{\frac{\eta -1 - \gamma}{\alpha} + 1 - \frac{1}{\alpha}}|v_1-v_2|^\gamma\rho^0(t-s,y-x)\\\left( \sum_{k=1}^m C^k(t-s)^{(k-1)\left(1 + \frac{\eta -1}{\alpha}\right)} \prod_{j=1}^{k-1} \BB\left(1 + \frac{\eta -1 - \gamma}{\alpha} + j \left(1 + \frac{\eta -1}{\alpha}\right), 1 - \frac{1}{\alpha}\right)\right).
	\end{multline}
	
	\item There exists $C>0$ such that for all $\gamma \in (0,1]$, $m \geq 1$, $\mu \in \PP$, $ 0 \leq s < t \leq T$, $x,y,v_1,v_2 \in \R^d$ \begin{equation}\label{linear_der_holder_v_Picard}
		\left\vert \del p_m(\mu,s,t,x,y)(v_1) - \del p_m(\mu,s,t,x,y)(v_2)\right\vert \leq C (t-s)^{ 1 - \frac{1+\gamma }{\alpha}}|v_1-v_2|^\gamma\rho^0(t-s,y-x).
	\end{equation}
	
	\item There exists $C>0$ such that for all $\gamma \in (0,1]$, $m \geq 1$, $\mu \in \PP$, $ 0 \leq s < t \leq T$, $x_1,x_2,y,v\in \R^d$ \begin{multline}\label{linear_der_holder_x_Picard}
		\left\vert \del p_m(\mu,s,t,x_1,y)(v) - \del p_m(\mu,s,t,x_2,y)(v)\right\vert \leq C (t-s)^{ 1 - \frac{1+\gamma}{\alpha}}|x_1-x_2|^\gamma\\ \left[\rho^0(t-s,y-x_1) + \rho^0(t-s,y-x_2)\right].
	\end{multline}
	
	\item For all $\gamma \in (0,1]\cap(0,\eta + \alpha -1)$, there exists $C>0$ such that for all $ m\geq 1$, $\mu \in \PP$, $ 0 \leq s < t \leq T$, $x_1,x_2,y,v\in \R^d$ \begin{multline}\label{gradient_linear_der_holder_x_Picard}
		\left\vert \pa_v\del p_m(\mu,s,t,x_1,y)(v) - \pa_v \del p_m(\mu,s,t,x_2,y)(v)\right\vert \leq  C(t-s)^{\frac{\eta -1 - \gamma}{\alpha} + 1 - \frac{1}{\alpha}}|x_1-x_2|^\gamma \\ \left[\rho^0(t-s,y-x_1) + \rho^0(t-s,y-x_2)\right].
	\end{multline}
	
	\item There exists $C>0$ such that for all $\gamma \in (0,1]$, $m \geq 1$, $\mu_1, \mu_2 \in \PP$, $ 0 \leq s < t \leq T$, $x,y \in \R^d$ \begin{equation}\label{density_holder_measure_Picard}
		\left\vert  p_m(\mu_1,s,t,x,y)-  p_m(\mu_2,s,t,x,y)\right\vert \leq C (t-s)^{ 1 - \frac{1+\gamma}{\alpha}}W_1^\gamma(\mu_1,\mu_2)\rho^0(t-s,y-x).
	\end{equation}
	
	\item For all $j \in \{0,1\}$ and $\gamma \in (0,1]$ with $\gamma \in (0, \alpha -1 + \eta)$ if $j=1$, there exists $C>0$ such that for all $m\geq 1$, $0\leq s <t \leq T$, $\mu_1,\mu_2 \in \PP$, $x,y \in \R^d$ 
	
	\begin{equation}\label{gradient_density_holder_measure_Picard}
		\left\vert \pa_x^j p_{m}(\mu_1,s,t,x,y) -  \pa_x^j p_{m}(\mu_2,s,t,x,y) \right\vert \leq C(t-s)^{1-\frac{1+\gamma +j}{\alpha}} W_1^\gamma (\mu_1,\mu_2) \rho^j(t-s,y-x).
	\end{equation}
	
	\item For all $\gamma \in (0,1]$, there exists $C>0$ such that for all $ m \geq 1$, $0 \leq s < t \leq T$, $\mu_1,\mu_2 \in \PP$, $x,y,v \in \R^d$
	
	\begin{multline}\label{linear_derivative_holder_measure_Picard}
		\left\vert \del p_m (\mu_1,s,t,x,y)(v) -  \del p_m (\mu_2,s,t,x,y)(v)\right\vert \leq (t-s)^{-\frac{\gamma}{\alpha} + 1 - \frac{1}{\alpha}}W_1^\gamma(\mu_1,\mu_2)\rho^0(t-s,y-x) \\ \sum_{k=1}^m C^k (t-s)^{(k-1)\left(1 - \frac{1}{\alpha}\right)} \prod_{j=1}^{k-1} \BB\left(1 - \frac{\gamma}{\alpha} + j\left(1 - \frac{1}{\alpha}\right), 1 - \frac{1}{\alpha}\right).
	\end{multline}
	
	\item For all $\gamma \in (0,1]$, there exists $C>0$ such that for all $ m \geq 1$, $0 \leq s < t \leq T$, $\mu_1,\mu_2 \in \PP$, $x,y,v \in \R^d$
	
	\begin{multline}\label{gradient_linear_derivative_holder_measure_Picard}
		\left\vert \pa_v\del p_m (\mu_1,s,t,x,y)(v) -  \pa_v\del p_m (\mu_2,s,t,x,y)(v)\right\vert \leq (t-s)^{\frac{\eta -\gamma -1}{\alpha} + 1 - \frac{1}{\alpha}}W_1^\gamma(\mu_1,\mu_2)\rho^0(t-s,y-x) \\ \sum_{k=1}^m C^k (t-s)^{(k-1)\left(1 + \frac{\eta-1}{\alpha}\right)} \prod_{j=1}^{k-1} \BB\left(1 + \frac{\eta -\gamma}{\alpha} + j\left(1 + \frac{\eta -1}{\alpha}\right), 1 - \frac{1}{\alpha}\right).
	\end{multline}

	\item For all $j \in \{0,1\}$, $ \gamma \in \left(0,1 - \frac{j}{\alpha}\right)$, there exists a constant $C>0$ such that for all $m\geq1$, $t \in (0,T],$ $s_1,s_2 \in [0,t)$, $\mu \in \PP$, $x,y \in \R^d$ 
	
	\begin{multline}\label{gradient_density_time_Holder_Picard}
		|\pa_x^j p_m(\mu,s_1,t,x,y) - \pa_x^jp_m(\mu,s_2,t,x,y)| \\ \leq C \left[\frac{|s_1-s_2|^\gamma}{(t - s_1)^{\gamma + \frac{j}{\alpha}}}\rho^{j}(t-s_1,y-x) + \frac{|s_1-s_2|^\gamma}{(t-s_2)^{\gamma+ \frac{j}{\alpha}}}\rho^{j}(t-s_2,y-x)\right].
	\end{multline}
	
	\item For all $\gamma \in \left(0,1\right)$, there exists a constant $C>0$ such that for all $m\geq 1$, $t \in (0,T]$, $s_1,s_2 \in[0,t)$, $\mu \in \PP$, $x,y,v \in \R^d$

	\begin{align}\label{lin_der_density_holder_time_Picard}
		\notag&\left\vert \del p_m(\mu,s_1,t,x,y)(v) - \del p_m(\mu,s_2,t,x,y)(v)\right\vert\\ &\leq \left[\frac{|s_1-s_2|^\gamma}{(t-s_1\vee s_2)^{\gamma + \frac{1}{\alpha}-1}} \rho^{0}(t-s_1\vee s_2,y-x) + \frac{|s_1-s_2|^\gamma }{(t-s_1\wedge s_2)^{\gamma +\frac{1}{\alpha}-1}} \rho^{0}(t-s_1\wedge s_2,y-x) \right] \\\notag&\hspace{3cm}\sum_{k=1}^{m} C^k(t-s_1\vee s_2)^{(k-1)\left(1 - \frac{1}{\alpha}\right)} \prod_{j=1}^{k-1}\BB\left(\left[2-\gamma - \frac{1}{\alpha} \right] \wedge 1 +(j-1)\left(1- \frac{1}{\alpha}\right), 1 - \frac{1}{\alpha}\right).
	\end{align}

	\item For all $\gamma \in \left(0,1+ \frac{\eta -1}{\alpha}\right)$, there exists a constant $C>0$ such that for all $m\geq 1$, $t \in (0,T]$, $s_1,s_2 \in[0,t)$, $\mu \in \PP$, $x,y,v \in \R^d$
	
	\begin{align}\label{gradient_lin_der_density_holder_time_Picard}
		\notag&\left\vert \pa_v\del p_m(\mu,s_1,t,x,y)(v) - \pa_v\del p_m(\mu,s_2,t,x,y)(v)\right\vert\\ &\leq \left[\frac{|s_1-s_2|^\gamma}{(t-s_1\vee s_2)^{\gamma + \frac{1}{\alpha}-1 + \frac{1-\eta}{\alpha}}} \rho^{0}(t-s_1\vee s_2,y-x) + \frac{|s_1-s_2|^\gamma }{(t-s_1\wedge s_2)^{\gamma +\frac{1}{\alpha}-1+ \frac{1-\eta}{\alpha}}} \rho^{0}(t-s_1\wedge s_2,y-x) \right] \\\notag&\hspace{2cm}\sum_{k=1}^{m} C^k(t-s_1\vee s_2)^{(k-1)\left(1 + \frac{\eta -1}{\alpha}\right)} \prod_{j=1}^{k-1}\BB\left(\left[2\left(1 + \frac{\eta-1}{\alpha}\right)-\gamma\right] \wedge 1 +(j-1)\left(1+ \frac{\eta -1}{\alpha}\right), 1 - \frac{1}{\alpha}\right).
	\end{align}
			
	\end{itemize}
\end{Prop}

\begin{Rq}
	Note that by the asymptotic behavior of the Beta function, we get that all the series appearing in the right-hand side members of the preceding inequalities are convergent. This provides controls which are uniform with respect to $m$.
\end{Rq}

\textbf{Step 2: Passage to the limit in the previous estimates.} We are going to take the limit $ m \to + \infty$ in all the estimates proved in Proposition \ref{Prop_density_Picard} to deduce that the transition density $p(\mu,s,t,x,y)$ of the McKean-Vlasov SDE \eqref{McKVSDE} satisfies the regularity properties and estimates of Theorem \ref{Thm_density_McKean}. This will be done along a converging subsequence given by the Arzelà-Ascoli theorem. Notice that all the partial sums of the series appearing in the all the upper-bounds of Proposition \ref{Prop_density_Picard} have a limit when $m \to + \infty$ using the asymptotic behavior of the Beta function.\\

First of all, note that Theorem \ref{Thm_wellposedness} yields 

\begin{equation*}
\sup_{r \in [s,t]} d_{TV} ([X^{s,\mu,(m)}_r],[X^{s,\mu}_r])  \underset{m \rightarrow + \infty}{\longrightarrow} 0.
\end{equation*}
It follows that for all $k \geq 1$, $\mu \in \PP$, $0 \leq s \leq r<t \leq T$, $x,y \in \R^d$ \[ \H_{m}^k(\mu,s,r,t,x,y) \underset{m \rightarrow + \infty}{\longrightarrow} \H^k(\mu,s,r,t,x,y), \]
where $\H(\mu,s,t,x,y)$ was defined in \eqref{defproxyMK}. We can thus let $m$ tend to infinity in the parametrix series \eqref{representationdensityparametrix_Picard} which yields the following pointwise convergence 

\begin{equation}\label{convergence_approximated_densities}
	p_m(\mu,s,t,x,y) \underset{m\rightarrow + \infty}{\longrightarrow} p(\mu,s,t,x,y)
\end{equation}
thanks to the parametrix expansion \eqref{param_rep_density_thm} of $p$. Let us fix $(t,y) \in (0,T] \times \R^d$ and $\KK$ a compact subset of $ \PP \times [0,t) \times \R^d$. Using \eqref{gradientdensity_holder_Picard}, \eqref{gradient_density_holder_measure_Picard} and \eqref{gradient_density_time_Holder_Picard}, we deduce that the sequence of functions $\left(p_m(\cdot,\cdot,t,\cdot,y)\right)_m \in \CC^0(\KK)^\N$ is uniformly equi-continuous on $\KK$. It is also uniformly bounded by \eqref{density_bound_Picard}. The Arzelà-Ascoli theorem ensures that we can extract a subsequence of $\left(p_m(\cdot,\cdot,t,\cdot,y)\right)_m$ which converges uniformly on $\KK$, necessarily towards $p(\cdot,\cdot,t,\cdot,y)$ by \eqref{convergence_approximated_densities}. This yields the continuity of $p(\cdot,\cdot,t,\cdot,y)$ on $\KK$ and thus, since $\KK$ is arbitrary, on  $\PP\times [0,t)\times \R^d$. Moreover, passing to the limit in \eqref{density_bound_Picard}, \eqref{gradientdensity_holder_Picard}, \eqref{gradient_density_holder_measure_Picard} and \eqref{gradient_density_time_Holder_Picard} for $j=1$, along the converging subsequence of $\left(p_m(\cdot,\cdot,t,\cdot,y)\right)_m$ previously obtained, we get that \eqref{density_bound_thm}, \eqref{gradientdensity_holder_thm}, \eqref{gradient_density_holder_measure_thm} and \eqref{gradient_density_time_Holder_thm} hold true for $j=0$. We now prove that $p(\mu,s,t,\cdot,y)$ is of class $\CC^1$ on $\R^d$ for any $\mu \in \PP$ and $s \in [0,t)$. To do this, we fix $R>0$. By \eqref{density_bound_Picard} and \eqref{gradientdensity_holder_Picard}, we can apply the Arzelà-Ascoli theorem to the sequence $\left(\pa_xp_m(\mu,s,t,\cdot,y)\right)_m \in \CC^0(B_R)^\N$, where $B_R$ denotes the open ball of $\R^d$ with radius $R$. Since $R$ is arbitrary, we can construct, using a diagonal extraction procedure, a continuous function on $\R^d$ which is the limit, uniformly on each compact subset of $\R^d$, of a subsequence of $\left(\pa_xp_m(\mu,s,t,\cdot,y)\right)_m$. This proves that $p_m(\mu,s,t,\cdot,y)$ is of class $\CC^1$. By \eqref{density_bound_Picard}, \eqref{gradientdensity_holder_Picard}, \eqref{gradient_density_holder_measure_Picard} and \eqref{gradient_density_time_Holder_Picard}, the continuity of $\pa_x p(\cdot,\cdot,t,\cdot,y)$ on $\PP\times [0,t) \times \R^d$ is again a consequence of the Arzelà-Ascoli theorem applied to $\left(\pa_xp_m(\cdot,\cdot,t,\cdot,y)\right)_m \in (\CC^0(\KK))^{\N}$, where $\KK$ is an arbitrary compact subset of $\PP \times [0,t) \times \R^d$. Taking the limit $m \to +\infty$ along a converging subsequence in \eqref{density_bound_Picard}, \eqref{gradientdensity_holder_Picard}, \eqref{gradient_density_holder_measure_Picard} and \eqref{gradient_density_time_Holder_Picard}, we deduce that \eqref{density_bound_thm}, \eqref{gradientdensity_holder_thm}, \eqref{gradient_density_holder_measure_thm} and \eqref{gradient_density_time_Holder_thm} are satisfied.\\

Let us now focus on the existence of the linear derivative of $p$. We fix $\tau< t$ and $\KK$ a closed and bounded subset of $\PP \times [0,\tau] \times (\R^d)^2$. Note that since $\beta >1$, $\KK$ is relatively compact in $\mathcal{P}_1(\R^d) \times [0,t) \times (\R^d)^2$ for the metric $d$ defined for all $\mu_1,\mu_2 \in \PP$, $s_1,s_2  \in [0,\tau]$, $x_1,x_2,v_1,v_2 \in \R^d$ by \[ d((\mu_1,s_1,x_1,v_1),(\mu_2,s_2,x_2,v_2)):= W_1(\mu_1,\mu_2) + |s_1-s_2| + |x_1-x_2| + |v_1-v_2|.\] Using \eqref{linear_der_holder_v_Picard}, \eqref{linear_der_holder_x_Picard}, \eqref{linear_derivative_holder_measure_Picard}, \eqref{lin_der_density_holder_time_Picard}, we deduce that the sequence of functions $\left(\del p_m(\cdot,\cdot,t,\cdot,y)(\cdot)\right)_m$ is uniformly equi-continuous on $\KK$ with respect to the metric $d$. Moreover \eqref{linear_der_bound_Picard} ensures that \begin{equation}\label{bound_lin_der_ascoli} \sup_{m\geq 1}\sup_{(\mu,s,x,v) \in \PP\times [0,\tau] \times (\R^d)^2} \left\vert \del p_m(\mu,s,t,x,y)(v) \right\vert < + \infty.\end{equation}
 Then, we apply the Arzelà-Ascoli theorem, which gives the existence of a subsequence of $\left(\del p_m(\cdot,\cdot,t,\cdot,y)(\cdot)\right)_m $ which converges uniformly on $\KK$ with respect to $d$. Since this is true for all $\tau <t$ and for every bounded subset $\KK$ of $\PP \times [0,\tau] \times (\R^d)^2$, we can use a diagonal extraction procedure. This yields the existence of function $g$ continuous with respect to $d$ on $\PP \times [0,t)\times (\R^d)^2$ such that, up to an extraction, $\left(\del p_m(\cdot,\cdot,t,\cdot,y)(\cdot)\right)_m $ converges towards $g$ uniformly on each compact subset of $\PP \times [0,t)\times (\R^d)^2$. Note that $g$ is also continuous with respect to the usual metric on $\PP \times [0,t) \times (\R^d)^2$ and that \eqref{bound_lin_der_ascoli} implies that for each $\tau \in [0,t)$, we have \begin{equation*}
 \sup_{(\mu,s,x,v) \in \PP\times [0,\tau] \times (\R^d)^2} \left\vert g(\mu,s,x,v) \right\vert < + \infty.
 \end{equation*}

 We now prove that $p(\cdot,s,t,x,y)$ admits a linear derivative given, for all $\mu \in \PP$, $s\in [0,t)$, $x,v\in \R^d$ by \begin{equation}\label{expresion_lin_der_ascoli}
	\del p(\mu,s,t,x,y)(v) = g(\mu,s,x,v),\end{equation}
 which is continuous on $\PP \times [0,t)\times (\R^d)^2$. For all $ \mu,\nu \in \PP$, one has \begin{equation*}
p_m(\mu,s,t,x,y) - p_m(\nu,s,t,x,y) = \int_0^1 \int_{\R^d} \del p_m(\lambda \mu + (1-\lambda)\nu,s,t,x,y)(v) \, d(\mu-\nu)(v) \, d\lambda.
\end{equation*}
We take the limit $m \to + \infty$ along the subsequence of $\left(\del p_m(\cdot,\cdot,t,\cdot,y)(\cdot)\right)_m $ converging towards $g$ that we have obtained above. By the dominated convergence theorem justified by \eqref{bound_lin_der_ascoli} and since $(p_m)_m$ converges pointwise towards $p$, we obtain that \begin{equation*}
	p(\mu,s,t,x,y) - p(\nu,s,t,x,y) = \int_0^1 \int_{\R^d} g(\lambda \mu + (1-\lambda)\nu,s,x,v) \, d(\mu-\nu)(v) \, d\lambda.
\end{equation*}
This proves \eqref{expresion_lin_der_ascoli}. Moreover, taking the limit $m \to + \infty$ in \eqref{linear_der_bound_Picard}, \eqref{linear_der_holder_v_Picard}, \eqref{linear_der_holder_x_Picard}, \eqref{linear_derivative_holder_measure_Picard} and \eqref{lin_der_density_holder_time_Picard} along the converging subsequence yields \eqref{linear_der_bound_thm}, \eqref{linear_der_holder_v_thm}, \eqref{linear_der_holder_x_thm}, \eqref{linear_derivative_holder_measure_thm} and \eqref{lin_der_density_holder_time_thm}. Using again the Arzelà-Ascoli theorem, we prove that for all $\mu \in \PP$, $s \in [0,t)$, $x,y \in \R^d$, the map $\del p(\mu,s,t,x,y)(\cdot)$ is of class $\CC^1$ on $\R^d$, that $\pa_v \del p(\cdot,\cdot,t,\cdot,y)(\cdot)$ is continuous on $\PP \times [0,t) \times (\R^d)^2$ and that it satisfies \eqref{gradient_linear_der_bound_thm}, \eqref{gradient_linear_der_holder_v_thm}, \eqref{gradient_linear_der_holder_x_thm}, \eqref{gradient_linear_derivative_holder_measure_thm} and \eqref{gradient_lin_der_density_holder_time_thm}.

\begin{Rq}\label{Rq_equi_cont_y}
	We have made the proof of the extraction of converging subsequences for a fixed $y \in \R^d$. However, following exactly the same lines as for the estimates of Hölder continuity with respect to $x$ \eqref{gradientdensity_holder_Picard}, \eqref{linear_der_holder_x_Picard} and \eqref{gradient_linear_der_holder_x_Picard}, we can prove similar estimates with respect to $y$. This ensures that the converging subsequences can also be assumed to converge uniformly on each compact subset of $\R^d$ with respect to $y$.
\end{Rq}

Let us prove \eqref{centering_prop_lin_der_density}. We fix $\eps >0$. Thanks to \eqref{linear_der_bound_Picard} and \eqref{gradient_linear_der_bound_Picard}, we can find a compact subset $\KK$ of $\R^d$ such that \begin{equation*}
	\int_{\KK^c} \sup_{m\geq 1} \left\{\left\vert \del p_m(\mu,s,t,x,y)(v)\right\vert +  \left\vert \del p_m(\mu,s,t,x,y)(v)\right\vert\right\} \, dy \leq \eps
\end{equation*}
By Remark \eqref{Rq_equi_cont_y}, up to extracting a subsequence, we can assume that the functions $\left(\del p_m(\mu,s,t,x,\cdot)(v)\right)_m$ and  $\left(\pa_v\del p_m(\mu,s,t,x,\cdot)(v)\right)_m$ converge uniformly on $\KK$ towards $\del p(\mu,s,t,x,\cdot)(v)$ and $\pa_v\del p(\mu,s,t,x,\cdot)(v)$. Noticing that \eqref{centering_prop_lin_der_density} is true for $p_m$ by \eqref{expression_linear_der}, \eqref{expression_linear_der_H}, \eqref{expression_gradient_linear_der} and \eqref{expression_gradient_linear_der_H}. We can thus write 

\begin{align*}
	&\left\vert \int_{\R^d} \del p(\mu,s,t,x,y)(v) \, dy \right\vert \\ &= \left\vert \int_{\R^d} \del p(\mu,s,t,x,y)(v) - \del p_m(\mu,s,t,x,y)(v)\, dy \right\vert \\ &\leq \int_{\KK} \left\vert\del p(\mu,s,t,x,y)(v) - \del p_m(\mu,s,t,x,y)(v)\right\vert\, dy + \eps.
\end{align*}
We conclude by letting $m$ tend to infinity and with a similar reasoning for $\pa_v \del p$. This ends the proof of Theorem \ref{Thm_density_McKean}. 

\begin{proof}[Proof of Theorem \ref{Thm_density_McKean_light}]
	
	 We first need to prove that for all $\mu \in \PP$, $t \in (0,T]$, $x,y\in \R^d$, the map $s \in [0,t) \mapsto p(\mu,\cdot,t,x,y)$ is of class $\CC^1$ on $[0,t)$, that $\pa_s p(\cdot,\cdot,t,\cdot,y)$ is continuous on $\PP \times [0,t) \times \R^d$ and satisfies \eqref{EDPKolmogorovfundamental}. Let us fix $h \in [0,s]$. By the well-posedness of the nonlinear martingale problem proved in Theorem \ref{Thm_wellposedness}, we deduce that the transition density satisfies the following Markov property \begin{equation}
p(\mu,s-h,t,x,y) = \E p([X^{s-h,\mu}_s],s,t,X^{s-h,x,\mu}_s,y).
\end{equation}
By \eqref{density_bound_thm}, \eqref{linear_der_bound_thm}, \eqref{gradient_linear_der_bound_thm}, \eqref{gradientdensity_holder_thm} and \eqref{gradient_linear_der_holder_v_thm} with $\gamma > \alpha -1$, we can apply Itô's formula of Proposition \ref{ito_formula} for the function $(\mu,x) \in \PP \times \R^d \mapsto p(\mu,s,t,x,y)$. Taking the expectation in Itô's formula, we obtain that 

\begin{equation*}
	\E p([X^{s-h,\mu}_s],s,t,X^{s-h,x,\mu}_s,y) = p(\mu,s,t,x,y) + \int_{s-h}^s \LL_r p(\cdot,s,t,\cdot,y)([X^{s-h,\mu}_r],X^{s-h,x,\mu}_r) \, dr,
\end{equation*}
where $\LL_r$ was defined in \eqref{def_operator_McKV}. We thus have \begin{equation*}
	\frac1h (p(\mu,s-h,t,x,y) - p(\mu,s,t,x,y)) = \frac1h \int_{s-h}^s \E \LL_rp(\cdot,s,t,\cdot,y)([X^{s-h,\mu}_r],X^{s-h,x,\mu}_r) \, dr.
\end{equation*}

Using the continuity and the boundedness of $b$ as well as the Hölder continuity and the bounds on $p(\cdot,\cdot,t,\cdot,y)$, $\pa_x p(\cdot,\cdot,t,\cdot,y)$, $\del p(\cdot,\cdot,t,\cdot,y)(\cdot)$ and $\pa_v \del p(\cdot,\cdot,t,\cdot,y)(\cdot)$ proved above, we find that 

\begin{equation*}
	\frac1h (p(\mu,s-h,t,x,y) - p(\mu,s,t,x,y)) \underset{h \to 0^+}{\longrightarrow}   \LL_s p(\cdot,s,t,\cdot,y)(\mu,x).
\end{equation*}

The map $s \in [0,t) \mapsto p(\mu,s,t,x,y)$ is thus left-differentiable on $[0,t)$. It also follows that the map $(\mu,s,x) \in \PP \times [0,t) \times \R^d \mapsto \LL_sp(\cdot,s,t,\cdot,y)(\mu,x)$ is continuous. This proves that $p(\mu,\cdot,t,x,y)$ is $\CC^1$ on $[0,t)$ and that it satisfies for all $\mu \in \PP$, $s \in [0,t)$, $x,y \in \R^d$ \begin{equation*}
	\pa_s p(\mu,s,t,x,y) = - \LL_sp(\cdot,s,t,\cdot,y)(\mu,x).
\end{equation*}

Let us now fix $ f : \R^d \rightarrow \R$ a bounded and uniformly continuous function. We fix $\eps >0$. There exists $\delta >0$ such that for all $x,y \in \R^d$ with $|x-y|\leq \delta$, we have $ |f(x) - f(y)| \leq \eps$. Using \eqref{density_bound_thm}, we obtain that \begin{align*}
	\sup_{x \in \R^d} \left\vert \int_{\R^d} f(y) p(\mu,s,t,x,y) \, dy - f(x) \right\vert& = 	\sup_{x \in \R^d} \left\vert \int_{\R^d} (f(y) - f(x)) p(\mu,s,t,x,y) \, dy \right\vert \\ &\leq \eps + C\|f\|_{\infty}\int_{|y| >\delta} \rho^0 (t-s,y) \, dy  \\ &\leq  \eps + C\|f\|_{\infty}\int_{|y| >\delta} (t-s)^{-\frac{d}{\alpha}}( 1 + (t-s)^{-\frac{1}{\alpha}}|y|)^{-d-\alpha} \, dy \\ &= \eps + C\|f\|_{\infty}\int_{|z| >(t-s)^{-\frac{1}{\alpha}}\delta} ( 1 + |z|)^{-d-\alpha} \, dz.
\end{align*}
We conclude taking the $\limsup$ when $s \rightarrow t$ in the preceding inequality that $p(\mu,s,t,x,\cdot) \underset{s\rightarrow t^-}{\longrightarrow} \delta_x$ in the weak sense.

The estimates \eqref{density_bound_thm_light} and \eqref{gradient_linear_der_bound_thm_light} have been proved in Theorem \ref{Thm_density_McKean} and the estimate \eqref{density_bound_laplacian_thm_light} is proved in Theorem \ref{Thmdensityparametrix}. It thus remains to prove \eqref{density_op_bound_thm}. Using  \eqref{linear_der_bound_thm}, \eqref{gradient_linear_der_bound_thm}, \eqref{gradient_linear_der_holder_v_thm} and the same reasoning as used in the proof of \eqref{timederivativedensity_bound} in Theorem \ref{Thmdensityparametrix}, we find that there exists a positive constant $C$ such that for all $\mu \in \PP$, $0 \leq s < t \leq T$, $x,y,v \in \R^d$ \begin{align*}
 & \left\vert\int_{\R^d} \left[ \del p(\mu,s,t,x,y)(v+z) - \del p(\mu,s,t,x,y)(v) - z\cdot \pa_v \del p(\mu,s,t,x,y)(v) \right] \, \frac{dz}{|z|^{d + \alpha}} \right\vert  \\ &\hspace{10cm} \leq C (t-s)^{-1+1 - \frac{1}{\alpha}}\rho^0(t-s,y-x).
\end{align*}
This concludes the proof of \eqref{density_op_bound_thm} and thus ends the proof of Theorem \ref{Thm_density_McKean_light}.

\end{proof}

\section{Backward Kolmogorov PDE on the space of measures}\label{section_PDE}
We prove Theorem \ref{Thm_EDP} in this section.\\

	\textbf{Step 1: Continuity of $\bm{U}$ on $\bm{[0,T]\times \PP}$.} Reasoning exactly as in the proof of \cite[Proposition $6.1$]{deraynal2021wellposedness}, the continuity of $U$ on $[0,T) \times \PP$ follows from the continuity of the map $(\mu,s,x) \in \PP \times [0,T) \times \R^d \mapsto p(\mu,s,t,x,y)$ and \eqref{density_bound_thm}. Let us prove it on $[0,T]\times \PP$. Let $(t_n)_n \in [0,T)^\N$, $(\mu_n)_n \in \PP^\N$ and $\mu \in \PP$ such that $|t_n-T| \underset{n \to + \infty}{\longrightarrow} 0$ and $W_{\beta}(\mu_n,\mu) \underset{n \to + \infty}{\longrightarrow} 0$. Since $\phi$ is continuous on $\PP$, it is enough to prove that $W_{\beta}([X^{t_n,\mu_n}_T],\mu) \underset{n \to + \infty}{\longrightarrow} 0$. We start with the weak convergence. Let us fix $f: \R^d \to \R$ a bounded and uniformly continuous function. 	We write 
	
	\begin{align*}
		\int_{\R^{2d}} f(y) p(\mu_n,t_n,T,x,y) \, dy \, d\mu_n(x) - \int_{\R^d}f(x) d\mu(x) & = \int_{\R^{2d}} (f(y) - f(x) ) p(\mu_n,t_n,T,x,y) \, dy \, d\mu_n(x) \\ &\quad + \int_{\R^d} f(x)d(\mu_n - \mu)(x) \\ &=: I_1 + I_2.
	\end{align*}

Since $W_{\beta}(\mu_n,\mu) \underset{n \to + \infty}{\longrightarrow} 0$, it is clear, by definition, that $I_2 \underset{n \to + \infty}{\longrightarrow} 0$. Let us now deal with $I_1$. We fix $\eps >0$. There exists $\delta >0$ such that for all $x,y \in \R^d$ with $|x-y|\leq \delta$, we have $ |f(x) - f(y)| \leq \eps$. Using \eqref{density_bound_thm}, we obtain that \begin{align*}
|I_1| &\leq \int_{\R^d} \left\vert \int_{\R^d} (f(y) -f(x)) p(\mu_n,t_n,T,x,y) \, dy  \right\vert\, d\mu_n(x) \\ &\leq \eps + C\|f\|_{\infty}\int_{|y| >\delta} \rho^0 (T-t_n,y) \, dy  \\ &\leq  \eps + C\|f\|_{\infty}\int_{|y| >\delta} (T-t_n)^{-\frac{d}{\alpha}}( 1 + (T-t_n)^{-\frac{1}{\alpha}}|y|)^{-d-\alpha} \, dy \\ &= \eps + C\|f\|_{\infty}\int_{|z| >(T-t_n)^{-\frac{1}{\alpha}}\delta} ( 1 + |z|)^{-d-\alpha} \, dz.
\end{align*}
We conclude that $I_1 \underset{n \to + \infty}{\longrightarrow} 0.$ It remains to show that \[ \int_{\R^{2d}} |y|^\beta p(\mu_n,t_n,T,x,y)  \, d\mu_n(x)\, dy \underset{n \to + \infty}{\longrightarrow} \int_{\R^d} |x|^\beta d\mu(x).\]
	To see this, we write 
	
	\begin{align*}
		\int_{\R^{2d}} |y|^\beta p(\mu_n,t_n,T,x,y) \, dy \, d\mu_n(x) - \int_{\R^d} |x|^\beta d\mu(x) & = \int_{\R^{2d}} (|y|^\beta -|x|^\beta ) p(\mu_n,t_n,T,x,y) \, dy \, d\mu_n(x) \\ &\quad + \int_{\R^d} |x|^\beta d(\mu_n - \mu)(x) \\ &=: J_1 + J_2.
 	\end{align*}
 Since $W_{\beta}(\mu_n,\mu) \underset{n \to + \infty}{\longrightarrow} 0$, we deduce that $J_2\underset{n \to + \infty}{\longrightarrow} 0$. For $J_1$, by the mean-value theorem, there exists a positive constant $C$ such that for all $x,y \in \R^d$ \[ ||y|^\beta - |x|^\beta| \leq C |y-x| (|y|^{\beta -1} + |x|^{\beta -1}).\] We obtain by \eqref{density_bound_thm} and the space-time inequality \eqref{scalingdensityref} that 
 
 \begin{align*}
 	|J_1| &  \leq C (T-t_n)^{\frac{1}{\alpha}} \int_{\R^{2d}}  (|y|^{\beta -1} + |x|^{\beta -1}) \rho^{-1}(T-t_n,y-x) \, dy \, d\mu_n(x)\\ &  \leq C (T-t_n)^{\frac{1}{\alpha}} \int_{\R^{2d}}  \rho^{-\beta}(T-t_n,y-x) + |x|^{\beta -1} \rho^{-1}(T-t_n,y-x) \, dy \, d\mu_n(x)\\&\leq C (T-t_n)^{\frac{1}{\alpha}} \left(1 + \int_{\R^d} |x|^{\beta-1} \, d\mu_n(x)\right) \\ &\leq  C (T-t_n)^{\frac{1}{\alpha}},
 \end{align*} 
since $\sup_{n \in \N} \int_{\R^d}|x|^{\beta -1} \, d\mu_n(x) < + \infty.$ We have thus proved that $W_{\beta}([X^{t_n,\mu_n}_T],\mu) \underset{n \to + \infty}{\longrightarrow} 0$, which concludes the first step.\\ 

\textbf{Step 2: Estimate \eqref{gradient_lin_der_bound_sol_PDE} and \eqref{gradient_lin_der_holder_sol_PDE}.} By Proposition \ref{Prop_lin_der_along_flow_density} and Proposition \ref{Prop_density_Picard}, we know that the map $U$ belongs to $\CC^1([0,T)\times \PP)$. Moreover, for any $t \in [0,T)$, $\mu \in \PP$ and $v \in \R^d$, we have \begin{align}\label{expression_lin_der_sol_EDP}
	\del U(t,\mu)(v) &= \int_{\R^d} \del \phi([X^{t,\mu}_T])(y) p(\mu,t,T,v,y) \, dy \\ \notag&\quad + \int_{\R^{2d}} \left( \del \phi([X^{t,\mu}_T])(y) - \del \phi([X^{t,\mu}_T])(x) \right) \del p(\mu,t,T,x,y)(v) \, dy  \, d\mu(x),
\end{align}
and
\begin{align}\label{expression_gradient_lin_deri_sol_PDE}
	\pa_v\del U(t,\mu)(v) &= \int_{\R^d} \left(\del \phi([X^{t,\mu}_T])(y) - \del \phi([X^{t,\mu}_T])(v) \right) \pa_xp(\mu,t,T,v,y) \, dy \\\notag &\quad + \int_{\R^{2d}} \left( \del \phi([X^{t,\mu}_T])(y) - \del \phi([X^{t,\mu}_T])(x) \right) \pa_v\del p(\mu,t,T,x,y)(v) \, dy  \, d\mu(x).
\end{align}

By \eqref{density_bound_thm}, \eqref{gradient_linear_der_bound_thm}, we similarly get that there exists a positive constant $C$ such that for all $t \in [0,T)$, $\mu \in \PP$ and $v \in \R^d$ \begin{align*}
	\left\vert \pa_v\del U(t,\mu)(v) \right\vert &\leq C (T-t)^{\frac{\delta-1}{\alpha}}  + C(T-t)^{\frac{\delta-1}{\alpha} +1 + \frac{\eta -1}{\alpha} } \\ &\leq C (T-t)^{\frac{\delta-1}{\alpha}}.
\end{align*}

We now prove \eqref{gradient_lin_der_holder_sol_PDE}. Let us first assume that $|v_1-v_2| \geq (T-t)^{\frac{1}{\alpha}}$. In this case, it follows from \eqref{gradient_lin_der_bound_sol_PDE} that \begin{align*}
	\left\vert \pa_v \del U(t,\mu)(v_1) -  \pa_v \del U(t,\mu)(v_2) \right\vert &\leq  	\left\vert \pa_v \del U(t,\mu)(v_1)\right\vert + 	\left\vert \pa_v \del U(t,\mu)(v_2) \right\vert \\ &\leq C (T-t)^{\frac{\delta -1}{\alpha}} \\ &\leq C (T-t)^{\frac{\delta - 1 - \gamma}{\alpha}}|v_1-v_2|^\gamma.
\end{align*}

Assume now that $|v_1 - v_2| < (T-t)^{\frac{1}{\alpha}}$. Using \eqref{expression_gradient_lin_deri_sol_PDE} and the fact that $\int_{\R^d} \pa_x p(\mu,t,T,v_1,y) \, dy = 0$, one has 

\begin{align*}
	&\pa_v\del U(t,\mu)(v_1) - \pa_v\del U(t,\mu)(v_2)\\ &= \int_{\R^d} \left(\del \phi([X^{t,\mu}_T])(y) - \del \phi([X^{t,\mu}_T])(v_2) \right) (\pa_xp(\mu,t,T,v_1,y) - \pa_xp(\mu,t,T,v_2,y)) \, dy \\ &\quad + \int_{\R^{2d}} \left( \del \phi([X^{t,\mu}_T])(y) - \del \phi([X^{t,\mu}_T])(x) \right) (\pa_v\del p(\mu,t,T,x,y)(v_1) - \pa_v\del p(\mu,t,T,x,y)(v_2)) \, dy  \, d\mu(x).
\end{align*}

Thanks to the uniform $\delta$-Hölder continuity of $\del \phi (\mu)(\cdot)$, \eqref{gradientdensity_holder_thm}, \eqref{gradient_linear_der_holder_v_thm} since $\gamma \in (0,(2\alpha-2)\wedge (\eta + \alpha -1))$ and \eqref{controldensityref} since $|v_1 - v_2| < (T-t)^{\frac{1}{\alpha}}$, we obtain that 

\begin{align*}
		&\left\vert\pa_v\del U(t,\mu)(v_1) - \pa_v\del U(t,\mu)(v_2)\right\vert\\ &\leq C \int_{\R^d} |y-v_2|^\delta (T-t)^{-\frac{\gamma +1}{\alpha}}|v_1-v_2|^\gamma \rho^1(T-t,y-v_2) \, dy \\ &\quad + C\int_{\R^{2d}} |y-x|^\delta (T-t)^{- \frac{\gamma +1}{\alpha} + 1 + \frac{\eta -1}{\alpha}} |v_1-v_2|^\gamma \rho^0(T-t,y-x) \, dy \, d\mu(x).
\end{align*}
Using the space-time inequality \eqref{scalingdensityref}, we deduce that 

\begin{align*}
	&\left\vert\pa_v\del U(t,\mu)(v_1) - \pa_v\del U(t,\mu)(v_2)\right\vert\\ &\leq C \ (T-t)^{\frac{\delta -\gamma -1}{\alpha}}|v_1-v_2|^\gamma.
\end{align*}

\textbf{Step 3: Time-derivative and backward Kolmogorov PDE \eqref{EDP_Kolmo_backward}.} Let us fix $t \in [0,T)$ and $h \in [0,t]$ such that $t-h \in [0,T)$. From the Markov property stemming from the well-posedness of the associated martingale problem related to \eqref{McKVSDE}, we obtain that \[ U(t-h,\mu) = \phi ([X^{t-h,\mu}_T]) = \phi ([X^{t,[X^{t-h,\mu}_t]}_T]) = U(t,[X^{t-h,\mu}_t]).\] Since $U \in \CC^{1}([0,T)\times \PP)$ and thanks to \eqref{gradient_lin_der_bound_sol_PDE} and \eqref{gradient_lin_der_holder_sol_PDE}, we can apply Itô's formula of Proposition \ref{ito_formula} for the function $U(t,\cdot)$. It yields

\begin{align*}
	U(t-h,\mu) &= U(t,\mu) + \int_{t-h}^t \E \left(\pa_v \del U(t,[X^{t-h,\mu}_s])(X^{t-h,\mu}_s)\cdot b(s,X^{t-h,\mu}_s,[X^{t-h,\mu}_s])\right) \, ds \\ & \quad + \int_{t-h}^t \E \int_{\R^d} \left( \del U(t,[X^{t-h,\mu}_s])(X^{t-h,\mu}_s + z) - \del U(t,[X^{t-h,\mu}_s])(X^{t-h,\mu}_s) \right. \\ &\left. \hspace{6cm}- \pa_v \del U(t,[X^{t-h,\mu}_s])(X^{t-h,\mu}_s) \cdot z  \right) \,\frac{dz}{|z|^{d+\alpha}} \, ds\\ &= U(t,\mu) + \int_{t-h}^t \mathscr{L}_sU(t,\cdot)([X^{t-h,\mu}_s]) \, ds,
\end{align*}
where $\mathscr{L}_s$ was defined in \eqref{def_operator_McKV_measure}. Using the continuity and the boundedness of $b$ as well as the regularity properties of the map $\mu \mapsto U(t,\mu)$ obtained above, we find that 

\begin{align*}
	\frac1h (	U(t-h,\mu) - U(t,\mu)) \underset{h \to 0^+}{\longrightarrow} \mathscr{L}_tU(t,\cdot)(\mu).
\end{align*}
This proves the left-differentiability of the map $t \in [0,T) \mapsto U(t,\mu)$. It also follows from the regularity of the drift $b$ and of $U$ on $[0,T) \times \PP$ that the map $(t,\mu) \in [0,T)\times \PP \mapsto \mathscr{L}_tU(t,\cdot)(\mu)$ is continuous. Thus, the map $t \in [0,T) \mapsto U(t,\mu)$ is $\CC^1$ and satisfies for all $t \in [0,T)$ and $\mu \in \PP$ \[ \pa_t U(t,\mu) + \mathscr{L}_tU(t,\cdot)(\mu) =0.\]
This shows that $U$ is solution to the backward Kolmogorov PDE \eqref{EDP_Kolmo_backward}.\\ 

\textbf{Step 4: Uniqueness of the solution to \eqref{EDP_Kolmo_backward}.} Let us consider $V \in \CC^0([0,T]\times \PP) \cap \CC^1([0,T)\times \PP)$ another solution to \eqref{EDP_Kolmo_backward} satisfying \eqref{gradient_lin_der_bound_sol_PDE} and \eqref{gradient_lin_der_holder_sol_PDE}. Let us fix $(t,\mu) \in [0,T) \times \PP$. For any $\tau \in [t,T)$, we can apply Itô's formula of Proposition \ref{ito_formula} for the map $(s,\mu) \in [t,\tau] \times \PP \mapsto V(s,\mu)$ which yields \begin{align*}
	V(\tau,[X^{t,\mu}_\tau]) = V(t,\mu) + \int_t^\tau \pa_s V(s,[X^{t,\mu}_s]) \, ds + \int_t^\tau \mathscr{L}_s V(s,\cdot)([X^{t,\mu}_s]) \, ds.
\end{align*}
Since $V$ solves \eqref{EDP_Kolmo_backward}, we obtain that $V(\tau,[X^{t,\mu}_\tau]) = V(t,\mu).$ We then use the continuity of the maps $(t,\mu) \in [0,T]\times \PP \mapsto V(t,\mu)$ and $\tau \in [t,T] \mapsto [X^{t,\mu}_\tau] \in \PP$ to let $\tau$ tends to $T$. This yields $$ \phi([X^{t,\mu}_T]) = U(t,\mu) = V(t,\mu).$$

\section{Quantitative weak propagation of chaos}\label{section_prop_chaos}

This section is dedicated to prove Theorem \ref{Thm_POC}. We first need to establish additional regularity properties on the solution to the PDE \eqref{EDP_Kolmo_backward} to prove our weak propagation of chaos result.

\begin{Prop}\label{prop_reg_sol_edp_for_prop_chaos}
	Let us fix $\delta \in (0,1]$, $L>0$ and $\gamma \in (0,1]\cap(0,(2\alpha -2) \wedge (\eta + \alpha -1))$. Then, under Assumption \ref{Assumption2}, there exists positive constant $C = C(d,T,\alpha,\beta,\ref{Assumption2},\delta,L,\gamma)$ such that for all $\phi \in \CC^{(2,\delta)}_L(\PP)$ (defined in Definition \ref{def_space_prop_chaos}), the solution $U$ of the backward Kolmogorov PDE \eqref{EDP_Kolmo_backward} with terminal condition $\phi$ at time $T$ satisfies the following properties. 
	
	\begin{itemize}
		
		\item For all $t \in [0,T)$, $\mu \in \PP$, $v \in \R^d$ 
		
		\begin{equation}\label{gradient_lin_der_bound_sol_PDE_POC}
			\left\vert \pa_v\del U (t,\mu)(v)\right\vert \leq C(T-t)^{\frac{\delta -1}{\alpha}}.
		\end{equation}
		
		\item For all $t \in [0,T)$, $\mu_1,\mu_2 \in \PP$, $v_1,v_2 \in \R^d$ 
		
		\begin{equation}\label{gradient_lin_der_holder_sol_PDE_POC}
			\left\vert \pa_v\del U (t,\mu_1)(v_1) - \pa_v\del U (t,\mu_2)(v_2)\right\vert \leq C(T-t)^{\frac{\delta -1 - \gamma}{\alpha}}\left(|v_1-v_2|^\gamma + W_1^\gamma (\mu_1,\mu_2)\right).
		\end{equation}
		
		\item For all $t \in [0,T)$, $\mu_1,\mu_2 \in \PP$, $v \in \R^d$ 
		
		\begin{equation}\label{lin_der_holder_sol_PDE_POC}
			\left\vert \del U (t,\mu_1)(v) - \del U (t,\mu_2)(v)\right\vert \leq C(T-t)^{\frac{\delta  - 1}{\alpha}} W_1 (\mu_1,\mu_2).
		\end{equation}
	\end{itemize}

\end{Prop}

\begin{proof}[Proof of Proposition \ref{prop_reg_sol_edp_for_prop_chaos}]
	First, note that \eqref{gradient_lin_der_bound_sol_PDE_POC} have been proved in Theorem \ref{Thm_EDP}.\\
	
	\noindent\textbf{Proof of \eqref{gradient_lin_der_holder_sol_PDE_POC}.} We write 
	
	\begin{align}\label{proof_lemme_POC_eq10}
		\pa_v \del U(t,\mu_1)(v_1) - \pa_v \del U(t,\mu_2)(v_2)  \notag&= \pa_v \del U(t,\mu_1)(v_1) - \pa_v \del U(t,\mu_1)(v_2) \\ \notag&\quad+ \pa_v \del U(t,\mu_1)(v_2) - \pa_v \del U(t,\mu_2)(v_2)\\ &=: I_1 + I_2.
	\end{align}
	Then, using \eqref{gradient_lin_der_holder_sol_PDE}, we obtain that 
	
	\begin{equation}\label{proof_lemme_POC_eq11}
		|I_1| \leq C(T-t)^{\frac{\delta -1 - \gamma}{\alpha}}|v_1-v_2|^\gamma.
	\end{equation}
	
	We now focus on $I_2$. Using \eqref{expression_gradient_lin_deri_sol_PDE}, we get

	\begin{align*}
		&\pa_v\del U(t,\mu_1)(v_2) - \pa_v\del U(t,\mu_2)(v_2)\\ &= \int_{\R^d} \left(\del \phi([X^{t,\mu_1}_T])(y) - \del \phi([X^{t,\mu_1}_T])(v_2) \right)  \pa_xp(\mu_1,t,T,v_2,y)\, dy \\ &\quad + \int_{\R^{2d}} \left( \del \phi([X^{t,\mu_1}_T])(y) - \del \phi([X^{t,\mu_1}_T])(x) \right)  \pa_v\del p(\mu_1,t,T,x,y)(v_2)\, dy  \, d\mu_1(x)\\ &\quad -\int_{\R^d} \left(\del \phi([X^{t,\mu_2}_T])(y) - \del \phi([X^{t,\mu_2}_T])(v_2) \right)  \pa_xp(\mu_2,t,T,v_2,y)\, dy \\ &\quad - \int_{\R^{2d}} \left( \del \phi([X^{t,\mu_2}_T])(y) - \del \phi([X^{t,\mu_2}_T])(x) \right)  \pa_v\del p(\mu_2,t,T,x,y)(v_2)\, dy  \, d\mu_2(x). \end{align*} 
	
	We decompose it in the following way \begin{align}\label{proof_lemme_POC_eq0}
		\notag&\pa_v\del U(t,\mu_1)(v_2) - \pa_v\del U(t,\mu_2)(v_2) \\\notag&= \int_{\R^d} \left(\del \phi([X^{t,\mu_1}_T])(y) - \del \phi([X^{t,\mu_1}_T])(v_2) - \left[\del \phi ([X^{t,\mu_2}_T])(y)  - \del \phi ([X^{t,\mu_2}_T])(v_2)\right]\right) \pa_x p(\mu_1,t,T,v_2,y) \, dy\\ \notag&\quad +   \int_{\R^d} \left(\del \phi ([X^{t,\mu_2}_T])(y) - \del \phi ([X^{t,\mu_2}_T])(v_2) \right)  \left(\pa_x p(\mu_1,t,T,v_2,y) - \pa_x p(\mu_2,t,T,v_2,y) \right) \, dy \\ \notag&\quad + \int_{\R^{2d}}  \left(\del \phi([X^{t,\mu_1}_T])(y) - \del \phi([X^{t,\mu_1}_T])(x) - \left[\del \phi ([X^{t,\mu_2}_T])(y)  - \del \phi ([X^{t,\mu_2}_T])(x)\right]\right) \\ \notag&\hspace{11cm} \pa_v \del p(\mu_1,t,T,x,y)(v_2) \, dy\, d\mu_1(x) \\ \notag&\quad +   \int_{\R^{2d}} \left(\del \phi ([X^{t,\mu_2}_T])(y) - \del \phi ([X^{t,\mu_2}_T])(x)\right) \\ \notag&\hspace{5cm} \left(\pa_v \del p(\mu_1,t,T,x,y)(v_2) - \pa_v \del p(\mu_2,t,T,x,y)(v_2) \right) \, dy \, d\mu_1(x)\\ \notag&\quad +   \int_{\R^{2d}} \left(\del \phi ([X^{t,\mu_2}_T])(y) - \del \phi ([X^{t,\mu_2}_T])(x)\right)\pa_v \del p(\mu_2,t,T,x,y)(v_2)\, dy \, d(\mu_1-\mu_2)(x) \\ &=: J_1 + J_2 + J_3 + J_4 + J_5.
	\end{align}
	
	Before estimating each term of the preceding decomposition, let us prove that for all $\gamma \in (0,1]$, there exists a positive constant $C$ such that for any $ \phi \in \CC^{(2,\delta)}_L(\PP)$, $t \in [0,T)$, $\mu_1,\mu_2, \in \PP$, $ x,y \in \R^d$ \begin{multline}\label{proof_lemme_POC_eq1}
		\left\vert\del \phi([X^{t,\mu_1}_T])(y) - \del \phi([X^{t,\mu_1}_T])(x) - \left[\del \phi ([X^{t,\mu_2}_T])(y)  - \del \phi ([X^{t,\mu_2}_T])(x)\right]\right\vert \\ \leq C(T-t)^{- \frac{\gamma}{\alpha}} |x-y|^\delta W_1^\gamma(\mu_1,\mu_2).
	\end{multline}
	
	To prove this, we write \begin{align*}
		& \del \phi([X^{t,\mu_1}_T])(y) - \del \phi([X^{t,\mu_1}_T])(x) - \left[\del \phi ([X^{t,\mu_2}_T])(y)  - \del \phi ([X^{t,\mu_2}_T])(x)\right] \\ &= \int_0^1\int_{\R^d} \left(\dell \phi(m_{\lambda})(y,v) - \dell \phi(m_{\lambda})(x,v) \right) (p(\mu_1,t,T,v) - p(\mu_2,s,T,v)) \, dv \, d\lambda \\ &= \int_0^1\int_{\R^{2d}} \left(\dell \phi(m_{\lambda})(y,v) - \dell \phi(m_{\lambda})(x,v) \right) p(\mu_1,t,T,x',v) \, d(\mu_1 - \mu_2)(x')\, dv \, d\lambda  \\ &\quad + \int_0^1\int_{\R^{2d}} \left(\dell \phi(m_{\lambda})(y,v) - \dell \phi(m_{\lambda})(x,v) \right) \left(p(\mu_1,t,T,x',v) - p(\mu_2,t,T,x',v) \right)\, d\mu_2(x')\, dv \, d\lambda  \\ &=: K_1 + K_2,
	\end{align*}
	where $m_\lambda := \lambda [X^{t,\mu_1}_T] + (1-\lambda)[X^{t,\mu_2}_T]$, for $\lambda \in [0,1]$. It follows from \eqref{gradient_density_holder_measure_thm} and the $\delta$-Hölder continuity of $\dell \phi (\mu)(\cdot,v)$ uniformly with respect to $\mu \in \PP$ and $v \in \R^d$ that
	
	\begin{equation}\label{proof_lemme_POC_eq2}
		|K_2| \leq C (T-t)^{1 - \frac{\gamma +1}{\alpha}}|y-x|^\delta W_1^\gamma (\mu_1,\mu_2).
	\end{equation}
	Concerning $K_1$, we set for $x' \in \R^d$ \[ F(x'):= \int_{\R^d} \left(\dell \phi(m_{\lambda})(y,v) - \dell \phi(m_{\lambda})(x,v) \right) p(\mu_1,t,T,x',v) \, dv.\]  Let us fix $x', x'' \in \R^d$. Thanks to \eqref{gradientdensity_holder_thm} and the $\delta$-Hölder continuity of $\dell \phi (\mu)(\cdot,v)$ uniformly with respect to $\mu \in \PP$ and $v \in \R^d$, one gets that
	
	\[|F(x') - F(x'')| \leq C |x-y|^\delta(T-t)^{-\frac{\gamma}{\alpha}} |x'-x''|^\gamma.\]
	Taking an optimal coupling between $\mu_1$ and $\mu_1$ for $W_1$ directly yields, with Jensen's inequality, \begin{equation}\label{proof_lemme_POC_eq3}
		|K_1| \leq  C |x-y|^\delta(T-t)^{-\frac{\gamma}{\alpha}} W_1^\gamma (\mu_1,\mu_2).
	\end{equation}
	Combining \eqref{proof_lemme_POC_eq2} and \eqref{proof_lemme_POC_eq3} concludes the proof of \eqref{proof_lemme_POC_eq1} since $\alpha \in (1,2)$. We can now turn to estimate each term of \eqref{proof_lemme_POC_eq0}. Using \eqref{proof_lemme_POC_eq1}, \eqref{density_bound_thm} and the space-time inequality \eqref{scalingdensityref}, one obtains that \begin{align}\label{proof_lemme_POC_eq4}
		|J_1| \notag& \leq  C\int_{\R^d} (T-t)^{-\frac{\gamma}{\alpha}} W_1^\gamma(\mu_1, \mu_2) |y-v_2|^\delta (T-t)^{-\frac{1}{\alpha}}\rho^1(T-t,y-v_2)  \, dy\\ &\leq C (T-t)^{\frac{\delta -\gamma -1}{\alpha}} W_1^\gamma(\mu_1, \mu_2) \int_{\R^d} \rho^{1-\delta}(T-t,y-v_2)  \, dy \\ \notag&\leq C (T-t)^{\frac{\delta - \gamma -1}{\alpha}}W_1^\gamma(\mu_1,\mu_2).
	\end{align} 
	
	For $J_2$, it follows from the $\delta$-Hölder continuity of $\del \phi(\mu)(\cdot)$ uniformly with respect to $\mu \in \PP$, \eqref{gradient_density_holder_measure_thm} since $\gamma < \eta  + \alpha -1$, and the space-time inequality \eqref{scalingdensityref} that 
	
	\begin{align}\label{proof_lemme_POC_eq5}
		|J_2| \notag &\leq C \int_{\R^d} |y-v_2|^\delta (T-t)^{1 - \frac{2 + \gamma}{\alpha}} W_1^\gamma (\mu_1,\mu_2) \rho^1(T-t,y-v_2) \, dy \\ &\leq  C  (T-t)^{1 + \frac{\delta -2 - \gamma}{\alpha}} W_1^\gamma (\mu_1,\mu_2) \\ \notag& \leq C  (T-t)^{\frac{\delta -1 - \gamma}{\alpha}} W_1^\gamma (\mu_1,\mu_2),
	\end{align}
	since $\alpha \in (1,2)$. For $J_3$, \eqref{proof_lemme_POC_eq1}, \eqref{gradient_linear_der_bound_thm} and the space-time inequality \eqref{scalingdensityref} yield 
	
	\begin{align}\label{proof_lemme_POC_eq6}
		|J_3| \notag &\leq C \int_{\R^{2d}} (T-t)^{ - \frac{\gamma}{\alpha}} |y-x|^\delta  W_1^\gamma(\mu_1,\mu_2) (T-t)^{\frac{\eta - 1}{\alpha} + 1 - \frac{1}{\alpha}}\rho^0(T-t,y-x) \, dy\, d\mu_1(x) \\ &\leq  C  (T-t)^{\frac{\delta -1 - \gamma}{\alpha} + 1 + \frac{\eta -1}{\alpha}} W_1^\gamma (\mu_1,\mu_2) \\ \notag& \leq C  (T-t)^{\frac{\delta -1 - \gamma}{\alpha}} W_1^\gamma (\mu_1,\mu_2).
	\end{align}
	
	It follows from the $\delta$-Hölder continuity of $\del \phi(\mu)(\cdot)$ uniformly with respect to $\mu \in \PP$, \eqref{gradient_linear_derivative_holder_measure_thm} and the space-time inequality \eqref{scalingdensityref} that 
	
	\begin{align}\label{proof_lemme_POC_eq7}
		|J_4| \notag &\leq C \int_{\R^{2d}} |y-x|^\delta (T-t)^{\frac{\eta -\gamma -1}{\alpha} + 1 - \frac{1}{\alpha}}W_1^\gamma(\mu_1,\mu_2)\rho^0(T-t,y-x) \, dy \, d\mu_1(x) \\ &\leq  C  (T-t)^{\frac{\delta -1 - \gamma}{\alpha} + 1 + \frac{\eta -1}{\alpha}} W_1^\gamma (\mu_1,\mu_2) \\ \notag& \leq C  (T-t)^{\frac{\delta -1 - \gamma}{\alpha}} W_1^\gamma (\mu_1,\mu_2).
	\end{align}

	We finally deal with $J_5$. We set for $x \in \R^d$ \[H(x):= \int_{\R^{d}} \left(\del \phi ([X^{t,\mu_2}_T])(y) - \del \phi ([X^{t,\mu_2}_T])(x)\right)\pa_v \del p(\mu_2,t,T,x,y)(v_2)\, dy. \]

	 Let us prove that for any $x,x' \in \R^d$ \begin{equation}\label{proof_lemme_POC_eq8}
		|H(x)-H(x')| \leq  C  (T-t)^{\frac{\delta-1-\gamma}{\alpha} + 1 + \frac{\eta -1}{\alpha}} |x-x'|^\gamma.
	\end{equation} Assume first that $|x-x'|>(T-t)^{\frac{1}{\alpha}}.$ The $\delta$-Hölder continuity of $\del \phi(\mu)(\cdot)$ uniformly with respect to $\mu \in \PP$, \eqref{gradient_linear_der_bound_thm} and the space-time inequality \eqref{scalingdensityref} ensure that \begin{align*}
		|H(x) - H(x')| &\leq |H(x)| + | H(x')| \\ &\leq C \int_{\R^d} |y-x|^\delta (T-t)^{\frac{\eta -1}{\alpha} +1 - \frac{1}{\alpha}} \rho^0(T-t,y-x) \, dy \\ &\leq  C  (T-t)^{\frac{\delta-1}{\alpha} + 1 + \frac{\eta -1}{\alpha}} \\ &\leq C  (T-t)^{\frac{\delta-1-\gamma}{\alpha} + 1 + \frac{\eta -1}{\alpha}} |x-x'|^\gamma.
	\end{align*}In the case where $|x-x'| \leq (T-t)^{\frac{1}{\alpha}}$, thanks to \eqref{centering_prop_lin_der_density}, we write  \begin{align*}
		&H(x) - H(x')\\ & =  \int_{\R^{d}} \left( \del \phi ([X^{t,\mu_2}_T])(x') - \del \phi ([X^{t,\mu_2}_T])(x) \right)\pa_v \del p(\mu_2,t,T,x,y)(v_2)\, dy \\ &\quad+ \int_{\R^d} \left(\del \phi ([X^{t,\mu_2}_T])(y) - \del \phi ([X^{t,\mu_2}_T])(x')\right) \left( \pa_v \del p(\mu_2,t,T,x,y)(v_2) - \pa_v \del p(\mu_2,t,T,x',y)(v_2)\right) \, dy\\ &= \int_{\R^d} \left(\del \phi ([X^{t,\mu_2}_T])(y) - \del \phi ([X^{t,\mu_2}_T])(x')\right) \left( \pa_v \del p(\mu_2,t,T,x,y)(v_2) - \pa_v \del p(\mu_2,t,T,x',y)(v_2)\right) \, dy.
	\end{align*}
	
	Thanks to the $\delta$-Hölder continuity of $\del \phi(\mu)(\cdot)$ uniformly with respect to $\mu \in \PP$ and \eqref{gradient_linear_der_holder_x_thm} since $\gamma < \eta + \alpha -1$, we directly get that  \begin{align*}
		|H(x)-H(x')| &\leq C \int_{\R^d} |y-x'|^\delta (T-t)^{\frac{\eta-1-\gamma}{\alpha} + 1 - \frac{1}{\alpha}}|x-x'|^\gamma  \left[\rho^0(T-t,y-x') + \rho^0(T-t,y-x)\right] \,dy.
	\end{align*}
	
	Since $|x-x'|\leq (T-t)^{\frac{1}{\alpha}}$, it follows from \eqref{controldensityref} and the space-time inequality \eqref{scalingdensityref} that \begin{align*}
		|H(x)-H(x')| & \leq C \int_{\R^d} |y-x'|^\delta (T-t)^{\frac{\eta-1-\gamma}{\alpha} + 1 - \frac{1}{\alpha}}|x-x'|^\gamma  \rho^0(T-t,y-x') \,dy \\&\leq C  (T-t)^{\frac{\delta-1-\gamma}{\alpha} + 1 + \frac{\eta -1}{\alpha}} |x-x'|^\gamma.
	\end{align*}
	
	Then, it follows from \eqref{proof_lemme_POC_eq8} that \begin{align}\label{proof_lemme_POC_eq9}
		|J_5|\notag &\leq  C  (T-t)^{\frac{\delta-1-\gamma}{\alpha} + 1 + \frac{\eta -1}{\alpha}} W_1^\gamma (\mu_1,\mu_2) \\ &\leq  C  (T-t)^{\frac{\delta-1-\gamma}{\alpha}} W_1^\gamma (\mu_1,\mu_2).
	\end{align}
	
	Coming back to \eqref{proof_lemme_POC_eq10} and using \eqref{proof_lemme_POC_eq0}, \eqref{proof_lemme_POC_eq4}, \eqref{proof_lemme_POC_eq5}, \eqref{proof_lemme_POC_eq6}, \eqref{proof_lemme_POC_eq7}, \eqref{proof_lemme_POC_eq9}, we have shown that \[ |I_2| \leq   C  (T-t)^{\frac{\delta-1-\gamma}{\alpha}} W_1^\gamma (\mu_1,\mu_2).\] 
	
	This estimate together with \eqref{proof_lemme_POC_eq11} ends the proof of \eqref{gradient_lin_der_holder_sol_PDE_POC}.\\
	
	\noindent\textbf{Proof of \eqref{lin_der_holder_sol_PDE_POC}.} It follows from \eqref{expression_gradient_lin_deri_sol_PDE} that

	\begin{align*}
		&\del U(t,\mu_1)(v) - \del U(t,\mu_2)(v)\\ &= \int_{\R^d} \del \phi([X^{t,\mu_1}_T])(y)  p(\mu_1,t,T,v,y)\, dy \\ &\quad + \int_{\R^{2d}} \left( \del \phi([X^{t,\mu_1}_T])(y) - \del \phi([X^{t,\mu_1}_T])(x) \right)  \del p(\mu_1,t,T,x,y)(v)\, dy  \, d\mu_1(x)\\ &\quad -\int_{\R^d} \del \phi([X^{t,\mu_2}_T])(y) p(\mu_2,t,T,v,y)\, dy \\ &\quad - \int_{\R^{2d}} \left( \del \phi([X^{t,\mu_2}_T])(y) - \del \phi([X^{t,\mu_2}_T])(x) \right) \del p(\mu_2,t,T,x,y)(v)\, dy  \, d\mu_2(x). \end{align*} 
	
	We decompose it in the following way \begin{align}\label{proof_lemme_POC_2eq0}
		\notag&\del U(t,\mu_1)(v) - \del U(t,\mu_2)(v) \\\notag&= \int_{\R^d} \left(\del \phi([X^{t,\mu_1}_T])(y)- \del \phi ([X^{t,\mu_2}_T])(y) \right)  p(\mu_1,t,T,v,y) \, dy\\ \notag&\quad +   \int_{\R^d} \del \phi ([X^{t,\mu_2}_T])(y) \left(p(\mu_1,t,T,v,y) -  p(\mu_2,t,T,v,y) \right) \, dy \\ \notag&\quad + \int_{\R^{2d}}  \left(\del \phi([X^{t,\mu_1}_T])(y) - \del \phi([X^{t,\mu_1}_T])(x) - \left[\del \phi ([X^{t,\mu_2}_T])(y)  - \del \phi ([X^{t,\mu_2}_T])(x)\right]\right) \\ \notag&\hspace{11cm}  \del p(\mu_1,t,T,x,y)(v_2) \, dy\, d\mu_1(x) \\ \notag&\quad +   \int_{\R^{2d}} \left(\del \phi ([X^{t,\mu_2}_T])(y) - \del \phi ([X^{t,\mu_2}_T])(x)\right) \\ \notag&\hspace{5cm} \left( \del p(\mu_1,t,T,x,y)(v_2) - \del p(\mu_2,t,T,x,y)(v_2) \right) \, dy \, d\mu_1(x)\\ \notag&\quad +   \int_{\R^{2d}} \left(\del \phi ([X^{t,\mu_2}_T])(y) - \del \phi ([X^{t,\mu_2}_T])(x)\right) \del p(\mu_2,t,T,x,y)(v_2)\, dy \, d(\mu_1-\mu_2)(x) \\ &=: J_1 + J_2 + J_3 + J_4 + J_5.
	\end{align}
	
	We first focus of $J_1$. Let us prove that there exists a positive constant $C$ such that for all $t \in [0,T)$, $\mu_1,\mu_2 \in \PP$, $y \in \R^d$

	\begin{equation}\label{proof_lemme_POC_2eq1}
		\left\vert \del \phi([X^{t,\mu_1}_T])(y)- \del \phi ([X^{t,\mu_2}_T])(y) \right\vert C(T-t)^{\frac{\delta -1}{\alpha}} W_1(\mu_1,\mu_2).
	\end{equation}
	
	To prove this, we write \begin{align*}
		&\del \phi([X^{t,\mu_1}_T])(y)- \del \phi ([X^{t,\mu_2}_T])(y) \\ &= \int_0^1 \int_{\R^d} \dell \phi(m_\lambda)(y,z) ( p(\mu_1,t,T,z) - p(\mu_2,t,T,z)) \, dz \, d\lambda \\ &=  \int_0^1 \int_{\R^{2d}} \dell \phi(m_\lambda)(y,z) ( p(\mu_1,t,T,x,z) - p(\mu_2,t,T,x,z))\, d\mu_1(x) \, dz \, d\lambda \\ &\quad +  \int_0^1 \int_{\R^{2d}} \dell \phi(m_\lambda)(y,z)p(\mu_2,t,T,x,z)\, d(\mu_1-\mu_2)(x) \, dz \, d\lambda  \\ &=: K_1 + K_2,
	\end{align*}
	where $m_\lambda := \lambda [X^{t,\mu_1}_T] + (1-\lambda)[X^{t,\mu_2}_T]$. We rewrite $K_1$ as \begin{equation*}
		K_1 =  \int_0^1 \int_{\R^{2d}} \dell \phi(m_\lambda)(y,z) \int_0^1  \int_{\R^d} \del p(M_r,t,T,x,z) (w) \, d(\mu_1 - \mu_2)(w) \, dr \, d\mu_1(x) \, dz \, d\lambda,
	\end{equation*}
	where $M_r = r \mu_1 + (1-r) \mu_2.$ We set for $w \in \R^d$ \[ F(w) := \int_{\R^d} \dell \phi(m_\lambda)(y,z) \del p(M_r,t,T,x,z) (w) \, dz,\]
	where $\lambda, r$ and $x$ are fixed. For $w_1,w_2 \in \R^d$, one has 
	
	\begin{align*}
		&F(w_1) - F(w_2)\\& = \int_{\R^d}  \dell \phi (m_\lambda)(y,z) \int_0^1\pa_v \del p(M_r,t,T,x,z)(sw_1 + (1-s)w_2) \cdot (w_1 - w_2) \, ds \, dz \\ &=\int_0^1\int_{\R^d}  \left(\dell \phi (m_\lambda)(y,z) - \dell \phi (m_\lambda)(y,x) \right) \pa_v \del p(M_r,t,T,x,z)(sw_1 + (1-s)w_2) \cdot (w_1 - w_2)  \, dz \,ds.
	\end{align*}
	
	It follows from the $\delta$-Hölder continuity of $\dell \phi(\mu)(y,\cdot)$ uniformly with respect to $\mu \in \PP$ and $y \in \R^d$, \eqref{gradient_linear_der_bound_thm} and the space-time inequality \eqref{scalingdensityref} that \[|F(w_1) - F(w_2)| \leq C (T-t)^{\frac{\delta-1}{\alpha } + 1 + \frac{\eta -1}{\alpha}}|w_1 - w_2|.\]
	
	This yields \begin{align}\label{proof_lemme_POC_2eq2}
		\notag|K_1| &\leq C (T-t)^{\frac{\delta-1}{\alpha} + 1 + \frac{\eta -1}{\alpha}} W_1(\mu_1,\mu_2)\\ &\leq   C (T-t)^{\frac{\delta-1}{\alpha}} W_1(\mu_1,\mu_2). \end{align}
	
	We control $K_2$ as for $K_1$ by studying the regularity with respect to $x$ of the function $G$ given by 
	
	\begin{equation*}
		G(x):=    \int_{\R^{d}} \dell \phi(m_\lambda)(y,z)p(\mu_2,t,T,x,z)\,  dz.
	\end{equation*}
	
	For $x_1,x_2 \in \R^d$, one has 
	
	\begin{align*} 
		G(x_1) - G(x_2)& =  \int_{\R^{d}} \dell \phi(m_\lambda)(y,z) \int_0^1 \pa_x p(\mu_2,t,T,rx_1 + (1-r)x_2,z) \cdot (x_1 - x_2)\,dr \,   dz \\ &=   \int_0^1 \int_{\R^{d}} \left(\dell \phi(m_\lambda)(y,z) - \dell \phi(m_\lambda)(y, rx_1 + (1-r)x_2) \right) \\ &\hspace{7cm} \pa_x p(\mu_2,t,T,rx_1 + (1-r)x_2,z) \cdot (x_1 - x_2) \,   dz \,dr.
	\end{align*}
	
	Using $\delta$-Hölder continuity of $\dell \phi(\mu)(y,\cdot)$ uniformly with respect to $\mu \in \PP$ and $y \in \R^d$, \eqref{density_bound_thm} and the space-time inequality \eqref{scalingdensityref}, we obtain that \[|G(x_1) - G(x_2)| \leq C (T-t)^{\frac{\delta-1}{\alpha }}|x_1 - x_2|.\]
	
	It proves that \begin{equation}\label{proof_lemme_POC_2eq3}
		|K_2| \leq   C (T-t)^{\frac{\delta-1}{\alpha}} W_1(\mu_1,\mu_2). \end{equation}
	
	Combining \eqref{proof_lemme_POC_2eq2} and \eqref{proof_lemme_POC_2eq3} concludes the proof of \eqref{proof_lemme_POC_2eq1}. We can now turn to estimate $J_1$ in \eqref{proof_lemme_POC_2eq0}. Using \eqref{proof_lemme_POC_2eq1}, one gets  \begin{align}\label{proof_lemme_POC_2eq4}
		|J_1| \notag &\leq  C\int_{\R^d} (T-t)^{\frac{\delta-1}{\alpha}} W_1(\mu_1,\mu_2) p(\mu_2,t,T,v,y) \, dy\\  &\leq  C (T-t)^{\frac{\delta-1}{\alpha}} W_1(\mu_1,\mu_2).
	\end{align} 
	
	For $J_2$, we rewrite it as \begin{align*}
		J_2 &= \int_{\R^d} \del \phi([X^{t,\mu_2}_T])(y) \int_0^1 \int_{\R^d} \del p(M_r,t,T,v,y)(z) \, d(\mu_1 - \mu_2)(z) \, dr\, dy,
	\end{align*}
	where $M_r:= r\mu_1 + (1-r)\mu_2.$ Following the same lines as for the proof of the above estimates on $K_1$ \eqref{proof_lemme_POC_2eq2}, we obtain that 
	
	\begin{equation}\label{proof_lemme_POC_2eq5}
		|J_2| \leq C (T-t)^{\frac{\delta- 1}{\alpha}}W_1(\mu_1,\mu_2).
	\end{equation}
	
	Following same lines as in the proof of \eqref{proof_lemme_POC_eq6}, \eqref{proof_lemme_POC_eq7} and \eqref{proof_lemme_POC_eq9} by using \eqref{linear_der_bound_thm}, \eqref{linear_derivative_holder_measure_thm} and \eqref{linear_der_holder_x_thm}, we can prove the following estimates

	\begin{align}\label{proof_lemme_POC_2eq6}
		|J_3| + |J_4| + |J_5| &\leq C  (T-t)^{\frac{\delta -1}{\alpha} +1 - \frac{1}{\alpha}} W_1 (\mu_1,\mu_2).
	\end{align}

	Gathering \eqref{proof_lemme_POC_2eq0} \eqref{proof_lemme_POC_2eq4}, \eqref{proof_lemme_POC_2eq5} and \eqref{proof_lemme_POC_2eq6}, we have proved \eqref{lin_der_holder_sol_PDE_POC}.

\end{proof}

\begin{proof}[Proof of Theorem \ref{Thm_POC}]
	
	Let us first introduce some notations. We can write for all $i \in \{ 1,\dots, N\}$ and for all $t \in [0,T]$ $$Z^i_t = \int_0^t \int_{\R^d} z \, \NNN^i(ds,dz),$$ where $\NN^i$ is the Poisson random measure associated with $Z^i$, and $\NNN^i$ is the associated compensated Poisson random measure. Then, we set for all $t \in [0,T]$ $$\textbf{Z}^N_t := \begin{pmatrix} Z^1_t \\ \vdots \\ Z^N_t \end{pmatrix} \in (\R^{d})^N.$$  
	As the Lévy processes $(Z^n)_n$ are independent, the process $(\textbf{Z}^N_t)_t$ is a Lévy process in $(\R^{d})^N.$ Il is a cylindrical $\alpha$-stable process. Its Poisson random measure $\bm{\NN}^N$ and its Lévy measure $\bm{\nu}^N$ are defined as follows. For all $\varphi : [0,T] \times (\R^{d})^N \rightarrow \R^+,$ one has \begin{equation}\label{Poissonmeasureparticles}
		\int_{0}^T\int_{\R^{Nd}} \varphi(s,\x) \, \bm{\NN}^N(ds,d\x) = \sum_{i=1}^N \int_0^t \int_{\R^d} \varphi(s,0, \dots,0,x_i,0,\dots,0) \, \NN^i(ds,dx_i).\end{equation} For all $\phi : \R^{Nd} \rightarrow \R^+,$ one has \begin{equation}\label{Levymeasureparticles}
		\int_{\R^{Nd}} \phi(\x) \, d\bm{\nu}^N(\x)= \sum_{i=1}^N  \int_{\R^d} \phi(0, \dots,0,x_i,0,\dots,0) \, d\nu(x_i).\end{equation}  Note that since $(Z^n)_n$ are independent processes, for all $t \in [0,T],$ the support of the random measure $\bm{\NN}^N(t,d\x)$ is contained in $$ \bigcup_{i=0}^{N-1} \{0\}_{\R^d}^i \times \R^d \times \{0\}_{\R^d}^{N-1-i} \subset (\R^{d})^N.$$

	Let us define for all $t \in [0,T],$ $\x=(x_1, \dots,x_N) \in (\R^d)^N$ \begin{equation*}
		\bm{b}^N(t,\x) := \begin{pmatrix} b(t,x_1,\muu^N_{\x}) \\ \vdots \\ b(t,x_N,\muu^N_{\x}) \end{pmatrix} \in (\R^{d})^N,
	\end{equation*}
	where $\muu^N_{\x} = \frac1N \displaystyle\sum_{j=1}^N \delta_{x_j}.$ Thus, writing for $t \in [0,T]$ and $N\geq 1$ $$\bm{X}^N_t = \begin{pmatrix}X^{1,N}_t \\ \vdots \\ X^{N,N}_t\end{pmatrix},$$  the SDE \eqref{edsparticles} defining the particle system can be rewritten as  \begin{equation}\label{edsparticulesvect}
		\left\{  \begin{array}{lll}
			&d\bm{X}^N_t = \bm{b}^N(t,\bm{X}^N_t) \,dt + d\textbf{Z}^N_t, \quad t \in [0,T], \\ &\bm{X}^N_0 = \begin{pmatrix} X^1_0 \\ \vdots \\ X^N_0 \end{pmatrix}.
		\end{array}\right.\end{equation}

	 \noindent\textbf{Proof of \eqref{eq1_thm_poc}.} Let us consider $U$ the solution to the backward Kolmogorov PDE \eqref{EDP_Kolmo_backward} with terminal condition $\phi$ at time $T$ given in Theorem \ref{Thm_EDP}. Using Lemma \ref{projempC1} and the fact that $U \in \CC^1([0,T)\times \PP)$, we obtain that the function $(t,\bm{x}) \in [0,T) \times (\R^{d})^N\mapsto U(t,\muu^N_{\bm{x}})$ belongs to $\CC^{1}([0,T)\times (\R^{d})^N).$ Moreover, for all $t\in [0,T)$ and $\bm{x} = (x_1, \dots, x_N)\in (\R^{d})^N,$ we have by Proposition \ref{projempC1} $$ \partial_{\bm{x}} U(t,\muu^N_{\bm{x}}) = \frac1N\begin{pmatrix}
	\partial_v \del U (t,\muu^N_{\bm{x}})(x_1) \\ \vdots \\ \partial_v \del U (t,\muu^N_{\bm{x}})(x_N)
\end{pmatrix}.$$ 

Let us fix $\gamma \in (0,1]\cap (\alpha -1, (\delta + \alpha -1)\wedge (2\alpha -2)\wedge (\eta + \alpha -1)).$ We easily see using \eqref{gradient_lin_der_holder_sol_PDE_POC} that the map $\bm{x}\mapsto \pa_{\bm{x}}U(t,\muu^N_{\bm{x}})$ is $\gamma$-Hölder continuous locally uniformly with respect to $t \in [0,T).$ This comes from the fact that for all $\bm{x}=(x_1,\dots,x_N), \bm{y}= (y_1,\dots,y_N) \in (\R^d)^N$, $$W_1^\gamma(\muu^N_{\bm{x}}, \muu^N_{\bm{y}}) \leq \frac{1}{N^\gamma} \sum_{k=1}^N|x_k - y_k|^\gamma.$$ We denote by $\bm{\NN}^N$ the Poisson random measure  associated with $\bm{Z}^N = (Z^1, \dots, Z^N)$ defined in \eqref{Poissonmeasureparticles} and by $\bm{\nu}^N$ its associated Lévy measure defined in \eqref{Levymeasureparticles}. Since $\gamma > \alpha -1$, we can apply the standard Itô formula for this function and the $(\R^{d})^N$-valued process $(\bm{X}^N_t)_t$. Noticing that $t\in[0,T] \mapsto U(t,\mu_t)$ is constant by the definition of $U$ and the well-posedness of the McKean-Vlasov SDE, we obtain that for all $t \in [0,T)$ \begin{align}\label{thmPCproof1}
	\notag&U(t,\muu^N_t) - U(t,\mu_t) - \left(U(0,\muu^N_0) - U(0,\mu_0)\right) \\\notag&= \int_0^t \partial_t U(s,\muu^N_s) \, ds \\ \notag&\quad+ \frac1N \sum_{i=1}^N \int_0^t \partial_v \del U (s,\muu^N_s)(X^{i,N}_s)\cdot b(s,X^{i,N}_s, \muu^N_s) \, ds \\ &\quad+ \int_0^t \int_{(\R^d)^N} \left[U(s,\muu^N_{\bm{X}^N_{s^-}+ \bm{z}}) - U(s,\muu^N_{\bm{X}^N_{s^-}}) - \partial_{\bm{x}}U(s,\muu^N_{\bm{X}^N_{s^-}})\cdot \bm{z}\right] \, d\bm{\nu}^N(\bm{z}) \, ds \\ \notag&\quad+ \int_0^t \int_{\{|\bm{z}| \geq 1\}} \left[U(s,\muu^N_{\bm{X}^N_{s^-}+ \bm{z}}) - U(s,\muu^N_{\bm{X}^N_{s^-}}) \right]\, \bm{\NNN}^N(ds,d\bm{z}) \\ \notag&\quad+ \int_0^t \int_{\{|\bm{z}|<1\}} \left[U(s,\muu^N_{\bm{X}^N_{s^-}+ \bm{z}}) - U(s,\muu^N_{\bm{X}^N_{s^-}})\right] \, \bm{\NNN}^N(ds,d\bm{z})\\ \notag&=: I_1 + I_2 + I_3 + I_4 + I_5.\end{align}
Note that the term $$\int_0^t \int_{\{|\bm{z}| \geq 1\}} \left[U(s,\muu^N_{\bm{X}^N_{s^-}+ \bm{z}}) - U(s,\muu^N_{\bm{X}^N_{s^-}}) \right]\, d\bm{\nu}^N(\bm{z}) \, ds$$ is well-defined. Indeed, if we set, for $h \in \R^d$, $\bm{\tilde{h_j}} := ( 0, \dots, 0,h,0, \dots, 0) \in (\R^d)^N,$ where $h$ appears in the $j$-th coordinate, one has using \eqref{gradient_lin_der_bound_sol_PDE_POC}  \begin{align*}
	&\int_0^t \int_{\{|\bm{z}| \geq 1\}} \left\vert U(s,\muu^N_{\bm{X}^N_{s^-}+ \bm{z}}) - U(s,\muu^N_{\bm{X}^N_{s^-}}) \right\vert\, d\bm{\nu}^N(\bm{z}) \, ds \\ &= \sum_{i=1}^N \int_0^t \int_{B_1^c} \left\vert U(s,\muu^N_{\bm{X}^N_{s^-}+ \bm{\tilde{z_i}}}) -U(s,\muu^N_{\bm{X}^N_{s^-}}) \right\vert\, d\nu(z) \, ds \\& \leq \frac{1}{N} \sum_{i=1}^N \int_0^t \int_{B_1^c} \int_0^1\left\vert \partial_v\del U(s,m^i_{s,z,w})(X^{i,N}_{s^-} + h z) \right\vert \vert z \vert\,dw \, d\nu(z) \, ds \\ &\leq C \int_0^t (T-s)^{\frac{\delta -1}{\alpha}} \, ds \int_{B_1^c} |z| \, d\nu(z),
\end{align*}
where $m^i_{s,z,w} := w \muu^N_{\bm{X}^N_{s^-}+ \bm{\tilde{z_i}}} + (1-w)\muu^N_{\bm{X}^N_{N,s^-}}$. We conclude since $\alpha \in (1,2)$.\\

 By \eqref{Levymeasureparticles}, we write \begin{align*}
	I_3 &= \sum_{i=1}^N \int_0^t \int_{(\R^d)^N} \left[U(s,\muu^N_{\bm{X}^N_{s^-}+ \bm{\tilde{z_i}}}) - U(s,\muu^N_{\bm{X}^N_{s^-}})  \right. \\ & \left.\hspace{5cm}-\frac1N  \partial_v \del U (s,\muu^N_{\bm{X}^N_{s^-}})(X^{i,N}_{s^-})\cdot z\right]\, d\nu(z) \, ds \\ &= \frac1N\sum_{i=1}^N \int_0^t \int_{\R^d} \int_0^1 \left[\del U(s, m^i_{s,z,w})(X^{i,N}_{s^{-}} + z)\right. \\ & \left.\hspace{3cm}- \del U(s,m^i_{s,z,w})(X^{i,N}_{s^{-}})  - \partial_v \del U (s,\muu^N_{\bm{X}^N_{s^-}})(X^{i,N}_{s^-})\cdot z\right]\, dw\, d\nu(z) \, ds,
\end{align*}
where $m^i_{s,z,w} = w \muu^N_{\bm{X}^N_{s^-}+ \bm{\tilde{z_i}}} + (1-w)\muu^N_{\bm{X}^N_{s^-}}.$ In order to make appear the backward Kolmogorov PDE \eqref{EDP_Kolmo_backward}, we decompose $I_3$ in the following way

\begin{align*}
	I_3 &=  \int_0^t \int_{\R^d} \int_{\R^d} \left[\del U(s, \muu^N_{\bm{X}^N_{s^-}})(x + z) - \del U(s,\muu^N_{\bm{X}^N_{s^-}})(x) -  \partial_v \del U (s,\muu^N_{\bm{X}^N_{s^-}})(x)\cdot z \right]\, d\nu(z) \,d\muu^N_{\bm{X}^N_{s^-}}(x) \, ds \\ &\quad + \frac1N\sum_{i=1}^N \int_0^t \int_{\R^d} \int_0^1 \left[\del U(s, m^i_{s,z,w})(X^{i,N}_{s^{-}} + z) -  \del U(s,m^i_{s,z,w})(X^{i,N}_{s^{-}}) \right. \\ & \hspace{5cm}+\left. \del U(s,\muu^N_{\bm{X}^N_{s^-}})(X^{i,N}_{s^{-}}) - \del U(s, \muu^N_{\bm{X}^N_{s^-}})(X^{i,N}_{s^{-}} + z) \right]\,dw\, d\nu(z) \, ds \\ &=: I_{3,A} + I_{3,B}.
\end{align*}

Since $(\bm{X}^N_{s})_{s \in [0,T]}$ is càd-làg, we deduce that almost surely for almost all $s \in [0,t]$ we have $ \muu^N_{\bm{X}^N_{s^{-}}} = \muu^N_{\bm{X}^N_{s}} = \muu^N_{s}.$ Owing to the backward Kolmogorov PDE \eqref{EDP_Kolmo_backward} in Theorem \ref{Thm_EDP}, one has \begin{align*}
	I_1 + I_2 + I_{3,A} &=  \int_0^t \partial_s U(s,\muu^N_{N,s}) + \mathscr{L}_s U(s,\cdot)(\muu^N_{s}) \, ds \\ &=0.
\end{align*}
Thus, we obtain the following decomposition, for any $t \in [0,T)$,\begin{align}\label{ThmPCproofeq1}
	\notag&U(t,\muu^N_{t}) - U(t,\mu_{t}) - \left(U(0,\muu^N_{0}) - U(0,\mu_{0})\right) \\ \notag&=  \frac1N\sum_{i=1}^N \int_0^t \int_{\R^d} \int_0^1 \left[\del U(s, m^i_{s,z,w})(X^{i,N}_{s^{-}} + z) -  \del U(s,m^i_{s,z,w})(X^{i,N}_{s^{-}}) \right. \\ \notag& \hspace{4cm}+\left. \del U(s,\muu^N_{\bm{X}^N_{s^-}})(X^{i,N}_{s^{-}}) - \del U(s, \muu^N_{\bm{X}^N_{s^-}})(X^{i,N}_{s^{-}} + z) \right]\,dw\, d\nu(z) \, ds  \\ &\quad + \int_0^t \int_{\{|\bm{z}| \geq 1\}} \left[U(s,\muu^N_{\bm{X}^N_{s^-}+ \bm{z}}) - U(s,\muu^N_{\bm{X}^N_{s^-}}) \right]\, \bm{\NNN}^N(ds,d\bm{z}) \\ \notag&\quad + \int_0^t \int_{\{|\bm{z}| < 1\}} \left[U(s,\muu^N_{\bm{X}^N_{s^-}+ \bm{z}}) - U(s,\muu^N_{\bm{X}^N_{s^-}})\right] \, \bm{\NNN}^N(ds,d\bm{z}) \\\notag &=  I_{3,B}  + I_4 + I_5.\end{align}

It follows that for all $t \in [0,T)$ \begin{equation}\label{thmPCproof3}
	\E|U(t,\muu^N_{t}) - U(t,\mu_{t})| \leq   \E |U(0,\muu^N_0) - U(0,\mu_0)| +  \E (|I_{3,B}| + |I_4| + |I_5|).
\end{equation}

We now treat each term separately. For $I_{3,B}$, one has 

\begin{align*}
	\E|I_{3,B}| &\leq \frac1N\sum_{i=1}^N \int_0^t \int_{B_1} \int_0^1 \left\vert\del U(s, m^i_{s,z,w})(X^{i,N}_{s^{-}} + z) -  \del U(s,m^i_{s,z,w})(X^{i,N}_{s^{-}}) \right. \\ \notag& \hspace{4cm}+\left. \del U(s,\muu^N_{\bm{X}^N_{s^-}})(X^{i,N}_{s^{-}}) - \del U(s, \muu^N_{\bm{X}^N_{s^-}})(X^{i,N}_{s^{-}} + z) \right\vert\,dw\, d\nu(z) \, ds \\ &\quad + \frac1N\sum_{i=1}^N \int_0^t \int_{B_1^c} \int_0^1 \left\vert \del U(s, m^i_{s,z,w})(X^{i,N}_{s^{-}} + z) -  \del U(s,m^i_{s,z,w})(X^{i,N}_{s^{-}}) \right. \\ \notag& \hspace{4cm}+\left. \del U(s,\muu^N_{\bm{X}^N_{s^-}})(X^{i,N}_{s^{-}}) - \del U(s, \muu^N_{\bm{X}^N_{s^-}})(X^{i,N}_{s^{-}} + z) \right\vert \,dw\, d\nu(z) \, ds \\ &\leq \frac1N\sum_{i=1}^N \int_0^t \int_{B_1} \int_0^1 \int_0^1 \left\vert \pa_v\del U(s, m^i_{s,z,w})(X^{i,N}_{s^{-}} + \lambda z) \right. \\ \notag& \hspace{6cm}\left. -\pa_v\del U(s, \muu^N_{\bm{X}^N_{s^-}})(X^{i,N}_{s^{-}} + \lambda z)\right\vert |z| \, d\lambda \,dw\, d\nu(z) \, ds \\ &\quad + \frac1N\sum_{i=1}^N \int_0^t \int_{B_1^c} \int_0^1 \left\vert \del U(s, m^i_{s,z,w})(X^{i,N}_{s^{-}} + z) - \del U(s, \muu^N_{\bm{X}^N_{s^-}})(X^{i,N}_{s^{-}} + z)  \right\vert \\ \notag& \hspace{4cm}+\left\vert \del U(s,\muu^N_{\bm{X}^N_{s^-}})(X^{i,N}_{s^{-}}) - \del U(s,m^i_{s,z,w})(X^{i,N}_{s^{-}}) \right\vert \,dw\, d\nu(z) \, ds.
\end{align*}

Recall that  $m^i_{s,z,w} = w \muu^N_{\bm{X}^N_{s^-}+ \bm{\tilde{z_i}}} + (1-w)\muu^N_{\bm{X}^N_{s^-}}.$ Moreover, one has 

\begin{align*}
	W_1 (m^i_{s,z,w}, \muu^N_{\bm{X}^N_{s^-}}) &\leq  	W_1 (\muu^N_{\bm{X}^N_{s^-}+ \bm{\tilde{z_i}}}, \muu^N_{\bm{X}^N_{s^-}})\\  &\leq \frac{|z|}{N}.
\end{align*}
Then, by \eqref{gradient_lin_der_holder_sol_PDE_POC} with $\gamma \in (0,1]\cap(\alpha -1, (\delta + \alpha -1)\wedge(2\alpha -2)\wedge(\eta + \alpha -1))$ and \eqref{lin_der_holder_sol_PDE_POC}, we deduce that 

\begin{align*}
	\E|I_{3,B}|  &\leq \frac{C}{N}\sum_{i=1}^N \int_0^t \int_{B_1} (T-s)^{\frac{\delta -1 - \gamma}{\alpha}} \frac{|z|^{1 + \gamma}}{N^\gamma} \, d\nu(z) \, ds \\ &\quad + \frac{C}{N}\sum_{i=1}^N \int_0^t \int_{B_1^c} (T-s)^{\frac{\delta -1}{\alpha}} \frac{|z|}{N} d\nu(z) \, ds.
\end{align*}
Since $\gamma < \delta +\alpha -1$, we have $\frac{\delta -1-\gamma}{\alpha} >-1$ and therefore the map $s \in [0,T) \mapsto (T-s)^{\frac{\delta -1 - \gamma}{\alpha}}$ is integrable and $z \in B_1 \mapsto |z|^{1 + \gamma} \in L^1(B_1,\nu)$. It yields \begin{align}\label{thmPCproof5}
	\E|I_{3,B}| &\leq \frac{C}{N^\gamma}.
\end{align}\\

Let us now focus on $I_4$. Recall that $m^i_{s,z,w} = w \muu^N_{\bm{X}^N_{s^-}+ \bm{\tilde{z_i}}} + (1-w)\muu^N_{\bm{X}^N_{s^-}}.$ We introduce $\Gamma \in (1,2]$ such that $\Gamma > \alpha$ and $\delta > 1- \alpha / \Gamma$, which is possible since $\delta >0$. Thanks to BDG's inequalities and the subadditivity of the map $|\cdot|^{\Gamma/2}$, there exists a constant $C>0$ independent of $N$ such that we have for all $t \in (0,T]$ \begin{align}\label{thmPCproof6}
	\notag\E |I_4|  & \leq  \left(\E \left\vert\int_0^t \int_{\{|\bm{z}| < 1\}} \left\vert U(s,\muu^N_{\bm{X}^N_{s^-}+ \bm{z}}) - U(s,\muu^N_{\bm{X}^N_{s^-}}) \right]\, \bm{\NNN}^N(ds,d\bm{z}) \right\vert^\Gamma \right)^{\frac{1}{\Gamma}} \\\notag & \leq C\left(\E \left[ \int_0^t \int_{\{ |\bm{z}| <1\}} \left\vert U(s,\muu^N_{\bm{X}^N_{s^-}+ \bm{z}}) - U(s,\muu^N_{\bm{X}^N_{s^-}}) \right\vert^2\, \bm{\NN}^N(ds,d\bm{z})\right]^{\frac{\Gamma}{2}} \right)^{\frac{1}{\Gamma}} \\& \leq C\left(\E \int_0^t \int_{\{|\bm{z}| <1\}} \left\vert U(s,\muu^N_{\bm{X}^N_{s^-}+ \bm{z}}) - U(s,\muu^N_{\bm{X}^N_{s^-}}) \right\vert^\Gamma\, \bm{\NN}^N(ds,d\bm{z}) \right)^{\frac{1}{\Gamma}}\\ &= C\left( \sum_{i=1}^N\E \int_0^t \int_{B_1} \left\vert U(s,\muu^N_{\bm{X}^N_{s^-}+ \bm{\tilde{z_i}}}) - U(s,\muu^N_{\bm{X}^N_{s^-}}) \right\vert^\Gamma\, d\nu(z) \, ds \right)^{\frac{1}{\Gamma}} \\ \notag&= C \left( \sum_{i=1}^N\E \int_0^t \int_{B_1} \left\vert \frac1N \int_{[0,1]^2}\partial_v\del U(s,m^i_{s,z,w})(X^{i,N}_{s^-}+ hz)\cdot z\,  dh \,dw \right\vert^\Gamma\, d\nu(z) \, ds \right)^{\frac{1}{\Gamma}} \\ \notag&\leq C \left( \sum_{i=1}^N\E \int_0^t \int_{B_1} \frac{1}{N^\Gamma}(T-s)^{\frac{\Gamma (\delta -1)}{\alpha}}  |z|^\Gamma\, d\nu(z) \, ds \right)^{\frac{1}{\Gamma}}\\ \notag&\leq \frac{C}{N^{1-\frac{1}{\Gamma}}}.
\end{align}

Indeed, the time integral is finite since $\Gamma (\delta - 1)/ \alpha >-1$ since $\delta > 1 - \alpha /\Gamma$ by choice of $\Gamma$ and $ z \in B_1 \mapsto |z|^\Gamma$ is integrable with respect to $\nu$ since $\gamma > \alpha$. Finally for $I_5,$ BDG's inequalities and the fact that $1<\beta< 2$ yields, for all $t \in [0,T)$, \begin{align}\label{ThmPCproofeq7}
	\notag\E |I_5| & \leq  \left(\E \left(\int_0^t \int_{\{|\bm{z}| >1\}} \left\vert U(s,\muu^N_{\bm{X}^N_{s^-}+ \bm{z}}) - U(s,\muu^N_{\bm{X}^N_{s^-}}) \right]\, \bm{\NNN}^N(ds,d\bm{z}) \right)^\beta \right)^{\frac{1}{\beta}}\\\notag & \leq C\left(\E \left[ \int_0^t \int_{\{|\bm{z}| >1\}} \left\vert U(s,\muu^N_{\bm{X}^N_{s^-}+ \bm{z}}) - U(s,\muu^N_{\bm{X}^N_{s^-}}) \right\vert^2\, \bm{\NN}^N(ds,d\bm{z})\right]^{\frac{\beta}{2}} \right)^{\frac{1}{\beta}} \\& \leq C\left(\E \int_0^t \int_{\{|\bm{z}| >1\}} \left\vert U(s,\muu^N_{\bm{X}^N_{s^-}+ \bm{z}}) - U(s,\muu^N_{\bm{X}^N_{s^-}}) \right\vert^\beta\, \bm{\NN}^N(ds,d\bm{z}) \right)^{\frac{1}{\beta}} \\ \notag&= C\left( \sum_{i=1}^N\E \int_0^t \int_{B_1^c} \left\vert \frac1N \int_{[0,1]^2}\partial_v\del U(s,m^i_{s,z,w})(X^{i,N}_{s^-}+ hz)\cdot z \,dh \, dw \right\vert^\beta\, d\nu(z) \, ds \right)^{\frac1\beta} \\ \notag&\leq  C\left( \sum_{i=1}^N\E \int_0^t \int_{B_1^c} \left\vert \frac1N (T-s)^{\frac{\delta -1}{\alpha}} \right\vert^\beta |z|^\beta\, d\nu(z) \, ds \right)^{\frac{1}{\beta}}\\ \notag&\leq \frac{C}{N^{1-\frac{1}{\beta}}}.
\end{align}

	Indeed, the time integral is finite since $1 < \beta < \alpha$ and thus $\beta(\delta -1)/\alpha > -1$, and the map $z \in B_1^c \mapsto |z|^\beta$ is integrable with respect to $\nu$ since $\beta < \alpha$. As a consequence of \eqref{ThmPCproofeq1}, \eqref{thmPCproof5}, \eqref{thmPCproof6},  \eqref{ThmPCproofeq7}, and the fact that $\gamma > 1 -\frac{1}{\beta}$ since $1<\beta < \alpha < 1 + \gamma$ and $1 - \frac{1}{\Gamma} > 1 - \frac{1}{\beta}$ since $\beta < \Gamma$, we have, for a constant $C>0$, that for all $ N\geq 1$ and $t \in [0,T)$ \begin{equation}\label{ThmPCproofeq8}
	\E|U(t,\muu^N_{t}) - U(t,\mu_{t})| \leq   \E |U(0,\muu^N_0) - U(0,\mu_0)| +\frac{C}{N^{1-\frac{1}{\beta}}}.
\end{equation}

It follows from \eqref{gradient_lin_der_bound_sol_PDE_POC} that the map $\del U(0,\mu)(\cdot)$ is Lipschitz uniformly with respect to $\mu \in \PP$. The Kantorovich-Rubinstein theorem ensures that \[\E |U(0,\muu^N_0) - U(0,\mu_0)| \leq CT^{\frac{\delta-1}{\alpha}} \E W_1(\muu^N_0,\mu_0).\]
 By Fatou's lemma, the continuity of $U$ on $[0,T] \times \PP$ and since $\muu^N_{T^-} = \muu^N_T$ almost surely, we can let $t$ tend to $T$ in \eqref{ThmPCproofeq8}, which concludes the proof of \eqref{eq1_thm_poc}.\\

\noindent\textbf{Proof of \eqref{eq3_thm_poc}.} Coming back to \eqref{ThmPCproofeq1} and using \eqref{thmPCproof6} and \eqref{ThmPCproofeq7}, we see that $I_4$ an $I_5$ are centered random variables because of the martingale property of compensated Poisson random integrals. We thus obtain that for any $t \in [0,T)$ $$ |\E (U(t,\muu^N_{t}) - U(t,\mu_{t}))| \leq |\E (U(0,\muu^N_0) - U(0,\mu_0))| + \E|I_{3,B}|.$$ The control of  $\E|I_{3,B}|$ has already been done in \eqref{thmPCproof5}. It follows that there exists a positive constant $C$ such that for all $t \in [0,T)$ \begin{equation}\label{ThmPCproofeq9}|\E(U(t,\muu^N_t) - U(t,\mu_t))| \leq  C\, |\E (U(0,\muu^N_0) - U(0,\mu_0))| + \frac{C}{N^{\gamma}}.\end{equation}
 Reasoning as in \cite{frikha2021backward}, after $(5.24)$, we can write using the exchangeability in law of the initial data $(\xi^i)_{i \in \{1,\dots,N\} }$ \begin{equation*}
\E(U(0,\muu^N_0) - U(0,\mu_0)) = \int_0^1 \E \left(  \del U(0,\tilde{\mu}^{N,\lambda_1}_0)(\tilde{\xi}) -  \del U(0,\mu^{N,\lambda_1}_0)(\tilde{\xi}) \right) d\lambda_1,
\end{equation*}
where $\mu_0^{N,\lambda_1}:= \lambda_1 \muu^N_0 + (1-\lambda_1)\mu_0$, $\tilde{\mu}_0^{N,\lambda_1}:= \lambda_1 \tilde{\mu}^N_0 + (1-\lambda_1)\mu_0$, $\tilde{\mu}_0^N := \muu^N_0 + \frac{1}{N} (\delta_{\tilde{\xi}} - \delta_{\xi^1})$, $(\xi^i)_i$ being an i.i.d.\ sequence of random variable with common distribution $\mu_0$ and $\tilde{\xi}$ a random variable independent of $(\xi^i)_i$ with distribution $\mu_0$. This ensures, using \eqref{lin_der_holder_sol_PDE_POC}, that there exists a positive constant $C$ such that for all $\phi \in \CC^{(2,\delta)}(\PP)$ and $N\geq 1$  \begin{align*}|\E (U(0,\muu^N_0) - U(0,\mu_0))|& \leq \int_0^1 CT^{\frac{\delta-1}{\alpha}} \E W_1(\tilde{\mu}^{N,\lambda_1}_0, \muu^{N,\lambda_1}_0) \, d\lambda_1 \\ & \leq \int_0^1 CT^{\frac{\delta-1}{\alpha}} \E W_1(\tilde{\mu}^N_0,\muu^N_0) \, d\lambda_1 \\ &\leq \int_0^1 \frac{C}{N}T^{\frac{\delta-1}{\alpha}} \E|\tilde{\xi}-\xi_1| \, d\lambda_1\\&\leq \frac{C}{N}T^{\frac{\delta-1}{\alpha}},\end{align*} 
the constant $C$ depending on $M_1(\mu_0)$ in \eqref{eq3_thm_poc}. 
We conclude by letting $t$ tend to $T$ in \eqref{ThmPCproofeq9} that \eqref{eq3_thm_poc} holds true.

\end{proof}

\section{Quantitative approximation of the distribution of one particle by the limiting McKean-Vlasov process at the level of densities}\label{section_poc_density}

This section is devoted to the proof of Theorem \ref{Thm_poc_density}. We start to show that the density bound \eqref{density_estimate_one_particle} holds. Let us fix $t \in (0,T]$. As in the proof of Theorem \ref{Thm_POC}, we easily see using Theorem \ref{Thm_density_McKean} that we can apply Itô's formula for the function $(s,\bm{x}) \in [0,t) \times (\R^d)^N \mapsto p(\muu^N_{\bm{x}},s,t,y) = \frac{1}{N} \sum_{k=1}^N p(\muu^N_{\bm{x}},s,t,x_k,y).$ Keeping the same notations as in the proof of Theorem \ref{Thm_POC}, we get that 

\begin{align}\label{proof_poc_density_dec}
	\notag p(\muu^N_s,s,t,y) - p(\muu^N_0,0,t,y)  &= \int_0^s \pa_r p(\muu^N_r,r,t,y)\,dr \\ \notag&\quad + \frac1N \sum_{i=1}^N \int_0^s  \pa_v \del p(\muu^N_r,r,t,y)(X^{i,N}_{r^-}) \cdot b(r,X^{i,N}_{r^-},\muu^N_r)  \, dr\\  &\quad+\sum_{i=1}^N \int_0^s \int_{(\R^d)^N} \left[p(\muu^N_{\bm{X}^N_{r^-}+ \bm{\tilde{z_i}}},r,t,y) -p(\muu^N_{\bm{X}^N_{r^-}},r,t,y)  \right. \\ \notag & \left.\hspace{5cm}-\frac1N  \partial_v \del p (\muu^N_{\bm{X}^N_{r^-}},r,t,y)(X^{i,N}_{r^-})\cdot z\right]\, d\nu(z) \, dr \\ \notag &\quad+ \int_0^s \int_{(\R^d)^N} p(\muu^N_{\bm{X}^N_{r^-}+ \bm{z}},r,t,y) - p(\muu^N_{\bm{X}^N_{r^-}},r,t,y)\, \bm{\NNN}^N(dr,d\bm{z}) \\ \notag &=: I_1 + I_2 + I_3 + I_4.
\end{align}

First, notice that reasoning as in the proof of Theorem \ref{Thm_POC}, by using \eqref{linear_derivative_holder_measure_thm} for the big jumps, and \eqref{gradient_linear_der_bound_thm} for the small jumps, we prove that $I_4$ is a true martingale and that \begin{equation}\label{proof_poc_density_mtgle}
	\E I_4 =0.
\end{equation}
By definition of the linear derivative, we can write \begin{align*}
	I_3  &= \frac1N\sum_{i=1}^N \int_0^s \int_{\R^d} \int_0^1 \left[\del p( m^i_{r,z,w},r,t,y)(X^{i,N}_{r^{-}} + z)\right. \\ & \left.\hspace{3cm}- \del p(m^i_{r,z,w},r,t,y)(X^{i,N}_{r^{-}})  - \partial_v \del p (\muu^N_{\bm{X}^N_{r^-}},r,t,y)(X^{i,N}_{r^-})\cdot z\right]\, dw\, d\nu(z) \, dr,
\end{align*}
where $m^i_{r,z,w} := w \muu^N_{\bm{X}^N_{r^-}+ \bm{\tilde{z_i}}} + (1-w)\muu^N_{\bm{X}^N_{r^-}}.$ We can now decompose $I_3$ in the following way

\begin{align*}
	I_3 &=  \int_0^s \int_{\R^d} \int_{\R^d} \left[\del p(\muu^N_{\bm{X}^N_{r^-}},r,t,y)(x + z) - \del p(\muu^N_{\bm{X}^N_{r^-}},r,t,y)(x) \right. \\ & \left. \hspace{7cm}-  \partial_v \del p (\muu^N_{\bm{X}^N_{r^-}},r,t,y)(x)\cdot z \right]\, d\nu(z) \,d\muu^N_{\bm{X}^N_{r^-}}(x) \, dr \\ &\quad + \frac1N\sum_{i=1}^N \int_0^s \int_{B_1^c} \int_0^1 \left[\del p( m^i_{r,z,w},r,t,y)(X^{i,N}_{r^{-}} + z) - \del p( \muu^N_{\bm{X}^N_{r^-}},r,t,y)(X^{i,N}_{r^{-}} + z)  \right] \\ & \hspace{5cm}-\left[  \del p(m^i_{r,z,w},r,t,y)(X^{i,N}_{r^{-}}) - \del p(\muu^N_{\bm{X}^N_{r^-}},r,t,y)(X^{i,N}_{r^{-}})  \right]\,dw\, d\nu(z) \, dr \\ &\quad + \frac1N\sum_{i=1}^N \int_0^s \int_{B_1} \int_0^1 \left[\del p( m^i_{r,z,w},r,t,y)(X^{i,N}_{r^{-}} + z) -  \del p(m^i_{r,z,w},r,t,y)(X^{i,N}_{r^{-}}) \right. \\ & \hspace{5cm}+\left. \del p(\muu^N_{\bm{X}^N_{r^-}},r,t,y)(X^{i,N}_{r^{-}}) - \del p( \muu^N_{\bm{X}^N_{r^-}},r,t,y)(X^{i,N}_{r^{-}} + z) \right]\,dw\, d\nu(z) \, dr  \\ &=: I_{3,A} + I_{3,B} + I_{3,C}.
\end{align*}

By taking the expectation in \eqref{proof_poc_density_dec}, we deduce thanks to \eqref{proof_poc_density_mtgle} and the backward Kolmogorov PDE \eqref{EDPKolmogorovfundamental_not_dec} satisfied by $p$ (see Theorem \ref{Thm_EDP_densité_McKV}) that \begin{equation}\label{proof_poc_density_EDP} \E p(\muu^N_s,s,t,y) = \E p(\muu^N_0,0,t,y) + \E I_{3,B} + \E I_{3,C}.\end{equation}
We first aim at controlling $I_{3,B}$. To do this, we write for $v \in \R^d$ \begin{align*}
 \del p(m^i_{r,z,w},r,t,y)(v)  &-	\del p(\muu^N_{\bm{X}^N_{r^-}},r,t,y)(v) \\& = p(m^i_{r,z,w},r,t,v,y) -  p(\muu^N_{\bm{X}^N_{r^-}},r,t,v,y) \\ &\quad + \int_{\R^d} \del p(m^i_{r,z,w},r,t,x,y)(v)\, w d(\muu^N_{\bm{X}^N_{r^-} + \bm{\tilde{z}_i}} - \muu^N_{\bm{X}^N_{r^-}})(x) \\ &\quad+ \int_{\R^d} \del p(m^i_{r,z,w},r,t,x,y)(v) - \del p(\muu^N_{\bm{X}^N_{r^-}},r,t,x,y)(v) \, d\muu^{N}_{\bm{X}^N_{r^-}}(x).
\end{align*}
Using \eqref{gradient_density_holder_measure_thm}, \eqref{linear_der_holder_x_thm} and \eqref{linear_derivative_holder_measure_thm} with $\gamma =1$, we obtain 

 \begin{multline*}
	\left\vert\del p(m^i_{r,z,w},r,t,y)(v)  -	\del p(\muu^N_{\bm{X}^N_{r^-}},r,t,y)(v)\right\vert \\ \leq  \frac{C}{N}(t-r)^{-\zeta} |z| \left[\rho^0(t-r,y-v) + \rho^0(t-r,y-X^{i,N}_{r^-} - z) +   \rho^0(t-r,y-X^{i,N}_{r^-}) + \frac1N \sum_{k=1}^N  \rho^0(t-r,y-X^{k,N}_{r^-})\right],
\end{multline*}
where we recall that $\zeta :=-\left( 1 - \frac{2+\gamma}{\alpha}\right).$ It follows that \begin{equation}\label{proof_poc_density_I3B1}
	\E |I_{3,B}| \leq \frac{C}{N} \int_0^s \int_{B_1^c} (t-r)^{-\zeta} |z| \E ( \rho^0(t-r,y-X^{i,N}_{r^-}) + \rho^0(t-r,y-X^{i,N}_{r^-} - z)) \, d\nu(z) \, dr.
\end{equation}
Let us now focus on $I_{3,C}$ which can be written in the following way 

\begin{align*}
	I_{3,C} &= \frac1N\sum_{i=1}^N \int_0^s \int_{B_1} \int_0^1\int_0^1 \left[\pa_v\del p( m^i_{r,z,w},r,t,y)(X^{i,N}_{r^{-}} + \lambda z)\right. \\ &\left. \hspace{5cm}-  \pa_v\del p(\muu^N_{\bm{X}^N_{r^-}},r,t,y)(X^{i,N}_{r^-} + \lambda z)\right] \cdot z \, d\lambda \,dw\, d\nu(z) \, dr.
\end{align*}
As previously, we need to control for $v \in \R^d$ \begin{align*}
	\pa_v\del p(m^i_{r,z,w},r,t,y)(v)  &-	\pa_v\del p(\muu^N_{\bm{X}^N_{r^-}},r,t,y)(v) \\& = \pa_x p(m^i_{r,z,w},r,t,v,y) -  \pa_x p(\muu^N_{\bm{X}^N_{r^-}},r,t,v,y) \\ &\quad + \int_{\R^d} \pa_v\del p(m^i_{r,z,w},r,t,x,y)(v)\, w d(\muu^N_{\bm{X}^N_{r^-} + \bm{\tilde{z}_i}} - \muu^N_{\bm{X}^N_{r^-}})(x) \\ &\quad+ \int_{\R^d} \pa_v\del p(m^i_{r,z,w},r,t,x,y)(v) - \pa_v \del p(\muu^N_{\bm{X}^N_{r^-}},r,t,x,y)(v) \, d\muu^{N}_{\bm{X}^N_{r^-}}(x).
\end{align*}
Using \eqref{gradient_density_holder_measure_thm}, \eqref{gradient_linear_der_holder_x_thm} and \eqref{gradient_linear_derivative_holder_measure_thm} with $\gamma$, we obtain 

\begin{multline*}
	\left\vert\pa_v\del p(m^i_{r,z,w},r,t,y)(v)  -	\pa_v\del p(\muu^N_{\bm{X}^N_{r^-}},r,t,y)(v)\right\vert  \leq  \frac{C}{N^\gamma}(t-r)^{1- \frac{2+\gamma}{\alpha}} |z|^\gamma\\  \left[\rho^0(t-r,y-v) + \rho^0(t-r,y-X^{i,N}_{r^-} - z)   +   \rho^0(t-r,y-X^{i,N}_{r^-}) + \frac1N \sum_{k=1}^N  \rho^0(t-r,y-X^{k,N}_{r^-})\right].
\end{multline*}
This yields, by definition of $\zeta$, \begin{equation}\label{proof_poc_density_I3B2}
	\E |I_{3,C}| \leq \frac{C}{N^\gamma} \int_0^s \int_{B_1} (t-r)^{-\zeta} |z|^{1+\gamma} \E ( \rho^0(t-r,y-X^{i,N}_{r^-}) + \rho^0(t-r,y-X^{i,N}_{r^-} - z)) \, d\nu(z) \, dr.
\end{equation}

Gathering \eqref{proof_poc_density_EDP}, \eqref{proof_poc_density_I3B1} and \eqref{proof_poc_density_I3B2}, we deduce that for any $s \in [0,t)$, $\mu_0 \in \PP$, $y \in \R^d$ and $N \geq 1$

\begin{align}\label{proof_poc_density_ineq_ito}
	\E p(\muu^N_s,s,t,y) &\leq \E p(\muu^N_0,0,t,y) \\ &\notag\quad+ \frac{C}{N^\gamma} \int_0^s \int_{\R^d} (t-r)^{-\zeta} f(z) \E ( \rho^0(t-r,y-X^{i,N}_{r^-}) + \rho^0(t-r,y-X^{i,N}_{r^-} - z)) \, d\nu(z) \, dr.
\end{align}
By \eqref{density_bound_thm}, one has \begin{align}\label{proof_poc_density_bound_initial}
\notag	\E p(\muu^N_0,0,t,y) & = \E\frac{1}{N} \sum_{i=1}^N \left(p(\muu^N_0,0,t,\xi^i,y)\right) \\ &\leq C q_0(\mu_0,0,t,y),
\end{align}
where we recall that $$q_0(\mu_0,0,t,y) := \int_{\R^d} \rho^0(t,y-x) \, d\mu_0(x).$$ 
The same reasoning as in the proof of \cite[Theorem $3.5$]{frikha2021backward} (see $(5.13)$ and below) shows that  \begin{equation}\label{proof_poc_density_limit}
	\lim_{s \to t^-} \E p(\muu^N_s,s,t,y) = p^{1,N}(\mu_0,0,t,y).
\end{equation}
This is justified by the parametrix expansion of the transition density $p$ of the McKean-Vlasov SDE (see \eqref{param_rep_density_thm}). Indeed, it allows to write, for $\mu_0 \in \PP$, $0 \leq s < t \leq T$, $x,y \in \R^d$, $ p(\mu_0,s,t,x,y) = \p (s,t,x,y) + R(\mu_0,s,t,x,y)$ with  $$ \vert R(\mu_0,s,t,x,y)\vert\leq  C(t-s)^{- \frac{1}{\alpha}} \rho^1(t-s,y-x).$$
Using \eqref{proof_poc_density_ineq_ito}, \eqref{proof_poc_density_bound_initial} and \eqref{proof_poc_density_limit}, we deduce that for any $t \in (0,T]$, $\mu_0 \in \PP$, $y \in \R^d$ and $N \geq 1$ \begin{multline*}
p^{1,N}(\mu_0,0,t,y) \leq C q_0(\mu_0,0,t,y)  + \frac{C}{N^\gamma} \int_0^t \int_{\R^d} \int_{\R^d} (t-r)^{-\zeta} f(z)[ \rho^0(t-r,y-w-z) \\ + \rho^0(t-r,y-w)] p^{1,N}(\mu_0,0,r,w) \, dw \, d\nu(z) \, dr.
\end{multline*}
Denoting, for $\mu \in \PP$, $0\leq s < t \leq T$ and $y \in \R^d$ \begin{equation}\label{proof_poc_density_def_g}
	g(\mu,s,t,y) := (t-s)^{-\zeta} \rho^0(t-s,y),
\end{equation}
the previous inequality can be rewritten as \begin{equation*}\label{proof_poc_density_ineq_implicit}
	p^{1,N}(\mu_0,0,t,y) \leq C q_0(\mu_0,0,t,y)  + \frac{C}{N^\gamma}  \int_{\R^d} [p^{1,N} \otimes g (\mu_0,0,t,y-z) + p^{1,N} \otimes g (\mu_0,0,t,y) ] f(z)\, d\nu(z),
\end{equation*}
where we recall that $$ p^{1,N} \otimes g (\mu,0,t,y)  := \int_0^t \int_{\R^d} p^{1,N}(\mu,0,r,w) g(\mu,r,t,y-w)\, dw \, dr.$$
Notice that $g$ yields a time-integrable singularity since $\zeta \in (0,1)$. Then, we can easily prove by induction that for any $M\geq 1$  \begin{align}\label{proof_poc_density_control_iterate}
	p^{1,N}(\mu_0,0,t,y) \notag&\leq C q_0(\mu_0,0,t,y) + \sum_{k=1}^M \frac{C^{k+1}}{N^{k\gamma}} \sum_{I \in P_k} \int_{(\R^d)^k} q_0\otimes g^{k} \left(\mu_0,0,t,y - \sum_{i\in I} z_i\right) \prod_{j=1}^k f(z_j) \, d\nu(z_j) \\ &\quad +  \frac{C^{M+1}}{N^{(M+1)\gamma}} \sum_{I \in P_{M+1}} \int_{(\R^d)^{M+1}} p^{1,N}\otimes g^{M+1} \left(\mu_0,0,t,y - \sum_{i\in I} z_i\right) \prod_{j=1}^{M+1} f(z_j) \, d\nu(z_j), 
\end{align}
where $P_k$ is the set of all subsets of $\{1,\dots,k\}$ and $g^k$ is defined by the recursive relation $g^{k+1} = g^k\otimes g.$ As in the proof of \eqref{H^k_bound}, we prove by induction that for any $\mu \in \PP$, $0 \leq s < t \leq T$ and $y \in \R^d$ \begin{equation} \label{proof_poc_density_g_iterate_bound}\vert g^k(\mu,s,t,y)\vert \leq C^{k-1} (t-s)^{-\zeta + (k-1)(1-\zeta)} \left( \prod_{j=1}^{k-1} \BB (j(1-\zeta),1-\zeta) \right) \rho^0(t-s,y).\end{equation}
Thus, sending $M$ to infinity in \eqref{proof_poc_density_control_iterate} and using \eqref{proof_poc_density_g_iterate_bound}, we deduce that
\begin{align*}
	p^{1,N}(\mu_0,0,t,y) &\leq C q_0(\mu_0,0,t,y) + \sum_{k=1}^{\infty} \frac{C^{k+1}}{N^{k\gamma}} \sum_{I \in P_k} \int_{(\R^d)^k} q_0\otimes g^k \left(\mu_0,0,t,y - \sum_{i\in I} z_i\right) \prod_{j=1}^k f(z_j) \, d\nu(z_j).
\end{align*}
By injecting \eqref{proof_poc_density_g_iterate_bound} in the previous inequality and using the convolution inequality \eqref{convolineqdensityref}, we get 

\begin{multline*}
	p^{1,N}(\mu_0,0,t,y) \leq C q_0(\mu_0,0,t,y) + \sum_{k=1}^{\infty} \frac{C^{k+1}}{N^{k\gamma}}   \left(\prod_{j=1}^{k-1}  \BB (j(1-\zeta),1-\zeta)\right)  \int_0^t (t-r)^{-\zeta + (k-1)(1-\zeta)} \, dr \\  \sum_{I \in P_k} \int_{(\R^d)^k} q_0\left(\mu_0,0,t,y - \sum_{i\in I} z_i\right) \prod_{j=1}^k f(z_j) \, d\nu(z_j).
\end{multline*}
This concludes the proof of \eqref{density_estimate_one_particle}.  \\

Before proving \eqref{eq_thm_poc_density}, we state and show the following Lemma.

\begin{Lemme}\label{Lemme_poc_density_initial}
	For any $\gamma \in [\alpha-1,1]$, there exists a constant $C>0$ such that for any $\mu_0 \in \PP$, $t \in (0,T]$, $y \in \R^d$ and $N \geq 1$ \begin{equation*}
		\left\vert \E p(\muu^N_0,0,t,y) - p(\mu_0,0,t,y) \right\vert \leq \frac{C}{N^\gamma}t^{1 - \frac{1+\gamma}{\alpha}} (1+M_\gamma(\mu_0))\int_{\R^d} (1+|x|^\gamma) \rho^0(t,y-x) \, d\mu_0(x).
	\end{equation*}
\end{Lemme}

\begin{proof}[Proof of Lemma \ref{Lemme_poc_density_initial}]
	Recall that $(\xi^i)_{i\geq1}$ is an i.i.d.\ sequence with common distribution $\mu_0 \in \PP$ and let us introduce $\tilde{\xi}$ a random variable with distribution $\mu_0$ independent of $(\xi^i)_{i\geq1}$. Reasoning as in Lemma $5.1$ in \cite{frikha2021backward} (see $(5.4)$ and below), we obtain that \begin{align*}
		\E p(\muu^N_0,0,t,y) - p(\mu_0,0,t,y) &= \frac1N \int_0^1 \E \left(\del p(\mu^{\lambda_1,N}_0,0,t,\xi^1,y)(\xi^1) - \del p(\mu^{\lambda_1,N}_0,0,t,\xi^1,y)(\tilde{\xi})\right) \, d\lambda_1 \\ &\quad+ \frac{N-1}{N} \int_0^1 \E \left( \del p(\tilde{\mu}^{\lambda_1,N}_0,0,t,\xi^1,y)(\tilde{\xi}) - \del p(\mu^{\lambda_1,N}_0,0,t,\xi^1,y)(\tilde{\xi})\right) \, d\lambda_1 \\ &=: I_1 + I_2,
	\end{align*}
	where $\mu^{\lambda_1,N}_0 := \lambda_1 \muu^N_0 + (1-\lambda_1) \mu_0$ and $\tilde{\mu}^{\lambda_1,N}_0 := \lambda_1 \tilde{\mu}^N_0 + (1-\lambda_1) \mu_0$ with $\tilde{\mu}^N_0 := \muu^N_0 + \frac1N (\delta_{\tilde{\xi}} - \delta_{\xi^2}).$ Using \eqref{linear_der_holder_v_thm}, we obtain that for some constant $C>0$ depending on $\gamma$ \begin{align*}|I_1| &\leq \frac{C}{N} \int_0^1 t^{1- \frac{1+\gamma}{\alpha}} \E (|\xi^1 - \tilde{\xi}|^\gamma \rho^0(t,y-\xi^1)) \, d\lambda_1 \\ & \leq \frac{C}{N} t^{1 - \frac{1+\gamma}{\alpha}} \left[ \E(|\xi^1|^\gamma \rho^0(t,y-\xi^1)) + \E(|\tilde{\xi}|^\gamma \rho^0(t,y-\xi^1)) \right] \\ &\leq \frac{C}{N} t^{1 - \frac{1+\gamma}{\alpha}}( 1 + M_\gamma(\mu_0)) \int_{\R^d} (1+|x|^\gamma) \rho^0(t,y-x) \, d\mu_0(x).\end{align*}
	Using now \eqref{linear_derivative_holder_measure_Picard}, one has \begin{align*}
		|I_2| &\leq C\frac{N-1}{N}t^{1-\frac{1+\gamma}{\alpha}} \int_0^1 \E (W_1^\gamma (\tilde{\mu}^{\lambda_1,N}_0,\mu^{\lambda_1,N}_0)\rho^0(t,y-\xi^1)) \, d\lambda_1.
	\end{align*}
	Since we have \begin{align*}
		W_1 (\tilde{\mu}^{\lambda_1,N}_0,\mu^{\lambda_1,N}_0) &\leq W_1(\tilde{\mu}^{N}_0,\muu^N_0) \\ &\leq \frac{1}{N}|\xi^2 - \tilde{\xi}|,
	\end{align*}
	we deduce that \begin{align*}
		|I_2| &\leq \frac{C}{N^\gamma}t^{1-\frac{1+\gamma}{\alpha}} \E (|\xi^2-\tilde{\xi}|\rho^0(t,y-\xi^1)) \,  \\ &\leq \frac{C}{N^\gamma} t^{1 - \frac{1+\gamma}{\alpha}}( 1 + M_\gamma(\mu_0)) \int_{\R^d} \rho^0(t,y-x) \, d\mu_0(x).
	\end{align*}
	This concludes the proof of Lemma \ref{Lemme_poc_density_initial}.
\end{proof}

Let us now prove the estimate \eqref{eq_thm_poc_density}. We come back to the identity \eqref{proof_poc_density_EDP} subtracting $p(\mu_0,0,t,y)$ from its both sides and use \eqref{proof_poc_density_I3B1} and \eqref{proof_poc_density_I3B2} to obtain that \begin{align}\label{proof_poc_density_ineq}
	\vert \E p(\muu^N_s,s,t,y) - p(\mu_0,0,t,y) \vert &\leq \vert \E p(\muu^N_0,0,t,y) - p(\mu_0,0,t,y)\vert + 	R^N(\mu_0,s,t,y).
\end{align}
where \begin{multline*}
	R^N(\mu_0,s,t,y) := \frac{C}{N^\gamma} \int_0^s \int_{(\R^d)^2} (t-r)^{-\zeta} f(z)\\ \left[ \rho^0(t-r,y-w-z) + \rho^0(t-r,y-w) \right] p^{1,N}(\mu_0,0,r,w) \, dw \, d\nu(z) \, dr 
\end{multline*} 
We can now use \eqref{density_estimate_one_particle} to get that \begin{align*} 
	R^N(\mu_0,s,t,y) &\leq \frac{C^2}{N^\gamma} \int_0^s \int_{(\R^d)^3}(t-r)^{-\zeta} f(z)\left[ \rho^0(t-r,y-w-z) + \rho^0(t-r,y-w) \right] \\ &\hspace{8cm}\rho^0(r,w-x) \, d\mu_0(x) \, dw \, d\nu(z) \, dr
	\\ &\quad + \sum_{k=1}^{\infty} \frac{C^{k+2}}{N^{(k+1)\gamma}} \left(\prod_{j=1}^{k-1} \BB(j(1-\zeta),1-\zeta)\right) \\ & \quad\quad \int_0^s \int_{(\R^d)^3} (t-r)^{-\zeta}\frac{r^{k(1-\zeta)}}{k(1-\zeta)}f(z_{k+1}) \left[ \rho^0(t-r,y-w-z_{k+1}) + \rho^0(t-r,y-w) \right]\\ &\hspace{3cm} \left[ \sum_{I \in P_k} \int_{(\R^d)^k} \rho^0 \left(r,w- x - \sum_{i\in I} z_i \right) \, \prod_{j=1}^{k} f(z_j) \, d\nu(z_j) \right]\, d\mu_0(x) \, dw\,d\nu(z_{k+1}) \, dr. 
\end{align*}
Thanks to the convolution inequality \eqref{convolineqdensityref}, we obtain that

\begin{align*}
 R^N(\mu,s,t,y) &\leq \sum_{k=1}^{\infty} \frac{C^{k+1}}{N^{k\gamma}} t^{k(1-\zeta)} \left( \prod_{j=1}^{k-2} \BB (j(1- \zeta),1-\zeta) \right)\BB(1 + (k-1)(1-\zeta), 1 -\zeta) \\ & \hspace{3cm}\sum_{I \in P_k} \int_{(\R^d)^k} q_0\left(\mu_0,0,t,y - \sum_{i\in I} z_i\right) \prod_{j=1}^k f(z_j) \, d\nu(z_j).
\end{align*}
Using the preceding inequality, Lemma \ref{Lemme_poc_density_initial} with $\gamma' \in [\alpha -1,1]$, and taking the limit $s \to t^-$ in \eqref{proof_poc_density_ineq} thanks to \eqref{proof_poc_density_limit} concludes the proof of \eqref{eq_thm_poc_density}.\\

The estimates  \eqref{eq_thm_poc_density_TV1} and  \eqref{eq_thm_poc_density_TV2} are direct consequences of \eqref{eq_thm_poc_density}. Indeed, we integrate \eqref{eq_thm_poc_density} over $y \in \R^d$ noticing that $$ \int_{\R^d} \int_{(\R^d)^k} \sum_{I \in P_k} q_0\left(\mu_0,t,y-\sum_{i\in I} z_i\right) \prod_{j=1}^k f(z_j)
\, d\nu(z_j) \, dy \leq 2^k\int_{\R^d} \rho^0(1,y) \, dy \left(\int_{\R^d} f(z) \, d\nu(z)\right)^k,$$ 
which is finite since $f \in L^1(\R^d,\nu)$. Using the asymptotics of the Beta function, we conclude the proof of \eqref{eq_thm_poc_density_TV1} (by taking  $\gamma ' = \gamma$ in \eqref{eq_thm_poc_density}) and \eqref{eq_thm_poc_density_TV2} (by taking $\gamma ' = \alpha -1$ in \eqref{eq_thm_poc_density}).

\section{Proof of Proposition \ref{Prop_density_Picard}}\label{section_proof_prop}
This section is dedicated to prove Proposition \ref{Prop_density_Picard}, which is rather long and technical. We proceed by induction on $m \geq 1$. The base case $m=1$ is immediate. It is enough to apply Theorem \ref{Thmdensityparametrix} since $p_1(\mu,s,t,x,y)$ does not depend on $\mu$. For the induction step, we assume that all the estimates of Proposition \ref{Prop_density_Picard}  are satisfied for $p_m$. For the sake of clarity, we denote by $K$ a positive constant depending only on $(d,\alpha,b,T)$ appearing in the induction step and independent of the induction assumption, which may change from line to line and which will determine the choice of the constant $C$ appearing in Proposition \ref{Prop_density_Picard}. Let us introduce a notation used in this section. If $f$ is a function defined on $\PP \times [0,t) \times (\R^d)^2$, we define $\Delta_{\mu_1,\mu_2} f(\cdot,s,x,v):= f(\mu_1,s,x,v)-f(\mu_2,s,x,v)$ for $\mu_1,\mu_2 \in \PP$. The same notations holds with respect to the variables $x$ and $v$. For the time variable, if $s_1,s_2 \in [0,t)$, we define $\Delta_{s_1,s_2} f(\mu,\cdot,x,v): = f(\mu,s_1\vee s_2,x,v) - f(\mu,s_1\wedge s_2,x,v)$.

\subsection{Preparatory technical results}

\begin{Lemme}\label{lemma_technical}
	
	\begin{itemize}
	\item There exists $K>0$ such that for all $0 \leq s \leq r  \leq T$, $\mu \in \PP$, $x,v \in \R^d$, one has \begin{equation}\label{b_linear_derivative_bound_Picard}
		\left\vert \del \left[ b(r,x,[X_r^{s,\mu,(m)}])\right] (v) \right\vert \leq K\left[1 + \sum_{k=1}^m C^k(r-s)^{k\left(1 - \frac{1}{\alpha}\right)} \prod_{j=1}^{k-1} \BB\left(1 + j \left(1 - \frac{1}{\alpha}\right), 1 - \frac{1}{\alpha}\right)\right].
		\end{equation}
	
		\item There exists $K>0$ such that for all $0 \leq s \leq r \leq T$, $\mu \in \PP$, $x,v \in \R^d$, one has \begin{multline}\label{b_gradient_linear_derivative_bound_Picard}
		\left\vert \pa_v  \del \left[ b(r,x,[X_r^{s,\mu,(m)}])\right] (v) \right\vert \\ \leq K(r-s)^{\frac{\eta -1}{\alpha}}\left[1 + \sum_{k=1}^m C^k(r-s)^{k\left(1 + \frac{\eta -1}{\alpha}\right)} \prod_{j=1}^{k-1} \BB\left(1 + \frac{\eta -1}{\alpha} + j \left(1 + \frac{\eta -1}{\alpha}\right), 1 - \frac{1}{\alpha}\right) \right].
	\end{multline}
	
		\item There exists $K>0$ such that for all $0 \leq s \leq r < t \leq T$, $\mu \in \PP$, $x,y,v \in \R^d$, one has \begin{multline}\label{H_linear_derivative_bound_Picard}
		\left\vert \del \H_{m+1}(\mu,s,r,t,x,y) (v) \right\vert \leq K (t-r)^{-\frac{1}{\alpha}}\rho^1(t-r,y-x) \\ \left[1+ \sum_{k=1}^m C^k(r-s)^{k\left(1 - \frac{1}{\alpha}\right)} \prod_{j=1}^{k-1} \BB\left(1 + j \left(1 - \frac{1}{\alpha}\right), 1 - \frac{1}{\alpha}\right)\right].
	\end{multline}

	\item There exists $K>0$ such that for all $0 \leq s \leq r < t \leq T$, $\mu \in \PP$, $x,y,v \in \R^d$, one has\begin{multline}\label{H_gradient_linear_derivative_bound_Picard}
		\left\vert \pa_v \del \H_{m+1}(\mu,s,r,t,x,y) (v) \right\vert  \leq K(r-s)^{\frac{\eta -1}{\alpha}} (t-r)^{-\frac{1}{\alpha}}\rho^1(t-r,y-x)\\ \left[1+ \sum_{k=1}^m C^k(r-s)^{k\left(1 + \frac{\eta -1}{\alpha}\right)} \prod_{j=1}^{k-1} \BB\left(1 + \frac{\eta -1}{\alpha} + j \left(1 + \frac{\eta -1}{\alpha}\right), 1 - \frac{1}{\alpha}\right) \right].\end{multline}

	\end{itemize}

\end{Lemme}

\begin{proof}[Proof of Lemma \ref{lemma_technical}]
Using the induction assumption, we deduce that for all $y \in \R^d$, the map $ \mu \in \PP \mapsto p_m(\mu,s,r,y)$ has a linear derivative given, for all $v \in \R^d$, by \[ \del p_m(\mu,s,t,y)(v) = p_m(\mu,s,r,v,y) + \int_{\R^d} \del p_m(\mu,s,r,x',y)(v) \, d\mu(x').\] Thus, the map $ \mu \in \PP \mapsto b(r,x,[X^{s,\mu,(m)}_r])$ has a linear derivative given by \begin{align}\label{lemma_technical_eq1} \del \left[b(r,x,[X^{s,\mu,(m)}_r])\right](v) &= \int_{\R^d} \del b(r,x,[X^{s,\mu,(m)}_r])(z) p_m(\mu,s,r,v,z) \, dz \\ \notag&\quad + \int_{\R^{2d}} \del b(r,x,[X^{s,\mu,(m)}_r])(z) \del p_m(\mu,s,r,x',z)(v) \, dz \, d\mu(x').\end{align} 

Using the boundedness of $\del b$ on $[0,T]\times \R^d \times \PP \times \R^d$, \eqref{density_bound_Picard}, and the induction assumption, we deduce that  \begin{align*}  &\left\vert \del\left[b(r,x,[X^{s,\mu,(m)}_r])\right](v)\right\vert \\&\leq K + K \int_{\R^{2d}} \left\vert \del p_m(\mu,s,r,x',z)(v)\right\vert \, dz \, d\mu(x')\\ & \leq K + K \int_{\R^{2d}} \left( \sum_{k=1}^m C^k(r-s)^{(k-1)\left(1 - \frac{1}{\alpha}\right)} \prod_{j=1}^{k-1} \BB\left(1 + j \left(1 - \frac{1}{\alpha}\right), 1 - \frac{1}{\alpha}\right)\right)\\ & \hspace{8cm}(r-s)^{1 - \frac{1}{\alpha}}\rho^0(r-s,z-x') \, dz \, d\mu(x') \\ & \leq K + K \sum_{k=1}^m C^k(r-s)^{k\left(1 - \frac{1}{\alpha}\right)} \prod_{j=1}^{k-1} \BB\left(1 + j \left(1 - \frac{1}{\alpha}\right), 1 - \frac{1}{\alpha}\right).
\end{align*} 
This proves \eqref{b_linear_derivative_bound_Picard}. Let us now focus on \eqref{b_gradient_linear_derivative_bound_Picard}. Note that for all $v,x' \in \R^d$, we have \[ \int_{\R^d} \pa_x p_m(\mu,s,t,v,z) \, dz = \int_{\R^d} \pa_v \del p_m(\mu,s,t,x',z)(v) \, dz = 0.\]
Thus, using the induction assumption and differentiating \eqref{lemma_technical_eq1} with respect to $v$, one has for all $v \in \R^d$ \begin{align}\label{b_expression_gradient_linr_der_Picard} \pa_v \del &\left[b(r,x,[X^{s,\mu,(m)}_r])\right](v) \\ \notag&= \int_{\R^d} \left[\del b(r,x,[X^{s,\mu,(m)}_r])(z) - \del b(r,x,[X^{s,\mu,(m)}_r])(v) \right] \pa_x p_m(\mu,s,r,v,z) \, dz \\\notag  &\quad + \int_{\R^{2d}} \left[\del b(r,x,[X^{s,\mu,(m)}_r])(z)  - \del b(r,x,[X^{s,\mu,(m)}_r])(x')\right] \pa_v\del p_m(\mu,s,r,x',z)(v) \, dz \, d\mu(x').\end{align} 

Using the $\eta$-Hölder continuity of $\del b$ with respect to $v$, which is uniform with respect to the other variables, \eqref{density_bound_Picard}, and the induction assumption, we deduce that  \begin{align*}  &\left\vert \pa_v\del\left[b(r,x,[X^{s,\mu,(m)}_r])\right](v)\right\vert \\&\leq K(r-s)^{\frac{\eta -1}{\alpha}} + K \int_{\R^{2d}} |z-x'|^\eta\left\vert \pa_v\del p_m(\mu,s,r,x',z)(v)\right\vert \, dz \, d\mu(x')\\ & \leq K(r-s)^{\frac{\eta -1}{\alpha}} + K \int_{\R^{2d}}  \left( \sum_{k=1}^m C^k(r-s)^{(k-1)\left(1 + \frac{\eta -1}{\alpha}\right)} \prod_{j=1}^{k-1} \BB\left(1 + \frac{\eta -1}{\alpha} + j \left(1 + \frac{\eta -1}{\alpha}\right), 1 - \frac{1}{\alpha}\right)\right)\\ &\hspace{6cm} (r-s)^{\frac{\eta -1}{\alpha} + 1 - \frac{1}{\alpha}} |z-x'|^\eta\rho^0(r-s,z-x')\, dz \, d\mu(x') \\ & \leq K(r-s)^{\frac{\eta -1}{\alpha}} + K   (r-s)^{\frac{\eta -1}{\alpha} }\sum_{k=1}^m C^k(r-s)^{k\left(1 + \frac{\eta -1}{\alpha}\right)} \prod_{j=1}^{k-1} \BB\left(1 + \frac{\eta -1}{\alpha} + j \left(1 + \frac{\eta -1}{\alpha}\right), 1 - \frac{1}{\alpha}\right).
\end{align*} 

It ends the proof of \eqref{b_gradient_linear_derivative_bound_Picard}. Using \eqref{b_linear_derivative_bound_Picard}, it is clear that $ \H_{m+1}$ has a linear derivative given for all $v \in \R^d$ by \begin{equation}\label{expression_linear_der_H} \del \H_{m+1}(\mu,s,r,t,x,y)(v) = \left(\del \left[ b(r,x,[X_r^{s,\mu,(m)}])\right] (v)\right) \cdot \pa_x \p(r,t,x,y),\end{equation}
 and that  \begin{equation}\label{expression_gradient_linear_der_H} \pa_v\del \H_{m+1}(\mu,s,r,t,x,y)(v) = \left(\pa_v\del \left[ b(r,x,[X_r^{s,\mu,(m)}])\right] (v)\right) \cdot \pa_x \p(r,t,x,y).\end{equation}
Thus, \eqref{H_linear_derivative_bound_Picard} and \eqref{H_gradient_linear_derivative_bound_Picard} follow directly from \eqref{b_linear_derivative_bound_Picard}, \eqref{b_gradient_linear_derivative_bound_Picard}, and \eqref{gradientestimatestable}.
\end{proof}
  \begin{Lemme}\label{lemma_technical1}
  	For any $k \geq 1$, $0 \leq s \leq r < t \leq T$, $x,y \in \R^d$, $\H_{m+1}^k(\cdot,s,r,t,x,y)$ has a linear derivative which is $\CC^1$ and satisfies the following properties.
  	
  	\begin{itemize}

  	\item There exists a positive constant $K_m$ depending on $m$ and a positive constant $C$ independent of $k$ and $m$ such that for all $k\geq 1$, $\mu \in \PP$, $0 \leq s \leq r < t \leq T$, $x,y,v \in \R^d$, one has  \begin{multline}\label{H_linear_derivative_bound_Picard_m}
  		\left\vert \del \H_{m+1}^k(\mu,s,r,t,x,y) (v) \right\vert \leq K_m kC^{k-1}(t-r)^{- \frac{1}{\alpha} +(k-1)\left(1 - \frac{1}{\alpha}\right)} \\ \prod_{j=1}^{k-1} \BB\left(j \left(1 - \frac{1}{\alpha}\right), 1 - \frac{1}{\alpha}\right)\rho^1(t-r,y-x),
  	\end{multline}
  and
  	\begin{multline}\label{H_gradient_linear_derivative_bound_Picard_m}
  		\left\vert \pa_v \del \H_{m+1}^k(\mu,s,r,t,x,y) (v) \right\vert  \leq K_mkC^{k-1}(r-s)^{\frac{\eta -1}{\alpha}}  (t-r)^{- \frac{1}{\alpha} +(k-1)\left(1 - \frac{1}{\alpha}\right)} \\ \prod_{j=1}^{k-1} \BB\left( j \left(1 - \frac{1}{\alpha}\right), 1 - \frac{1}{\alpha}\right)  \rho^1(t-r,y-x).
  	\end{multline}
  
  \item  There exists a positive constant $K_m$ depending on $m$ and a positive constant $C$ independent of $k$ and $m$ such that for all $k\geq 1$, $\mu \in \PP$, $0 \leq s < t \leq T$, $x,y,v \in \R^d$, one has  \begin{multline}\label{potimesH_linear_derivative_bound_Picard_m}
  	\left\vert \p\otimes\del \H_{m+1}^k(\mu,s,t,x,y) (v) \right\vert \leq K_m kC^{k-1}(t-s)^{1- \frac{1}{\alpha} +(k-1)\left(1 - \frac{1}{\alpha}\right)} \BB\left(k\left(1 - \frac{1}{\alpha}\right),1\right) \\ \prod_{j=1}^{k-1} \BB\left(j \left(1 - \frac{1}{\alpha}\right), 1 - \frac{1}{\alpha}\right)\rho^0(t-s,y-x),
  \end{multline}
  and
  \begin{multline}\label{potimesH_gradient_linear_derivative_bound_Picard_m}
  	\left\vert \p\otimes\pa_v \del \H_{m+1}^k(\mu,s,t,x,y) (v) \right\vert  \leq K_mkC^{k-1}(t-s)^{\frac{\eta -1}{\alpha} + k \left(1 - \frac{1}{\alpha}\right)} \BB\left(k\left(1-\frac{1}{\alpha}\right),1 + \frac{\eta -1}{\alpha}\right)  \\ \prod_{j=1}^{k-1} \BB\left( j \left(1 - \frac{1}{\alpha}\right), 1 - \frac{1}{\alpha}\right)  \rho^0(t-s,y-x).
  	\end{multline}

\end{itemize}

 \end{Lemme}
\begin{proof}[Proof of Lemma \ref{lemma_technical1}] We proceed by induction on $k$ to prove \eqref{H_linear_derivative_bound_Picard_m} for all $m \geq 1$. The base case $k=1$ is a direct consequence of \eqref{H_linear_derivative_bound_Picard}. Assume now that \eqref{H_linear_derivative_bound_Picard_m} holds for $\H_{m+1}^k$ and let us prove it for $\H_{m+1}^{k+1}$. By definition, we have \[\H^{k+1}_{m+1}(\mu,s,r,t,x,y) = \int_r^t\int_{\R^d} \H_{m+1}(\mu,s,r,r',x,z) \H_{m+1}^{k}(\mu,s,r',t,z,y) \, dz\, dr'.\]
	
	Using Lemma \ref{lemma_technical} and the induction assumption, we deduce that $\H_{m+1}^{k+1}(\cdot,s,r,t,x,y)$ has a linear derivative given by \[ \del \H_{m+1}^{k+1} = \del \H_{m+1} \otimes \H_{m+1}^k + \H_{m+1} \otimes \del \H_{m+1}^k.\]
	
	Using \eqref{H^k_m_bound}, the induction assumption \eqref{H_linear_derivative_bound_Picard_m} and the convolution inequality \eqref{convolineqdensityref}, we deduce that \begin{align*}
		&\left\vert\del \H_{m+1}^{k+1} (\mu,s,r,t,x,y)(v) \right\vert \\ &\leq K_m \int_r^t \int_{\R^d}(r'-r)^{- \frac{1}{\alpha}} \rho^1(r'-r,z-x) C^{k} (t-r')^{-\frac{1}{\alpha} + (k-1)\left( 1 - \frac{1}{\alpha} \right)}\\ &\hspace{6cm} \prod_{j=1}^{k-1} \BB \left( j\left( 1 - \frac{1}{\alpha} \right),1 - \frac{1}{\alpha} \right) \rho^1(t-r',y-z) \, dz\,dr' \\ &\quad + \int_r^t \int_{\R^d} K (r'-r)^{-\frac{1}{\alpha}} \rho^1(r'-r,z-x)  K_m kC^{k-1}(t-r')^{(k-1)\left(1 - \frac{1}{\alpha}\right)} \\ &\hspace{4cm} \prod_{j=1}^{k-1} \BB\left( j \left(1 - \frac{1}{\alpha}\right), 1 - \frac{1}{\alpha}\right)(t-r')^{ - \frac{1}{\alpha}}\rho^1(t-r',y-z) \, dz\,dr'\\ &\leq  K_m C^{k}(t-r)^{- \frac{1}{\alpha} +k\left(1 - \frac{1}{\alpha}\right)} \prod_{j=1}^{k} \BB\left( j \left(1 - \frac{1}{\alpha}\right), 1 - \frac{1}{\alpha}\right) \rho^1(t-r,y-x) \\ &\quad +  K K_m kC^{k-1}(t-r)^{- \frac{1}{\alpha} +k\left(1 - \frac{1}{\alpha}\right)} \prod_{j=1}^{k} \BB\left( j \left(1 - \frac{1}{\alpha}\right), 1 - \frac{1}{\alpha}\right) \rho^1(t-r,y-x) \\ &\leq  K_m (k+1)C^{k}(t-r)^{- \frac{1}{\alpha} +k\left(1 - \frac{1}{\alpha}\right)} \prod_{j=1}^{k-1} \BB\left(j \left(1 - \frac{1}{\alpha}\right), 1 - \frac{1}{\alpha}\right) \rho^1(t-r,y-x),
	\end{align*}
if we choose $C \geq K.$ Following the same lines, we prove that for any $k \geq 1$, $\del \H_{m+1}^{k}$ is $\CC^1$ with respect to $v \in \R^d$ and that it satisfies \[ \pa_v\del \H_{m+1}^{k+1} = \pa_v \del \H_{m+1} \otimes \H_{m+1}^k + \H_{m+1} \otimes \pa_v \del \H_{m+1}^k.\]

Using \eqref{H^k_m_bound}, the induction assumption \eqref{H_gradient_linear_derivative_bound_Picard_m} and the convolution inequality \eqref{convolineqdensityref}, we deduce that \begin{align*}
	&\left\vert\pa_v\del \H_{m+1}^{k+1} (\mu,s,r,t,x,y)(v) \right\vert \\ &\leq K_m \int_r^t \int_{\R^d}(r-s)^{\frac{\eta-1}{\alpha}}(r'-r)^{- \frac{1}{\alpha}} \rho^1(r'-r,z-x) C^{k} (t-r')^{-\frac{1}{\alpha} + (k-1)\left( 1 - \frac{1}{\alpha} \right)}\\ &\hspace{6cm} \prod_{j=1}^{k-1} \BB \left( j\left( 1 - \frac{1}{\alpha} \right),1 - \frac{1}{\alpha} \right) \rho^1(t-r',y-z) \, dz\,dr' \\ &\quad + \int_r^t \int_{\R^d} K (r'-r)^{-\frac{1}{\alpha}} \rho^1(r'-r,z-x)  K_m kC^{k-1}(r'-s)^{\frac{\eta-1}{\alpha}}(t-r')^{(k-1)\left(1 - \frac{1}{\alpha}\right)} \\ &\hspace{4cm} \prod_{j=1}^{k-1} \BB\left( j \left(1 - \frac{1}{\alpha}\right), 1 - \frac{1}{\alpha}\right)(t-r')^{ - \frac{1}{\alpha}}\rho^1(t-r',y-z) \, dz\,dr\\ &\leq  K_m C^{k}(r-s)^{\frac{\eta-1}{\alpha}}(t-r)^{- \frac{1}{\alpha} +k\left(1 - \frac{1}{\alpha}\right)} \prod_{j=1}^{k} \BB\left( j \left(1 - \frac{1}{\alpha}\right), 1 - \frac{1}{\alpha}\right) \rho^1(t-r,y-x) \\ &\quad +  K K_m kC^{k-1}(r-s)^{\frac{\eta-1}{\alpha}}(t-r)^{- \frac{1}{\alpha} +k\left(1 - \frac{1}{\alpha}\right)} \prod_{j=1}^{k} \BB\left( j \left(1 - \frac{1}{\alpha}\right), 1 - \frac{1}{\alpha}\right) \rho^1(t-r,y-x) \\ &\leq  K_m (k+1)C^{k}(r-s)^{\frac{\eta-1}{\alpha}}(t-r)^{- \frac{1}{\alpha} +k\left(1 - \frac{1}{\alpha}\right)} \prod_{j=1}^{k-1} \BB\left(j \left(1 - \frac{1}{\alpha}\right), 1 - \frac{1}{\alpha}\right) \rho^1(t-r,y-x),
\end{align*}
if we chose $C \geq K.$ It proves \eqref{H_gradient_linear_derivative_bound_Picard_m}. The estimates \eqref{potimesH_linear_derivative_bound_Picard_m} and \eqref{potimesH_gradient_linear_derivative_bound_Picard_m} follow immediately from \eqref{H_linear_derivative_bound_Picard_m}, \eqref{H_gradient_linear_derivative_bound_Picard_m}, Lemma \ref{Lemmegradientestimatestable} and the convolution inequality \eqref{convolineqdensityref}.

\end{proof}

\begin{Lemme}\label{lemma_technical5}
	\begin{itemize}
		\item For any $0 \leq s < t \leq T$, $\mu \in \PP,$ $x\in \R^d$, we have \begin{equation}\label{b_time_derivative}
			\left\vert\pa_s \left[ b(t,x,[X^{s,\mu,(m)}_t])\right] \right\vert \leq K \int_{\R^{2d}} (1 \wedge |x'-y|^\eta) | \pa_s p_m(\mu,s,t,x',y)| \, dy \, d\mu(x').
		\end{equation}
		
		\item  For any $0 \leq s <r < t \leq T$, $\mu \in \PP,$ $x,y\in \R^d$, we have
		
		\begin{equation}\label{H_time_derivarive_bound}
			\left\vert \pa_s\H_{m+1} (\mu,s,r,t,x,y)\right\vert \leq K (t-r)^{-\frac{1}{\alpha}} \rho^1(t-r,y-x) \int_{\R^{2d}} (1 \wedge |x'-y|^\eta) |\pa_s p_m(\mu,s,r,x',y)| \, dy \, d\mu(x').
		\end{equation}
	\end{itemize}
\end{Lemme}

\begin{proof}[Proof of Lemma \ref{lemma_technical5}]
	\noindent\textbf{Proof of \eqref{b_time_derivative}.} 
	
	By the induction assumption at step $m$, we see by the dominated convergence theorem that the map $s \in [0,t) \mapsto b(t,x,[X_t^{s,\mu,(m)}])$ is differentiable and that \begin{align*}
		\pa_s \left[b(t,x,[X_t^{s,\mu,(m)}])\right] &= \int_{\R^{2d}} \left( \del b(t,x,[X_t^{s,\mu,(m)}]) (y) -  \del b(t,x,[X_t^{s,\mu,(m)}]) (x') \right) \pa_s p_m(\mu_,s,t,x',y)\, d\mu(x') \, dy. \end{align*}
	
	The boundedness of $\del b $ and the $\eta$-Hölder continuity of $\del b(t,x,\mu)(\cdot)$ allows to conclude.\\
	
	\noindent\textbf{Proof of \eqref{H_time_derivarive_bound}.}
	
	Because of the expression \eqref{def_proxy_kernel_Picard} of $\H_{m+1}$ and since $ \p(\mu,s,r,t,x,y)$ does not depend on $s$, \eqref{H_time_derivarive_bound} follows from \eqref{b_time_derivative} and \eqref{gradientestimatestable}.
\end{proof}

\begin{Lemme}\label{lemma_technical6}
	\begin{itemize}
		\item There exists a positive constant $K_m$ depending on $m$ and a positive constant $C$ such that for all $k \geq 1$, $m\geq 1$, $\mu\in \PP$, $0 \leq s <r <t \leq T$, $x,y \in \R^d$ 
		
		\begin{multline}\label{H^k_time_derivative_bound}
			\left\vert \pa_s \H^{k}_{m+1}(\mu,s,r,t,x,y) \right\vert \leq K_mkC^{k-1} (r-s)^{\frac{\eta}{\alpha}-1}(t-r)^{-\frac{1}{\alpha} + (k-1)\left(1 - \frac{1}{\alpha}\right)}\rho^1(t-r,y-x)\\ \prod_{j=1}^{k-1} \BB\left(j\left(1 - \frac{1}{\alpha}\right), 1 - \frac{1}{\alpha}\right).
		\end{multline}
		
		\item There exists a positive constant $K_m$ depending on $m$ such that for all  $k \geq 1$, $m\geq 1$, $\mu\in \PP$, $0 \leq s <r <t \leq T$, $x,y \in \R^d$ 
		
		\begin{equation}\label{Phi_time_derivative_bound}
			\left\vert \pa_s\Phi_{m+1}(\mu,s,r,t,x,y) \right\vert \leq K_m (r-s)^{\frac{\eta}{\alpha}-1}(t-r)^{-\frac{1}{\alpha}}\rho^1(t-r,y-x).
		\end{equation}
	\end{itemize}
\end{Lemme}
\begin{proof}[Proof of Lemma \ref{lemma_technical6}.] 
	
	\noindent\textbf{Proof of \eqref{H^k_time_derivative_bound}.} 
	We proceed by induction on $k \geq 1$. The base case $k=1$ comes from \eqref{H_time_derivarive_bound} and the induction assumption \eqref{density_time_der_bound_Picard} for $\pa_s p_m$. Indeed, the induction assumption and the space-time inequality \eqref{scalingdensityref} ensure that \[ \int_{\R^{2d}} (1 \wedge |x'-y|^\eta) |\pa_s p_m(\mu,s,r,x',y)| \, dy \, d\mu(x') \leq K_m(r-s)^{\frac{\eta}{\alpha}-1} \int_{\R^{2d}} \rho^{-\tilde{\eta}-\eta}(r-s,y-x') \, dy \ d\mu(x').\] We conclude since $\eta + \tilde{\eta} < \alpha$, the map $\rho^{-\tilde{\eta}-\eta}(r-s,\cdot)$ belongs to $L^1(\R^d)$ and $\int_{\R^d} \rho^{-\tilde{\eta}-\eta} (t-s,y)\,dy$ is equal to a constant independent of $s$ and $t.$ For the induction step, we assume that \eqref{H^k_time_derivative_bound} holds at step $k$. By definition, we have \begin{align*}
		\H_{m+1}^{k+1} (\mu,s,r,t,x,y) = \int_r^t \int_{\R^d} \H_{m+1}(\mu,s,r,r',x,z) \H_{m+1}^k(\mu,s,r',t,z,y) \, dz \, dr'.
	\end{align*}
	
	It follows from the induction assumption and the dominated convergence theorem that the map $s \in [0,t) \mapsto \H_{m+1}^{k+1}(\mu,s,r,t,x,y)$ is differentiable and that \begin{align*}
		\pa_s \H_{m+1}^{k+1}(\mu,s,r,t,x,y) &= \int_r^t\int_{\R^d} \pa_s \H_{m+1} (\mu,s,r,r',x,z) \H_{m+1}^k(\mu,s,r',t,z,y) \, dz\,dr' \\ &\quad+   \int_r^t\int_{\R^d}  \H_{m+1} (\mu,s,r,r',x,z) \pa_s\H_{m+1}^k(\mu,s,r',t,z,y) \, dz\,dr' \\ &=: I_1 + I_2.
	\end{align*}
	
	Using the base case $k=1$ and \eqref{H^k_m_bound}, we deduce that 
	
	\begin{align*}
		|I_1| &\leq \int_{r}^t \int_{\R^d} K_m (r-s)^{\frac{\eta}{\alpha}-1} (r'-r)^{-\frac{1}{\alpha}} \rho^1(r'-r,z-x) C^k (t-r')^{-\frac{1}{\alpha} + (k-1)\left(1 - \frac{1}{\alpha}\right)} \\ &\hspace{8cm}\prod_{j=1}^{k-1} \BB\left(j\left(1 - \frac{1}{\alpha}\right), 1 - \frac{1}{\alpha}\right) \rho^1(t-r',y-z) \, dr'\, dz \\ &\leq K_mC^k (r-s)^{\frac{\eta}{\alpha} -1}(t-r)^{-\frac{1}{\alpha} + k \left(1- \frac{1}{\alpha}\right)} \prod_{j=1}^{k} \BB\left(j\left(1 - \frac{1}{\alpha}\right), 1 - \frac{1}{\alpha}\right) \rho^1(t-r,y-x).
	\end{align*}
	
	For $I_2$, from the induction assumption \eqref{H^k_time_derivative_bound}, \eqref{H^k_m_bound}, the convolution inequality \eqref{convolineqdensityref} and since $\frac{\eta}{\alpha} < 1$, we obtain that
	
	\begin{align*}
		|I_2| &\leq \int_r^t \int_{\R^d} K (r'-r)^{-\frac{1}{\alpha}}\rho^1(r'-r,z-x)  K_mkC^{k-1} (r'-s)^{\frac{\eta}{\alpha}-1}(t-r')^{-\frac{1}{\alpha} + (k-1)\left(1 - \frac{1}{\alpha}\right)} \\ &\hspace{7cm}\prod_{j=1}^{k-1} \BB\left(j\left(1 - \frac{1}{\alpha}\right), 1 - \frac{1}{\alpha}\right)\rho^1(t-r',y-z)\, dz\,dr' \\ &\leq (r-s)^{\frac{\eta}{\alpha}-1}\int_r^t \int_{\R^d} K (r'-r)^{-\frac{1}{\alpha}}\rho^1(r'-r,z-x)  K_mkC^{k-1} (t-r')^{-\frac{1}{\alpha} + (k-1)\left(1 - \frac{1}{\alpha}\right)} \\ &\hspace{7cm}\prod_{j=1}^{k-1} \BB\left(j\left(1 - \frac{1}{\alpha}\right), 1 - \frac{1}{\alpha}\right)\rho^1(t-r',y-z)\, dz\,dr' \\ &\leq K_mkC^k(r-s)^{\frac{\eta}{\alpha}-1} (t-r)^{-\frac{1}{\alpha} + k\left(1 - \frac{1}{\alpha}\right)}\rho^1(t-r,y-x)  \prod_{j=1}^{k} \BB\left(j\left(1 - \frac{1}{\alpha}\right), 1 - \frac{1}{\alpha}\right),
	\end{align*}
	provided that we choose $C \geq K$ in \eqref{H^k_time_derivative_bound}. This concludes the induction step for \eqref{H^k_time_derivative_bound}.\\
	
	\noindent\textbf{Proof of \eqref{Phi_time_derivative_bound}.} Using the definition of $\Phi_{m+1}$ \eqref{defvolterraMK_Picard} and \eqref{H^k_time_derivative_bound}, we obtain by the dominated convergence theorem that $s \in [0,r) \mapsto \phi_{m+1}(\mu,s,r,t,x,y)$ is continuously differentiable with \[ \pa_s \Phi_{m+1} (\mu,s,r,t,x,y) = \sum_{k=1}^{\infty} \pa_s \H^k_{m+1} (\mu,s,r,t,x,y).\] 
	
	Then, \eqref{Phi_time_derivative_bound} follows immediately from \eqref{H^k_time_derivative_bound}.

\end{proof}

\subsection{First part of the proof of the induction step}

We split the proof of the induction step into different parts for the sake of clarity. We start by proving the estimates \eqref{linear_der_bound_Picard}, \eqref{gradient_linear_der_bound_Picard} and \eqref{density_time_der_bound_Picard}.\\

\noindent\textbf{Proof of \eqref{linear_der_bound_Picard} and \eqref{gradient_linear_der_bound_Picard}.} We start by showing that for all $ 0 \leq s < t \leq T$, $x,y \in \R^d$, the map $p_{m+1}(\cdot,s,t,x,y)$ admits a linear derivative given which is $\CC^1$ with respect to $v\in \R^d$ and such that for all $v \in \R^d$ \begin{equation}\label{expression_linear_der}
	\del p_{m+1}(\mu,s,t,x,y)(v) = \sum_{k=1}^{\infty} \p \otimes \del \H_{m+1}^k(\mu,s,t,x,y)(v),
\end{equation}
 and \begin{equation}\label{expression_gradient_linear_der}
 	\pa_v\del p_{m+1}(\mu,s,t,x,y)(v) = \sum_{k=1}^{\infty} \p \otimes \pa_v\del \H_{m+1}^k(\mu,s,t,x,y)(v).
 \end{equation}
where the series are absolutely convergent. Moreover, we also have the following representation formulas \begin{equation}\label{representation_linear_der}
	\del p_{m+1}(\mu,s,t,x,y)(v) = \sum_{k=0}^{\infty} p_{m+1} \otimes \del \H_{m+1} \otimes \H_{m+1}^k(\mu,s,t,x,y)(v),
\end{equation}
and \begin{equation}\label{representation_gradient_linear_der}
	\pa_v \del p_{m+1}(\mu,s,t,x,y)(v) = \sum_{k=0}^{\infty} p_{m+1} \otimes \pa_v\del \H_{m+1} \otimes \H^k_{m+1}(\mu,s,t,x,y)(v).
\end{equation}

To prove \eqref{expression_linear_der}, we fix $\mu_1,\mu_2 \in \PP$ and we write thanks to Lemma \ref{lemma_technical1} and Fubini's theorem

\begin{align*}
	&p_{m+1}(\mu_1,s,t,x,y) - p_{m+1}(\mu_2,s,t,x,y) \\&= \sum_{k=1}^{\infty} \p\otimes \H_{m+1}^k(\mu_1,s,t,x,y) - \p\otimes \H_{m+1}^k(\mu_2,s,t,x,y) \\ &= \sum_{k=1}^{\infty} \int_s^t \int_{\R^d} \p (s,r,x,z)(\H_{m+1}^k(\mu_1,s,r,t,z,y) - \H_{m+1}^k(\mu_2,s,r,t,z,y))  \,dz \,dr  \\ &= \sum_{k=1}^{\infty} \int_s^t \int_{\R^d} \p (s,r,x,z)\int_0^1 \int_{\R^d} \del\H_{m+1}^k(\lambda\mu_1 + (1-\lambda)\mu_2,s,r,t,z,y)(v) \, d(\mu_1 - \mu_2)(v)\, d\lambda  \,dz \,dr \\ &=\int_0^1 \int_{\R^d} \sum_{k=1}^{\infty} \int_s^t \int_{\R^d} \p (s,r,x,z) \del\H_{m+1}^k(\lambda\mu_1 + (1-\lambda)\mu_2,s,r,t,z,y)(v) \\ &\hspace{8cm}\, d(\mu_1 - \mu_2)(v)\, d\lambda  \,dz \,dr  \, d(\mu_1 - \mu_2)(v)\, d\lambda \\ &= \int_0^1 \int_{\R^d} \sum_{k=1}^{\infty} \p \otimes \del\H_{m+1}^k (\lambda\mu_1 + (1-\lambda)\mu_2,s,t,x,y)(v) \, d(\mu_1 - \mu_2)(v)\, d\lambda.
\end{align*}
Note that owing to \eqref{potimesH_linear_derivative_bound_Picard_m} the series \eqref{expression_linear_der} is absolutely convergent, locally uniformly with respect to $(\mu,s,x,v)\in \PP\times [0,t)\times \R^d \times \R^d$. This concludes the proof of \eqref{expression_linear_der}. Moreover, we have proved that the map $(\mu,s,x,v) \mapsto \del p_{m+1}(\mu,s,t,x,y)(v)$ is continuous. By differentiation under the integral using \eqref{H_gradient_linear_derivative_bound_Picard_m}, we obtain \eqref{expression_gradient_linear_der}. The dominated convergence theorem yields the continuity of the map $(\mu,s,x,v) \mapsto \pa_v\del p_{m+1}(\mu,s,t,x,y)(v)$.\\

Let us now focus on the representation formula \eqref{representation_linear_der}. Using the parametrix expansion \eqref{representationdensityparametrix_Picard} of $p_{m+1}$, one has for all $\mu \in \PP$, $0 \leq s <t\leq T$, $x,y \in \R^d$ \[p_{m+1}(\mu,s,t,x,y) = \p(s,t,x,y) + p_{m+1}\otimes \H_{m+1}(\mu,s,t,x,y).\] 
 We can easily see by induction thanks to \eqref{H_linear_derivative_bound_Picard_m} that for any $k \geq 1$, one has \[ \del \H_{m+1}^k(\mu,s,r,t,x,y)(v) = \sum_{j=1}^{k} \H_{m+1}^{k-j} \otimes \del \H_{m+1} \otimes \H_{m+1}^{j-1} (\mu,s,r,t,x,y)(v).\]
 We plug this expression into \eqref{expression_linear_der}, and since the series is absolutely convergent, we obtain setting $l=k-j$ and $i=j-1$ and by Fubini's theorem that \begin{align*}
 	\del p_{m+1}(\mu,s,t,x,y)(v) &= \sum_{k=1}^{\infty} \sum_{j=1}^k \p \otimes \H_{m+1}^{k-j} \otimes \del \H_{m+1} \otimes \H_{m+1}^{j-1}(\mu,s,t,x,y)(v) \\ &= \sum_{l=0}^{\infty} \sum_{i=0}^{\infty} \p \otimes \H_{m+1}^{l} \otimes \del \H_{m+1} \otimes \H_{m+1}^{i}(\mu,s,t,x,y)(v) \\ &= \sum_{i=0}^{\infty} p_{m+1}\otimes \del \H_{m+1} \otimes \H_{m+1}^i(\mu,s,t,x,y)(v).
 \end{align*}
 
 This is exactly \eqref{representation_linear_der}.
Let us now prove that the estimate \eqref{linear_der_bound_Picard} is still true at step $m+1$. It follows from \eqref{H_linear_derivative_bound_Picard}, \eqref{density_bound_Picard} and the convolution inequality \eqref{convolineqdensityref} that for some positive constant $K$ independent of $m$ one has \begin{align*}
	&\left\vert \p \otimes \del \H_{m+1}(\mu,s,t,x,y) (v) \right\vert\\ &\leq K \int_s^t\int_{\R^d} \left[1+ \sum_{k=1}^m C^k(r-s)^{k\left(1 - \frac{1}{\alpha}\right)} \prod_{j=1}^{k-1} \BB\left(1 + j \left(1 - \frac{1}{\alpha}\right), 1 - \frac{1}{\alpha}\right)\right]\\ &\hspace{7cm} \rho^0(r-s,z-x) (t-r)^{-\frac{1}{\alpha}}\rho^1(t-r,y-z) \, dz\,dr  \\ &\leq K (t-s)^{1-\frac{1}{\alpha}} + K\sum_{k=1}^m C^k(t-s)^{(k+1)\left(1 - \frac{1}{\alpha}\right)} \prod_{j=1}^{k} \BB\left(1 + j \left(1 - \frac{1}{\alpha}\right), 1 - \frac{1}{\alpha}\right)\rho^0(t-s,y-x).
	\end{align*}
Using this inequality, the bound \eqref{H^k_m_bound} for $\H_{m+1}^k$, the convolution inequality \eqref{convolineqdensityref} and summing over $k \geq 0$, we find that there exists a constant $K$ independent of $m$ such that 

\begin{align*}
	&\left\vert \del p_{m+1}(\mu,s,t,x,y)(v)\right\vert \\ &\leq \sum_{k=0}^{\infty} \left\vert p_{m+1}\otimes \del \H_{m+1} \otimes \H_{m+1}^k (\mu,s,t,x,y)(v) \right\vert\\ &\leq K \left( (t-s)^{1-\frac{1}{\alpha}} + \sum_{k=1}^m C^k(t-s)^{(k+1)\left(1 - \frac{1}{\alpha}\right)} \prod_{j=1}^{k} \BB\left(1 + j \left(1 - \frac{1}{\alpha}\right), 1 - \frac{1}{\alpha}\right)\right)\rho^0(t-s,y-x)\\ &\leq (t-s)^{1 - \frac{1}{\alpha}}\rho^0(t-s,y-x)  \left( \sum_{k=1}^{m+1} C^{k+1}(t-s)^{k\left(1 - \frac{1}{\alpha}\right)} \prod_{j=1}^{k-1} \BB\left(1 + j \left(1 - \frac{1}{\alpha}\right), 1 - \frac{1}{\alpha}\right)\right),
\end{align*}

provided that we choose $C\geq K$ in \eqref{linear_der_bound_Picard}. It ends the proof of the induction step for \eqref{linear_der_bound_Picard}.\\

 Notice that by differentiating under the integral \eqref{representation_linear_der} using \eqref{density_bound_Picard}, \eqref{H_gradient_linear_derivative_bound_Picard} and \eqref{H^k_m_bound}, we obtain the representation formula \eqref{representation_gradient_linear_der} for $\pa_v \del p_{m+1}$. Let us now prove that the estimate \eqref{gradient_linear_der_bound_Picard} is verified at step $m+1$. It follows from \eqref{H_gradient_linear_derivative_bound_Picard}, \eqref{density_bound_Picard} and the convolution inequality \eqref{convolineqdensityref} that for some positive constant $K$ independent of $m$ one has \begin{align*}
 	&\left\vert \p \otimes \pa_v\del \H_{m+1}(\mu,s,t,x,y) (v) \right\vert\\ &\leq K \int_s^t\int_{\R^d} (r-s)^{\frac{\eta-1}{\alpha}}\left[1+ \sum_{k=1}^m C^k(r-s)^{k\left(1 +\frac{\eta-1}{\alpha}\right)}  \prod_{j=1}^{k-1} \BB\left(1 + \frac{\eta -1}{\alpha} + j \left(1 + \frac{\eta -1}{\alpha}\right), 1 - \frac{1}{\alpha}\right)\right]\\ &\hspace{7cm} \rho^0(r-s,z-x) (t-r)^{-\frac{1}{\alpha}}\rho^1(t-r,y-z) \, dz\,dr  \\ &\leq K (t-s)^{\frac{\eta-1}{\alpha} + 1-\frac{1}{\alpha}} \left[ 1+ \sum_{k=1}^m C^k(t-s)^{k\left(1 +\frac{\eta-1}{\alpha}\right)} \prod_{j=1}^{k} \BB\left(1 + \frac{\eta -1}{\alpha} + j \left(1 + \frac{\eta -1}{\alpha}\right), 1 - \frac{1}{\alpha}\right)\right]\rho^0(t-s,y-x).
 \end{align*}
 Using this inequality, the bound \eqref{H^k_m_bound} for $\H_{m+1}^k$, the convolution inequality \eqref{convolineqdensityref} and summing over $k \geq 0$, we find that there exists a constant $K$ independent of $m$ such that 
 
 \begin{align*}
 	&\left\vert \pa_v\del p_{m+1}(\mu,s,t,x,y)(v)\right\vert \\ &\leq \sum_{k=0}^{\infty} \left\vert p_{m+1}\otimes \pa_v \del \H_{m+1} \otimes \H_{m+1}^k (\mu,s,t,x,y)(v) \right\vert\\ &\leq K (t-s)^{\frac{\eta-1}{\alpha} + 1-\frac{1}{\alpha}} \left[ 1+ \sum_{k=1}^m C^k(t-s)^{k\left(1 +\frac{\eta-1}{\alpha}\right)} \prod_{j=1}^{k} \BB\left(1 + \frac{\eta -1}{\alpha} + j \left(1 + \frac{\eta -1}{\alpha}\right), 1 - \frac{1}{\alpha}\right)\right]\\ &\hspace{11cm} (t-s)^{\frac{\eta-1}{\alpha} + 1-\frac{1}{\alpha}}\rho^0(t-s,y-x)\\ &\leq (t-s)^{\frac{\eta-1}{\alpha} + 1 - \frac{1}{\alpha}}\rho^0(t-s,y-x)  \left( \sum_{k=1}^{m+1} C^{k}(t-s)^{(k-1)\left(1 + \frac{\eta-1}{\alpha}\right)} \prod_{j=1}^{k-1} \BB\left(1 + \frac{\eta -1}{\alpha} + j \left(1 + \frac{\eta -1}{\alpha}\right), 1 - \frac{1}{\alpha}\right)\right),
 \end{align*}
 
 provided that we choose $C\geq K$ in \eqref{gradient_linear_der_bound_Picard}. It ends the proof of the induction step for \eqref{gradient_linear_der_bound_Picard}.\\ 
 
 \noindent\textbf{Proof of \eqref{density_time_der_bound_Picard}.} Let us first prove the following representation formula for $\pa_s \Phi_{m+1}$ \begin{multline}\label{representation_time_der_Phi}
 	\pa_s \Phi_{m+1} (\mu,s,r,t,x,y) = \left[\pa_s \H_{m+1} + \pa_s \H_{m+1} \otimes \Phi_{m+1}\right](\mu,s,r,t,x,y) \\+ \Phi_{m+1} \otimes \left[\pa_s \H_{m+1} + \pa_s \H_{m+1} \otimes \Phi_{m+1}\right](\mu,s,r,t,x,y).
 \end{multline}
 By differentiating the relation $\Phi_{m+1} = \H_{m+1} + \H_{m+1}\otimes \Phi_{m+1}$, we obtain that 
 \begin{multline*}
 	\pa_s \Phi_{m+1}(\mu,s,r,t,x,y) = \pa_s\H_{m+1}(\mu,s,r,t,x,y) +\pa_s \H_{m+1} \otimes \Phi_{m+1} (\mu,s,r,t,x,y) \\+ \H_{m+1}\otimes \pa_s \Phi_{m+1}(\mu,s,r,t,x,y).
 \end{multline*}
 Notice that by \eqref{H^k_time_derivative_bound} for $k=1$ and \eqref{Phi_m_bound}, we get that \begin{equation*}
 	\left\vert \left[\pa_s\H_{m+1} +\pa_s \H_{m+1} \otimes \Phi_{m+1}\right] (\mu,s,r,t,x,y) \right\vert \leq K_m(r-s)^{\frac{\eta}{\alpha}-1}(t-r)^{-\frac{1}{\alpha}} \rho^1(t-r,y-x).
 \end{equation*}
 The kernel $\left[\pa_s\H_{m+1} +\pa_s \H_{m+1} \otimes \Phi_{m+1}\right] (\mu,s,r,t,x,y)$ yields a time-integrable singularity of order $(r-t)^{-\frac{1}{\alpha}}$. We can thus iterate the previous relation to obtain that \begin{align*}
 	&\pa_s \Phi_{m+1} (\mu,s,r,t,x,y) \\&= \sum_{k=0}^{\infty} \H_{m+1}^k \otimes \left[ \pa_s \H_{m+1} + \pa_s \H_{m+1} \otimes \Phi_{m+1}\right](\mu,s,r,t,x,y) \\ &=  \left[\pa_s \H_{m+1} + \pa_s \H_{m+1} \otimes \Phi_{m+1}\right](\mu,s,r,t,x,y) + \Phi_{m+1} \otimes \left[\pa_s \H_{m+1} + \pa_s \H_{m+1} \otimes \Phi_{m+1}\right](\mu,s,r,t,x,y).
 \end{align*}
 
 In order to deal with the differentiability of the map $s\in[0,t) \mapsto p_{m+1}(\mu,s,t,x,y)$, keeping in mind the parametrix expansion \eqref{representationdensityparametrix_Picard}, we first study the differentiability of the map \[s\in[0,r) \mapsto \int_{\R^d}\p(s,r,x,z) \Phi_{m+1}(\mu,s,r,t,z,y) \, dz.\] 
 
 Since $\int_{\R^d}\pa_s \p(s,r,x,z) \, dz =0$, we deduce by the dominated convergence theorem that \begin{align}\label{proof_time_der_bound_Picard_eq1}
 	\notag\pa_s \left(\int_{\R^d}\p(s,r,x,z) \Phi_{m+1}(\mu,s,r,t,z,y) \, dz\right) &= \int_{\R^d} \pa_s \p(s,r,x,z) \left(\Phi_{m+1}(\mu,s,r,t,z,y) - \Phi_{m+1}(\mu,s,r,t,x,y)\right)\,dz \\ &\quad+ \int_{\R^d} \p(s,r,x,z) \pa_s\Phi_{m+1}(\mu,s,r,t,z,y) \, dz \\ \notag&=: A_1 +A_2,
 \end{align}
 which is continuous with respect to $(\mu,s,x) \in \PP \times [0,r) \times \R^d.$ We now control $A_1$ and $A_2$. For $A_1$, it follows from \eqref{timederivativedensity_bound}, \eqref{Phi_m_Holder} with $\gamma = \tilde{\eta}$ and the space-time inequality \eqref{scalingdensityref} that \begin{align*}
 	|A_1| &\leq K\int_{\R^d}  (r-s)^{-1} \rho^{0}(r-s,z-x) (t-r)^{-\frac{\tilde{\eta}+1}{\alpha}} |z-x|^{\tilde{\eta}} \left[\rho^1(t-r,y-z) + \rho^1(t-r,y-x)\right] \, dz \\ &\leq K \int_{\R^d}  (r-s)^{\frac{\tilde{\eta}}{\alpha}-1} \rho^{-\tilde{\eta}}(r-s,z-x) (t-r)^{-\frac{\tilde{\eta}+1}{\alpha}}  \left[\rho^{-\tilde{\eta}}(t-r,y-z) + \rho^{-\tilde{\eta}}(t-r,y-x)\right] \,dz.
 \end{align*}
 
 Note that $\int_{\R^d}\rho^{ -\tilde{\eta}}(r-s,z-x) \, dz$ is a constant independent of $s$ and $r$ . By the definition of $\rho^{-\tilde{\eta}}$ \eqref{defdensityref}, the fact that $r >s$ and the convolution inequality \eqref{convolineqdensityref}, one has 
 \begin{equation}\label{proof_time_der_Picard_eq2}
 	|A_1| \leq K  (r-s)^{\frac{\tilde{\eta}}{\alpha}-1}(t-r)^{-\frac{\tilde{\eta}+1}{\alpha}}  \left[\rho^{-\tilde{\eta}}(t-s,y-x) + (t-r)^{-\frac{d}{\alpha}}(1+(t-s)^{-\frac{1}{\alpha}}|y-x|)^{-d-\alpha +\tilde{\eta}}\right]. 
 \end{equation}
 
 Concerning $A_2$, it follows from \eqref{Phi_time_derivative_bound} and \eqref{gradientestimatestable} that 
 
 \begin{align*}
 	|A_2| &\leq \int_{\R^d} K \rho^0(r-s,z-x) K_m(r-s)^{\frac{\eta}{\alpha}-1}(t-r)^{-\frac{1}{\alpha}} \rho^1(t-r,y-z)\,dz \\ &\leq K_m(r-s)^{\frac{\eta}{\alpha}-1}(t-r)^{-\frac{1}{\alpha}} \rho^0(t-s,y-x).
 \end{align*}
 
 By the dominated convergence theorem justified by the controls previously obtained on $A_1$ and $A_2$, we obtain that the map $s \in [0,t) \mapsto \p \otimes \Phi_{m+1}(\mu,s,t,x,y)$ is differentiable with 
 
 \begin{align*}
 	\pa_s \left(\p \otimes \Phi_{m+1} \right)(\mu,s,t,x,y) &= - \Phi_{m+1}(\mu,s,t,x,y) + \pa_s \p \otimes \Phi_{m+1}(\mu,s,t,x,y)  + \p \otimes \pa_s \Phi_{m+1}(\mu,s,t,x,y)  \\ &=  - \Phi_{m+1}(\mu,s,t,x,y) \\ &\quad+ \int_s^t \int_{\R^d}\pa_s \p(s,r,x,z) \left(\Phi_{m+1}(\mu,s,r,t,z,y) - \Phi_{m+1}(\mu,s,r,t,x,y)\right)\,dz\,dr \\ &\quad + \int_s^t\int_{\R^d} \p(s,r,x,z) \pa_s \Phi_{m+1}(\mu,s,r,t,z,y)\,dz\,dr.
 \end{align*}
 
 Thus, the map $s \in [0,t) \mapsto p_{m+1}(\mu,s,t,x,y)$ is differentiable with \begin{multline}\label{expression_time_derivative_density_Picard}
 	\pa_s p_{m+1}(\mu,s,t,x,y) = \pa_s \p(s,t,x,y)  - \Phi_{m+1}(\mu,s,t,x,y)\\ + \pa_s \p \otimes \Phi_{m+1}(\mu,s,t,x,y)  + \p \otimes \pa_s \Phi_{m+1}(\mu,s,t,x,y).
 \end{multline}
 
 Then, plugging \eqref{representation_time_der_Phi} into \eqref{expression_time_derivative_density_Picard}, we obtain that \begin{align*}
 	\pa_s p_{m+1}(\mu,s,t,x,y) &= \pa_s \p(s,t,x,y) - \Phi_{m+1}(\mu,s,t,x,y) + \pa_s\p \otimes \Phi_{m+1}(\mu,s,t,x,y) \\ &\quad+ \p \otimes \left(\pa_s \H_{m+1} + \pa_s \H_{m+1} \otimes \Phi_{m+1}\right)(\mu,s,t,x,y) \\ &\quad +  \p\otimes\Phi_{m+1} \otimes \left(\pa_s \H_{m+1} + \pa_s \H_{m+1} \otimes \Phi_{m+1}\right)(\mu,s,t,x,y).
 \end{align*}
 
 By the parametrix expansion \eqref{representationdensityparametrix_Picard} of $p_{m+1}$, we obtain the following representation formula \begin{align}\label{representation_time_der_density_Picard}
 	\notag\pa_s p_{m+1}(\mu,s,t,x,y) & = \pa_s \p(s,t,x,y)-\Phi_{m+1}(\mu,s,t,x,y) \\ \notag&\quad + \pa_s \p \otimes \Phi_{m+1}(\mu,s,t,x,y) \\ \notag&\quad+ p_{m+1}\otimes \left(\pa_s \H_{m+1} + \pa_s \H_{m+1} \otimes \Phi_{m+1}\right) (\mu,s,t,x,y)\\ \notag&=: I_1 + I_2 + I_3.
 \end{align}
 
 We can now prove that the estimate \eqref{density_time_der_bound_Picard} is still true at step $m+1$ for some choice of the constant $C$ which appears in \eqref{density_time_der_bound_Picard}.\\ 
 
 For $I_1$, we note that \eqref{eqproofproxy7} and \eqref{Phi_m_bound} yield \[ |I_1| \leq K (t-s)^{-1}\rho^{0}(t-s,y-x)  \leq K (t-s)^{-1}\rho^{-\tilde{\eta}}(t-s,y-x),\] since $\alpha \in (1,2)$ and $\rho^{0}(t-s,y-x) \leq \rho^{-\tilde{\eta}}(t-s,y-x).$

 Concerning $I_2$, it can be rewritten as \[ I_2 = \int_s^t \int_{\R^d} \pa_s \p(s,r,x,z) \left( \Phi_{m+1}(\mu,s,r,t,z,y) - \Phi_{m+1}(\mu,s,r,t,x,y)\right) \, dz \,dr.\]

 Owing to the bound \eqref{proof_time_der_Picard_eq2} obtained for $A_1$, which was defined in \eqref{proof_time_der_bound_Picard_eq1}, we deduce that 
 
 \begin{align*}
 	|I_2| &\leq K\int_s^t  (r-s)^{\frac{\tilde{\eta}}{\alpha}-1}(t-r)^{-\frac{\tilde{\eta}+1}{\alpha}}  \left[\rho^{-\tilde{\eta}}(t-s,y-x) + (t-r)^{-\frac{d}{\alpha}}(1+(t-s)^{-\frac{1}{\alpha}}|y-x|)^{-d-\alpha +\tilde{\eta}}\right]\, dr \\ &\leq K (t-s)^{-\frac{1}{\alpha}} \rho^{ - \tilde{\eta}}(t-s,y-x) + K (t-s)^{- \frac{1}{\alpha} - \frac{d}{\alpha}} (1+(t-s)^{-\frac{1}{\alpha}}|y-x|)^{-d-\alpha +\tilde{\eta}} \\ &\leq K (t-s)^{-\frac{1}{\alpha}} \rho^{ - \tilde{\eta}}(t-s,y-x) \\ &\leq K (t-s)^{-1} \rho^{ - \tilde{\eta}}(t-s,y-x). 
 \end{align*}
 
 We now focus on $I_3$. Using the induction assumption \eqref{density_time_der_bound_Picard}, \eqref{H_time_derivarive_bound}, the fact that $\eta+ \tilde{\eta} <\alpha$ and the space-time inequality  \eqref{scalingdensityref}, we get that 
 
 \begin{align*}
 	&\left\vert \pa_s \H_{m+1}(\mu,s,r,t,x,y) \right\vert \\& \leq  K (t-r)^{-\frac{1}{\alpha}} \rho^1(t-r,y-x) \int_{\R^{2d}} |x'-y|^\eta (r-s)^{-1} \rho^{-\tilde{\eta}}(r-s,y-x') \\ &\hspace{3cm} \sum_{k=1}^m C^k (r-s)^{(k-1)\left(1 + \frac{\eta-1}{\alpha}\right)} \prod_{j=1}^{k-1} \BB\left(\frac{\eta}{\alpha} + (j-1)\left(1 + \frac{\eta -1}{\alpha}\right),1 - \frac{1}{\alpha}\right) \, dy \, d\mu(x') \\ &\leq   K (r-s)^{\frac{\eta}{\alpha}-1} (t-r)^{-\frac{1}{\alpha}} \rho^1(t-r,y-x)   \\ &\hspace{3cm} \sum_{k=1}^m C^k (r-s)^{(k-1)\left(1 + \frac{\eta-1}{\alpha}\right)} \prod_{j=1}^{k-1} \BB\left(\frac{\eta}{\alpha} + (j-1)\left(1 + \frac{\eta -1}{\alpha}\right),1 - \frac{1}{\alpha}\right). 
 \end{align*}
 
 Since the kernel $\Phi_{m+1}$ yields a time-integrable singularity by \eqref{Phi_m_bound}, we deduce with the preceding inequality that 
 
 \begin{align*}
 	&\left\vert \left[\pa_s \H_{m+1} + \pa_s \H_{m+1} \otimes \Phi_{m+1}\right](\mu,s,r,t,x,y) \right\vert \\ &\leq   K (r-s)^{\frac{\eta}{\alpha}-1} (t-r)^{-\frac{1}{\alpha}} \rho^1(t-r,y-x)   \\ &\hspace{3cm} \sum_{k=1}^m C^k (r-s)^{(k-1)\left(1 + \frac{\eta-1}{\alpha}\right)} \prod_{j=1}^{k-1} \BB\left(\frac{\eta}{\alpha} + (j-1)\left(1 + \frac{\eta -1}{\alpha}\right),1 - \frac{1}{\alpha}\right). 
 \end{align*}
 
 It follows from this inequality, \eqref{density_bound_Picard} and the convolution inequality \eqref{convolineqdensityref} that 
 
 \begin{align*}
 	|I_3| &\leq \int_{s}^t \int_{\R^d} \rho^0(r-s,z-x)  K (r-s)^{\frac{\eta}{\alpha}-1} (t-r)^{-\frac{1}{\alpha}} \rho^1(t-r,y-z)   \\ &\hspace{3cm} \sum_{k=1}^m C^k (r-s)^{(k-1)\left(1 + \frac{\eta-1}{\alpha}\right)} \prod_{j=1}^{k-1} \BB\left(\frac{\eta}{\alpha} + (j-1)\left(1 + \frac{\eta -1}{\alpha}\right),1 - \frac{1}{\alpha}\right)\, dz\,dr \\ &\leq  K (t-s)^{-1}\rho^{0}(t-s,y-x)\sum_{k=1}^m C^k (r-s)^{k\left(1 + \frac{\eta-1}{\alpha}\right)} \prod_{j=1}^{k} \BB\left(\frac{\eta}{\alpha} + (j-1)\left(1 + \frac{\eta -1}{\alpha}\right),1 - \frac{1}{\alpha}\right).
 \end{align*}

Gathering the previous estimates on $I_1$, $I_2$ and $I_3$, we have \begin{align*}
	&|\pa_s p_{m+1}(\mu,s,t,x,y) | \\&\leq   (t-s)^{-1}\rho^{-\tilde{\eta}}(t-s,y-x)K\left[1 + \sum_{k=1}^m C^k (r-s)^{k\left(1 + \frac{\eta-1}{\alpha}\right)} \prod_{j=1}^{k} \BB\left(\frac{\eta}{\alpha} + (j-1)\left(1 + \frac{\eta -1}{\alpha}\right),1 - \frac{1}{\alpha}\right)\right] \\ &\leq (t-s)^{-1} \rho^{-\tilde{\eta}}(t-s,y-x) \sum_{k=1}^{m+1} C^k (t-s)^{(k-1)\left(1 + \frac{\eta-1}{\alpha}\right)}  \prod_{j=1}^{k-1} \BB\left(\frac{\eta}{\alpha} + (j-1)\left(1 + \frac{\eta -1}{\alpha}\right),1 - \frac{1}{\alpha}\right)
\end{align*}
 provided that we choose $C\geq K$ in \eqref{density_time_der_bound_Picard}. This concludes the induction step.

\subsection{Preparatory technical results}

\begin{Lemme}\label{lemma_technical2}
	
	\begin{itemize} \item For all $\gamma \in (0,1]\cap (0,(2\alpha -2)\wedge (\eta + \alpha -1))$, there exists a positive constant $K$ independent of $m$ and such that for all $\mu \in \PP$, $0 \leq s \leq r  \leq T$, $x,v_1,v_2 \in \R^d$, one has \begin{align}\label{b_gradient_linear_derivative_holder_v_Picard}
			&\left\vert \pa_v \del \left[b(r,x,[X^{s,\mu,(m)}_r])\right] (v_1) - \pa_v \del \left[b(r,x,[X^{s,\mu,(m)}_r])\right](v_2) \right\vert  \\ &\notag\leq  K (r-s)^{\frac{\eta -1 - \gamma}{\alpha}}|v_1-v_2|^\gamma \left[ 1+ \sum_{k=1}^m C^k(r-s)^{k\left(1 + \frac{\eta -1}{\alpha}\right)}  \prod_{j=1}^{k-1} \BB\left(1 + \frac{\eta -1 - \gamma}{\alpha} + j \left(1 + \frac{\eta -1}{\alpha}\right), 1 - \frac{1}{\alpha}\right)\right].\end{align}
		
		\item For all $\gamma \in (0,1]\cap(0,(2\alpha -2)\wedge (\eta + \alpha -1))$, there exists a positive constant $K$ independent of $m$ and such that for all $\mu \in \PP$, $0 \leq s \leq r < t \leq T$, $x,y,v_1,v_2 \in \R^d$, one has \begin{align}\label{H_gradient_linear_derivative_holder_v_Picard}
			&\left\vert \pa_v \del \H_{m+1}(\mu,s,r,t,x,y) (v_1) - \pa_v \del \H_{m+1}(\mu,s,r,t,x,y) (v_2) \right\vert  \\ &\notag\leq  K (r-s)^{\frac{\eta -1 - \gamma}{\alpha}}(t-r)^{-\frac{1}{\alpha}}|v_1-v_2|^\gamma \rho^1(t-r,y-x) \\ \notag&\hspace{3cm}\left[ 1+ \sum_{k=1}^m C^k(r-s)^{k\left(1 + \frac{\eta -1}{\alpha}\right)}  \prod_{j=1}^{k-1} \BB\left(1 + \frac{\eta -1 - \gamma}{\alpha} + j \left(1 + \frac{\eta -1}{\alpha}\right), 1 - \frac{1}{\alpha}\right)\right].\end{align}
		
	\end{itemize}
	
\end{Lemme}

\begin{proof}[Proof of Lemma \ref{lemma_technical2}]
	
	First, note that we only need to show \eqref{b_gradient_linear_derivative_holder_v_Picard} since it implies \eqref{H_gradient_linear_derivative_holder_v_Picard} because of \eqref{expression_linear_der_H}. We can write \begin{align*}
		&\Delta_{v_1,v_2}\pa_v \del \left[ b(r,x,[X^{s,\mu,(m)}_r])\right](\cdot)  \\& =  \int_{\R^d} \del b (r,x,[X^{s,\mu,(m)}_r]) (z) \Delta_{v_1,v_2}\pa_x p_m(\mu,s,r,\cdot,z)\,dz \\& \quad+ \int_{\R^{2d}} \del b (r,x,[X^{s,\mu,(m)}_r]) (z)  \Delta_{v_1,v_2} \pa_v \del p_m(\mu,s,r,x',z)(\cdot)  \, dz \, d\mu(x')   \\ & =  \int_{\R^d} \left(\del b (r,x,[X^{s,\mu,(m)}_r]) (z) - \del b (r,x,[X^{s,\mu,(m)}_r]) (v_2) \right) \Delta_{v_1,v_2}\pa_x p_m(\mu,s,r,\cdot,z)\,dz \\ & \quad+ \int_{\R^{2d}} \left( \del b (r,x,[X^{s,\mu,(m)}_r]) (z) - \del b (r,x,[X^{s,\mu,(m)}_r]) (x') \right) \\ & \hspace{5cm} \Delta_{v_1,v_2} \pa_v \del p_m(\mu,s,r,x',z)(\cdot)  \, dz \, d\mu(x')  \\ &=: I_1 + I_2.
	\end{align*}
	Let us first focus on $I_1$. We start by assuming that $|v_1-v_2| \leq (r-s)^{\frac{1}{\alpha}}$. Using the $\eta$-Hölder continuity of $\del b(r,x,\mu)(\cdot)$, \eqref{gradientdensity_holder_Picard}, \eqref{controldensityref} since $|v_1-v_2| \leq (r-s)^{\frac{1}{\alpha}}$ and \eqref{scalingdensityref}, we get that \begin{align*}
		|I_1| &\leq K \int_{\R^d}|z-v_2|^\eta (r-s)^{-\frac{\gamma +1}{\alpha}} |v_1-v_2|^\gamma \rho^1(r-s,z-v_2)\, dz \\ &\leq K (r-s)^{\frac{\eta - \gamma -1}{\alpha}} |v_1 - v_2|^\gamma.
	\end{align*}
	
	Now assume that $|v_1 - v_2|> (r-s)^{\frac{1}{\alpha}}$. The gradient estimate \eqref{density_bound_Picard} yields \begin{align*}
		|I_1| & = \left\vert \int_{\R^d} \left(\del b (r,x,[X^{s,\mu,(m)}_r]) (z) - \del b (r,x,[X^{s,\mu,(m)}_r]) (v_1) \right) \pa_x p_m(\mu,s,r,v_1,z)\,dz \right. \\ &\left.\quad + \int_{\R^d} \left(\del b (r,x,[X^{s,\mu,(m)}_r]) (z) - \del b (r,x,[X^{s,\mu,(m)}_r]) (v_2) \right) \pa_x p_m(\mu,s,r,v_2,z)\,dz \right\vert\\ &\leq K (r-s)^{\frac{\eta -1}{\alpha}} \\ &\leq K (r-s)^{\frac{\eta-1-\gamma}{\alpha}}|v_1-v_2|^\gamma. 
	\end{align*}
	We have thus shown that \[ |I_1| \leq  K (r-s)^{\frac{\eta-1-\gamma}{\alpha}}|v_1-v_2|^\gamma. \]
	
	Then, we have \begin{align*}
		|I_2| &\leq K \int_{\R^{2d}} |z-x'|^\eta (r-s)^{\frac{\eta -1 - \gamma}{\alpha} + 1 - \frac{1}{\alpha}}|v_1-v_2|^\gamma\rho^0(r-s,z-x')\\ & \hspace{2cm}\left( \sum_{k=1}^m C^k(r-s)^{(k-1)\left(1 + \frac{\eta -1}{\alpha}\right)} \prod_{j=1}^{k-1} \BB\left(1 + \frac{\eta -1 - \gamma}{\alpha} + j \left(1 + \frac{\eta -1}{\alpha}\right), 1 - \frac{1}{\alpha}\right)\right)\, dz \, d\mu(x')\\ &\leq  K  (r-s)^{\frac{\eta -1 - \gamma}{\alpha} + 1 + \frac{\eta -1}{\alpha}}|v_1-v_2|^\gamma\\ & \hspace{2cm}\left( \sum_{k=1}^m C^k(r-s)^{(k-1)\left(1 + \frac{\eta -1}{\alpha}\right)} \prod_{j=1}^{k-1} \BB\left(1 + \frac{\eta -1 - \gamma}{\alpha} + j \left(1 + \frac{\eta -1}{\alpha}\right), 1 - \frac{1}{\alpha}\right)\right).
	\end{align*}

\end{proof}

\subsection{Second part of the induction step}
We prove here that the estimates  \eqref{gradient_linear_der_holder_v_Picard}, \eqref{linear_der_holder_v_Picard}, \eqref{linear_der_holder_x_Picard}, \eqref{gradient_linear_der_holder_x_Picard} and \eqref{density_holder_measure_Picard} hold true.\\

 \noindent\textbf{Proof of \eqref{gradient_linear_der_holder_v_Picard}.} It follows from \eqref{H_gradient_linear_derivative_holder_v_Picard}, \eqref{density_bound_Picard} and the convolution inequality \eqref{convolineqdensityref} that for some positive constant $K$ independent of $m$ one has \begin{align*}
	&\left\vert \Delta_{v_1,v_2}\left(\p \otimes \pa_v\del \H_{m+1}\right)(\mu,s,t,x,y) (\cdot) \right\vert\\ &\leq K \int_{s}^t \int_{\R^d}(r-s)^{\frac{\eta -1 - \gamma}{\alpha}}(t-r)^{-\frac{1}{\alpha}}|v_1-v_2|^\gamma \rho^1(t-r,y-x) \rho^0(r-s,z-x) \\ \notag&\hspace{3cm}\left[ 1+ \sum_{k=1}^m C^k(r-s)^{k\left(1 + \frac{\eta -1}{\alpha}\right)}  \prod_{j=1}^{k-1} \BB\left(1 + \frac{\eta -1 - \gamma}{\alpha} + j \left(1 + \frac{\eta -1}{\alpha}\right), 1 - \frac{1}{\alpha}\right)\right] \, dz\,dr \\ &\leq K (t-s)^{\frac{\eta - 1 - \gamma}{\alpha} + 1 - \frac{1}{\alpha}}|v_1-v_2|^\gamma \rho^0(t-s,y-x)   \left[ 1+\sum_{k=1}^m C^k(r-s)^{k\left(1 + \frac{\eta -1}{\alpha}\right)} \right.\\ &\left. \hspace{6cm} \prod_{j=1}^{k} \BB\left(1 + \frac{\eta -1 - \gamma}{\alpha} + j \left(1 + \frac{\eta -1}{\alpha}\right), 1 - \frac{1}{\alpha}\right)\right].
\end{align*}

Starting from the representation formula \eqref{representation_gradient_linear_der}, using this inequality, the bound \eqref{H^k_m_bound} for $\H_{m+1}^k$, the convolution inequality \eqref{convolineqdensityref} and finally summing over $k \geq 0$, we find that there exists a constant $K$ independent of $m$ such that 

\begin{align*}
	&\left\vert \Delta_{v_1,v_2}\pa_v\del p_{m+1}(\mu,s,t,x,y)(\cdot)\right\vert \\ &\leq \sum_{k=0}^{\infty} \left\vert \Delta_{v_1,v_2} \left(p_{m+1}\otimes \pa_v \del \H_{m+1}  \otimes \H_{m+1}^k\right) (\mu,s,t,x,y)(\cdot) \right\vert\\ &\leq K (t-s)^{\frac{\eta - 1 - \gamma}{\alpha} + 1 - \frac{1}{\alpha}}|v_1-v_2|^\gamma \rho^0(t-s,y-x)   \left[ 1+\sum_{k=1}^m C^k(r-s)^{k\left(1 + \frac{\eta -1}{\alpha}\right)} \right.\\ &\left. \hspace{6cm} \prod_{j=1}^{k} \BB\left(1 + \frac{\eta -1 - \gamma}{\alpha} + j \left(1 + \frac{\eta -1}{\alpha}\right), 1 - \frac{1}{\alpha}\right)\right]\\ &\leq (t-s)^{\frac{\eta -1 - \gamma}{\alpha} + 1 - \frac{1}{\alpha}}|v_1-v_2|^\gamma\rho^0(t-s,y-x)\\ &\hspace{2cm}\left( \sum_{k=1}^{m+1} C^k(t-s)^{(k-1)\left(1 + \frac{\eta -1}{\alpha}\right)} \prod_{j=1}^{k-1} \BB\left(1 + \frac{\eta -1 - \gamma}{\alpha} + j \left(1 + \frac{\eta -1}{\alpha}\right), 1 - \frac{1}{\alpha}\right)\right).
\end{align*}
provided that we choose $C\geq K$ in \eqref{gradient_linear_der_holder_v_Picard}. It ends the proof of the induction step for \eqref{gradient_linear_der_holder_v_Picard}.\\

\noindent\textbf{Proof of \eqref{linear_der_holder_v_Picard}.} Let us first assume that $|v_1-v_2| \geq (t-s)^{\frac{1}{\alpha}}.$ Using \eqref{linear_der_bound_Picard} and the fact the series appearing in this bound has a limit when $m$ tends to infinity, one has in this case \begin{align*}
	\left\vert \Delta_{v_1,v_2} \del p_m(\mu,s,t,x,y)(\cdot)\right\vert &\leq C (t-s)^{1 - \frac{1}{\alpha}} \rho^0(t-s,y-x) \\ &\leq C(t-s)^{- \frac{\gamma}{\alpha}+1 - \frac{1}{\alpha}} |v_1-v_2|^\gamma \rho^0(t-s,y-x).
\end{align*}
In the case where $|v_1-v_2|<(t-s)^{\frac{1}{\alpha}}$, the mean value theorem and \eqref{gradient_linear_der_bound_Picard} yield \begin{align*}
	\left\vert \Delta_{v_1,v_2} \del p_m(\mu,s,t,x,y)(\cdot)\right\vert &\leq C (t-s)^{\frac{\eta-1}{\alpha}+1 - \frac{1}{\alpha}}|v_1-v_2| \rho^0(t-s,y-x) \\ &\leq C(t-s)^{\frac{\eta}{\alpha}- \frac{\gamma}{\alpha}+1 - \frac{1}{\alpha}} |v_1-v_2|^\gamma \rho^0(t-s,y-x) \\ &\leq C(t-s)^{- \frac{\gamma}{\alpha}+1 - \frac{1}{\alpha}} |v_1-v_2|^\gamma \rho^0(t-s,y-x).
\end{align*}

\noindent\textbf{Proof of \eqref{linear_der_holder_x_Picard}.} We first note that since \eqref{linear_der_bound_Picard} has been proved before, \eqref{H_linear_derivative_bound_Picard} ensures that there exists a positive constant $C$ such that for all $m \geq 1$ \begin{equation*}
	\left\vert \del \H_{m} (\mu,s,r,t,x,y)(v) \right\vert \leq C (t-r)^{-\frac{1}{\alpha}} \rho^1(t-r,y-x).
\end{equation*}
Using this inequality, \eqref{gradientdensity_holder_Picard} and the convolution inequality \eqref{convolineqdensityref}, we obtain that \begin{align*}
	&\left\vert \Delta_{x_1,x_2} \left[p_m \otimes \del \H_{m}\right](\mu,s,r,t,\cdot,y)(v) \right\vert \\ &= \left\vert \int_r^t\int_{\R^d} \Delta_{x_1,x_2}p_m(\mu,r,r',\cdot,z) \del \H_m (\mu,s,r',t,z,y)\,dz\,dr'\right\vert \\ &\leq C\int_r^t \int_{\R^d} (r'-r)^{-\frac{\gamma}{\alpha}}|x_1-x_2|^\gamma \left[\rho^0(r'-r,z-x_1) + \rho^0(r'-r,z-x_2)\right] (t-r')^{-\frac{1}{\alpha}} \rho^1(t-r',y-z)\, dz\, dr' \\ &\leq C (t-r)^{1 - \frac{1+\gamma }{\alpha}} |x_1-x_2|^\gamma \left[\rho^0(t-r,y-x_1) + \rho^0(t-r,y-x_2)\right].
\end{align*}

We conclude by the representation formula \eqref{representation_linear_der} and \eqref{H^k_m_bound} (since the series appearing is convergent) that we have \begin{align*}
	\left\vert \Delta_{x_1,x_2} \del p_m(\mu,s,t,\cdot,y)(v) \right\vert &\leq \sum_{k=0}^{\infty} \left\vert \Delta_{x_1,x_2} \left[ p_m \otimes \del \H_m \otimes \H_m^k \right] (\mu,s,t,\cdot,y)(v)\right\vert \\ &\leq C (t-s)^{1 - \frac{1+\gamma }{\alpha}} |x_1-x_2|^\gamma \left[\rho^0(t-s,y-x_1) + \rho^0(t-s,y-x_2)\right].
\end{align*}

\noindent\textbf{Proof of \eqref{gradient_linear_der_holder_x_Picard}.} We first note that since \eqref{gradient_linear_der_bound_Picard} has been proved before, \eqref{H_gradient_linear_derivative_bound_Picard} ensures that there exists a positive constant $C$ such that for all $m \geq 1$ \begin{equation*}
	\left\vert \pa_v\del \H_{m} (\mu,s,r,t,x,y)(v) \right\vert \leq C (r-s)^{\frac{\eta-1}{\alpha}}(t-r)^{-\frac{1}{\alpha}} \rho^1(t-r,y-x).
\end{equation*}
Using this inequality, \eqref{gradientdensity_holder_Picard}, the convolution inequality \eqref{convolineqdensityref} and the fact that $\gamma < \alpha - 1 - \eta$, we obtain that \begin{align*}
	&\left\vert \Delta_{x_1,x_2} \left[p_m \otimes \pa_v\del \H_{m}\right](\mu,s,s,t,\cdot,y)(v) \right\vert \\ &= \left\vert \int_s^t\int_{\R^d} \Delta_{x_1,x_2}p_m(\mu,s,r,\cdot,z) \pa_v\del \H_m (\mu,s,r,t,z,y)\,dz\,dr\right\vert \\ &\leq C\int_s^t \int_{\R^d} (r-s)^{-\frac{\gamma}{\alpha}}|x_1-x_2|^\gamma \left[\rho^0(r-s,z-x_1) + \rho^0(r-s,z-x_2)\right] (r-s)^{\frac{\eta-1}{\alpha}}(t-r)^{-\frac{1}{\alpha}} \rho^1(t-r,y-z)\, dz\, dr \\ &\leq C(t-s)^{\frac{\eta-1-\gamma}{\alpha}+1 - \frac{1}{\alpha}}  |x_1-x_2|^\gamma \left[\rho^0(t-r,y-x_1) + \rho^0(t-r,y-x_2)\right].
\end{align*}

We conclude by the representation formula \eqref{representation_gradient_linear_der} and \eqref{H^k_m_bound} that we have \begin{align*}
	\left\vert \Delta_{x_1,x_2} \pa_v\del p_m(\mu,s,t,\cdot,y)(v) \right\vert &\leq \sum_{k=0}^{\infty} \left\vert \Delta_{x_1,x_2} \left[ p_m \otimes \pa_v\del \H_m \otimes \H_m^k \right] (\mu,s,t,\cdot,y)(v)\right\vert \\ &\leq C(t-s)^{\frac{\eta-1-\gamma}{\alpha}+1 - \frac{1}{\alpha}}  |x_1-x_2|^\gamma \left[\rho^0(t-r,y-x_1) + \rho^0(t-r,y-x_2)\right].
\end{align*}

\noindent\textbf{Proof of \eqref{density_holder_measure_Picard}.} Let us first assume that $W_1(\mu_1,\mu_2) \leq (t-s)^{\frac{1}{\alpha}}$. In this case, by definition of the linear derivative, the Kantorovich-Rubinstein theorem, and \eqref{gradient_linear_der_bound_Picard}, we get 

\begin{align*}
	\left\vert\Delta_{\mu_1,\mu_2} p_m(\cdot,s,t,x,y)\right\vert &= \left\vert\int_0^1 \int_{\R^d} \del p_m (\lambda \mu_1 + (1-\lambda)\mu_2,s,t,x,y)(v) \, d(\mu_1 - \mu_2)(v) \, d\lambda \right\vert \\ &\leq \sup_{\lambda \in [0,1], \, v\in \R^d} \left\vert \pa_v \del p_m  (\lambda \mu_1 + (1-\lambda)\mu_2,s,t,x,y)(v)\right\vert  W_1(\mu_1,\mu_2)\\ &\leq C (t-s)^{\frac{\eta-1}{\alpha} + 1 - \frac{1}{\alpha}}W_1(\mu_1,\mu_2) \rho^0(t-s,y-x)\\ & \leq C(t-s)^{1 - \frac{1+\gamma }{\alpha}}W_1^\gamma(\mu_1,\mu_2)\rho^0(t-s,y-x).
\end{align*}

In the case where $W_1(\mu_1,\mu_2)> (t-s)^\frac{1}{\alpha}$, the parametrix expansion \eqref{representationdensityparametrix_Picard}, \eqref{gradientestimatestable} and \eqref{Phi_m_bound} yield

\begin{align*}
	\left\vert \Delta_{\mu_1,\mu_2} p_m(\cdot,s,t,x,y) \right\vert &\leq  \left\vert \Delta_{\mu_1,\mu_2} \left[\p \otimes \Phi_m\right](\cdot,s,t,x,y) \right\vert \\ &\leq  \left\vert  \p \otimes \Phi_m(\mu_1,s,t,x,y) \right\vert +  \left\vert  \p \otimes \Phi_m(\mu_2,s,t,x,y) \right\vert \\ &\leq C (t-s)^{
		1 - \frac{1}{\alpha}}\rho^0(t-s,y-x) \\ &\leq C(t-s)^{1 - \frac{1+\gamma }{\alpha}}W_1^\gamma(\mu_1,\mu_2)\rho^0(t-s,y-x).
\end{align*}
\subsection{Preparatory technical results}

\begin{Lemme}\label{lemma_technical3}
	\begin{itemize}
	\item For any $\gamma \in [\eta ,1]$, there exists a positive constant $K$ such that for all $m\geq 1$, $\mu_1,\mu_2 \in \PP$, $ 0 \leq s <t \leq T$, $ x \in \R^d$
	\begin{equation}\label{b_holder_mu_Picard}
		|b(t,x,[X^{s,\mu_1,(m)}_t]) - b(t,x,[X^{s,\mu_2,(m)}_t])| \leq K (t-s)^{\frac{\eta - \gamma}{\alpha}} W_1^\gamma(\mu_1,\mu_2).
	\end{equation}
\item For any $\gamma \in [\eta ,1]$, there exists a positive constant $K$ such that for all $m\geq 1$, $\mu_1,\mu_2 \in \PP$, $ 0 \leq s\leq r <t \leq T$, $ x,y \in \R^d$
\begin{equation}\label{H_holder_mu_Picard}
	|\H_m(\mu_1,s,r,t,x,y) - \H_m(\mu_2,s,r,t,x,y)| \leq K (r-s)^{\frac{\eta - \gamma}{\alpha}}(t-r)^{-\frac{1}{\alpha}} W_1^\gamma(\mu_1,\mu_2)\rho^{1}(t-r,y-x).
\end{equation}

\item For any $\gamma \in [\eta ,1]$, there exists a positive constant $K$ such that for all $m\geq 1$, $\mu_1,\mu_2 \in \PP$, $ 0 \leq s\leq r <t \leq T$, $ x,y \in \R^d$
\begin{equation}\label{Phi_holder_mu_Picard}
	|\Phi_m(\mu_1,s,r,t,x,y) - \Phi_m(\mu_2,s,r,t,x,y)| \leq K (r-s)^{\frac{\eta - \gamma}{\alpha}}(t-r)^{-\frac{1}{\alpha}} W_1^\gamma(\mu_1,\mu_2)\rho^{1}(t-r,y-x).
\end{equation}
	\end{itemize}
\end{Lemme}

\begin{proof}[Proof of Lemma \ref{lemma_technical3}]
	\noindent\textbf{Proof of \eqref{b_holder_mu_Picard} and \eqref{H_holder_mu_Picard}.} By definition of the linear derivative and denoting by $m_l := l [X^{s,\mu_1,(m)}_t] + (1- l) [X^{s,\mu_2,(m)}_t]$, we have 
	
	\begin{align*}
		&b(t,x,[X^{s,\mu_1,(m)}_t]) - b(t,x,[X^{s,\mu_2,(m)}_t])\\ &= \int_0^1 \int_{\R^d} \del b (t,x,m_l)(y') [p_m(\mu_1,s,t,y') - p_m(\mu_2,s,t,y')] \, dy'\, dl \\ &= \int_0^1 \int_{\R^{2d}} \del b (t,x,m_l)(y') p_m(\mu_1,s,t,x',y') \, dy'\,d(\mu_1-\mu_2)(x') dl \\ &\quad + \int_0^1 \int_{\R^{2d}} \del b (t,x,m_l)(y') [p_m(\mu_1,s,t,x',y') - p_m(\mu_2,s,t,x',y')] \, dy'\,d\mu_2(x') dl \\ &=: I_1 + I_2.
	\end{align*}

For $I_1$, we need to control, for $x',x'' \in \R^d$ \[\int_{\R^{d}} \del b (t,x,m_l)(y') [p_m(\mu_1,s,t,x',y') - p_m(\mu_1,s,t,x'',y')] \, dy'. \]

In the case where $|x'-x''| \leq (t-s)^{\frac{1}{\alpha}}$, we write

\begin{align*}
& \int_{\R^{d}} \del b (t,x,m_l)(y') [p_m(\mu_1,s,t,x',y') - p_m(\mu_1,s,t,x'',y')] \, dy' \\ & = \int_0^1\int_{\R^{d}}  \left[\del b (t,x,m_l)(y') - \del b (t,x,m_l)(l'x' + (1- l')x'')\right]\\ & \hspace{4cm}\pa_x p_m(\mu_1,s,t,l'x' + (1- l')x'',y') \cdot (x'-x'') \, dy' \, dl'.
\end{align*}

Using the $\eta$-Hölder continuity of $\del b$ and \eqref{density_bound_Picard}, we obtain that 

\begin{align*}
	& \left\vert\int_{\R^{d}} \del b (t,x,m_l)(y') [p_m(\mu_1,s,t,x',y') - p_m(\mu_1,s,t,x'',y')] \, dy' \right\vert\\ & \leq K(t-s)^{\frac{\eta-1}{\alpha}}|x'-x''| \\ &\leq K (t-s)^{\frac{\eta - \gamma}{\alpha}}|x'-x''|^{\gamma}.
\end{align*}

Assume now that $|x'-x''|>(t-s)^{\frac{1}{\alpha}}$. One can write \begin{align*}
	&\int_{\R^{d}} \del b (t,x,m_l)(y') [p_m(\mu_1,s,t,x',y') - p_m(\mu_1,s,t,x'',y')] \, dy' \\&= \del b (t,x,m_l)(x') - \del b (t,x,m_l)(x'') \\ &\quad + 	\int_{\R^{d}} \left[\del b (t,x,m_l)(y') - \del b (t,x,m_l)(x')  \right] p_m(\mu_1,s,t,x',y') \, dy' \\ &\quad +  	\int_{\R^{d}} \left[\del b (t,x,m_l)(y') - \del b (t,x,m_l)(x'')  \right] p_m(\mu_1,s,t,x'',y') \, dy'
\end{align*}

It follows from the $\eta$-Hölder continuity of $\del b$, \eqref{density_bound_Picard} and the space-time inequality \eqref{scalingdensityref} that \begin{align*}
	&\int_{\R^{d}} \del b (t,x,m_l)(y') [p_m(\mu_1,s,t,x',y') - p_m(\mu_1,s,t,x'',y')] \, dy' \\& \leq K(|x'-x''|^\eta + (t-s)^{\frac{\eta}{\alpha}}) \\ &\leq K(|x'-x''|^\eta + (t-s)^{\frac{\eta}{\alpha}}) \\&\leq K(t-s)^{\frac{\eta-\gamma}{\alpha}}|x'-x''|^\gamma,
\end{align*}
since $\gamma \geq \eta.$ Jensen's inequality yields \[|I_1| \leq K(t-s)^{\frac{\eta - \gamma}{\alpha}}W_1^\gamma (\mu_1,\mu_2).\]

It remains to study $I_2$ which can be rewritten as \begin{multline}
I_2 =\int_0^1 \int_{\R^{2d}} \left[\del b (t,x,m_l)(y') -\del b (t,x,m_l)(x')\right] [p_m(\mu_1,s,t,x',y') - p_m(\mu_2,s,t,x',y')] \, dy'\,d\mu_2(x') dl.\end{multline}
 Thanks to \eqref{density_holder_measure_Picard}, \eqref{density_bound_Picard} and the space-time inequality \eqref{scalingdensityref}, one has since $\alpha \in (1,2)$ \begin{align*}
	|I_2| &\leq  K (t-s)^{ 1 - \frac{1+\gamma}{\alpha} + \frac{\eta}{\alpha}}W_1^\gamma(\mu_1,\mu_2) \\ &\leq K(t-s)^{ \frac{\eta-\gamma}{\alpha}}W_1^\gamma(\mu_1,\mu_2).
\end{align*}
This concludes the proof of \eqref{b_holder_mu_Picard}. The proof of \eqref{H_holder_mu_Picard} immediately follows from the definition \eqref{def_proxy_kernel_Picard} of $\H_m$ and \eqref{gradientestimatestable}.\\

\noindent\textbf{Proof of \eqref{Phi_holder_mu_Picard}.} We prove by induction of $k\geq 2$ that there exists a constant $K$ independent of $k$ such that for any $k \geq 2$, $m \geq 1$, $\mu_1,\mu_2 \in \PP,$ $0 \leq s \leq r < t \leq T$, $x,y \in \R^d$ \begin{multline}\label{eq_proof_phi_holder_mu_picard}
	|\Delta_{\mu_1,\mu_2}\H_m^k(\cdot,s,r,t,x,y)| \leq (r-s)^{\frac{\eta- \gamma}{\alpha}}(t-r)^{-\frac{1}{\alpha}}  W_1^\gamma(\mu_1,\mu_2)\rho^1(t-r,y-x) \\ \sum_{i=1}^{k-1} C^i (t-r)^{i\left(1 - \frac{1}{\alpha}\right)} \prod_{j=1}^i \BB \left(j \left( 1 - \frac{1}{\alpha}\right), 1 - \frac{1}{\alpha}\right).
\end{multline}
We do not prove the base case $k=2$ since it relies on the same computations as the induction step using \eqref{H_holder_mu_Picard}. For the induction step, assume that \eqref{eq_proof_phi_holder_mu_picard} holds for $\H_m^k$ and let us prove it for $\H_m^{k+1}$. One has 

\begin{align*}
	|\Delta_{\mu_1,\mu_2} \H_m^{k+1}(\cdot,s,r,t,x,y)| &\leq \left\vert \int_r^t \int_{\R^d} \H_m(\mu_1,s,r,r',x,z) \Delta_{\mu_1,\mu_2} \H_m^k(\cdot,s,r',t,z,y) \, dz\,dr'\right\vert \\ &\quad+ \left\vert  \int_r^t \int_{\R^d} \Delta_{\mu_1,\mu_2}\H_m(\cdot,s,r,r',x,z) \H_m^k(\mu_2,s,r',t,z,y) \, dz\,dr'\right\vert \\ &=: I_1 + I_2.
\end{align*}
By \eqref{H_holder_mu_Picard} and the bound of $\H_m^k$ \eqref{H^k_m_bound} (note that the series appearing in this bound in convergent and thus can be bounded independently of $k$ and $m$), we have \begin{align*}
	I_2 &\leq K \int_r^t \int_{\R^d} (r-s)^{\frac{\eta - \gamma}{\alpha}} (r'-r)^{-\frac{1}{\alpha}} W_1^\gamma(\mu_1,\mu_2) \rho^1(r'-r,z-x) (t-r')^{-\frac{1}{\alpha}}\rho^1(t-r',y-z) \, dz\,dr' \\ &\leq K (r-s)^{\frac{\eta - \gamma}{\alpha}}(t-r)^{-\frac{1}{\alpha} +1 - \frac{1}{\alpha}} W_1^\gamma(\mu_1,\mu_2) \rho^1(t-r,y-x)\BB \left(1 - \frac{1}{\alpha}, 1 - \frac{1}{\alpha}\right).
\end{align*}
By induction and using \eqref{H^k_m_bound}, one has \begin{align*}
	I_1 &\leq K\int_r^t \int_{\R^d} (r'-r)^{-\frac{1}{\alpha}} \rho^1(r'-r,z-x) (r'-s)^{\frac{\eta- \gamma}{\alpha}}(t-r')^{-\frac{1}{\alpha}}  W_1^\gamma(\mu_1,\mu_2)\rho^1(t-r',y-z) \\ &\hspace{6cm} \sum_{i=1}^{k-1} C^i (t-r')^{i\left(1 - \frac{1}{\alpha}\right)} \prod_{j=1}^i \BB \left(j \left( 1 - \frac{1}{\alpha}\right), 1 - \frac{1}{\alpha}\right) \, dz\,dr' \\  &\leq K (r-s)^{\frac{\eta- \gamma}{\alpha}}\int_r^t \int_{\R^d} (r'-r)^{-\frac{1}{\alpha}} \rho^1(r'-r,z-x) (t-r')^{-\frac{1}{\alpha}}  W_1^\gamma(\mu_1,\mu_2)\rho^1(t-r',y-z) \\ &\hspace{6cm} \sum_{i=1}^{k-1} C^i (t-r')^{i\left(1 - \frac{1}{\alpha}\right)} \prod_{j=1}^i \BB \left(j \left( 1 - \frac{1}{\alpha}\right), 1 - \frac{1}{\alpha}\right) \, dz\,dr' \\ &\leq K (r-s)^{\frac{\eta- \gamma}{\alpha}} (t-r)^{-\frac{1}{\alpha}} W_1^\gamma(\mu_1,\mu_2)\rho^1(t-r,y-x) \sum_{i=1}^{k-1} C^i (t-r)^{(i+1)\left(1 - \frac{1}{\alpha}\right)} \prod_{j=1}^{i+1} \BB \left(j \left( 1 - \frac{1}{\alpha}\right), 1 - \frac{1}{\alpha}\right)
\end{align*}

Summing the bounds obtained for $I_1$ and $I_2$, we have proved that \begin{multline*}
	|\Delta_{\mu_1,\mu_2}\H_m^{k+1}(\cdot,s,r,t,x,y)| \leq (r-s)^{\frac{\eta- \gamma}{\alpha}}(t-r)^{-\frac{1}{\alpha}}  W_1^\gamma(\mu_1,\mu_2)\rho^1(t-r,y-x) \\ \sum_{i=1}^{k} C^i (t-r)^{i\left(1 - \frac{1}{\alpha}\right)} \prod_{j=1}^i \BB \left(j \left( 1 - \frac{1}{\alpha}\right), 1 - \frac{1}{\alpha}\right),
\end{multline*}
if we choose $C\geq K$ since the beginning. Finally, summing \eqref{eq_proof_phi_holder_mu_picard} over $k \geq 1$ yields \eqref{Phi_holder_mu_Picard}.

\end{proof}

\begin{Lemme}\label{lemma_technical4}
	\begin{itemize}
		
		\item  For all $\gamma \in [\eta,1]$, there exists $K>0$ such that for all $m \geq 1$, $\mu_1,\mu_2 \in \PP$, $ 0 \leq s < t \leq T$, $x,v \in \R^d$

		\begin{equation}\label{b_lin_der_holder_mu_Picard}
			\left\vert  \del  b(t,x,[X^{s,\mu_1,(m)}_t])(v) -  \del  b(t,x,[X^{s,\mu_2,(m)}_t])(v) \right\vert \leq K(t-s)^{\frac{\eta - \gamma}{\alpha}} W_1^\gamma (\mu_1,\mu_2).
		\end{equation}

		\item For all $\gamma \in (0,1]$, there exists $K>0$ such that for all $m \geq 1$, $\mu_1,\mu_2 \in \PP$, $ 0 \leq s < t \leq T$, $x,v \in \R^d$

		\begin{multline}\label{b_lin_der_holder_mu_Picard_2}
			\left\vert  \del \left[ b(t,x,[X^{s,\mu_1,(m)}_t])\right](v) -  \del \left[ b(t,x,[X^{s,\mu_2,(m)}_t])\right](v) \right\vert \leq K(t-s)^{-\frac{ \gamma}{\alpha}}W_1^\gamma (\mu_1,\mu_2) \\ +K \int_{\R^{2d}}  \left\vert \Delta_{\mu_1,\mu_2} \del p_m(\cdot,s,t,x',y')\right\vert \, dy'\,d\mu_2(x').
		\end{multline}
	
		\item For all $\gamma \in (0,1]\cap(0,\eta + \alpha -1)$, there exists $K>0$ such that for all $m \geq 1$, $\mu_1,\mu_2 \in \PP$, $ 0 \leq s < t \leq T$, $x,v \in \R^d$ 
	
	\begin{multline}\label{b_gradient_lin_der_holder_mu_Picard_2}
		\left\vert  \pa_v\del \left[ b(t,x,[X^{s,\mu_1,(m)}_t])\right](v) -  \pa_v\del \left[ b(t,x,[X^{s,\mu_2,(m)}_t])\right](v) \right\vert \leq K(t-s)^{\frac{\eta - \gamma-1}{\alpha}}W_1^\gamma(\mu_1,\mu_2) \\ +K \int_{\R^{2d}} (1 \wedge |x'-y'|^\eta) \left\vert \Delta_{\mu_1,\mu_2} \pa_v\del p_m(\cdot,s,t,x',y')\right\vert \, dy'\,d\mu_2(x').
	\end{multline}

	\item For all $\gamma \in (0,1]$, there exists $K>0$ such that for all $m \geq 1$, $\mu_1,\mu_2 \in \PP$, $ 0 \leq s < t \leq T$, $x,v \in \R^d$

\begin{multline}\label{H_lin_der_holder_mu_Picard}
	\left\vert  \del \H_{m+1}(\mu_1,s,r,t,x,y)(v) -  \del \H_{m+1}(\mu_2,s,r,t,x,y)(v) \right\vert \leq K(r-s)^{-\frac{ \gamma}{\alpha}}(t-r)^{-\frac{1}{\alpha}}W_1^\gamma (\mu_1,\mu_2)\rho^1(t-r,y-x) \\ +K(t-r)^{-\frac{1}{\alpha}} \int_{\R^{2d}}  \left\vert \Delta_{\mu_1,\mu_2} \del p_m(\cdot,s,r,x',y')\right\vert \, dy'\,d\mu_2(x') \rho^1(t-r,y-x).
\end{multline}

\item For all $\gamma \in (0,1]\cap(0,\eta + \alpha -1)$, there exists $K>0$ such that for all $m \geq 1$, $\mu_1,\mu_2 \in \PP$, $ 0 \leq s < t \leq T$, $x,v \in \R^d$ 

\begin{align}\label{H_gradient_lin_der_holder_mu_Picard}
	\notag&\left\vert  \pa_v\del\H_{m+1}(\mu_1,s,r,t,x,y)(v)- \pa_v\del\H_{m+1}(\mu_1,s,r,t,x,y)(v) \right\vert \\ &\leq K(r-s)^{\frac{\eta - \gamma-1}{\alpha}}(t-r)^{-\frac{1}{\alpha}}W_1^\gamma(\mu_1,\mu_2) \rho^1(t-r,y-x) \\\notag &\quad +K(t-r)^{-\frac{1}{\alpha}} \int_{\R^{2d}} (1 \wedge |x'-y'|^\eta) \left\vert \Delta_{\mu_1,\mu_2} \pa_v\del p_m(\cdot,s,r,x',y')\right\vert \, dy'\,d\mu_2(x') \rho^1(t-r,y-x).
\end{align}
	\end{itemize}
\end{Lemme}

\begin{proof}[Proof of Lemma \ref{lemma_technical4}]
	
	\noindent\textbf{Proof of \eqref{b_lin_der_holder_mu_Picard}.} By definition of the linear derivative and denoting by $m_l := l [X^{s,\mu_1,(m)}_t] + (1- l) [X^{s,\mu_2,(m)}_t]$, we have 
	
	\begin{align*}
		&\del b(t,x,[X^{s,\mu_1,(m)}_t])(v) - \del b(t,x,[X^{s,\mu_2,(m)}_t])(v)\\ &= \int_0^1 \int_{\R^d} \dell b (t,x,m_l)(v,y') [p_m(\mu_1,s,t,y') - p_m(\mu_2,s,t,y')] \, dy'\, dl \\ &= \int_0^1 \int_{\R^{2d}} \dell b (t,x,m_l)(v,y') p_m(\mu_1,s,t,x',y') \, dy'\,d(\mu_1-\mu_2)(x') dl \\ &\quad + \int_0^1 \int_{\R^{2d}} \dell b (t,x,m_l)(v,y') [p_m(\mu_1,s,t,x',y') - p_m(\mu_2,s,t,x',y')] \, dy'\,d\mu_2(x') dl \\ &=: I_1 + I_2.
	\end{align*}

	For $I_1$, we need to control, for $x',x'' \in \R^d$ \[\int_{\R^{d}} \dell b (t,x,m_l)(v,y') [p_m(\mu_1,s,t,x',y') - p_m(\mu_1,s,t,x'',y')] \, dy'. \]
	
	In the case where $|x'-x''| \leq (t-s)^{\frac{1}{\alpha}}$, we write
	
	\begin{align*}
		& \int_{\R^{d}} \dell b (t,x,m_l)(v,y') [p_m(\mu_1,s,t,x',y') - p_m(\mu_1,s,t,x'',y')] \, dy' \\ & = \int_0^1\int_{\R^{d}}  \left[\dell b (t,x,m_l)(v,y') - \dell b (t,x,m_l)(v,l'x' + (1- l')x'')\right]\\ & \hspace{4cm}\pa_x p_m(\mu_1,s,t,l'x' + (1- l')x'',y') \cdot (x'-x'') \, dy' \, dl'.
	\end{align*}
	
	Using the $\eta$-Hölder continuity of $\dell b(t,x,\mu)(v,\cdot)$ and \eqref{density_bound_Picard}, we obtain that 
	
	\begin{align*}
		& \left\vert\int_{\R^{d}} \dell b (t,x,m_l)(v,y') [p_m(\mu_1,s,t,x',y') - p_m(\mu_1,s,t,x'',y')] \, dy' \right\vert\\ & \leq K(t-s)^{\frac{\eta-1}{\alpha}}|x'-x''| \\ &\leq K (t-s)^{\frac{\eta - \gamma}{\alpha}}|x'-x''|^{\gamma}.
	\end{align*}
	
	Assume now that $|x'-x''|>(t-s)^{\frac{1}{\alpha}}$. One can write \begin{align*}
		&\int_{\R^{d}} \dell b (t,x,m_l)(v,y') [p_m(\mu_1,s,t,x',y') - p_m(\mu_1,s,t,x'',y')] \, dy' \\&= \dell b (t,x,m_l)(v,x') - \dell b (t,x,m_l)(v,x'') \\ &\quad + 	\int_{\R^{d}} \left[\dell b (t,x,m_l)(v,y') - \dell b (t,x,m_l)(v,x')  \right] p_m(\mu_1,s,t,x',y') \, dy' \\ &\quad +  	\int_{\R^{d}} \left[\dell b (t,x,m_l)(v,y') - \dell b (t,x,m_l)(v,x'')  \right] p_m(\mu_1,s,t,x'',y') \, dy'
	\end{align*}
	
	It follows from the $\eta$-Hölder continuity of $\dell b(t,x,\mu)(v,\cdot)$, \eqref{density_bound_Picard} and the space-time inequality \eqref{scalingdensityref} that \begin{align*}
		&\int_{\R^{d}} \dell b (t,x,m_l)(v,y') [p_m(\mu_1,s,t,x',y') - p_m(\mu_1,s,t,x'',y')] \, dy' \\& \leq K(|x'-x''|^\eta + (t-s)^{\frac{\eta}{\alpha}}) \\ &\leq K(|x'-x''|^\eta + (t-s)^{\frac{\eta}{\alpha}}) \\&\leq K(t-s)^{\frac{\eta-\gamma}{\alpha}}|x'-x''|^\gamma,
	\end{align*}
	since $\gamma \geq \eta.$ Jensen's inequality yields \[|I_1| \leq K(t-s)^{\frac{\eta - \gamma}{\alpha}}W_1^\gamma (\mu_1,\mu_2).\]
	
	It remains to study $I_2$ which can be rewritten as \begin{multline}
		I_2 =\int_0^1 \int_{\R^{2d}} \left[\dell b (t,x,m_l)(v,y') -\dell b (t,x,m_l)(v,x')\right] [p_m(\mu_1,s,t,x',y') - p_m(\mu_2,s,t,x',y')] \, dy'\,d\mu_2(x') dl.\end{multline}
	Thanks to \eqref{density_holder_measure_Picard}, \eqref{density_bound_Picard} and the space-time inequality \eqref{scalingdensityref}, one has since $\alpha \in (1,2)$ \begin{align*}
		|I_2| &\leq  K (t-s)^{ 1 - \frac{1+\gamma}{\alpha} + \frac{\eta}{\alpha}}W_1^\gamma(\mu_1,\mu_2) \\ &\leq K(t-s)^{ \frac{\eta-\gamma}{\alpha}}W_1^\gamma(\mu_1,\mu_2).
	\end{align*}

This concludes the proof of \eqref{b_lin_der_holder_mu_Picard}.\\

\noindent\textbf{Proof of \eqref{b_lin_der_holder_mu_Picard_2}.} We begin to treat the case where $W_1(\mu_1,\mu_2) \geq (t-s)^{\frac{1}{\alpha}}$. Using \eqref{b_linear_derivative_bound_Picard} and noting that the series appearing in the bound is convergent, one has

\begin{align*}
	\left\vert \Delta_{\mu_1,\mu_2} \del \left[ b(t,x,[X^{s,\cdot,(m)}_t])\right](v) \right\vert &\leq  \left\vert  \del \left[ b(t,x,[X^{s,\mu_1,(m)}_t])\right](v) \right\vert + \left\vert  \del \left[ b(t,x,[X^{s,\mu_2,(m)}_t])\right](v) \right\vert \\ &\leq K \\ &\leq K(t-s)^{-\frac{\gamma}{\alpha}} W_1(\mu_1,\mu_2)^\gamma.
\end{align*}

We now focus on the case where $W_1(\mu_1,\mu_2)< (t-s)^\frac{1}{\alpha}$. Using \eqref{lemma_technical_eq1}, one has 
	
	\begin{align*}
		&\Delta_{\mu_1,\mu_2} \del \left[ b(t,x,[X^{s,\cdot,(m)}_t])\right](v)  \\& =  \int_{\R^d} \Delta_{\mu_1,\mu_2}\del b (t,x,[X^{s,\cdot,(m)}_t]) (y) p_m(\mu_1,s,t,v,y)\,dy \\ & \quad +\int_{\R^{d}} \del b (t,x,[X^{s,\mu_2,(m)}_t]) (y) \Delta_{\mu_1,\mu_2} p_m(\cdot,s,t,v,y)\,dy \\ &\quad +\int_{\R^{2d}} \Delta_{\mu_1,\mu_2}\del b (t,x,[X^{s,\cdot,(m)}_t]) (y) \del p_m(\mu_1,s,t,x',y)(v) \, dy \, d\mu_1(x') \\ & \quad +\int_{\R^{2d}} \del b (t,x,[X^{s,\mu_2,(m)}_t]) (y) \del p_m(\mu_1,s,t,x',y)(v)\,dy\, d(\mu_1-\mu_2)(x') \\ & \quad +\int_{\R^{2d}} \del b (t,x,[X^{s,\mu_2,(m)}_t]) (y) \Delta_{\mu_1,\mu_2} \del p_m(\cdot,s,t,x',y)(v)\,dy\, d\mu_2(x') \\ &=: I_1 + I_2 + I_3 + I_4 + I_5.
	\end{align*}

Thanks to \eqref{density_bound_Picard} and \eqref{b_lin_der_holder_mu_Picard} which is true for all $\gamma \in (0,1]$ since $W_1(\mu_1,\mu_2)< (t-s)^\frac{1}{\alpha}$, we obtain that \begin{align*}
|I_1| & \leq K (t-s)^{\frac{\eta - \gamma}{\alpha}}W_1^\gamma (\mu_1,\mu_2)\\ &\leq  K (t-s)^{-\frac{ \gamma}{\alpha}}W_1^\gamma (\mu_1,\mu_2).
\end{align*}

It follows from the boundedness of $\del b$ and \eqref{density_holder_measure_Picard} that \begin{align*}
	|I_2| & \leq K (t-s)^{ 1 - \frac{1+ \gamma}{\alpha}} W_1^\gamma (\mu_1,\mu_2)\\ &\leq   K (t-s)^{-\frac{\gamma}{\alpha}} W_1^\gamma (\mu_1,\mu_2).
\end{align*}

Owing to \eqref{b_lin_der_holder_mu_Picard} which is true for all $\gamma \in (0,1]$ since $W_1(\mu_1,\mu_2)< (t-s)^\frac{1}{\alpha}$ and \eqref{linear_der_bound_Picard} (the series appearing in this estimates is convergent), one has \begin{align*}
|I_3|  &\leq K(t-s)^{\frac{\eta - \gamma}{\alpha} + 1 - \frac{1}{\alpha}}W_1^\gamma (\mu_1,\mu_2) \\ &\leq K(t-s)^{-\frac{ \gamma}{\alpha}}W_1^\gamma (\mu_1,\mu_2) . 
\end{align*}

Let us now deal with $I_4.$ We need to control \[ J:= \int_{\R^{d}} \del b (t,x,[X^{s,\mu_2,(m)}_t]) (y) \left[\del p_m(\mu_1,s,t,x',y)(v) - \del p_m(\mu_1,s,t,x'',y)(v) \right]\,dy.\] 

Thanks to \eqref{linear_der_holder_x_Picard} and the boundedness $\del b$, we obtain that \begin{align*}
	|J| &\leq K\int_{\R^d} (t-s)^{1-\frac{1+ \gamma}{\alpha}}|x'-x''|^\gamma \left[\rho^0(t-s,y-x') +\rho^0(t-s,y-x'') \right] \, dy \\ &\leq   K (t-s)^{1-\frac{1+ \gamma}{\alpha}}|x'-x''|^\gamma.
\end{align*}

Jensen's inequality implies that \begin{align*} |I_4| &\leq K(t-s)^{1-\frac{1+ \gamma}{\alpha}} W_1^\gamma (\mu_1,\mu_2) \\ &\leq  K (t-s)^{-\frac{\gamma}{\alpha}} W_1^\gamma (\mu_1,\mu_2).  \end{align*}

Finally, for $I_5$, one has 

\begin{align*}
	|I_5| &= \left\vert \int_{\R^{2d}} \del b (t,x,[X^{s,\mu_2,(m)}_t]) (y)  \Delta_{\mu_1,\mu_2} \del p_m(\cdot,s,t,x',y)(v)\,dy\, d\mu_2(x') \right\vert \\ &\leq K \int_{\R^{2d}} \left\vert \Delta_{\mu_1,\mu_2} \del p_m(\cdot,s,t,x',y) \right\vert \, dy \, d\mu_2(x'). 
\end{align*}

Gathering all the previous estimates, we have proved that \begin{multline*}
	\left\vert  \del \left[ b(t,x,[X^{s,\mu_1,(m)}_t])\right](v) -  \del \left[ b(t,x,[X^{s,\mu_2,(m)}_t])\right](v) \right\vert \leq K(t-s)^{-\frac{ \gamma}{\alpha}}W_1^\gamma (\mu_1,\mu_2) \\ +K \int_{\R^d}\left\vert \Delta_{\mu_1,\mu_2} \del p_m(\cdot,s,t,x',y')\right\vert \, dy\,d\mu_2(x').
\end{multline*}

\noindent\textbf{Proof of \eqref{b_gradient_lin_der_holder_mu_Picard_2}.}
We begin to treat the case where $W_1(\mu_1,\mu_2) \geq (t-s)^{\frac{1}{\alpha}}$. Using \eqref{b_gradient_linear_derivative_bound_Picard} and noting that the series appearing in the bound is convergent, one has

\begin{align*}
	\left\vert \Delta_{\mu_1,\mu_2} \pa_v\del \left[ b(t,x,[X^{s,\cdot,(m)}_t])\right](v) \right\vert &\leq  \left\vert  \pa_v\del \left[ b(t,x,[X^{s,\mu_1,(m)}_t])\right](v) \right\vert + \left\vert  \pa_v\del \left[ b(t,x,[X^{s,\mu_2,(m)}_t])\right](v) \right\vert \\ &\leq K(t-s)^{\frac{\eta -1}{\alpha}} \\ &\leq K(t-s)^{\frac{\eta -1-\gamma}{\alpha}} W_1(\mu_1,\mu_2)^\gamma.
\end{align*}

We now focus on the case where $W_1(\mu_1,\mu_2)< (t-s)^\frac{1}{\alpha}$. In this case, by \eqref{b_expression_gradient_linr_der_Picard}, we have 

\begin{align*}
	&\Delta_{\mu_1,\mu_2} \pa_v\del \left[ b(t,x,[X^{s,\cdot,(m)}_t])\right](v)  \\& =  \int_{\R^d} \Delta_{\mu_1,\mu_2}\del b (t,x,[X^{s,\cdot,(m)}_t]) (y) \pa_x p_m(\mu_1,s,t,v,y)\,dy \\ & \quad +\int_{\R^{d}} \del b (t,x,[X^{s,\mu_2,(m)}_t]) (y) \Delta_{\mu_1,\mu_2} \pa_xp_m(\cdot,s,t,v,y)\,dy \\ &\quad +\int_{\R^{2d}} \Delta_{\mu_1,\mu_2}\del b (t,x,[X^{s,\cdot,(m)}_t]) (y) \pa_v \del p_m(\mu_1,s,t,x',y)(v) \, dy \, d\mu_1(x') \\ & \quad +\int_{\R^{2d}} \del b (t,x,[X^{s,\mu_2,(m)}_t]) (y) \pa_v \del p_m(\mu_1,s,t,x',y)(v)\,dy\, d(\mu_1-\mu_2)(x') \\ & \quad +\int_{\R^{2d}} \del b (t,x,[X^{s,\mu_2,(m)}_t]) (y) \Delta_{\mu_1,\mu_2} \pa_v \del p_m(\cdot,s,t,x',y)(v)\,dy\, d\mu_2(x') \\ &=: I_1 + I_2 + I_3 + I_4 + I_5.
\end{align*}

Thanks to \eqref{density_bound_Picard} and \eqref{b_lin_der_holder_mu_Picard} which is true for all $\gamma \in (0,1]$ since $W_1(\mu_1,\mu_2)< (t-s)^\frac{1}{\alpha}$, we obtain that \begin{align*}
	|I_1| & \leq K (t-s)^{\frac{\eta - \gamma-1}{\alpha}}W_1^\gamma (\mu_1,\mu_2).
\end{align*}

Concerning $I_2$, it can be rewritten as \[ I_2 = \int_{\R^{d}} \left(\del b (t,x,[X^{s,\mu_2,(m)}_t]) (y) - \del b (t,x,[X^{s,\mu_2,(m)}_t]) (v)\right) \Delta_{\mu_1,\mu_2}\pa_x p_m(\cdot,s,t,v,y)\,dy. \]

It follows from the $\eta$-Hölder continuity of $\del b(t,x,\mu)(\cdot)$, \eqref{gradient_density_holder_measure_Picard} since $\gamma \in (0,\eta + \alpha -1)$ and the space-time inequality \eqref{scalingdensityref} that \begin{align*}
	|I_2| & \leq K(t-s)^{\frac{\eta}{\alpha}+1-\frac{1+\gamma +1}{\alpha}} W_1^\gamma (\mu_1,\mu_2)\\ & \leq K (t-s)^{\frac{\eta - \gamma-1}{\alpha}} W_1^\gamma (\mu_1,\mu_2).
\end{align*}

Owing to \eqref{b_lin_der_holder_mu_Picard} which is true for all $\gamma \in (0,1]$ since $W_1(\mu_1,\mu_2)< (t-s)^\frac{1}{\alpha}$ and \eqref{gradient_linear_der_bound_Picard} (the series appearing in this estimates is convergent), one has \begin{align*}
	|I_3|  &\leq K(t-s)^{\frac{\eta - \gamma}{\alpha} + \frac{\eta -1}{\alpha}+ 1 - \frac{1}{\alpha}}W_1^\gamma(\mu_1,\mu_2) \\ &\leq K(t-s)^{\frac{\eta - \gamma-1}{\alpha}}W_1^\gamma(\mu_1,\mu_2). 
\end{align*}

Let us now deal with $I_4.$ We need to control \[ J:= \int_{\R^{d}} \del b (t,x,[X^{s,\mu_2,(m)}_t]) (y) \left[\pa_v\del p_m(\mu_1,s,t,x',y)(v) - \pa_v\del p_m(\mu_1,s,t,x'',y)(v) \right]\,dy.\] 

Let us prove that \[|J| \leq K(t-s)^{\frac{\eta - \gamma-1}{\alpha} + 1 + \frac{\eta-1}{\alpha}}|x'-x''|^\gamma.\]

Assume first that $|x'-x''|\leq (t-s)^{\frac{1}{\alpha}}.$ In this case using \eqref{representation_gradient_linear_der} and \eqref{expression_linear_der_H}, we have $\int_{\R^d} \pa_v \del p_m(\mu_1,s,t,x',y)(v) \, dy = 0 $. We can thus rewrite \begin{align*}
	J	&=\int_{\R^{d}} \left[\del b (t,x,[X^{s,\mu_2,(m)}_t]) (y) - \del b (t,x,[X^{s,\mu_2,(m)}_t]) (x'') \right]\\ &\hspace{5cm}\left[\pa_v \del p_m(\mu_1,s,t,x',y)(v) - \pa_v \del p_m(\mu_1,s,t,x'',y)(v) \right] \, dy \\
\end{align*}

Thanks to \eqref{gradient_linear_der_holder_x_Picard}, \eqref{controldensityref} since $|x'-x''|\leq (t-s)^{\frac{1}{\alpha}}$ and the $\eta$-Hölder continuity of $\del b (t,x,\mu)(\cdot)$, we obtain that \begin{align*}
	|J| &\leq K\int_{\R^d}|y-x''|^\eta (t-s)^{\frac{\eta -1 - \gamma}{\alpha} + 1 - \frac{1}{\alpha}}|x'-x''|^\gamma \rho^1(t-s,y-x'') \, dy.
\end{align*}
The space-time inequality \eqref{scalingdensityref} yields \begin{align*}
	|J| &\leq K(t-s)^{\frac{\eta}{\alpha} + \frac{\eta -1 - \gamma}{\alpha} + 1 - \frac{1}{\alpha}}|x'-x''|^\gamma.
\end{align*}

Assume now that $|x'-x''|> (t-s)^{\frac{1}{\alpha}}.$ In this case, we rewrite \begin{align*}
	J &= \int_{\R^{d}} \left[\del b (t,x,[X^{s,\mu_2,(m)}_t]) (y) - \del b (t,x,[X^{s,\mu_2,(m)}_t]) (x')  \right] \pa_v \del p_m(\mu_1,s,t,x',y)(v) \, dy \\ &\quad + \int_{\R^{d}} \left[\del b (t,x,[X^{s,\mu_2,(m)}_t]) (y) - \del b (t,x,[X^{s,\mu_2,(m)}_t]) (x'')  \right] \pa_v \del p_m(\mu_1,s,t,x'',y)(v) \, dy.
\end{align*}

Then, we use \eqref{gradient_linear_der_bound_Picard} (the series appearing in the bound being convergent) and the $\eta$-Hölder continuity of $\del b(t,x,\mu)(\cdot)$ which yield \begin{align*}
	|J| &\leq K (t-s)^{\frac{\eta}{\alpha} + \frac{\eta -1}{\alpha} + 1 - \frac{1}{\alpha}} \\ &\leq K (t-s)^{\frac{\eta}{\alpha} + \frac{\eta -\gamma-1}{\alpha} + 1 - \frac{1}{\alpha}}|x'-x''|^\gamma.
\end{align*}

Jensen's inequality implies that \begin{align*} |I_4| &\leq K (t-s)^{\frac{\eta - \gamma-1}{\alpha} + 1 + \frac{\eta-1}{\alpha}} W_1^\gamma (\mu_1,\mu_2) \\ &\leq  K (t-s)^{\frac{\eta - \gamma-1}{\alpha}} W_1^\gamma (\mu_1,\mu_2).  \end{align*}

Finally, for $I_5$, one has 

\begin{align*}
	|I_5| &= \left\vert \int_{\R^{2d}} \left(\del b (t,x,[X^{s,\mu_2,(m)}_t]) (y)  - \del b (t,x,[X^{s,\mu_2,(m)}_t]) (x') \right) \Delta_{\mu_1,\mu_2} \pa_v \del p_m(\cdot,s,t,x',y)(v)\,dy\, d\mu_2(x') \right\vert \\ &\leq K \int_{\R^{2d}} (1 \wedge |x'-y|^\eta) \left\vert \Delta_{\mu_1,\mu_2}\pa_v \del p_m(\cdot,s,t,x',y) \right\vert \, dy \, d\mu_2(x'). 
\end{align*}

Gathering all the previous estimates, we have proved that \begin{multline*}
	\left\vert  \pa_v\del \left[ b(t,x,[X^{s,\mu_1,(m)}_t])\right](v) -  \pa_v\del \left[ b(t,x,[X^{s,\mu_2,(m)}_t])\right](v) \right\vert \leq K(t-s)^{\frac{\eta - \gamma-1}{\alpha}}W_1^\gamma (\mu_1,\mu_2) \\ +K \int_{\R^d} (1 \wedge |x'-y|^\eta) \left\vert \Delta_{\mu_1,\mu_2} \pa_v \del p_m(\cdot,s,t,x',y')\right\vert \, dy\,d\mu_2(x').
\end{multline*}

\noindent\textbf{Proof of \eqref{H_lin_der_holder_mu_Picard} and \eqref{H_gradient_lin_der_holder_mu_Picard}.} Both estimates are immediate consequences of \eqref{b_lin_der_holder_mu_Picard_2} and \eqref{b_gradient_lin_der_holder_mu_Picard_2} using the expression of $\del \H_{m+1}$ given by \eqref{expression_linear_der_H} and \eqref{density_bound_Picard}.

\end{proof}

\subsection{Third part of the induction step}

We prove here that the estimates \eqref{gradient_density_holder_measure_Picard},
\eqref{linear_derivative_holder_measure_Picard}, and 
\eqref{gradient_linear_derivative_holder_measure_Picard} hold true.\\

\noindent\textbf{Proof of \eqref{gradient_density_holder_measure_Picard}.} Assume first that $W_1(\mu_1,\mu_2)> (t-s)^{\frac{1}{\alpha}}$. In this case, the parametrix expansion \eqref{representationdensityparametrix_Picard} (which can be differentiated with respect to $x$), \eqref{gradientestimatestable} and \eqref{Phi_m_bound} yield

\begin{align*}
	\left\vert \Delta_{\mu_1,\mu_2} \pa_x^jp_m(\cdot,s,t,x,y) \right\vert &\leq  \left\vert \Delta_{\mu_1,\mu_2} \left[\pa_x^j\p \otimes \Phi_m\right](\cdot,s,t,x,y) \right\vert \\ &\leq  \left\vert  \pa_x^j\p \otimes \Phi_m(\mu_1,s,t,x,y) \right\vert +  \left\vert  \pa_x^j\p \otimes \Phi_m(\mu_2,s,t,x,y) \right\vert \\ &\leq C (t-s)^{
		1 - \frac{1+j}{\alpha}}\rho^j(t-s,y-x) \\ &\leq C(t-s)^{1 - \frac{1+j+\gamma }{\alpha}}W_1^\gamma(\mu_1,\mu_2)\rho^j(t-s,y-x).
\end{align*}

We now focus on the case where $W_1(\mu_1,\mu_2) \leq (t-s)^{\frac{1}{\alpha}}.$ Let us prove that the following representation formula holds true \begin{equation}\label{representation_delta_mu_density}
	\Delta_{\mu_1,\mu_2} \pa_x^j p_{m+1}(\cdot,s,t,x,y) = \sum_{k=0}^\infty \left(\pa_x^j p_{m+1} \otimes \Delta_{\mu_1,\mu_2} \H_{m+1} \right)\otimes \H^k_{m+1}(\mu_2,s,t,x,y).
\end{equation}
Indeed, using the representation formula \eqref{representationdensityparametrix_Picard} which can be differentiated with respect to $x \in \R^d$, we get that \begin{equation*}
\Delta_{\mu_1,\mu_2} \pa_x^j p_{m+1}(\cdot,s,t,x,y) = \pa_x^j p_{m+1} \otimes \Delta_{\mu_1,\mu_2} \H_{m+1} (\mu_2,s,t,x,y) + \Delta_{\mu_1,\mu_2} \pa_x^j p_{m+1} \otimes \H_{m+1} (\mu_1,s,t,x,y).
\end{equation*}
Hence, we easily prove by induction that for all $n \geq $1, one has \begin{multline}\label{representation_del_mu_density_lemme}
	\Delta_{\mu_1,\mu_2} \pa_x^j p_{m+1}(\cdot,s,t,x,y) = \sum_{k=0}^n \left(\pa_x^j p_{m+1} \otimes \Delta_{\mu_1,\mu_2} \H_{m+1} \right) \otimes \H_{m+1}^k(\mu_2,s,t,x,y) \\+ \Delta_{\mu_1,\mu_2} \pa_x^j p_{m+1}\otimes \H_{m+1}^{n+1} (\mu_1,s,t,x,y).
\end{multline}

Thanks to \eqref{density_bound_Picard} and \eqref{H^k_m_bound}, we deduce that $\Delta_{\mu_1,\mu_2} \pa_x^j p_{m+1} \otimes \H_{m+1}\otimes \H_{m+1}^{n+1} (\mu_1,s,t,x,y)$ converges to $0$ as $n$ tends to infinity. Then, it follows from \eqref{density_bound_Picard}, \eqref{H^k_m_bound} and the convolution inequality \eqref{convolineqdensityref} that \begin{align*}
	&\left\vert \pa_x^j p_{m+1} \otimes \Delta_{\mu_1,\mu_2} \H_{m+1}(\mu_2,s,t,x,y) \right\vert \\&\leq \int_s^t \int_{\R^d} \left\vert\pa_x^j p_{m+1}(\mu_2,s,r,x,z) \Delta_{\mu_1,\mu_2} \H_{m+1}(\cdot,s,r,t,z,y)\right\vert \, dz\,dr \\ &\leq K \int_s^t \int_{\R^d} (r-s)^{-\frac{j}{\alpha}} \rho^j(r-s,z-x)  (t-r)^{-\frac{1}{\alpha}} \rho^1(t-r,y-z) \, dz\,dr\\ &\leq K (t-s)^{-\frac{j}{\alpha} + 1 - \frac{1}{\alpha}} \rho^j(t-s,y-x).
\end{align*}

From this estimate and \eqref{H^k_m_bound}, we deduce that for any $k \geq 1$ 

\begin{multline*}
	\left\vert\left(\pa_x^j p_{m+1} \otimes \Delta_{\mu_1,\mu_2} \H_{m+1} \right) \otimes \H_{m+1}^k(\mu_2,s,t,x,y)\right\vert \leq K (t-s)^{-\frac{j}{\alpha} + 1 - \frac{1}{\alpha} + k\left(1 - \frac{1}{\alpha}\right)}\rho^j(t-s,y-x) \\ \BB \left(k \left(1 - \frac{1}{\alpha}\right), 1 - \frac{j}{\alpha} + 1 - \frac{1}{\alpha} \right) \prod_{i=1}^{k-1} \BB\left(i \left(1 - \frac{1}{\alpha}\right), 1 - \frac{1}{\alpha}\right).
	\end{multline*}
Thanks to the asymptotic behavior of the Beta function, we deduce that the series appearing in \eqref{representation_del_mu_density_lemme} is absolutely convergent, which yields the representation formula \eqref{representation_delta_mu_density}. Now, since $\gamma \in (0,1]$ if $j=0$ and $\gamma \in (0,\eta + \alpha -1)$ if $j=1$, we can chose $\tilde{\gamma} \geq \gamma$ such that $\tilde{\gamma} \in [\eta,1]$ if $j=0$ and $\tilde{\gamma} \in [\eta, \eta + \alpha -1)$ if $j=1.$ It follows from \eqref{density_bound_Picard}, \eqref{H_holder_mu_Picard}, the convolution inequality \eqref{convolineqdensityref} and since $\tilde{\gamma} \in [\eta, \eta + \alpha -1)$ when $j=1$ that
 \begin{align*}
	&\left\vert \pa_x^j p_{m+1} \otimes \Delta_{\mu_1,\mu_2} \H_{m+1}(\mu_2,s,t,x,y) \right\vert \\&\leq \int_s^t \int_{\R^d} \left\vert\pa_x^j p_{m+1}(\mu_2,s,r,x,z) \Delta_{\mu_1,\mu_2} \H_{m+1}(\cdot,s,r,t,z,y)\right\vert \, dz\,dr \\ &\leq K \int_s^t \int_{\R^d} (r-s)^{-\frac{j}{\alpha}} \rho^j(r-s,z-x) (r-s)^{\frac{\eta - \tilde{\gamma}}{\alpha}} (t-r)^{-\frac{1}{\alpha}} W_1^{\tilde{\gamma}}(\mu_1,\mu_2)\rho^1(t-r,y-z) \, dz\,dr\\ &\leq K (t-s)^{\frac{\eta - \tilde{\gamma} -j}{\alpha} + 1 - \frac{1}{\alpha}} W_1^{\tilde{\gamma}} (\mu_1,\mu_2) \rho^j(t-s,y-x).
\end{align*}
From this estimate and \eqref{H^k_m_bound}, we deduce that for any $k \geq 1$ 

\begin{multline*}
	\left\vert\left(\pa_x^j p_{m+1} \otimes \Delta_{\mu_1,\mu_2} \H_{m+1} \right) \otimes \H_{m+1}^k(\mu_2,s,t,x,y)\right\vert \leq K (t-s)^{\frac{\eta - \tilde{\gamma} -j}{\alpha} + 1 - \frac{1}{\alpha} + k\left(1 - \frac{1}{\alpha}\right)} W_1^{\tilde{\gamma}} (\mu_1,\mu_2)\rho^j(t-s,y-x) \\ \BB \left(k \left(1 - \frac{1}{\alpha}\right), 1 + \frac{\eta - \tilde{\gamma} -j}{\alpha} + 1 - \frac{1}{\alpha} \right) \prod_{i=1}^{k-1} \BB\left(i \left(1 - \frac{1}{\alpha}\right), 1 - \frac{1}{\alpha}\right).
\end{multline*}
Summing over $k \geq 0$, we find that, since $\gamma \leq \tilde{\gamma}$ and $W_1(\mu_1,\mu_2) < (t-s)^{\frac{1}{\alpha}}$, we have
\begin{align*}
	\left\vert \Delta_{\mu_1,\mu_2} \pa_x^j p_{m+1} (\cdot,s,t,x,y) \right\vert &\leq  K (t-s)^{\frac{\eta - \tilde{\gamma} -j}{\alpha} +1 - \frac{1}{\alpha}} W_1^{\tilde{\gamma}}(\mu_1,\mu_2)\rho^j(t-s,y-x) \\ &\leq  K (t-s)^{\frac{\eta - \gamma -j}{\alpha}+1 - \frac{1}{\alpha}} W_1^{\gamma}(\mu_1,\mu_2) \rho^j(t-s,y-x)\\ &\leq K(t-s)^{1-\frac{1+\gamma +j}{\alpha}} W_1^{\gamma}(\mu_1,\mu_2)\rho^j(t-s,y-x).
\end{align*}

This concludes the proof of \eqref{gradient_density_holder_measure_Picard}.\\

\noindent\textbf{Proof of \eqref{linear_derivative_holder_measure_Picard}.} We first assume that $W_1(\mu_1,\mu_2) > (t-s)^{\frac{1}{\alpha}}$. In this case, using \eqref{linear_der_bound_Picard}, one has \begin{align*}
	&\left\vert \Delta_{\mu_1,\mu_2} \del p_{m+1}(\cdot,s,t,x,y)(v)\right\vert \\&\leq K (t-s)^{1 - \frac{1}{\alpha}} \rho^0(t-s,y-x) \\ &\leq  K (t-s)^{-\frac{\gamma}{\alpha}+1 - \frac{1}{\alpha}}W_1^\gamma(\mu_1,\mu_2) \rho^0(t-s,y-x) \\ &\leq  (t-s)^{-\frac{\gamma}{\alpha} + 1 - \frac{1}{\alpha}}W_1^\gamma(\mu_1,\mu_2)\rho^0(t-s,y-x) \\ &\hspace{4cm} \sum_{k=1}^{m+1} C^k (t-s)^{(k-1)\left(1 - \frac{1}{\alpha}\right)} \prod_{j=1}^{k-1} \BB\left(1 - \frac{\gamma}{\alpha} + (j-1)\left(1 - \frac{1}{\alpha}\right), 1 - \frac{1}{\alpha}\right),
\end{align*}
provided that we choose $C \geq K$ in \eqref{linear_derivative_holder_measure_Picard}. We now treat the case $W_1(\mu_1,\mu_2) \leq (t-s)^{\frac{1}{\alpha}}$. By the representation formula \eqref{representation_linear_der}, we have the following decomposition \begin{align*}
	\Delta_{\mu_1,\mu_2} \del p_{m+1}(\cdot,s,t,x,y)(v) &= p_{m+1}\otimes \Delta_{\mu_1,\mu_2} \del\H_{m+1} (\mu_2,s,t,x,y)(v) \\ &\quad + \Delta_{\mu_1,\mu_2} p_{m+1}\otimes \del \H_{m+1} (\mu_1,s,t,x,y)(v) \\ &\quad+\left( \Delta_{\mu_1,\mu_2} \left[p_{m+1} \otimes \del \H_{m+1}\right] \right)\otimes \Phi_{m+1}(\mu_1,s,t,x,y)(v) \\ &\quad+ \left[p_{m+1} \otimes \del \H_{m+1}\right] \otimes \Delta_{\mu_1,\mu_2}\Phi_{m+1}(\mu_2,s,t,x,y)(v)\\ &=: I_1 + I_2 + I_3 + I_4.
\end{align*}

For $I_1$, using \eqref{density_bound_Picard}, \eqref{H_lin_der_holder_mu_Picard}, the induction assumption \eqref{linear_derivative_holder_measure_Picard} and the convolution inequality \eqref{convolineqdensityref}, one has \begin{align*}
|I_1| &\leq K \int_s^t \int_{\R^d} \rho^0(r-s,z-x) \left[(r-s)^{-\frac{ \gamma}{\alpha}}(t-r)^{-\frac{1}{\alpha}}W_1^\gamma (\mu_1,\mu_2)\rho^1(t-r,y-z) \right.\\ &\left.\hspace{4cm} +(t-r)^{-\frac{1}{\alpha}} \int_{\R^{2d}}  \left\vert \Delta_{\mu_1,\mu_2} \del p_m(\cdot,s,r,x',y')\right\vert \, dy'\,d\mu_2(x') \rho^1(t-r,y-z)\right]\, dz\,dr \\ &\leq K (t-s)^{-\frac{\gamma}{\alpha} +1 - \frac{1}{\alpha}} W_1^\gamma(\mu_1,\mu_2)\rho^0(t-s,y-x) \\ &\quad+ K (t-s)^{-\frac{\gamma}{\alpha} + 1 - \frac{1}{\alpha}}W_1^\gamma (\mu_1,\mu_2)\rho^0(t-s,y-x)  \sum_{k=1}^m C^k (t-s)^{k\left(1 - \frac{1}{\alpha}\right)} \prod_{j=1}^{k} \BB\left(1 - \frac{\gamma}{\alpha} + j\left(1 - \frac{1}{\alpha}\right), 1 - \frac{1}{\alpha}\right).
\end{align*}

For $I_2$, it follows from \eqref{density_holder_measure_Picard}, \eqref{H_linear_derivative_bound_Picard} (the series appearing in the bound being convergent) and the convolution inequality \eqref{convolineqdensityref} that 

\begin{align*}
|I_2| &\leq K \int_s^t \int_{\R^d} (r-s)^{1 - \frac{1 + \gamma}{\alpha}}W_1^\gamma(\mu_1,\mu_2) \rho^0(r-s,z-x) (t-r)^{-\frac{1}{\alpha}} \rho^1(t-r,y-z) \, dz\,dr \\ &\leq K (t-s)^{1 - \frac{1 + \gamma}{\alpha} + 1 - \frac{1}{\alpha}} W_1^\gamma (\mu_1,\mu_2) \rho^0(t-s,y-x) \\ &\leq K (t-s)^{- \frac{ \gamma}{\alpha} + 1 - \frac{1}{\alpha}} W_1^\gamma (\mu_1,\mu_2) \rho^0(t-s,y-x). 
\end{align*}

Concerning $I_3$, we note that it writes \[I_3 = \left(I_1 + I_2\right) \otimes \Phi_{m+1}(\mu_1,s,t,x,y)(v).\] Thus, using the preceding bounds obtained for $I_1$ and $I_2$, \eqref{Phi_m_bound} and the convolution inequality \eqref{convolineqdensityref}, we get that \begin{align*}
	|I_3| & \leq K (t-s)^{-\frac{\gamma}{\alpha} +1 - \frac{1}{\alpha} +1 - \frac{1}{\alpha}} W_1^\gamma(\mu_1,\mu_2)\rho^0(t-s,y-x) \\ &\hspace{5cm}\quad\left[ 1 + \sum_{k=1}^m C^k (t-s)^{k\left(1 - \frac{1}{\alpha}\right)} \prod_{j=1}^{k} \BB\left(1 - \frac{\gamma}{\alpha} + j\left(1 - \frac{1}{\alpha}\right), 1 - \frac{1}{\alpha}\right) \right]\\ &\leq K (t-s)^{-\frac{\gamma}{\alpha} +1 - \frac{1}{\alpha}} W_1^\gamma(\mu_1,\mu_2)\rho^0(t-s,y-x) \\ &\hspace{5cm}\quad\left[ 1 + \sum_{k=1}^m C^k (t-s)^{k\left(1 - \frac{1}{\alpha}\right)} \prod_{j=1}^{k} \BB\left(1 - \frac{\gamma}{\alpha} + j\left(1 - \frac{1}{\alpha}\right), 1 - \frac{1}{\alpha}\right) \right].
\end{align*}

We finally deal with $I_4$. The convolution inequality \eqref{convolineqdensityref}, \eqref{density_bound_Picard} and \eqref{H_linear_derivative_bound_Picard} (the series appearing in this bound being convergent) yield \[\left\vert p_{m+1} \otimes \del \H_{m+1} (\mu_2,s,t,x,y)(v) \right\vert \leq K (t-s)^{1 - \frac{1}{\alpha}} \rho^0(t-s,y-x).\] 
Then, it follows from \eqref{Phi_holder_mu_Picard}, which is valid for all $\gamma \in (0,1]$ since $W_1(\mu_1,\mu_2)<(t-s)^{\frac{1}{\alpha}}$ that 

\begin{align*}
	|I_4| &\leq K \int_s^t \int_{\R^d} K (r-s)^{1 - \frac{1}{\alpha}} \rho^0(r-s,z-x) (r-s)^{\frac{\eta - \gamma}{\alpha}}(t-r)^{-\frac{1}{\alpha}}W_1^\gamma(\mu_1,\mu_2) \rho^1(t-r,y-z) \, dz \,dr \\ &\leq K (t-s)^{\frac{\eta - \gamma}{\alpha} + 1 - \frac{1}{\alpha} + 1 - \frac{1}{\alpha}}W_1^\gamma (\mu_1,\mu_2) \rho^0(t-s,y-x) \\ &\leq K (t-s)^{\frac{\eta - \gamma}{\alpha} + 1 - \frac{1}{\alpha}}W_1^\gamma (\mu_1,\mu_2) \rho^0(t-s,y-x).
\end{align*}

Gathering all the previous estimates, we have proved that 
\begin{align*}
&\left\vert \Delta_{\mu_1,\mu_2} \del p_{m+1}(\cdot,s,t,x,y)(v) \right\vert \\ & \leq  K (t-s)^{-\frac{\gamma}{\alpha} +1 - \frac{1}{\alpha}} W_1^\gamma(\mu_1,\mu_2)\rho^0(t-s,y-x) \\ &\hspace{2cm}\left[ 1 + \sum_{k=1}^m C^k (t-s)^{k\left(1 - \frac{1}{\alpha}\right)} \prod_{j=1}^{k} \BB\left(1 - \frac{\gamma}{\alpha} + j\left(1 - \frac{1}{\alpha}\right), 1 - \frac{1}{\alpha}\right)\right] \\ &\leq  (t-s)^{-\frac{\gamma}{\alpha} +1 - \frac{1}{\alpha}} W_1^\gamma(\mu_1,\mu_2)\rho^0(t-s,y-x) \sum_{k=1}^{m+1}C^k (t-s)^{(k-1)\left(1 - \frac{1}{\alpha}\right)} \prod_{j=1}^{k-1} \BB\left(1 - \frac{\gamma}{\alpha} + j\left(1 - \frac{1}{\alpha}\right), 1 - \frac{1}{\alpha}\right),
\end{align*}
provided that we choose $C\geq K$ in \eqref{linear_derivative_holder_measure_Picard}. This ends the proof of the induction step for \eqref{linear_derivative_holder_measure_Picard}.\\

\noindent\textbf{Proof of \eqref{gradient_linear_derivative_holder_measure_Picard}.} Assume first that $W_1(\mu_1,\mu_2) > (t-s)^{\frac{1}{\alpha}}$. In this case,  \eqref{linear_der_bound_Picard} yields \begin{align*}
	&\left\vert \Delta_{\mu_1,\mu_2} \pa_v\del p_{m+1}(\cdot,s,t,x,y)(v)\right\vert \\&\leq K (t-s)^{\frac{\eta -1}{\alpha}+1 - \frac{1}{\alpha}} \rho^0(t-s,y-x) \\ &\leq  K (t-s)^{\frac{\eta -\gamma -1}{\alpha}+1 - \frac{1}{\alpha}}W_1^\gamma(\mu_1,\mu_2) \rho^0(t-s,y-x) \\ &\leq  (t-s)^{\frac{\eta - \gamma-1}{\alpha} + 1 - \frac{1}{\alpha}}W_1^\gamma(\mu_1,\mu_2)\rho^0(t-s,y-x) \\ &\hspace{4cm} \sum_{k=1}^{m+1} C^k (t-s)^{(k-1)\left(1 + \frac{\eta -1}{\alpha}\right)} \prod_{j=1}^{k-1} \BB\left(1 + \frac{\eta-\gamma}{\alpha} + (j-1)\left(1 +\frac{\eta-1}{\alpha}\right), 1 - \frac{1}{\alpha}\right),
\end{align*}

provided that we choose $C \geq K$ in \eqref{linear_derivative_holder_measure_Picard}. We now treat the case $W_1(\mu_1,\mu_2) \leq (t-s)^{\frac{1}{\alpha}}$. By the representation formula \eqref{representation_gradient_linear_der}, we have the following decomposition \begin{align*}
	\Delta_{\mu_1,\mu_2} \pa_v\del p_{m+1}(\cdot,s,t,x,y)(v) &= p_{m+1}\otimes \Delta_{\mu_1,\mu_2} \pa_v\del\H_{m+1} (\mu_2,s,t,x,y)(v) \\ &\quad + \Delta_{\mu_1,\mu_2} p_{m+1}\otimes \pa_v\del \H_{m+1} (\mu_1,s,t,x,y)(v) \\ &\quad+\left( \Delta_{\mu_1,\mu_2} \left[p_{m+1} \otimes \pa_v\del \H_{m+1}\right] \right)\otimes \Phi_{m+1}(\mu_1,s,t,x,y)(v) \\ &\quad+ \left[p_{m+1} \otimes \pa_v \del \H_{m+1}\right] \otimes \Delta_{\mu_1,\mu_2}\Phi_{m+1}(\mu_2,s,t,x,y)(v)\\ &=: I_1 + I_2 + I_3 + I_4.
\end{align*}

For $I_1$, using \eqref{density_bound_Picard}, \eqref{H_gradient_lin_der_holder_mu_Picard}, the induction assumption \eqref{gradient_linear_derivative_holder_measure_Picard}, the space-time inequality \eqref{scalingdensityref} and the convolution inequality \eqref{convolineqdensityref}, one has since $\gamma \in (0, \eta + \alpha -1)$ \begin{align*}
	|I_1| &\leq K \int_s^t \int_{\R^d} \rho^0(r-s,z-x) \left[(r-s)^{\frac{\eta- \gamma-1}{\alpha}}(t-r)^{-\frac{1}{\alpha}}W_1^\gamma (\mu_1,\mu_2)\rho^1(t-r,y-z) \right.\\ &\left.\hspace{2cm} +(t-r)^{-\frac{1}{\alpha}} \int_{\R^{2d}}  (1 \wedge |x'-y'|^\eta)\left\vert \Delta_{\mu_1,\mu_2} \del p_m(\cdot,s,r,x',y')\right\vert \, dy'\,d\mu_2(x') \rho^1(t-r,y-z)\right]\, dz\,dr \\ &\leq K (t-s)^{\frac{\eta -\gamma-1}{\alpha} +1 - \frac{1}{\alpha}} W_1^\gamma(\mu_1,\mu_2)\rho^0(t-s,y-x) \\ &\hspace{3cm}\left[ 1+  \sum_{k=1}^m C^k (t-s)^{k\left(1 + \frac{\eta-1}{\alpha}\right)} \prod_{j=1}^{k} \BB\left(1 + \frac{\eta-\gamma}{\alpha} + j\left(1 + \frac{\eta-1}{\alpha}\right), 1 - \frac{1}{\alpha}\right)\right].
\end{align*}

For $I_2$, it follows from \eqref{density_holder_measure_Picard}, \eqref{H_gradient_linear_derivative_bound_Picard} and the convolution inequality \eqref{convolineqdensityref} that 

\begin{align*}
	|I_2| &\leq K \int_s^t \int_{\R^d} (r-s)^{1 - \frac{1 + \gamma}{\alpha}}W_1^\gamma(\mu_1,\mu_2) \rho^0(r-s,z-x)(r-s)^{\frac{\eta -1}{\alpha}} (t-r)^{-\frac{1}{\alpha}} \rho^1(t-r,y-z) \, dz\,dr \\ &\leq K (t-s)^{\frac{\eta -1}{\alpha}+1 - \frac{1 + \gamma}{\alpha} + 1 - \frac{1}{\alpha}} W_1^\gamma (\mu_1,\mu_2) \rho^0(t-s,y-x) \\ &\leq K (t-s)^{ \frac{\eta - \gamma-1}{\alpha} + 1 - \frac{1}{\alpha}} W_1^\gamma (\mu_1,\mu_2) \rho^0(t-s,y-x). 
\end{align*}

Concerning $I_3$, we note that it writes \[I_3 = \left(I_1 + I_2\right) \otimes \Phi_{m+1}(\mu_1,s,t,x,y)(v).\] Thus, using the preceding bounds obtained for $I_1$ and $I_2$, \eqref{Phi_m_bound} and the convolution inequality \eqref{convolineqdensityref}, we get that \begin{align*}
	|I_3| & \leq K (t-s)^{\frac{\eta-\gamma-1}{\alpha} +1 - \frac{1}{\alpha} +1 - \frac{1}{\alpha}} W_1^\gamma(\mu_1,\mu_2)\rho^0(t-s,y-x) \\ &\hspace{5cm}\quad\left[ 1 + \sum_{k=1}^m C^k (t-s)^{k\left(1+ \frac{\eta-1}{\alpha}\right)} \prod_{j=1}^{k} \BB\left(1 + \frac{\eta-\gamma}{\alpha} + j\left(1 + \frac{\eta-1}{\alpha}\right), 1 - \frac{1}{\alpha}\right) \right]\\ &\leq K (t-s)^{\frac{\eta-\gamma-1}{\alpha} +1 - \frac{1}{\alpha}} W_1^\gamma(\mu_1,\mu_2)\rho^0(t-s,y-x) \\ &\hspace{5cm}\quad\left[ 1 + \sum_{k=1}^m C^k (t-s)^{k\left(1 +\frac{\eta-1}{\alpha}\right)} \prod_{j=1}^{k} \BB\left(1 + \frac{\eta-\gamma}{\alpha} + j\left(1 + \frac{\eta-1}{\alpha}\right), 1 - \frac{1}{\alpha}\right) \right].
\end{align*}

We finally deal with $I_4$. The convolution inequality \eqref{convolineqdensityref}, \eqref{density_bound_Picard} and \eqref{H_gradient_linear_derivative_bound_Picard} (the series appearing in this bound being convergent) yield \[\left\vert p_{m+1} \otimes \pa_v \del \H_{m+1} (\mu_2,s,t,x,y)(v) \right\vert \leq K (t-s)^{\frac{\eta-1}{\alpha}+1 - \frac{1}{\alpha}} \rho^0(t-s,y-x).\] 
Then, it follows from \eqref{Phi_holder_mu_Picard}, which is valid for all $\gamma \in (0,1]$ since $W_1(\mu_1,\mu_2)<(t-s)^{\frac{1}{\alpha}}$ that 

\begin{align*}
	|I_4| &\leq K \int_s^t \int_{\R^d} K (r-s)^{\frac{\eta-1}{\alpha}+1 - \frac{1}{\alpha}} \rho^0(r-s,z-x) (r-s)^{\frac{\eta - \gamma}{\alpha}}(t-r)^{-\frac{1}{\alpha}}W_1^\gamma(\mu_1,\mu_2) \rho^1(t-r,y-z) \, dz \,dr \\ &\leq K (t-s)^{\frac{\eta-1}{\alpha}+\frac{\eta - \gamma}{\alpha} + 1 - \frac{1}{\alpha} + 1 - \frac{1}{\alpha}}W_1^\gamma (\mu_1,\mu_2) \rho^0(t-s,y-x) \\ &\leq K (t-s)^{\frac{\eta - \gamma-1}{\alpha} + 1 - \frac{1}{\alpha}}W_1^\gamma (\mu_1,\mu_2) \rho^0(t-s,y-x).
\end{align*}

Gathering all the previous estimates, we have proved that 
\begin{align*}
	&\left\vert \Delta_{\mu_1,\mu_2}\pa_v \del p_{m+1}(\cdot,s,t,x,y)(v) \right\vert \\ & \leq  K (t-s)^{\frac{\eta-\gamma-1}{\alpha} +1 - \frac{1}{\alpha}} W_1^\gamma(\mu_1,\mu_2)\rho^0(t-s,y-x) \\ &\hspace{2cm}\left[ 1 + \sum_{k=1}^m C^k (t-s)^{k\left(1 + \frac{\eta-1}{\alpha}\right)} \prod_{j=1}^{k} \BB\left(1+\frac{\eta-\gamma}{\alpha} + j\left(1 + \frac{\eta-1}{\alpha}\right), 1 - \frac{1}{\alpha}\right)\right] \\ &\leq  (t-s)^{\frac{\eta-\gamma-1}{\alpha} +1 - \frac{1}{\alpha}} W_1^\gamma(\mu_1,\mu_2)\rho^0(t-s,y-x) \\ &\hspace{3cm}\sum_{k=1}^{m+1}C^k (t-s)^{(k-1)\left(1 + \frac{\eta-1}{\alpha}\right)} \prod_{j=1}^{k-1} \BB\left(1 + \frac{\eta-\gamma}{\alpha} + j\left(1 + \frac{\eta-1}{\alpha}\right), 1 - \frac{1}{\alpha}\right),
\end{align*}

provided that we choose $C\geq K$ in \eqref{gradient_linear_derivative_holder_measure_Picard}. This ends the proof of the induction step for \eqref{gradient_linear_derivative_holder_measure_Picard}.\\

\subsection{Preparatory technical results}

\begin{Lemme}\label{lemma_technical_7}
	\begin{itemize}
	\item For any $\gamma \in (0,1]$, $\tilde{\eta} \in (0,\eta \wedge (\alpha -1))$, there exists a positive constant $K$ such that for all $m \geq 1$, $t \in (0,T]$, $s_1,s_2 \in [0,t)$, $\mu \in \PP$, $x,y, \in \R^d$ \begin{equation}\label{density_Holder_s_Picard_lemme}
		\left\vert \Delta_{s_1,s_2} p_m(\mu,\cdot,t,x,y) \right\vert \leq K  \left[\frac{|s_1-s_2|^\gamma}{(t-s_1)^{\gamma}} \rho^{-\tilde{\eta}}(t-s_1,y-x) + \frac{|s_1-s_2|^\gamma}{(t-s_2)^{\gamma}} \rho^{-\tilde{\eta}}(t-s_2,y-x) \right].
	\end{equation}

\item For any $\gamma \in (0,1]$, there exists a positive constant $K$ such that for all $m \geq 1$, $t \in (0,T]$, $s_1,s_2 \in [0,t)$, $\mu \in \PP$, $x \in \R^d$ \begin{equation}\label{b_Holder_s_Picard_lemme}
	\left\vert \Delta_{s_1,s_2}b(t,x,[X^{\cdot,\mu,(m)}_t]) \right\vert \leq K |s_1-s_2|^\gamma \left[(t-s_1)^{\frac{\eta}{\alpha}-\gamma} + (t-s_2)^{\frac{\eta}{\alpha}-\gamma}\right].
\end{equation}

\item For any $\gamma \in (0,1]$, there exists a positive constant $K$ such that for all $m \geq 1$, $0 \leq r <t \leq T  $, $s_1,s_2 \in [0,r)$, $\mu \in \PP$, $x,y \in \R^d$ \begin{equation}\label{H_Holder_s_Picard_lemme}
	\left\vert \Delta_{s_1,s_2}\H_{m+1}(\mu,\cdot,r,t,x,y) \right\vert \leq K (t-r)^{- \frac{1}{\alpha}} \rho^1(t-r,y-x)|s_1-s_2|^\gamma \left[(r-s_1)^{\frac{\eta}{\alpha}-\gamma} + (r-s_2)^{\frac{\eta}{\alpha}-\gamma}\right],
\end{equation}

and \begin{equation}\label{Phi_Holder_s_Picard_lemme}
	\left\vert \Delta_{s_1,s_2}\Phi_{m+1}(\mu,\cdot,r,t,x,y) \right\vert \leq K (t-r)^{- \frac{1}{\alpha}} \rho^1(t-r,y-x)|s_1-s_2|^\gamma \left[(r-s_1)^{\frac{\eta}{\alpha}-\gamma} + (r-s_2)^{\frac{\eta}{\alpha}-\gamma}\right].
\end{equation}
	
	\end{itemize}
\end{Lemme}

\begin{proof}[Proof of Lemma \eqref{lemma_technical_7}]
	
	\noindent\textbf{Proof of \eqref{density_Holder_s_Picard_lemme}.} Assume first that $|s_1-s_2| > t-s_1\vee s_2$. Then, using \eqref{density_bound_Picard}, we get that 
	
	\begin{align*}
		&|\Delta_{s_1,s_2}p_m(\mu,\cdot,t,x,y)| \\&\leq K\left[\rho^0(t-s_1\vee s_2, y-x) + \rho^0(t-s_1\wedge s_2,y-x)\right] \\ &\leq K \left[ \frac{|s_1-s_2|^\gamma}{(t-s_1\vee s_2)^\gamma} \rho^0(t-s_1\vee s_2,y-x) +  \frac{(t-s_1 \wedge s_2)^\gamma  +|s_1-s_2|^\gamma}{(t-s_1 \wedge s_2)^\gamma} \rho^0(t-s_1\wedge s_2,y-x)  \right]\\ &\leq K \left[ \frac{|s_1-s_2|^\gamma}{(t-s_1\vee s_2)^\gamma} \rho^0(t-s_1\vee s_2,y-x) +  \frac{|s_1-s_2|^\gamma}{(t-s_1 \wedge s_2)^\gamma} \rho^0(t-s_1\wedge s_2,y-x)  \right].
	\end{align*}

We now focus on the case $|s_1-s_2| \leq t-s_1\vee s_2.$ For $\lambda \in [0,1]$, we set $ s_\lambda := \lambda s_1\vee s_2 + (1-\lambda)s_1\wedge s_2$. We have thus by \eqref{density_time_der_bound_Picard} (the series appearing being convergent) \begin{align*}
	&| \Delta_{s_1,s_2}p_m(\mu,\cdot,t,x,y)| \\ &\leq \int_0^1 |\pa_s p_m(\mu,s_\lambda,t,x,y)|\,|s_1-s_2| \, d\lambda \\ &\leq K|s_1-s_2| \int_0^1 (t-s_\lambda)^{-1} \rho^{-\tilde{\eta}}(t-s_\lambda,y-x) \, d\lambda \\ &\leq C|s_1-s_2| \int_0^1 (t-s_\lambda)^{-1-\frac{d}{\alpha}} (1+(t-s_\lambda)^{-\frac{1}{\alpha}} |y-x|)^{-d-\alpha + \tilde{\eta}}\, d\lambda \\ &\leq K |s_1-s_2|^\gamma (t-s_1\vee s_2)^{1-\gamma}(t-s_1\vee s_2)^{-1-\frac{d}{\alpha}} \left[(1+(t-s_1\vee s_2)^{-\frac{1}{\alpha}} |y-x|)^{-d-\alpha + \tilde{\eta}}\right. \\ &\left. \hspace{9cm} + (1+(t-s_1\wedge s_2)^{-\frac{1}{\alpha}} |y-x|)^{-d-\alpha + \tilde{\eta}}\right],
\end{align*}
for some $\tilde{\eta} \in (0,\eta \wedge(\alpha-1))$. Since $|s_1-s_2| \leq t-s_1\vee s_2$, we easily check that $(t-s_1\vee s_2)^{-1} \leq 2 (t-s_1\wedge s_2)^{-1}$. It follows that 

\begin{align*}
	&| \Delta_{s_1,s_2}p_m(\mu,\cdot,t,x,y) | \\ &\leq K |s_1-s_2|^\gamma\left[ (t-s_1\vee s_2)^{-\gamma-\frac{d}{\alpha}} (1+(t-s_1\vee s_2)^{-\frac{1}{\alpha}} |y-x|)^{-d-\alpha + \tilde{\eta}}\right. \\ &\left. \hspace{4cm} +  (t-s_1\wedge s_2)^{-\gamma-\frac{d}{\alpha}}(1+(t-s_1\wedge s_2)^{-\frac{1}{\alpha}} |y-x|)^{-d-\alpha +\tilde{\eta}}\right] \\ &  \leq K \left[ \frac{|s_1 - s_2|^\gamma}{(t-s_1\wedge s_2)^{\gamma}} \rho^{-\tilde{\eta}}(t-s_1 \wedge s_2, y-x) + \frac{|s_1 - s_2|^\gamma}{(t-s_1\vee s_2)^{\gamma }} \rho^{-\tilde{\eta}}(t-s_1 \vee s_2, y-x) \right].
\end{align*}
This concludes the proof of \eqref{density_Holder_s_Picard_lemme}. \\

\noindent\textbf{Proof of \eqref{b_Holder_s_Picard_lemme}.} By definition of the linear derivative and a centering argument, one can write setting $m_\lambda := \lambda [X^{s_1\vee s_2,\mu,(m)}_t] + (1-\lambda)  [X^{s_1\wedge s_2,\mu,(m)}_t]$ for $\lambda \in [0,1]$

\begin{align*}
	\Delta_{s_1,s_2} b(t,x,[X^{\cdot,\mu,(m)}_t]) = \int_{\R^{2d}} \left(\del b(t,x,m_\lambda)(y) - \del b(t,x,m_\lambda)(x')\right) \Delta_{s_1,s_2}p_m(\mu,\cdot,t,x',y) \, d\mu(x') \, dy \, d\lambda.
\end{align*}

We deduce using \eqref{density_Holder_s_Picard_lemme} for some $\tilde{\eta} \in (0,\eta \wedge (\alpha-1))$, the boundedness of $\del b$ and the $\eta$-Hölder continuity of $\del b(t,x,\mu)(\cdot)$ that  

\begin{align*}
	&\left\vert 	\Delta_{s_1,s_2} b(t,x,[X^{\cdot,\mu,(m)}_t]) \right\vert \\ &\leq K\int_{\R^{2d}} |y-x'|^\eta  \left[ \frac{|s_1 - s_2|^\gamma}{(t-s_1\wedge s_2)^{\gamma}} \rho^{-\tilde{\eta}}(t-s_1 \wedge s_2, y-x') \right.\\ &\left. \hspace{6cm}+ \frac{|s_1 - s_2|^\gamma}{(t-s_1\vee s_2)^{\gamma }} \rho^{-\tilde{\eta}}(t-s_1 \vee s_2, y-x') \right] \, dy \, d\mu(x').
\end{align*}
Notice that because $\eta + \tilde{\eta}< \alpha$, $\rho^{-\tilde{\eta} -\eta}(t-s,\cdot)$ is integrable. By the space-time inequality \eqref{scalingdensityref} and since $\int_{\R^d}\rho^{-\tilde{\eta} -\eta}(t-s,y) \, dy$ is a finite constant independent of $s$ and $t$, we obtain \eqref{b_Holder_s_Picard_lemme}.\\

\noindent\textbf{Proof of \eqref{H_Holder_s_Picard_lemme} and \eqref{Phi_Holder_s_Picard_lemme}.} By the definition \eqref{def_proxy_kernel_Picard} of $\H_{m+1}$ and since $\p(r,t,x,y)$ does not depend on $\mu$ and $s$, we immediately deduce \eqref{H_Holder_s_Picard_lemme} from \eqref{b_Holder_s_Picard_lemme} and \eqref{gradientestimatestable}.\\ 

Concerning the proof of \eqref{Phi_Holder_s_Picard_lemme}, we start from the Volterra integral equation \eqref{Volterra_eq_Picard} which yields \begin{align}\label{representation_delta_time_phi} \Delta_{s_1,s_2}\Phi_{m+1}(\mu,\cdot,r,t,x,y) \notag&= \Delta_{s_1,s_2}\H_{m+1} (\mu,\cdot,r,t,x,y) + \Delta_{s_1,s_2}\H_{m+1}\otimes \Phi_{m+1}(\mu,s_1\vee s_2,r,t,x,y) \\&+ \H_{m+1}\otimes \Delta_{s_1,s_2} \Phi_{m+1} (\mu,s_1\wedge s_2,r,t,x,y). \end{align}

Using \eqref{H_Holder_s_Picard_lemme} and \eqref{Phi_m_bound}, we deduce that 

\begin{align*}
	&\left\vert \Delta_{s_1,s_2}\H_{m+1} \otimes \Phi_{m+1}(\mu,s_1\vee s_2,r,t,x,y) \right\vert \\ &\leq K\int_r^t \int_{\R^d}   (r'-r)^{- \frac{1}{\alpha}} \rho^1(r'-r,y-x)|s_1-s_2|^\gamma \left[(r'-s_1)^{\frac{\eta}{\alpha}-\gamma} + (r'-s_2)^{\frac{\eta}{\alpha}-\gamma}\right] (t-r')^{-\frac{1}{\alpha}} \rho^1(t-r',y-z)\, dz \, dr'\\ &\leq K  (t-r)^{- \frac{1}{\alpha} + 1 - \frac{1}{\alpha}} \left[(r-s_1)^{\frac{\eta}{\alpha}-\gamma} + (r-s_2)^{\frac{\eta}{\alpha}-\gamma}\right]  \rho^1(t-r,y-x).
\end{align*}

This inequality and \eqref{H_Holder_s_Picard_lemme} ensure that \begin{multline}\label{Phi_holder_times_control_kernel_rep}
\left\vert \left[ \Delta_{s_1,s_2}\H_{m+1} + \Delta_{s_1,s_2}\H_{m+1} \otimes \Phi_{m+1}\right](\mu,s_1\vee s_2,r,t,x,y) \right\vert \\  \leq K  (t-r)^{- \frac{1}{\alpha}} \left[(r-s_1)^{\frac{\eta}{\alpha}-\gamma} + (r-s_2)^{\frac{\eta}{\alpha}-\gamma}\right]  \rho^1(t-r,y-x).
\end{multline}

As $\Phi_{m+1}$ yields a time-integrable singularity, we can iterate the implicit representation formula \eqref{representation_delta_time_phi}. We thus obtain \begin{multline*}
\Delta_{s_1,s_2} \Phi_{m+1} (\mu,\cdot,t,x,y) = \left[ \Delta_{s_1,s_2}\H_{m+1} + \Delta_{s_1,s_2}\H_{m+1} \otimes \Phi_{m+1}\right] (\mu,s_1\vee s_2,t,x,y)\\ + \sum_{k=1}^{\infty}\H_{m+1}^k \otimes  \left[ \Delta_{s_1,s_2}\H_{m+1} + \Delta_{s_1,s_2}\H_{m+1} \otimes \Phi_{m+1}\right](\mu,s_1\wedge s_2,t,x,y).
\end{multline*}
After standard computations using \eqref{Phi_holder_times_control_kernel_rep} and \eqref{H^k_m_bound} that we omit, we conclude that \eqref{Phi_Holder_s_Picard_lemme} holds true.

\end{proof}

\subsection{Fourth part of the induction step}
	
We prove here that the estimate \eqref{gradient_density_time_Holder_Picard} holds true.\\

\noindent\textbf{Proof of \eqref{gradient_density_time_Holder_Picard}.} We start by assuming that $|s_1-s_2| > t- s_1\vee s_2$. In this case, \eqref{density_bound_Picard} directly yields \begin{align*}
	&\left\vert \Delta_{s_1,s_2}\pa_x^j p_m(\mu,\cdot,t,x,y)\right\vert \\&\leq  	\left\vert \pa_x^j p_m(\mu,s_1,t,x,y)\right\vert + 	\left\vert \pa_x^j p_m(\mu,s_2,t,x,y)\right\vert \\ &\leq K \left[(t-s_1\vee s_2)^{-\frac{j}{\alpha}}\rho^j(t-s_1\vee s_2,y-x) + (t-s_1\wedge s_2)^{-\frac{j}{\alpha}}\rho^j(t-s_1\wedge s_2,y-x)\right]\\ &\leq K \left[\frac{|s_1-s_2|^\gamma}{(t-s_1\vee s_2)^{\gamma + \frac{j}{\alpha}}} \rho^j(t-s_1\vee s_2,y-x) + \frac{|s_1-s_2|^\gamma + (t-s_1\vee s_2)^\gamma}{(t-s_1\wedge s_2)^\gamma(t-s_1\wedge s_2)^{\frac{j}{\alpha}}} \rho^j(t-s_1\wedge s_2,y-x) \right]\\ &\leq   K \left[\frac{|s_1-s_2|^\gamma}{(t-s_1\vee s_2)^{\gamma + \frac{j}{\alpha}}} \rho^j(t-s_1\vee s_2,y-x) + \frac{|s_1-s_2|^\gamma }{(t-s_1\wedge s_2)^{\gamma +\frac{j}{\alpha}}} \rho^j(t-s_1\wedge s_2,y-x) \right],
\end{align*}
which proves \eqref{gradient_density_time_Holder_Picard}. \\

We now turn to the case $|s_1 - s_2|\leq t-s_1\vee s_2$. Differentiating with respect to $x$ the parametrix expansion \eqref{representationdensityparametrix_Picard}, we get \[\pa_x^j p_{m} (\mu,s,t,x,y) = \pa_x^j \p(s,t,x,y) + \pa_x^j \p \otimes \Phi_m(\mu,s,t,x,y).\] 
We are going to use the following decomposition 

\begin{align*}
	\Delta_{s_1,s_2} \pa_x^jp_m(\mu,\cdot,t,x,y) &= \Delta_{s_1,s_2} \pa_x^j \p(\cdot,t,x,y) \\ &\quad+ \int_{s_1\vee s_2}^t \int_{\R^d} \Delta_{s_1,s_2} \pa_x^j\p(\cdot,r,x,z) \Phi_{m}(\mu,s_1\vee s_2,r,t,z,y) \, dz \,dr \\ &\quad+ \int_{s_1\vee s_2}^t \int_{\R^d}  \pa_x^j\p(s_1\wedge s_2,r,x,z) \Delta_{s_1,s_2}\Phi_{m}(\mu,\cdot,r,t,z,y) \, dz \,dr\\ &\quad - \int_{s_1\wedge s_2} ^{s_1\vee s_2} \int_{\R^d} \pa_x^j \p (s_1\wedge s_2,r,x,z) \Phi_m(\mu,s_1\wedge s_2,r,t,z,y) \, dz \,dr \\ &=: I_1 + I_2 + I_3 + I_4.
\end{align*}

By \eqref{gradientHoldertimestable}, we obtain that \begin{align*}
	|I_1| &\leq K \left[ \frac{|s_1 - s_2|^\gamma}{(t-s_1\wedge s_2)^{\gamma + \frac{j}{\alpha}}} \rho^{j}(t-s_1 \wedge s_2, y-x) + \frac{|s_1 - s_2|^\gamma}{(t-s_1\vee s_2)^{\gamma + \frac{j}{\alpha}}} \rho^{j}(t-s_1 \vee s_2, y-x) \right].
\end{align*}

We now focus on $I_2$. Thanks to \eqref{gradientHoldertimestable}, \eqref{Phi_m_bound} and the convolution inequality \eqref{convolineqdensityref}, one has \begin{align*}
	|I_2| &\leq K \int_{s_1\vee s_2}^t \int_{\R^d} \left[ \frac{|s_1 - s_2|^\gamma}{(r-s_1\wedge s_2)^{\gamma + \frac{j}{\alpha}}} \rho^{j}(r-s_1 \wedge s_2, z-x) + \frac{|s_1 - s_2|^\gamma}{(r-s_1\vee s_2)^{\gamma + \frac{j}{\alpha}}} \rho^{j}(r-s_1 \vee s_2, z-x) \right]\\ &\hspace{12cm}(t-r)^{-\frac{1}{\alpha}}\rho^1(t-r,y-z) \, dz\,dr \\ &\leq  K \int_{s_1\wedge s_2}^t \int_{\R^d} \frac{|s_1 - s_2|^\gamma}{(r-s_1\wedge s_2)^{\gamma + \frac{j}{\alpha}}} \rho^{j}(r-s_1 \wedge s_2, z-x) (t-r)^{-\frac{1}{\alpha}}\rho^1(t-r,y-z) \, dz\,dr\\ & \quad+ K \int_{s_1\vee s_2}^t \int_{\R^d}   \frac{|s_1 - s_2|^\gamma}{(r-s_1\vee s_2)^{\gamma + \frac{j}{\alpha}}} \rho^{j}(r-s_1 \vee s_2, z-x)(t-r)^{-\frac{1}{\alpha}}\rho^1(t-r,y-z) \, dz\,dr \\ &\leq K \left[ \frac{|s_1 - s_2|^\gamma}{(t-s_1\wedge s_2)^{\gamma + \frac{j}{\alpha}-1 + \frac{1}{\alpha}}} \rho^{j}(t-s_1 \wedge s_2, y-x) + \frac{|s_1 - s_2|^\gamma}{(t-s_1\vee s_2)^{\gamma + \frac{j}{\alpha}-1 + \frac{1}{\alpha}}} \rho^{j}(t-s_1 \vee s_2, y-x) \right] \\ &\leq K \left[ \frac{|s_1 - s_2|^\gamma}{(t-s_1\wedge s_2)^{\gamma + \frac{j}{\alpha}}} \rho^{j}(t-s_1 \wedge s_2, y-x) + \frac{|s_1 - s_2|^\gamma}{(t-s_1\vee s_2)^{\gamma + \frac{j}{\alpha}}} \rho^{j}(t-s_1 \vee s_2, y-x) \right].
\end{align*}

For $I_3$, it follows from \eqref{gradientestimatestable}, \eqref{Phi_Holder_s_Picard_lemme} and the convolution inequality \eqref{convolineqdensityref} that \begin{align*}
|I_3| &\leq K \int_{s_1 \vee s_2}^t \int_{\R^d} (r-s_1 \wedge s_2)^{-\frac{j}{\alpha}}\rho^j(r-s_1 \wedge s_2,z-x) (t-r)^{- \frac{1}{\alpha}} \rho^1(t-r,y-z) \\ &\hspace{8cm}|s_1-s_2|^\gamma \left[\frac{1}{(r-s_1\vee s_2)^{\gamma}} + \frac{1}{(r-s_1\wedge s_2)^{\gamma}}\right] \, dz\,dr \\ &\leq K |s_1-s_2|^\gamma \int_{s_1 \vee s_2}^t  (r-s_1 \vee s_2)^{-\frac{j}{\alpha}- \gamma} (t-r)^{- \frac{1}{\alpha}} \, dr \left[\rho^{j}(t-s_1\vee s_2,y-x) + \rho^{j}(t-s_1\wedge s_2,y-x)\right] \\ &\leq  K |s_1-s_2|^\gamma  (t-s_1 \vee s_2)^{-\frac{j}{\alpha}- \gamma + 1 - \frac{1}{\alpha}}  \left[\rho^{j}(t-s_1\vee s_2,y-x) + \rho^{j}(t-s_1\wedge s_2,y-x)\right] \\ &\leq K\left[ \frac{|s_1-s_2|^\gamma}{(t-s_1 \vee s_2)^{\frac{j}{\alpha}+ \gamma }} \rho^{j}(t-s_1\vee s_2,y-x) + \frac{|s_1-s_2|^\gamma}{(t-s_1 \wedge s_2)^{\frac{j}{\alpha}+ \gamma }} \rho^{j}(t-s_1\wedge s_2,y-x) \right], 
\end{align*}
where the last inequality comes from the fact that $t-s_1\wedge s_1 \leq 2(t-s_1\vee s_2)$. We finally deal with $I_4$. Owing to \eqref{gradientestimatestable} and \eqref{Phi_m_bound}, we have 

\begin{align*}
	|I_4| &\leq K  \int_{s_1\wedge s_2}^{s_1\vee s_2} \int_{\R^d}(r-s_1\wedge s_2)^{-\frac{j}{\alpha}} \rho^j(r-s_1\wedge s_2, z-x) (t-r)^{-\frac{1}{\alpha}} \rho^1(t-r,y-z) \, dz\,dr \\ &\leq  K  \int_{s_1\wedge s_2}^{s_1\vee s_2} (r-s_1\wedge s_2)^{-\frac{j}{\alpha}} (t-r)^{-\frac{1}{\alpha}} \,dr  \rho^j(t-s_1\wedge s_2, y-x) \\ &\leq   K (t-s_1\vee s_2)^{-\frac{1}{\alpha}} \int_{s_1\wedge s_2}^{s_1\vee s_2} (r-s_1\wedge s_2)^{-\frac{j}{\alpha}}  \,dr  \rho^j(t-s_1\wedge s_2, y-x) \\ &\leq   K (t-s_1\vee s_2)^{-\frac{1}{\alpha}}  |s_1 - s_2|^{1-\frac{j}{\alpha}}  \left[ \rho^j(t-s_1\vee s_2, y-x) + \rho^j(t-s_1\wedge s_2, y-x)\right].
\end{align*}
Since $\gamma \leq 1 - \frac{j}{\alpha}$ by assumption, we obtain \begin{align*}
	|I_4| &\leq    K (t-s_1\vee s_2)^{-\frac{1}{\alpha} - \gamma +1 - \frac{j}{\alpha}}  |s_1 - s_2|^{\gamma}  \left[ \rho^j(t-s_1\vee s_2, y-x) + \rho^j(t-s_1\wedge s_2, y-x)\right] \\&\leq K\left[ \frac{|s_1-s_2|^\gamma}{(t-s_1 \vee s_2)^{\frac{j}{\alpha}+ \gamma }} \rho^{j}(t-s_1\vee s_2,y-x) + \frac{|s_1-s_2|^\gamma}{(t-s_1 \wedge s_2)^{\frac{j}{\alpha}+ \gamma }} \rho^{j}(t-s_1\wedge s_2,y-x) \right], 
\end{align*}
where the last inequality comes from the fact that $t-s_1\wedge s_1 \leq 2(t-s_1\vee s_2)$ and $1 - \frac{1}{\alpha} >0$. This concludes the proof of \eqref{gradient_density_time_Holder_Picard}.\\

\subsection{Preparatory technical results}

\begin{Lemme}\label{lemma_technical8}	
	\begin{itemize}
		
		\item For any $\gamma \in (0,1]$, there exists a positive constant $K$ such that for all $m \geq 1$, $t \in (0,T]$, $s_1,s_2 \in [0,t)$, $\mu \in \PP$, $x,v \in \R^d$ \begin{equation}\label{b_lin_der_holder_s_Picard_lemme}
			\left\vert \Delta_{s_1,s_2}\delb(t,x,[X^{\cdot,\mu,(m)}_t])(v) \right\vert \leq K |s_1-s_2|^\gamma \left[(t-s_1)^{\frac{\eta}{\alpha}-\gamma} + (t-s_2)^{\frac{\eta}{\alpha}-\gamma}\right].
		\end{equation} 
		
		\item For any $\gamma \in (0,1]$, there exists a positive constant $K$ such that for all $m \geq 1$, $t \in (0,T]$, $s_1,s_2 \in [0,t)$, $\mu \in \PP$, $x,v \in \R^d$ \begin{equation}\label{b_Holder_lin_der_s_Picard_lemme2}
			\left\vert \Delta_{s_1,s_2}\del \left[b(t,x,[X^{\cdot,\mu,(m)}_t])\right](v) \right\vert \leq K \frac{|s_1-s_2|^\gamma}{(t-s_1\vee s_2)^{\gamma}}  + K \int_{\R^{2d}} \left\vert \Delta_{s_1,s_2} \del p_m(\mu,\cdot,t,x',y) \right\vert\,dy \, d\mu(x').
		\end{equation}
		
		\item For any $\gamma \in (0,1]$, there exists a positive constant $K$ such that for all $m \geq 1$, $t \in (0,T]$, $s_1,s_2 \in [0,t)$, $r\in [s_1\vee s_2,t)$, $\mu \in \PP$, $x,y,v \in \R^d$ \begin{multline}\label{H_Holder_lin_der_s_Picard_lemme2}
			\left\vert \Delta_{s_1,s_2}\del \H_{m+1}(\mu,\cdot,r,t,x,y)(v)\right\vert \leq K (t-r)^{-\frac{1}{\alpha}} \rho^{1}(t-r,y-x)\\\left[\frac{|s_1-s_2|^\gamma}{(r-s_1\vee s_2)^{\gamma}}  + \int_{\R^{2d}} \left\vert \Delta_{s_1,s_2} \del p_m(\mu,\cdot,r,x',y) \right\vert\,dy \, d\mu(x')\right].
		\end{multline}

		\item For any $\gamma \in (0,1]$, there exists a positive constant $K$ such that for all $m \geq 1$, $t \in (0,T]$, $s_1,s_2 \in [0,t)$, $\mu \in \PP$, $x,v \in \R^d$ \begin{multline}\label{b_Holder_gradient_lin_der_s_Picard_lemme2}
			\left\vert \Delta_{s_1,s_2}\pa_v\del \left[b(t,x,[X^{\cdot,\mu,(m)}_t])\right](v) \right\vert \leq K \frac{|s_1-s_2|^\gamma}{(t-s_1\vee s_2)^{\gamma + \frac{1-\eta}{\alpha}}} \\ + K \int_{\R^{2d}} (1\wedge |y-x'|^\eta) \left\vert \Delta_{s_1,s_2} \pa_v\del p_m(\mu,\cdot,t,x',y) \right\vert\,dy \, d\mu(x').
		\end{multline}
		
		\item For any $\gamma \in (0,1]$, there exists a positive constant $K$ such that for all $m \geq 1$, $t \in (0,T]$, $s_1,s_2 \in [0,t)$, $r\in [s_1\vee s_2,t)$, $\mu \in \PP$, $x,y,v \in \R^d$ \begin{multline}\label{H_Holder_gradient_lin_der_s_Picard_lemme2}
			\left\vert \Delta_{s_1,s_2}\pa_v\del \H_{m+1}(\mu,\cdot,r,t,x,y)(v)\right\vert \leq K (t-r)^{-\frac{1}{\alpha}} \rho^{1}(t-r,y-x)\\ \left[\frac{|s_1-s_2|^\gamma}{(r-s_1\vee s_2)^{\gamma + \frac{1-\eta}{\alpha}}}  + \int_{\R^{2d}} (1\wedge |y-x'|^\eta) \left\vert \Delta_{s_1,s_2} \pa_v\del p_m(\mu,\cdot,r,x',y) \right\vert\,dy \, d\mu(x')\right].
		\end{multline}
		
	\end{itemize}
\end{Lemme}

\begin{proof}[Proof of Lemma \ref{lemma_technical8}]
	
	\noindent\textbf{Proof of \eqref{b_lin_der_holder_s_Picard_lemme}.} By definition of the linear derivative and a centering argument, one can write, setting $m_\lambda := \lambda [X^{s_1\vee s_2,\mu,(m)}_t] + (1-\lambda)  [X^{s_1\wedge s_2,\mu,(m)}_t]$, for $\lambda \in [0,1]$,
	
	\begin{multline*}
		\Delta_{s_1,s_2} \del b(t,x,[X^{\cdot,\mu,(m)}_t])(v) = \int_{\R^{2d}} \left(\dell b(t,x,m_\lambda)(v,y) - \dell b(t,x,m_\lambda)(v,x')\right) \\\Delta_{s_1,s_2}p_m(\mu,\cdot,t,x',y) \, d\mu(x') \, dy \, d\lambda.
	\end{multline*}
	
	We deduce using \eqref{density_Holder_s_Picard_lemme} and uniform the $\eta$-Hölder continuity of $\dell b(t,x,\mu)(v,\cdot)$ that 
	
	\begin{align*}
		&\left\vert 	\Delta_{s_1,s_2} \del b(t,x,[X^{\cdot,\mu,(m)}_t])(v) \right\vert \\ &\leq K\int_{\R^{2d}} |y-x'|^\eta  \left[ \frac{|s_1 - s_2|^\gamma}{(t-s_1\wedge s_2)^{\gamma}} \rho^{-\tilde{\eta}}(t-s_1 \wedge s_2, y-x') \right.\\ &\left. \hspace{6cm}+ \frac{|s_1 - s_2|^\gamma}{(t-s_1\vee s_2)^{\gamma }} \rho^{-\tilde{\eta}}(t-s_1 \vee s_2, y-x') \right] \, dy \, d\mu(x').
	\end{align*}
	Notice that because $\eta + \tilde{\eta} < \alpha  $, $\rho^{-\tilde{\eta} - \eta}(t-s,\cdot)$ is integrable. We conclude by the space-time inequality \eqref{scalingdensityref} and since $\int_{\R^d}\rho^{-\tilde{\eta}-\eta}(t-s,y) \, dy$ is a finite constant independent of $s$ and $t$.\\
	
	\noindent\textbf{Proof of \eqref{b_Holder_lin_der_s_Picard_lemme2} and \eqref{H_Holder_lin_der_s_Picard_lemme2}.} Coming back to \eqref{lemma_technical_eq1}, we have the following decomposition 
	
	\begin{align*}
		&\Delta_{s_1,s_2}\del \left[b(t,x,[X^{\cdot,\mu,(m)}_t])\right](v)\\ &=\int_{\R^d} \Delta_{s_1,s_2} \del b(t,x,[X^{\cdot,\mu,(m)}_t])(y) p_m(\mu,s_1\vee s_2,t,v,y) \, dy \\ &\quad+ \int_{\R^d} \del b (t,x,[X^{s_1 \wedge s_2,\mu,(m)}_t])(y) \Delta_{s_1,s_2}p_m(\mu,\cdot,t,v,y) \,dy \\ &\quad+ \int_{\R^{2d}} \Delta_{s_1,s_2} \del b(t,x,[X^{\cdot,\mu,(m)}_t])(y) \del p_m(\mu,s_1\vee s_2,t,x',y)(v)\, d\mu(x') \, dy \\ &\quad+ \int_{\R^{2d}} \del b (t,x,[X^{s_1 \wedge s_2,\mu,(m)}_t])(y) \Delta_{s_1,s_2}\del p_m(\mu,\cdot,t,x',y)(v) \, d\mu(x') \,dy \\ &=: I_1 + I_2 + I_3 + I_4.
	\end{align*}
	
	It immediately follows from \eqref{b_lin_der_holder_s_Picard_lemme} and \eqref{density_bound_Picard} that \begin{align*}
		|I_1| &\leq K \left[\frac{|s_1-s_2|^\gamma}{(t-s_1)^{\gamma - \frac{\eta}{\alpha}}} + \frac{|s_1-s_2|^\gamma}{(t-s_2)^{\gamma - \frac{\eta}{\alpha}}} \right] \\ &\leq K \frac{|s_1-s_2|^\gamma}{(t-s_1\vee s_2)^{\gamma}}.
	\end{align*}
	Then, the boundedness of $\del b$ and \eqref{gradient_density_time_Holder_Picard} yield 
	
	\begin{align*}
		|I_2| &\leq K\int_{\R^d}  \frac{|s_1-s_2|^\gamma}{(t - s_1)^{\gamma}}\rho^{-\tilde{\eta}}(t-s_1,y-v) + \frac{|s_1-s_2|^\gamma}{(t-s_2)^{\gamma}}\rho^{-\tilde{\eta}}(t-s_2,y-v) \,dy \\ &\leq K \frac{|s_1-s_2|^\gamma}{(t - s_1\vee s_2)^{\gamma}}.
	\end{align*}
	
	Similarly to $I_1$ and using \eqref{linear_der_bound_Picard}, we deduce that  \begin{align*}
		|I_3| &\leq K \frac{|s_1-s_2|^\gamma}{(t-s_1\vee s_2)^{\gamma}}.
	\end{align*}
	
	Gathering the previous estimates on $I_1,$ $I_2$, $I_3$ and noting that $\del b$ is bounded, we conclude that \eqref{b_Holder_lin_der_s_Picard_lemme2} holds true. Remark that \eqref{H_Holder_lin_der_s_Picard_lemme2} follows directly from the expression of $\del \H_{m+1}$ given in \eqref{expression_linear_der_H} and from \eqref{b_Holder_lin_der_s_Picard_lemme2}. \\
	
	\noindent\textbf{Proof of \eqref{b_Holder_gradient_lin_der_s_Picard_lemme2}.} We use \eqref{b_expression_gradient_linr_der_Picard}, which yields the following decomposition 
	
	\begin{align*}
		&\Delta_{s_1,s_2}\pa_v\del \left[b(t,x,[X^{\cdot,\mu,(m)}_t])\right](v)\\ &=\int_{\R^d} \Delta_{s_1,s_2} \del b(t,x,[X^{\cdot,\mu,(m)}_t])(y) \pa_x p_m(\mu,s_1\vee s_2,t,v,y) \, dy \\ &\quad+ \int_{\R^d} \del b (t,x,[X^{s_1 \wedge s_2,\mu,(m)}_t])(y) \Delta_{s_1,s_2}\pa_x p_m(\mu,\cdot,t,v,y) \,dy \\ &\quad+ \int_{\R^{2d}} \Delta_{s_1,s_2} \del b(t,x,[X^{\cdot,\mu,(m)}_t])(y) \pa_v\del p_m(\mu,s_1\vee s_2,t,x',y)(v)\, d\mu(x') \, dy \\ &\quad+ \int_{\R^{2d}} \del b (t,x,[X^{s_1 \wedge s_2,\mu,(m)}_t])(y) \Delta_{s_1,s_2}\pa_v\del p_m(\mu,\cdot,t,x',y)(v) \, d\mu(x') \,dy \\ &=: I_1 + I_2 + I_3 + I_4.
	\end{align*}
	
	It immediately follows from \eqref{b_lin_der_holder_s_Picard_lemme} and \eqref{density_bound_Picard} that \begin{align*}
		|I_1| &\leq K \left[\frac{|s_1-s_2|^\gamma}{(t-s_1)^{\gamma - \frac{\eta}{\alpha}}} + \frac{|s_1-s_2|^\gamma}{(t-s_2)^{\gamma - \frac{\eta}{\alpha}}} \right](t-s_1\vee s_2)^{-\frac{1}{\alpha}} \\ &\leq K \frac{|s_1-s_2|^\gamma}{(t-s_1\vee s_2)^{\gamma + \frac{1 - \eta}{\alpha}}}.
	\end{align*}
	
	We rewrite $I_2$ in the following form 
	
	\begin{equation*}
		I_2 = \int_{\R^d} \left(\del b (t,x,[X^{s_1 \wedge s_2,\mu,(m)}_t])(y) - \del b (t,x,[X^{s_1 \wedge s_2,\mu,(m)}_t])(v) \right) \Delta_{s_1,s_2}\pa_x p_m(\mu,\cdot,t,v,y) \,dy. 
	\end{equation*}
	Then, the $\eta$-Hölder continuity of $\del b (t,x,\mu)(\cdot)$ and \eqref{gradient_density_time_Holder_Picard} yield 
	
	\begin{align*}
		|I_2| &\leq K\int_{\R^d}  \frac{|s_1-s_2|^\gamma}{(t - s_1)^{\gamma + \frac{1}{\alpha}}}|y-v|^\eta\rho^{j}(t-s_1,y-v) + \frac{|s_1-s_2|^\gamma}{(t-s_2)^{\gamma + \frac{1}{\alpha}}}|y-v|^\eta\rho^{j}(t-s_2,y-v) \,dy.
	\end{align*}
	The space-time inequality \eqref{scalingdensityref} ensures that 
	
	\begin{align*}
		|I_2| &\leq K \frac{|s_1-s_2|^\gamma}{(t - s_1\vee s_2)^{\gamma + \frac{1-\eta}{\alpha}}}.
	\end{align*}
	
	By using \eqref{b_lin_der_holder_s_Picard_lemme} and \eqref{gradient_linear_der_bound_Picard}, we deduce that  \begin{align*}
		|I_3| &  \leq K\left[ \frac{|s_1-s_2|^\gamma}{(t-s_1)^{\gamma - \frac{\eta}{\alpha}}} + \frac{|s_1-s_2|^\gamma}{(t-s_2)^{\gamma - \frac{\eta}{\alpha}}}\right](t-s_1\vee s_2)^{\frac{\eta -1}{\alpha} + 1- \frac{1}{\alpha}}\\ &\leq K \frac{|s_1-s_2|^\gamma}{(t-s_1\vee s_2)^{\gamma + \frac{1-\eta}{\alpha}}}.
	\end{align*}
	
	Note that since $\int_{\R^d} \pa_v \del p_m(\mu,s,t,x,y)(v) \, dy =0$, one can rewrite $I_4$ as 
	
	\begin{equation*}
		I_4 = \int_{\R^{2d}} \left(\del b (t,x,[X^{s_1 \wedge s_2,\mu,(m)}_t])(y) - \del b (t,x,[X^{s_1 \wedge s_2,\mu,(m)}_t])(x') \right)\Delta_{s_1,s_2}\pa_v\del p_m(\mu,\cdot,t,x',y)(v) \, d\mu(x') \,dy.
	\end{equation*}
	The boundedness of $\del b$ and the $\eta$-Hölder continuity of $\del b(t,x,\mu)(\cdot)$ yield
	
	\begin{align*}
		|I_4| &\leq \int_{\R^{2d}} (1\wedge |y-x'|^\eta) \left\vert \Delta_{s_1,s_2} \pa_v\del p_m(\mu,\cdot,r,x',y) \right\vert\,dy \, d\mu(x').
	\end{align*}
	
	Gathering the previous estimates on $I_1,$ $I_2$, $I_3$ and $I_4$, we conclude that \eqref{b_Holder_gradient_lin_der_s_Picard_lemme2} holds true. Finally, notice that \eqref{H_Holder_gradient_lin_der_s_Picard_lemme2} follows directly from the expression of $\pa_v \del \H_{m+1}$ obtained by differentiating  \eqref{expression_linear_der_H} with respect to $v$ and from \eqref{b_Holder_gradient_lin_der_s_Picard_lemme2}.

\end{proof}

\subsection{Fifth part of the induction step} We prove here that the estimates \eqref{lin_der_density_holder_time_Picard} and \eqref{gradient_lin_der_density_holder_time_Picard} hold true.\\

\noindent\textbf{Proof of \eqref{lin_der_density_holder_time_Picard}.} We separate the proof of the induction step into two disjoint cases. First, we assume that $|s_1-s_2| > t- s_1\vee s_2$. By \eqref{linear_der_bound_Picard} (the series appearing in the bound being convergent), we can write \begin{align*}
	&\left\vert \Delta_{s_1,s_2}\del p_m(\mu,\cdot,t,x,y)(v)\right\vert \\&\leq  	\left\vert \del p_m(\mu,s_1,t,x,y)(v)\right\vert + 	\left\vert \del p_m(\mu,s_2,t,x,y)(v)\right\vert \\ &\leq K \left[(t-s_1\vee s_2)^{1-\frac{1}{\alpha}}\rho^0(t-s_1\vee s_2,y-x) + (t-s_1\wedge s_2)^{1-\frac{1}{\alpha}}\rho^0(t-s_1\wedge s_2,y-x)\right]\\ &\leq K \left[\frac{|s_1-s_2|^\gamma}{(t-s_1\vee s_2)^{\gamma +\frac{1}{\alpha} -1}} \rho^0(t-s_1\vee s_2,y-x) + \frac{|s_1-s_2|^\gamma + (t-s_1\vee s_2)^\gamma}{(t-s_1\wedge s_2)^\gamma(t-s_1\wedge s_2)^{\frac{1}{\alpha}-1}} \rho^0(t-s_1\wedge s_2,y-x) \right]\\ &\leq   K \left[\frac{|s_1-s_2|^\gamma}{(t-s_1\vee s_2)^{\gamma + \frac{1}{\alpha}-1}} \rho^{0}(t-s_1\vee s_2,y-x) + \frac{|s_1-s_2|^\gamma }{(t-s_1\wedge s_2)^{\gamma +\frac{1}{\alpha}-1}} \rho^{0}(t-s_1\wedge s_2,y-x) \right],
\end{align*} This shows that \eqref{lin_der_density_holder_time_Picard} holds true at step $m+1$ provided that we choose $C\geq K$ in \eqref{lin_der_density_holder_time_Picard}. We now turn to the case $|s_1 - s_2|\leq t-s_1\vee s_2$. Using the representation formula \eqref{representation_linear_der}, we get that 

\begin{align*}
	\Delta_{s_1,s_2} \del p_{m+1}(\mu,\cdot,t,x,y)(v) &= \Delta_{s_1,s_2} \left(p_{m+1}\otimes \del \H_{m+1} \right)(\mu,\cdot,t,x,y)(v) \\ &\quad+ \Delta_{s_1,s_2} \left(p_{m+1}\otimes\del \H_{m+1} \otimes \Phi_{m+1}\right)(\mu,\cdot,t,x,y)(v) \\ &=: I_1 + I_2.
\end{align*}

Then, we decompose $I_1$ in the following way 

\begin{align*}
	I_1 &= \int_{s_1\vee s_2}^t \int_{\R^d} \Delta_{s_1,s_2} p_{m+1}(\mu,\cdot,r,x,z) \del \H_{m+1}(\mu,s_1\vee s_2,r,t,z,y) \, dz \,dr \\ &\quad+ \int_{s_1\vee s_2}^t \int_{\R^d}  p_{m+1}(\mu,s_1 \wedge s_2,r,x,z) \Delta_{s_1,s_2}\del \H_{m+1}(\mu,\cdot,r,t,z,y) \, dz \,dr\\ &\quad - \int_{s_1\wedge s_2}^{s_1\vee s_2} \int_{\R^d} p_{m+1}(\mu,s_1\wedge s_2,r,x,z) \del \H_{m+1}(\mu,s_1\wedge s_2,r,t,z,y) \, dz \,dr \\ &=: I_{1,1} + I_{1,2} + I_{1,3}.
\end{align*}

It follows from \eqref{gradient_density_time_Holder_Picard}, \eqref{H_linear_derivative_bound_Picard} (the series appearing in the bound being convergent) and the convolution inequality \eqref{convolineqdensityref} that 

\begin{align*}
	|I_{1,1}| &\leq K\int_{s_1 \vee s_2}^t \int_{\R^d} \frac{|s_1-s_2|^\gamma}{(r - s_1\vee s_2)^{\gamma}}\left[\rho^{0}(r-s_1\vee s_2,y-x) +\rho^{0}(r-s_1 \wedge s_2,y-x)\right] \\ & \hspace{12cm}(t-r)^{-\frac{1}{\alpha}} \rho^1(t-r,y-z)\, dz\,dr \\ &\leq  K\frac{|s_1-s_2|^\gamma}{(t- s_1\vee s_2)^{\gamma+\frac{1}{\alpha}-1}}\left[\rho^{0}(t-s_1\vee s_2,y-x) \rho^{0}(t-s_1 \wedge s_2,y-x)\right].
\end{align*}
We now deal with $I_{1,2}$. Using the induction assumption and \eqref{H_Holder_lin_der_s_Picard_lemme2}, we deduce that

\begin{align}\label{proof_lin_der_holder_times_eq4}
	&\notag\left\vert \Delta_{s_1,s_2}\del \H_{m+1}(\mu,\cdot,r,t,x,y)(v)\right\vert\\ \notag& \leq K (t-r)^{-\frac{1}{\alpha}} \rho^{1}(t-r,y-x)\left[\frac{|s_1-s_2|^\gamma}{(r-s_1\vee s_2)^{\gamma}}  + \int_{\R^{2d}} \left\vert \Delta_{s_1,s_2} \del p_m(\mu,\cdot,r,x',y) \right\vert\,dy \, d\mu(x')\right] \\ &\leq   K (t-r)^{-\frac{1}{\alpha}} \rho^{1}(t-r,y-x)\frac{|s_1-s_2|^\gamma}{(r-s_1\vee s_2)^{\gamma}}  \\ \notag&\quad +  K (t-r)^{-\frac{1}{\alpha}} \rho^{1}(t-r,y-x) \left[ \frac{|s_1-s_2|^\gamma}{(r-s_1\vee s_2)^{\gamma + \frac{1}{\alpha}-1}} + \frac{|s_1-s_2|^\gamma}{(r-s_1 \wedge s_2)^{\gamma + \frac{1}{\alpha}-1}} \right]\\\notag &\hspace{3cm} \sum_{k=1}^{m} C^k(r-s_1\vee s_2)^{(k-1)\left(1 - \frac{1}{\alpha}\right)} \prod_{j=1}^{k-1}\BB\left(\left[2-\gamma - \frac{1}{\alpha} \right] \wedge 1 +(j-1)\left(1- \frac{1}{\alpha}\right), 1 - \frac{1}{\alpha}\right).
\end{align} 

Notice that \[ \begin{cases}
	&(r-s_1 \wedge s_2)^{\gamma + \frac{1}{\alpha}-1} \geq (r-s_1 \vee s_2)^{\gamma + \frac{1}{\alpha}-1}, \quad\text{if $ \gamma + \frac{1}{\alpha}-1 \geq 0$,}\\  &(r-s_1 \wedge s_2)^{\gamma + \frac{1}{\alpha}-1} \geq (t-s_1 \wedge s_2)^{\gamma + \frac{1}{\alpha}-1} \geq 2^{\gamma + \frac{1}{\alpha}-1} (t-s_1\vee s_2)^{\gamma + \frac{1}{\alpha}-1},\quad \text{if $ \gamma + \frac{1}{\alpha}-1 < 0$,}
\end{cases}\]
since $t-s_1\wedge s_2 \leq 2 (t-s_1 \vee s_2)$. Using these inequalities, one has 

\begin{align*}
	&\int_{s_1\vee s_2}^t (t-r)^{-\frac{1}{\alpha}}  \left[ \frac{|s_1-s_2|^\gamma}{(r-s_1\vee s_2)^{\gamma + \frac{1}{\alpha}-1}} + \frac{|s_1-s_2|^\gamma}{(r-s_1 \wedge s_2)^{\gamma + \frac{1}{\alpha}-1}} \right](r-s_1\vee s_2)^{(k-1)\left(1 - \frac{1}{\alpha}\right)} \,dr\\ &\leq K \frac{|s_1-s_2|^\gamma}{(t-s_1\vee s_2)^{\gamma + \frac{1}{\alpha}-1}}(t-s_1\vee s_2)^{k\left(1 - \frac{1}{\alpha}\right)} \BB \left(1 - \gamma + k \left(1 - \frac{1}{\alpha}\right), 1 - \frac{1}{\alpha}\right) \\ &\quad+ K \frac{|s_1-s_2|^\gamma}{(t-s_1\vee s_2)^{\gamma + \frac{1}{\alpha}-1}}(t-s_1\vee s_2)^{k\left(1 - \frac{1}{\alpha}\right)} \BB \left(1 + (k-1) \left(1 - \frac{1}{\alpha}\right), 1 - \frac{1}{\alpha}\right) \\ &\leq  K \frac{|s_1-s_2|^\gamma}{(t-s_1\vee s_2)^{\gamma + \frac{1}{\alpha}-1}}(t-s_1\vee s_2)^{k\left(1 - \frac{1}{\alpha}\right)} \BB \left( \left[2 - \gamma - \frac{1}{\alpha} \right]\wedge 1 + (k-1) \left(1 - \frac{1}{\alpha}\right), 1 - \frac{1}{\alpha}\right), 
\end{align*}
since the Beta function is decreasing with respect to its first argument. From this inequality, \eqref{proof_lin_der_holder_times_eq4} and \eqref{density_bound_Picard}, we obtain 

\begin{align}\label{proof_lin_der_holder_times_eq5}
	|I_{1,2}| &\leq K  \frac{|s_1-s_2|^\gamma}{(t-s_1\vee s_2)^{\gamma + \frac{1}{\alpha}-1}} \rho^0 (t-s_1\wedge s_2 ,y-x)  \\ \notag& \hspace{2cm}\left[1+  \sum_{k=1}^{m} C^k(t-s_1\vee s_2)^{k\left(1 - \frac{1}{\alpha}\right)} \prod_{j=1}^k\BB\left(\left[2-\gamma - \frac{1}{\alpha} \right] \wedge 1 +(j-1)\left(1- \frac{1}{\alpha}\right), 1 - \frac{1}{\alpha}\right) \right].
\end{align}
We now turn to estimate $I_{1,3}$. Thanks to \eqref{density_bound_Picard} and \eqref{H_linear_derivative_bound_Picard}, one has 

\begin{align*}
	|I_{1,3}| & \leq \int_{s_1 \wedge s_2}^{s_1 \vee s_2} \rho^0(r-s_1\wedge s_2,z-x) (t-r)^{-\frac{1}{\alpha}}\rho^1(t-r,y-z) \, dz \,dr \\ &\leq K \frac{|s_1-s_2|}{(t-s_1\vee s_2)^{\frac{1}{\alpha}}} \rho^0(t-s_1\wedge s_2,y-x) \\ &\leq   K \frac{|s_1-s_2|^\gamma}{(t-s_1\vee s_2)^{\gamma +\frac{1}{\alpha}-1}} \rho^0(t-s_1\wedge s_2,y-x).
\end{align*}
Gathering the preceding estimates on $I_{1,1}$ and $I_{1,2}$ and $I_{1,3}$, we have proved that \begin{align}\label{proof_lin_der_holder_times_eq2}
	|I_1| &\leq K  \frac{|s_1-s_2|^\gamma}{(t-s_1\vee s_2)^{\gamma + \frac{1}{\alpha}-1}} \left[\rho^{0} (t-s_1 \vee s_2 ,y-x) +  \rho^{0}(t-s_1 \wedge s_2 ,y-x)\right] \\\notag & \hspace{2cm}\left[1+  \sum_{k=1}^{m} C^k(t-s_1\vee s_2)^{k\left(1 - \frac{1}{\alpha}\right)} \prod_{j=1}^k\BB\left(\left[2-\gamma - \frac{1}{\alpha} \right] \wedge 1 +(j-1)\left(1- \frac{1}{\alpha}\right), 1 - \frac{1}{\alpha}\right) \right].
\end{align}

As done before for $I_1$, we decompose $I_2$ in the following way 

\begin{align*}
	I_2 &= \int_{s_1\vee s_2}^t \int_{\R^d} \Delta_{s_1,s_2} \left(p_{m+1}\otimes \del \H_{m+1}\right)(\mu,\cdot,r,x,z) \Phi_{m+1}(\mu,s_1\vee s_2,r,t,z,y) \, dz \,dr \\ &\quad+ \int_{s_1\vee s_2}^t \int_{\R^d}  \left(p_{m+1} \otimes \del \H_{m+1}\right)(\mu,s_1 \wedge s_2,r,x,z) \Delta_{s_1,s_2}\Phi_{m+1}(\mu,\cdot,r,t,z,y) \, dz \,dr\\ &\quad - \int_{s_1\wedge s_2}^{s_1\vee s_2} \int_{\R^d} \left(p_{m+1}\otimes \del \H_{m+1}\right)(\mu,s_1\wedge s_2,r,x,z)  \Phi_{m+1}(\mu,s_1\wedge s_2,r,t,z,y) \, dz \,dr \\ &=: I_{2,1} + I_{2,2} + I_{2,3}.
\end{align*}
We follow the same lines of reasoning as for $I_1$. Using the bound \eqref{proof_lin_der_holder_times_eq2} previously obtained for $I_1$ and \eqref{Phi_m_bound}, we have that \begin{align*}
	|I_{2,1}| & \leq  K  \frac{|s_1-s_2|^\gamma}{(t-s_1\vee s_2)^{\gamma + \frac{1}{\alpha}-1}} (t-s_1 \vee s_2)^{ 1 - \frac{1}{\alpha}} \left[\rho^{0} (t-s_1 \vee s_2 ,y-x) +  \rho^{0}(t-s_1 \wedge s_2 ,y-x)\right] \\ & \hspace{2cm}\left[1+  \sum_{k=1}^{m} C^k(t-s_1\vee s_2)^{k\left(1 - \frac{1}{\alpha}\right)} \prod_{j=1}^k\BB\left(\left[2-\gamma - \frac{1}{\alpha} \right] \wedge 1 +(j-1)\left(1- \frac{1}{\alpha}\right), 1 - \frac{1}{\alpha}\right) \right].
\end{align*}
Notice that \eqref{density_bound_Picard}, \eqref{H_linear_derivative_bound_Picard} and the convolution inequality \eqref{convolineqdensityref} yield \begin{equation}\label{proof_lin_deriv_holder_time_eq1} \left\vert p_{m+1}\otimes \del \H_{m+1}(\mu,s,t,x,y)(v)\right\vert \leq K (t-s)^{1 - \frac{1}{\alpha}} \rho^0(t-s,y-x).\end{equation}
Using this inequality and \eqref{Phi_Holder_s_Picard_lemme}, we deduce that \begin{align*}
	|I_{2,2}| &\leq K \int_{s_1\vee s_2}^t \int_{\R^d} (t-s_1\wedge s_2)^{1 - \frac{1}{\alpha}} \rho^0(r-s_1 \wedge s_2,z-x) \frac{|s_1 - s_2|^\gamma }{(r-s_1\vee s_2)^\gamma} (t-r)^{-\frac{1}{\alpha}} \rho^1(t-r,y-z) \, dz \,dr.
\end{align*}

Since $t-s_1\wedge s_2 \leq 2 (t-s_1\vee s_2)$, we get 

\begin{align*}
	|I_{2,2}| & \leq  K \frac{|s_1 - s_2|^\gamma }{(t-s_1\vee s_2)^{\gamma + \frac{1}{\alpha}-1}}(t-s_1\vee s_2)^{1 - \frac{1}{\alpha}} \rho^0(t-s_1 \wedge s_2,y-x).
\end{align*}
For $I_{2,3}$, \eqref{proof_lin_deriv_holder_time_eq1} and \eqref{Phi_m_bound} yield 

\begin{align*}
	|I_{2,3}| & \leq \int_{s_1 \wedge s_2}^{s_1 \vee s_2} (r-s_1\wedge s_2)^{1 - \frac{1}{\alpha}}\rho^0(r-s_1\wedge s_2,z-x) (t-r)^{-\frac{1}{\alpha}}\rho^1(t-r,y-z) \, dz \,dr \\ &\leq K \frac{|s_1-s_2|^{1 + 1 - \frac{1}{\alpha}}}{(t-s_1\vee s_2)^{\frac{1}{\alpha}}} \rho^0(t-s_1\wedge s_2,y-x) \\ &\leq   K \frac{|s_1-s_2|^\gamma}{(t-s_1\vee s_2)^{\gamma +\frac{1}{\alpha}-1}} (t-s_1 \vee s_2)^{1 - \frac{1}{\alpha}} \rho^0(t-s_1\wedge s_2,y-x),
\end{align*}
since $|s_1 - s_2| \leq t - s_1\vee s_2$. Gathering the estimates obtained, we have proved that \begin{align}\label{proof_lin_der_holder_times_eq3}
	|I_2| \notag&\leq K  \frac{|s_1-s_2|^\gamma}{(t-s_1\vee s_2)^{\gamma + \frac{1}{\alpha}-1}}(t -s_1\vee s_2)^{1 - \frac{1}{\alpha}} \left[\rho^{0} (t-s_1 \vee s_2 ,y-x) +  \rho^{0}(t-s_1 \wedge s_2 ,y-x)\right] \\ & \hspace{1cm}\left[1+  \sum_{k=1}^{m} C^k(t-s_1\vee s_2)^{k\left(1 - \frac{1}{\alpha}\right)} \prod_{j=1}^k\BB\left(\left[2-\gamma - \frac{1}{\alpha} \right] \wedge 1 +(j-1)\left(1- \frac{1}{\alpha}\right), 1 - \frac{1}{\alpha}\right) \right].
\end{align}

We conclude by \eqref{proof_lin_der_holder_times_eq2}, \eqref{proof_lin_der_holder_times_eq3} and the fact that $t-s_1 \vee s_2 \leq t-s_1\wedge s_2 \leq 2(t-s_1 \vee s_2)$ that 

\begin{align*}
	&\left\vert \Delta_{s_1,s_2} \del p_{m+1}(\mu,\cdot,t,x,y)(v) \right\vert\\ &\leq K  \frac{|s_1-s_2|^\gamma}{(t-s_1\vee s_2)^{\gamma + \frac{1}{\alpha}-1}} \left[\rho^{0} (t-s_1 \vee s_2 ,y-x) +  \rho^{0}(t-s_1 \wedge s_2 ,y-x)\right] \\\notag & \hspace{3cm}\left[1+  \sum_{k=1}^{m} C^k(t-s_1\vee s_2)^{k\left(1 - \frac{1}{\alpha}\right)} \prod_{j=1}^k\BB\left(\left[2-\gamma - \frac{1}{\alpha} \right] \wedge 1 +(j-1)\left(1- \frac{1}{\alpha}\right), 1 - \frac{1}{\alpha}\right) \right] \\ &\leq    \left[\frac{|s_1-s_2|^\gamma}{(t-s_1\vee s_2)^{\gamma + \frac{1}{\alpha}-1}}\rho^{0} (t-s_1 \vee s_2 ,y-x) +  \frac{|s_1-s_2|^\gamma}{(t-s_1\wedge s_2)^{\gamma + \frac{1}{\alpha}-1}} \rho^{0}(t-s_1 \wedge s_2 ,y-x)\right] \\\notag & \hspace{3cm}  \sum_{k=1}^{m+1} C^k(t-s_1\vee s_2)^{(k-1)\left(1 - \frac{1}{\alpha}\right)} \prod_{j=1}^{k-1}\BB\left(\left[2-\gamma - \frac{1}{\alpha} \right] \wedge 1 +(j-1)\left(1- \frac{1}{\alpha}\right), 1 - \frac{1}{\alpha}\right),
\end{align*}
provided that we choose $C \geq K$ in \eqref{lin_der_density_holder_time_Picard}. This ends the proof of the induction step for \eqref{lin_der_density_holder_time_Picard}.\\

\noindent\textbf{Proof of \eqref{gradient_lin_der_density_holder_time_Picard}.} We start by treating the case $|s_1-s_2| > t- s_1\vee s_2$. Reasoning as before by using \eqref{gradient_linear_der_bound_Picard} (the series appearing in the bound being convergent), we get that \begin{align*}
	&\left\vert \Delta_{s_1,s_2}\pa_v\del p_m(\mu,\cdot,t,x,y)(v)\right\vert  \\ &\leq K \left[(t-s_1\vee s_2)^{\frac{\eta-1}{\alpha}+1-\frac{1}{\alpha}}\rho^0(t-s_1\vee s_2,y-x) + (t-s_1\wedge s_2)^{\frac{\eta-1}{\alpha} +1-\frac{1}{\alpha}}\rho^0(t-s_1\wedge s_2,y-x)\right]\\ &\leq   K \left[\frac{|s_1-s_2|^\gamma}{(t-s_1\vee s_2)^{\gamma + \frac{1}{\alpha}-1 + \frac{1-\eta}{\alpha}}} \rho^{0}(t-s_1\vee s_2,y-x) + \frac{|s_1-s_2|^\gamma }{(t-s_1\wedge s_2)^{\gamma +\frac{1}{\alpha}-1+ \frac{1-\eta}{\alpha}}} \rho^{0}(t-s_1\wedge s_2,y-x) \right].
\end{align*}

This shows that \eqref{gradient_lin_der_density_holder_time_Picard} holds true at step $m+1$ provided that we choose $C\geq K$ in \eqref{gradient_lin_der_density_holder_time_Picard}. We now turn to the case $|s_1 - s_2|\leq t-s_1\vee s_2$. By the representation formula \eqref{representation_gradient_linear_der}, we write

\begin{align*}
	\Delta_{s_1,s_2} \pa_v\del p_{m+1}(\mu,\cdot,t,x,y)(v) &= \Delta_{s_1,s_2} \left(p_{m+1}\otimes \pa_v\del \H_{m+1} \right)(\mu,\cdot,t,x,y)(v) \\ &\quad+ \Delta_{s_1,s_2} \left(p_{m+1}\otimes\pa_v\del \H_{m+1} \otimes \Phi_{m+1}\right)(\mu,\cdot,t,x,y)(v) \\ &=: I_1 + I_2.
\end{align*}

Then, we decompose $I_1$ in the following way 

\begin{align*}
	I_1 &= \int_{s_1\vee s_2}^t \int_{\R^d} \Delta_{s_1,s_2} p_{m+1}(\mu,\cdot,r,x,z) \pa_v \del \H_{m+1}(\mu,s_1\vee s_2,r,t,z,y) \, dz \,dr \\ &\quad+ \int_{s_1\vee s_2}^t \int_{\R^d}  p_{m+1}(\mu,s_1 \wedge s_2,r,x,z) \Delta_{s_1,s_2}\pa_v\del \H_{m+1}(\mu,\cdot,r,t,z,y) \, dz \,dr\\ &\quad - \int_{s_1\wedge s_2}^{s_1\vee s_2} \int_{\R^d} p_{m+1}(\mu,s_1\wedge s_2,r,x,z) \pa_v\del \H_{m+1}(\mu,s_1\wedge s_2,r,t,z,y) \, dz \,dr \\ &=: I_{1,1} + I_{1,2} + I_{1,3}.
\end{align*}

It follows from \eqref{gradient_density_time_Holder_Picard}, \eqref{H_gradient_linear_derivative_bound_Picard} (the series appearing in the bound being convergent), the convolution inequality \eqref{convolineqdensityref} and since $\gamma < 1 +\frac{\eta -1}{\alpha}$ that 

\begin{align*}
	|I_{1,1}| &\leq K\int_{s_1 \vee s_2}^t \int_{\R^d} \frac{|s_1-s_2|^\gamma}{(r - s_1\vee s_2)^{\gamma}}\left[\rho^{0}(r-s_1\vee s_2,y-x) +\rho^{0}(r-s_1 \wedge s_2,y-x)\right] \\ & \hspace{7cm}(r-s_1\vee s_2)^{\frac{\eta -1}{\alpha}}(t-r)^{-\frac{1}{\alpha}} \rho^1(t-r,y-z)\, dz\,dr \\ &\leq  K\frac{|s_1-s_2|^\gamma}{(t- s_1\vee s_2)^{\gamma+\frac{1}{\alpha}-1 + \frac{1-\eta}{\alpha}}}\left[\rho^{0}(t-s_1\vee s_2,y-x) + \rho^{0}(t-s_1 \wedge s_2,y-x)\right].
\end{align*}

We now deal with $I_{1,2}$. Using the induction assumption and \eqref{H_Holder_gradient_lin_der_s_Picard_lemme2}, we deduce that

\begin{align*}
	&\left\vert \Delta_{s_1,s_2}\pa_v\del \H_{m+1}(\mu,\cdot,r,t,x,y)(v)\right\vert\\ & \leq K (t-r)^{-\frac{1}{\alpha}} \rho^{1}(t-r,y-x)\left[\frac{|s_1-s_2|^\gamma}{(r-s_1\vee s_2)^{\gamma + \frac{1-\eta}{\alpha}}}  + \int_{\R^{2d}}(1\wedge |y-x'|^\eta) \left\vert \Delta_{s_1,s_2} \pa_v\del p_m(\mu,\cdot,r,x',y) \right\vert\,dy \, d\mu(x')\right] \\ &\leq   K (t-r)^{-\frac{1}{\alpha}} \rho^{1}(t-r,y-x)\frac{|s_1-s_2|^\gamma}{(r-s_1\vee s_2)^{\gamma+ \frac{1-\eta}{\alpha}}}   \\ \notag&\quad + K (t-r)^{-\frac{1}{\alpha}} \rho^{1}(t-r,y-x)\left[ \frac{|s_1-s_2|^\gamma  (r-s_1\vee s_2)^{\frac{\eta}{\alpha}} }{(r-s_1\vee s_2)^{\gamma + \frac{1}{\alpha}-1 + \frac{1-\eta}{\alpha}}} + \frac{|s_1-s_2|^\gamma (r-s_1\wedge s_2)^{\frac{\eta}{\alpha}} }{(r-s_1\wedge s_2)^{\gamma + \frac{1}{\alpha}-1 + \frac{1-\eta}{\alpha}}} \right]\\ &\notag \hspace{2cm}\sum_{k=1}^{m} C^k(r-s_1\vee s_2)^{(k-1)\left(1 + \frac{\eta-1}{\alpha}\right)} \prod_{j=1}^{k-1}\BB\left(\left[2\left(1 + \frac{\eta-1}{\alpha}\right)-\gamma\right] \wedge 1 +(j-1)\left(1+ \frac{\eta -1}{\alpha}\right), 1 - \frac{1}{\alpha}\right).
\end{align*} 
Note that \begin{equation}\label{proof_gradient_lin_der_holder_times_disjonction}
	\begin{cases}
		&(r-s_1 \wedge s_2)^{\gamma + \frac{1}{\alpha}-1 + \frac{1-2\eta}{\alpha}} \geq (r-s_1 \vee s_2)^{\gamma + \frac{1}{\alpha}-1+ \frac{1-2\eta}{\alpha}}, \quad\text{if $ \gamma + \frac{1}{\alpha}-1+ \frac{1-2\eta}{\alpha} \geq 0$, }\\  &(r-s_1 \wedge s_2)^{\gamma + \frac{1}{\alpha}-1+ \frac{1-2\eta}{\alpha}}  \geq 2^{\gamma + \frac{1}{\alpha}-1 + \frac{1-2\eta}{\alpha}} (t-s_1\vee s_2)^{\gamma + \frac{1}{\alpha}-1+ \frac{1-2\eta}{\alpha}},\quad \text{if $ \gamma + \frac{1}{\alpha}-1+ \frac{1-2\eta}{\alpha} < 0$,}
\end{cases}\end{equation}
since $t-s_1\wedge s_2 \leq 2 (t-s_1 \vee s_2)$. It yields

\begin{align*}
	&\int_{s_1\vee s_2}^t (t-r)^{-\frac{1}{\alpha}}  \left[ \frac{|s_1-s_2|^\gamma}{(r-s_1\vee s_2)^{\gamma + \frac{1}{\alpha}-1+ \frac{1-2\eta}{\alpha}}} + \frac{|s_1-s_2|^\gamma}{(r-s_1 \wedge s_2)^{\gamma + \frac{1}{\alpha}-1+ \frac{1-2\eta}{\alpha}}} \right](r-s_1\vee s_2)^{(k-1)\left(1 + \frac{\eta-1}{\alpha}\right)} \,dr\\ &\leq K \frac{|s_1-s_2|^\gamma}{(t-s_1\vee s_2)^{\gamma + \frac{1}{\alpha}-1+ \frac{1-\eta}{\alpha}}}(t-s_1\vee s_2)^{k\left(1 + \frac{\eta-1}{\alpha}\right)} \BB \left(2 \left(1 + \frac{\eta-1}{\alpha}\right)-\gamma + (k-1) \left(1 - \frac{1}{\alpha}\right), 1 - \frac{1}{\alpha}\right) \\ &\quad+ K \frac{|s_1-s_2|^\gamma}{(t-s_1\vee s_2)^{\gamma + \frac{1}{\alpha}-1+ \frac{1-\eta}{\alpha}}}(t-s_1\vee s_2)^{k\left(1 + \frac{\eta-1}{\alpha}\right)} \BB \left(1 + (k-1) \left(1 - \frac{1}{\alpha}\right), 1 - \frac{1}{\alpha}\right) \\ &\leq  K \frac{|s_1-s_2|^\gamma}{(t-s_1\vee s_2)^{\gamma + \frac{1}{\alpha}-1+ \frac{1-\eta}{\alpha}}}(t-s_1\vee s_2)^{k\left(1 + \frac{\eta-1}{\alpha}\right)} \\ &\hspace{4cm}\BB\left(\left[2\left(1 + \frac{\eta-1}{\alpha}\right)-\gamma\right] \wedge 1 +(j-1)\left(1+ \frac{\eta -1}{\alpha}\right), 1 - \frac{1}{\alpha}\right), 
\end{align*}
since the Beta function is decreasing with respect to its first argument. From this inequality, \eqref{proof_lin_der_holder_times_eq4} and \eqref{density_bound_Picard}, we deduce that 

\begin{align*}
	|I_{1,2}| &\leq K  \frac{|s_1-s_2|^\gamma}{(t-s_1\vee s_2)^{\gamma + \frac{1}{\alpha}-1 + \frac{1-\eta}{\alpha}}} \rho^0 (t-s_1\wedge s_2 ,y-x)  \\ & \hspace{1cm}\left[1+  \sum_{k=1}^{m} C^k(t-s_1\vee s_2)^{k\left(1 + \frac{\eta-1}{\alpha}\right)} \prod_{j=1}^k\BB\left(\left[2\left(1 + \frac{\eta-1}{\alpha}\right)-\gamma\right] \wedge 1 +(j-1)\left(1+ \frac{\eta -1}{\alpha}\right), 1 - \frac{1}{\alpha}\right) \right].
\end{align*}

We now turn to estimate $I_{1,3}$. Thanks to \eqref{density_bound_Picard} and \eqref{H_gradient_linear_derivative_bound_Picard}, one has 

\begin{align*}
	|I_{1,3}| & \leq \int_{s_1 \wedge s_2}^{s_1 \vee s_2} \rho^0(r-s_1\wedge s_2,z-x)(r-s_1\wedge s_2)^{\frac{\eta-1}{\alpha}} (t-r)^{-\frac{1}{\alpha}}\rho^1(t-r,y-z) \, dz \,dr \\ &\leq K \frac{|s_1-s_2|^{1 + \frac{\eta -1}{\alpha}}}{(t-s_1\vee s_2)^{\frac{1}{\alpha}}} \rho^0(t-s_1\wedge s_2,y-x) \\ &\leq   K \frac{|s_1-s_2|^\gamma}{(t-s_1\vee s_2)^{\gamma +\frac{1}{\alpha}-1 + \frac{1-\eta}{\alpha}}} \rho^0(t-s_1\wedge s_2,y-x),
\end{align*}
since $\gamma < 1 + \frac{\eta -1}{\alpha}$. Gathering the estimates obtained, we have proved that \begin{align}\label{proof_gradient_lin_der_holder_times_eq2}
	|I_1| &\leq K  \frac{|s_1-s_2|^\gamma}{(t-s_1\vee s_2)^{\gamma + \frac{1}{\alpha}-1 + \frac{1 - \eta}{\alpha}}} \left[\rho^{0} (t-s_1 \vee s_2 ,y-x) +  \rho^{0}(t-s_1 \wedge s_2 ,y-x)\right] \\\notag & \hspace{1cm}\left[1+  \sum_{k=1}^{m} C^k(t-s_1\vee s_2)^{k\left(1 + \frac{\eta-1}{\alpha}\right)} \prod_{j=1}^k\BB\left(\left[2\left(1 + \frac{\eta-1}{\alpha}\right)-\gamma\right] \wedge 1 +(j-1)\left(1+ \frac{\eta -1}{\alpha}\right), 1 - \frac{1}{\alpha}\right) \right].
\end{align}

Then, we decompose $I_2$ in the following way 

\begin{align*}
	I_2 &= \int_{s_1\vee s_2}^t \int_{\R^d} \Delta_{s_1,s_2} \left(p_{m+1}\otimes \pa_v\del \H_{m+1}\right)(\mu,\cdot,r,x,z) \Phi_{m+1}(\mu,s_1\vee s_2,r,t,z,y) \, dz \,dr \\ &\quad+ \int_{s_1\vee s_2}^t \int_{\R^d}  \left(p_{m+1} \otimes \pa_v\del \H_{m+1}\right)(\mu,s_1 \wedge s_2,r,x,z) \Delta_{s_1,s_2}\Phi_{m+1}(\mu,\cdot,r,t,z,y) \, dz \,dr\\ &\quad - \int_{s_1\wedge s_2}^{s_1\vee s_2} \int_{\R^d} \left(p_{m+1}\otimes \pa_v\del \H_{m+1}\right)(\mu,s_1\wedge s_2,r,x,z)  \Phi_{m+1}(\mu,s_1\wedge s_2,r,t,z,y) \, dz \,dr \\ &=: I_{2,1} + I_{2,2} + I_{2,3}.
\end{align*}

We follow the same lines of reasoning as for $I_1$. Using the bound \eqref{proof_gradient_lin_der_holder_times_eq2} previously obtained for $I_1$ and \eqref{Phi_m_bound}, we show that \begin{align*}
	|I_{2,1}| & \leq  K  \frac{|s_1-s_2|^\gamma}{(t-s_1\vee s_2)^{\gamma + \frac{1}{\alpha}-1 + \frac{1 -\eta}{\alpha}}} (t-s_1 \vee s_2)^{ 1 - \frac{1}{\alpha}} \left[\rho^{0} (t-s_1 \vee s_2 ,y-x) +  \rho^{0}(t-s_1 \wedge s_2 ,y-x)\right] \\ & \hspace{1cm}\left[1+  \sum_{k=1}^{m} C^k(t-s_1\vee s_2)^{k\left(1 + \frac{\eta-1}{\alpha}\right)} \prod_{j=1}^k\BB\left(\left[2\left(1 + \frac{\eta-1}{\alpha}\right)-\gamma\right] \wedge 1 +(j-1)\left(1+ \frac{\eta -1}{\alpha}\right), 1 - \frac{1}{\alpha}\right) \right].
\end{align*}

Notice that it follows from \eqref{density_bound_Picard}, \eqref{H_gradient_linear_derivative_bound_Picard} and the convolution inequality \eqref{convolineqdensityref} that \begin{equation}\label{proof_gradient_lin_deriv_holder_time_eq1} \left\vert p_{m+1}\otimes \pa_v\del \H_{m+1}(\mu,s,t,x,y)(v)\right\vert \leq K (t-s)^{\frac{\eta -1}{\alpha}+1 - \frac{1}{\alpha}} \rho^0(t-s,y-x).\end{equation}
Using this inequality and \eqref{Phi_Holder_s_Picard_lemme}, we deduce that \begin{align*}
	|I_{2,2}| &\leq K \int_{s_1\vee s_2}^t \int_{\R^d} (r-s_1\wedge s_2)^{\frac{\eta -1}{\alpha}+1 - \frac{1}{\alpha}} \rho^0(r-s_1 \wedge s_2,z-x) \frac{|s_1 - s_2|^\gamma }{(r-s_1\vee s_2)^\gamma} (t-r)^{-\frac{1}{\alpha}} \rho^1(t-r,y-z) \, dz \,dr.
\end{align*}

Reasoning as in \eqref{proof_gradient_lin_der_holder_times_disjonction} to control $(r-s_1\wedge s_2)^{\frac{\eta-1}{\alpha} + 1 - \frac{1}{\alpha}}$ since $t - s_1 \vee s_2 \leq t-s_1\wedge s_2 \leq 2 (t-s_1\vee s_2)$, we get that 

\begin{align*}
	|I_{2,2}| & \leq  K \frac{|s_1 - s_2|^\gamma }{(t-s_1\vee s_2)^{\gamma + \frac{1}{\alpha}-1 + \frac{1 - \eta}{\alpha}}}(t-s_1\vee s_2)^{1 - \frac{1}{\alpha}} \rho^0(t-s_1 \wedge s_2,y-x).
\end{align*}

For $I_{2,3}$, note that \eqref{proof_gradient_lin_deriv_holder_time_eq1}, \eqref{Phi_m_bound} and the convolution inequality \eqref{convolineqdensityref} yield 

\begin{align*}
	|I_{2,3}| & \leq \int_{s_1 \wedge s_2}^{s_1 \vee s_2} (r-s_1\vee s_2)^{\frac{\eta-1}{\alpha}+1 - \frac{1}{\alpha}}\rho^0(r-s_1\wedge s_2,z-x) (t-r)^{-\frac{1}{\alpha}}\rho^1(t-r,y-z) \, dz \,dr \\ &\leq K \frac{|s_1-s_2|^{1 + \frac{\eta -1}{\alpha}+1 - \frac{1}{\alpha}}}{(t-s_1\vee s_2)^{\frac{1}{\alpha}}} \rho^0(t-s_1\wedge s_2,y-x) \\ &\leq   K \frac{|s_1-s_2|^\gamma}{(t-s_1\vee s_2)^{\gamma +\frac{1}{\alpha}-1 + \frac{1- \eta}{\alpha}}} (t-s_1 \vee s_2)^{1 - \frac{1}{\alpha}} \rho^0(t-s_1\wedge s_2,y-x),
\end{align*}
since $|s_1 - s_2| \leq t - s_1\vee s_2$. Gathering the estimates obtained, we have proved that \begin{align}\label{proof_gradient_lin_der_holder_times_eq3}
	|I_2| &\leq K  \frac{|s_1-s_2|^\gamma}{(t-s_1\vee s_2)^{\gamma + \frac{1}{\alpha}-1 + \frac{1-\eta}{\alpha}}}(t -s_1\vee s_2)^{1 - \frac{1}{\alpha}} \left[\rho^{0} (t-s_1 \vee s_2 ,y-x) +  \rho^{0}(t-s_1 \wedge s_2 ,y-x)\right] \\ \notag& \quad\left[1+  \sum_{k=1}^{m} C^k(t-s_1\vee s_2)^{k\left(1 + \frac{\eta -1}{\alpha}\right)} \prod_{j=1}^k\BB\left(\left[2\left(1 + \frac{\eta-1}{\alpha}\right)-\gamma\right] \wedge 1 +(j-1)\left(1+ \frac{\eta -1}{\alpha}\right), 1 - \frac{1}{\alpha}\right) \right].
\end{align}

We conclude from \eqref{proof_gradient_lin_der_holder_times_eq2}, \eqref{proof_gradient_lin_der_holder_times_eq3} and the fact that $t-s_1\wedge s_2 \leq 2(t-s_1 \vee s_2)$ that 

\begin{align*}
	&\left\vert \Delta_{s_1,s_2} \pa_v\del p_{m+1}(\mu,\cdot,t,x,y)(v) \right\vert\\ &\leq K  \frac{|s_1-s_2|^\gamma}{(t-s_1\vee s_2)^{\gamma + \frac{1}{\alpha}-1+ \frac{1-\eta}{\alpha}}} \left[\rho^{0} (t-s_1 \vee s_2 ,y-x) +  \rho^{0}(t-s_1 \wedge s_2 ,y-x)\right] \\\notag & \hspace{1cm}\left[1+  \sum_{k=1}^{m} C^k(t-s_1\vee s_2)^{k\left(1 + \frac{\eta -1}{\alpha}\right)} \prod_{j=1}^k\BB\left(\left[2\left(1 + \frac{\eta-1}{\alpha}\right)-\gamma\right] \wedge 1 +(j-1)\left(1+ \frac{\eta -1}{\alpha}\right), 1 - \frac{1}{\alpha}\right)  \right] \\ &\leq    \left[\frac{|s_1-s_2|^\gamma}{(t-s_1\vee s_2)^{\gamma + \frac{1}{\alpha}-1 + \frac{1-\eta}{\alpha}}}\rho^{0} (t-s_1 \vee s_2 ,y-x) +  \frac{|s_1-s_2|^\gamma}{(t-s_1\wedge s_2)^{\gamma + \frac{1}{\alpha}-1}} \rho^{0}(t-s_1 \wedge s_2 ,y-x)\right] \\\notag & \hspace{1cm}  \sum_{k=1}^{m+1} C^k(t-s_1\vee s_2)^{(k-1)\left(1 + \frac{\eta-1}{\alpha}\right)} \prod_{j=1}^{k-1}\BB\left(\left[2\left(1 + \frac{\eta-1}{\alpha}\right)-\gamma\right] \wedge 1 +(j-1)\left(1+ \frac{\eta -1}{\alpha}\right), 1 - \frac{1}{\alpha}\right),
\end{align*}

provided that we choose $C \geq K$ in \eqref{gradient_lin_der_density_holder_time_Picard}. This ends the proof of the induction step for \eqref{gradient_lin_der_density_holder_time_Picard}.

\appendix

\section{Differential calculus for functions of a measure variable}\label{section_appendix_diff_calculus}

Let us fix $\beta \in [0,2]$. We use the following convention $\mathcal{P}_0(\R^d):= \PPP,$ endowed with the weak topology. 

\begin{Def}[Linear derivative] \label{deflinearderivative}A function $u : \PP\rightarrow \R$ is said to have a linear derivative if there exists a function $\del u \in \CC^0(\PP\times \R^d;\R)$ satisfying the following properties.
	\begin{enumerate}
		\item For all compact subset $\KK \subset \PP,$ there exists a constant $C_{\KK}>0$ such that $$ \forall \mu \in \KK,\, \forall v \in \R^d, \, \left\vert \delu (\mu)(v) \right\vert \leq C_{\KK} (1+|v|^{\beta}).$$
		\item For all $\mu,\nu \in \PP,$ we have $$ u(\mu)-u(\nu) = \int_0^1 \int_{\R^d} \del u (t \mu + (1-t)\nu)(v)\, d(\mu-\nu)(v) \, dt.$$
		
	\end{enumerate}

The function $u$ is said to have a linear derivative of order two if for all $v \in \R^d$, the map $\del u(\cdot)(v)$ admits a linear derivative $\dell(\cdot)(v,\cdot)$ such that $\dell u$ is continuous on $\PP \times \R^d \times \R^d$ and for all compact subset $\KK \subset \PP,$ there exists a constant $C_{\KK}>0$ such that $$ \forall \mu \in \KK,\, \forall v,v' \in \R^d, \, \left\vert \dell u (\mu)(v,v') \right\vert \leq C_{\KK} (1+|v|^{\beta}+|v'|^\beta).$$
\end{Def}

\begin{Def}\label{def_C1_space}
	We define the space $\CC^{1}([0,T] \times \R^d \times \PP)$ as the set of continuous functions $u : [0,T] \times \R^d \times \PP \to \R$ satisfying the following properties. \begin{enumerate} \item For any $\mu \in \PP$, the map $u(\cdot,\cdot,\mu)$ belongs to $\CC^{1}([0,T]\times \R^d)$ with $\pa_t u$ and $\pa_x u$ continuous on $[0,T] \times \R^d \times \PP$.
		\item For any $(t,x) \in [0,T] \times \R^d$, the map $u(t,x,\cdot)$ admits a linear derivative $(\mu,v) \in \PP \times \R^d \mapsto \del u(t,x,\mu)(v)$ such that $\del u$ is continuous on  $[0,T] \times \R^d \times \PP$. 
		\item For any $(t,x,\mu) \in [0,T] \times \R^d\times \PP$, the map $\del u(t,x,\mu)$ is of class $\CC^1$ on $\R^d$ and $\pa_v \del u$ is continuous on  $[0,T] \times \R^d \times \PP\times \R^d$. 
	\end{enumerate}
\end{Def}
We now introduce the notion of empirical projection.

\begin{Def}[Empirical projection]\label{def_empirical_proj}
	Fix $u : \PP \rightarrow \R$. For all $N \geq 1$, the empirical projection $u^N$ of $u$ is defined, for all $\bm{x} = (x_1, \dots,x_N) \in (\R^d)^N$, by $$ u^N(\bm{x}) = u (\muu^N_{\bm{x}}),$$
	where $\muu^N_{\x} = \frac1N \displaystyle\sum_{j=1}^N \delta_{x_j}.$
\end{Def}
The following proposition is the analogue of \cite[Proposition $5.91$]{CarmonaProbabilisticTheoryMean2018} where $\beta =2$, so we don't give the proof.

\begin{Prop}\label{projempC1}
	Let $ u : \PP\rightarrow \R$ be a function admitting a linear derivative $ \del u $ such that all $\mu \in \PP$, $\del u (\mu)(\cdot) \in \CC^1(\R^d)$ and $\partial_v \del u$ is continuous on $\PP \times \R^d$. Then, for all $N \geq 1$, the empirical projection $u^N$ of $u$ is of class $\CC^1$. Moreover for all $\bm{x}=(x_1, \dots, x_N) \in (\R^{d})^N$  $$\partial_{x_i} u^N(x_1, \dots, x_N) = \frac{1}{N}  \partial_v \del u (\overline{\mu}^N_{\bm{x}})(x_i).$$
\end{Prop}
The next proposition illustrates how a smooth flow of measures admitting a transition density can regularize a function defined on $\PP$. It is clearly reminiscent of \cite[Proposition $2.3$]{frikha2021backward}. We don't prove it since it can be done in a completely analogous manner.

\begin{Prop}[Regularization by a smooth flow of density functions]\label{Prop_lin_der_along_flow_density}Let us fix $\phi : \PP \rightarrow \R$ a function admitting a linear derivative and consider a map $(\mu,s,x) \in \PP \times [0,T) \times \R^d \mapsto p(\mu,s,T,x,y)$, where $ T>0$ is fixed and such that $p(\mu,s,T,x,\cdot)$ is a density function. We define the measure-valued map $\Theta : (s,\mu) \in [0,T) \times \PP \mapsto \Theta(s,\mu)(dy):= \left(\int_{\R^d}p(\mu,s,T,x,y) \, d\mu(x) \right)\, dy \in \PPP.$  We assume that the following properties hold true.
	
	\begin{enumerate}
		\item For any compact subset $\KK$ of $[0,T) \times \PP$, $$\int_{\R^d} \sup_{(s,\mu) \in \KK} |y|^\beta \Theta(s,\mu)(dy) < +\infty.$$
		\item For all $y \in \R^d$, the map $(\mu,s,x) \in \PP \times [0,T) \times \R^d\mapsto p(\mu,s,T,x,y)$ belongs to $\CC^{1}(\PP \times [0,T) \times \R^d)$.
		
		\item  For any compact subset $\KK$ of $\PP \times [0,T) \times \R^d \times \R^d$ and for any $j \in \{0,1\}$  \[  \int_{\R^d} \sup_{(\mu,s,x,v)\in \KK}\left\{ \left\vert \pa_v^j\del p(\mu,s,T,x,y)(v) \right\vert +  \left\vert \pa_x^j p(\mu,s,T,x,y) \right\vert  +   \left\vert\pa_s p(\mu,s,T,x,y) \right\vert  \right\}\, dy<+\infty.\]
		\item For any compact subset $\KK$ of $\PP \times [0,T)$, there exists a positive constant $C$ such that  \[ \int_{\R^d} (1+|y|^\beta)\sup_{(\mu,s)\in \KK}  \left\vert p(\mu,s,T,v,y) \right\vert \, dy \leq C(1 + |v|^\beta),\] \[  \int_{\R^d} (1+|y|^\beta)\sup_{(\mu,s)\in \KK} \left\vert \del p(\mu,s,T,x,y)(v) \right\vert \, dy \leq C(1+|x|^\beta)(1 + |v|^\beta),\] and  \[  \sup_{x\in\R^d}\int_{\R^d} \sup_{(\mu,s)\in \KK} \left\vert \del p(\mu,s,T,x,y)(v) \right\vert \, dy \leq C(1 + |v|^\beta).\]
	\end{enumerate}
	
	Then, the function $(s,\mu) \in [0,T) \times \PP \mapsto \phi(\Theta(s,\mu))$ belongs to $\CC^{1}([0,T) \times \PP).$ Moreover, we have 
	
	  \begin{equation*}
		\pa_s \left[\phi(\Theta(s,\mu))\right] = \int_{\R^{2d}} \left(\del \phi (\Theta(s,\mu))(y) -\del \phi (\Theta(s,\mu))(x) \right) \pa_s p(\mu,s,T,x,y) \, dy\, d\mu(x),
	\end{equation*}

\begin{align*}
	\del \left[\phi(\Theta(s,\mu))\right](v) &= \int_{\R^d} \del \phi(\Theta(s,\mu))(y) p(\mu,s,T,v,y) \, dy \\ &\quad + \int_{\R^{2d}} \left( \del \phi(\Theta(s,\mu))(y) - \del \phi(\Theta(s,\mu))(x) \right) \del p(\mu,s,T,x,y)(v) \, dy  \, d\mu(x),
\end{align*}

and\begin{align*}
	\pa_v\del \left[\phi(\Theta(s,\mu))\right](v) &= \int_{\R^d} \left(\del \phi(\Theta(s,\mu))(y) - \del \phi(\Theta(s,\mu))(v) \right) \pa_xp(\mu,s,T,v,y) \, dy \\ &\quad + \int_{\R^{2d}} \left( \del \phi(\Theta(s,\mu))(y) - \del \phi(\Theta(s,\mu))(x) \right) \pa_v\del p(\mu,s,T,x,y)(v) \, dy  \, d\mu(x).
\end{align*}

\end{Prop}

We now focus on Itô's formula along the flow of probability measures associated with a jump process. Let us fix $Z^1=(Z^1_t)_t$ and $Z^2 = (Z^2_t)_t$ two $\alpha$-stable processes on $\R^d$ with $\alpha \in (1,2).$ Their associated Poisson random measures are respectively denoted by $\NN^1$ and $\NN^2$, their compensated Poisson random measures by $\NNN^1$ and $\NNN^2$ and their Lévy measures by $\nu^1$ and $\nu^2$. Since $\alpha \in (1,2)$, we can write for all $t \geq 0$ \[ Z^1_t = \int_0^t\int_{\R^d} z \, \NNN^1(ds,dz) \quad \text{and} \quad Z^2_t = \int_0^t\int_{\R^d} z \, \NNN^2(ds,dz). \] 

We fix $\beta \in (1,\alpha)$ and $\gamma \in (0,1]$ such that $\gamma > \alpha -1$. We consider two jump processes $X=(X_t)_{t \in [0,T]}$ and $Y=(Y_t)_{t \in [0,T]}$ defined for all $t \in [0,T]$ by \begin{equation}\label{process}
	X_t := X_0 + \int_0^t b_s \, ds + Z^1_t, \quad \text{and} \quad  Y_t = Y_0 + \int_0^t \eta_s \, ds + Z^2_t,\end{equation} where $X_0,Y_0 \in L^\beta(\Omega,\FF_0)$, $b,\eta: [0,T] \times \Omega \rightarrow \R^d$ are bounded predictable processes. The distribution of $X_t$ is denoted by $\mu_t.$\\

We state in the next proposition Itô's formula, which is deduced from \cite[Theorem $2$]{Cavallazzi_Ito_jump} for the specific type of processes that are considered in the present work.

\begin{Prop}[Itô's formula]\label{ito_formula}

	Let $u:[0,T]\times \R^d \times \PP \rightarrow \R$ be a continuous function satisfying the following properties. \begin{enumerate}
		\item  The function $u$ belongs to $\CC^1([0,T]\times \R^d \times \PP)$ and $\partial_x u(t,\cdot,\mu)$ is $\gamma$-Hölder continuous uniformly with respect to $t$ and $\mu$.
		\item For all compact $\KK \subset \R^d \times\PP,$ there exists $C_{\KK} >0$ such that  $$ \forall t \in [0,T], \, \forall (x,\mu) \in \KK,\, \forall v\in \R^d, \, \left\vert \del u (t,x,\mu)(v) \right\vert \leq C_{\KK}(1+|v|^\beta).$$
		
		\item If $\gamma>0,$ for any compact $\KK \subset \R^d \times\PP,$ there exists $C_{\KK}>0$ such that $$  \forall t\in [0,T],\, \forall (x,\mu) \in \KK, \, \forall v,v' \in \R^d,\, \displaystyle\left\vert \partial_v\del u(t,x,\mu)(v) - \partial_v\del u(t,x,\mu)(v')  \right\vert \leq C_{\KK} |v-v'|^{\gamma}.$$

		\item For any compact $\KK \subset \R^d \times \PP,$ we have  $$ \sup_{t \in[0,T]}\sup_{(x,\mu) \in \KK} \int_{\R^d} \left\vert \partial_v \del u (t,x,\mu)(v) \right\vert^{\beta'} \, d\mu(v) < + \infty.$$
	\end{enumerate}

	Then, the function $(t,x) \in [0,T]\times \R^d \mapsto u(t,x,\mu_t)$ is of class $\CC^{1}$, with $\partial_x u(t,\cdot,\mu_t)$ $\gamma$-Hölder continuous uniformly with respect to $t$. Moreover, we have almost surely for all $t \in [0,T]$ \begin{align}\label{itoformula2}
		\notag&u(t,Y_t,\mu_t) - u(0,Y_0,\mu_0) \\ \notag&=  \int_0^t\partial_t u(s,Y_s,\mu_s) \, ds  + \int_0^t \overline{\E} \left( \partial_v \del u (s,Y_s,\mu_s)(\overline{X}_s) \cdot \overline{b}_s \right)\, ds  \\\notag &\quad + \int_0^t \int_{\R^d} \overline{\E} \left[ \del u(s,Y_s,\mu_s)(\overline{X}_{s^-}+ z) - \del u(s,Y_s,\mu_s)(\overline{X}_{s^-}) \right.\\ &\left. \hspace{7cm}- \partial_v \del u(s,Y_s,\mu_s)(\overline{X}_{s^-})\cdot z\right] \, d\nu^1(z) \,ds \\ &\notag \quad+\int_0^t \partial_x u(s,Y_s,\mu_s)\cdot \eta_s\,ds + \int_0^t \int_{\R^d} u(s,Y_{s^-} + z,\mu_s) - u(s,Y_{s^-},\mu_s) \, \NNN^2(ds,dz) \\  &\notag \quad+ \int_0^t \int_{\R^d} \left[ u(s,Y_{s^-} + z,\mu_s) - u(s,Y_{s^-},\mu_s) - \partial_x u(s,Y_{s^-},\mu_s)\cdot z \right]\, d\nu^2(z) \, ds
	\end{align}
	where $(\overline{\Omega},\overline{\FF},\overline{\P})$ is an independent copy of $(\Omega,\FF,\P)$ and $(\overline{b},\overline{X})$ is a copy of $(b,X)$.
	
\end{Prop}

\section{Parametrix expansion for stable-driven SDEs}\label{section_appendix_parametrix}

Let us fix $Z=(Z_t)_t$ a rotationally invariant $\alpha$-stable process on $\R^d$ with $\alpha \in (1,2).$ Its associated Poisson random measure is denoted by $\NN$, the compensated Poisson random measure by $\NNN$. Since $\alpha \in (1,2)$, we can write for all $t \geq 0$ \[ Z_t = \int_0^t\int_{\R^d} z \, \NNN(ds,dz). \] The Lévy measure $\nu$ of $Z$ is given by \[ d\nu(z) := \frac{d z}{|z|^{d+\alpha}}.\] We consider a function $b:[0,T]\times \R^d \to \R^d$ satisfying the following properties.

\begin{enumerate}
	\item The function $b$ is jointly continuous and globally bounded on $[0,T]\times \R^d$.
	\item The function $b$ is $\eta$-Hölder continuous on $\R^d$ uniformly in time, with $\eta \in ( 0, 1]$, i.e.\ there exists $C>0$ such that for all $t\in [0,T]$ and $x_1,x_2 \in \R^d$ \[|b(t,x_1)-b(t,x_2)| \leq C |x_1 - x_2|^\eta.\]
\end{enumerate}

We fix $s\in [0,T)$ and we consider the following stable-driven SDE  \begin{equation}\label{SDEparametrix}
\begin{cases}
dX^{s,x}_t= b(t,X_t^{s,x})\, dt +dZ_t,\quad t \in [s,T], \\ X_s^{s,x}=x \in \R^d.
\end{cases}
\end{equation}

It is well-posed in the weak sense by \cite{Mikulevicius_MP_stable}. The density of $Z_t$ is denoted by $q(t,\cdot)$. We denote by $\Delta^{\frac{\alpha}{2}}$ the fractional Laplacian associated with $Z$ defined for all $f \in \CC^{1+\gamma}_b(\R^d;\R)$, with $\gamma > \alpha -1$ (i.e.\ $f$ belongs to $\CC^1_b(\R^d;\R)$ and $\nabla f$ is $\gamma$-Hölder continuous) and for all $x \in \R^d$ by \begin{equation}\label{stableoperator}
	\Delta^{\frac{\alpha}{2}} f(x):= \int_{\R^d} (f(x+z) - f(x) - \nabla f (x) \cdot z )\, d\nu(z).
\end{equation}

We define for all $s\in [0,T)$, $t \in (s,T]$ and $x,y \in \R^d$ \begin{align}
& \p(s,t,x,y) := q(t-s,y-x), \\ &\notag \H(s,t,x,y) := b(s,x) \cdot\pa_x \p(s,t,x,y).
\end{align}
Note that the proxy $\p(s,t,x,\cdot)$ is the density at time $t >s$ of the solution to \begin{equation}\label{SDEproxy}
\begin{cases}
d\widehat{X}^{s,x}_t= dZ_t, \quad  t \in [s,T], \\ \widehat{X}^{s,x}_s=x \in \R^d,
\end{cases}
\end{equation}
\\
and $\H$ is the associated parametrix kernel. We also define the space-time convolution operator between to functions $f$ and $g$ by \begin{equation}\label{defconvolop} f \otimes g (s,t,x,y) := \int_s^t \int_{\R^d} f(s,r,x,z) g(r,t,z,y) \, dz \, dr,
\end{equation}
when it is well-defined. The space-time convolution iterates $\H^k$ of $\H$ are defined recursively by $\H^1 := \H$ and $\H^{k+1} := \H \otimes \H^{k}.$ By convention $f \otimes \H^0$ is equal to $f$. Finally, we denote by $\Phi$ the solution to the following Volterra integral equation \[ \Phi(s,t,x,y) = \H(s,t,x,y) + \H \otimes \Phi(s,t,x,y),\] which is given by the uniform convergent series \begin{equation}\label{defvolterra}
\Phi(s,t,x,y) = \sum_{k=1}^{\infty} \H^k(s,t,x,y).
\end{equation}

Let us also define, for $k > - \alpha $ the function $\rho^k$ by \begin{equation}\label{defdensityref}
	\forall t>0, \, x \in \R^d,\, \rho^k(t,x):= t^{-\frac{d}{\alpha}}(1+t^{-\frac{1}{\alpha}}|x|)^{-d-\alpha -k}.
\end{equation}

\begin{Thm}\label{Thmdensityparametrix} For any $s \in [0,T)$, $s<t \leq T$ and $x \in \R^d$, the distribution of $X^{s,x}_t$ has a density with respect to the Lebesgue measure denoted by $p(s,t,x,\cdot)$ and given by the absolutely convergent parametrix series \begin{align}\label{representationdensityparametrix}
		\notag p(s,t,x,y) &= \p(s,t,x,y) + \sum_{k=1}^{\infty} \p \otimes \H^k(s,t,x,y)\\ &= \p(s,t,x,y) + \p \otimes \Phi(s,t,x,y).
\end{align}

For any $t \in(0,T]$ and $y \in \R^d$, $p(\cdot,t,\cdot,y)$ is of class $\CC^1$ on $[0,t)\times \R^d$ and $p(\cdot,t,\cdot,y)$, $\pa_s p(\cdot,t,\cdot,y)$ and $\pa_x p(\cdot,t,\cdot,y)$ are continuous on $[0,t) \times \R^d$. The function $p(\cdot,t,\cdot,y)$ is solution to the following backward Kolmogorov PDE \begin{equation}\label{EDPKolmogorovparametrix}
	\begin{cases}
		\pa_s p(s,t,x,y) + b(s,x)\cdot \pa_x p(s,t,x,y) + \Delta^{\frac{\alpha}{2}} p(s,t,\cdot,y)(x) =0, \quad \forall (s,x) \in [0,t) \times \R^d,\\ 
		p(s,t,x,\cdot) \underset{s\to t^-}{\longrightarrow} \delta_x,
	\end{cases}
\end{equation}where $p(s,t,x,\cdot) \underset{s\to t^-}{\longrightarrow} \delta_x$ means that for all function $f : \R^d \rightarrow \R$ bounded and uniformly continuous, one has \[ \sup_{x\in \R^d}\left\vert \int_{\R^d} f(y) p(s,t,x,y) \,dy - f(x)\right\vert \underset{s \rightarrow t^-}{\to} 0.\]Moreover, $p$ satisfies the following estimates. \begin{itemize}
\item There exists $C>0$ such that for all $j \in \{0,1\}$, $0 \leq s < t \leq T$ and $x,y \in \R^d$

\begin{equation}\label{density_bound}
	|\pa_x^j p(s,t,x,y) |
\leq C (t-s)^{-\frac{j}{\alpha}} \rho^{j}(t-s,y-x).
\end{equation}

\begin{align}\label{timederivativedensity_bound}
|\Delta^{\frac{\alpha}{2}} p(s,t,\cdot,y)(x) | \leq C (t-s)^{-1} \rho^{0}(t-s,y-x).
\end{align}
\item For all $j \in \{ 0,1\}$ and $\gamma \in (0, 1]$ with $\gamma \in (0, (2\alpha -2) \wedge (\eta + \alpha -1))$ if $n=1$, there exists $C>0$ such that for all $0 \leq s < t \leq T$ and $x_1,x_2,y \in \R^d$

\begin{equation}\label{gradientdensity_holder}
	|\pa_x^j p(s,t,x_1,y) -  \pa^j_x p(s,t,x_2,y)| \leq C (t-s)^{-\frac{ j+\gamma}{\alpha}} |x_1-x_2|^\gamma \left[\rho^{j}(t-s,y-x_1) +\rho^{j}(t-s,y-x_2)\right].
\end{equation}

\end{itemize}

\end{Thm}

Before proving Theorem \ref{Thmdensityparametrix}, we recall some properties satisfied by the functions $\rho^k$.

\begin{Lemme}\label{Lemmedensityreferencecontrol}
	
	The functions $\rho^k$ satisfy the following properties.
	
	\begin{itemize}
		
		\item For all $k >- \alpha$ and $\gamma \in [0,1]$ with $k - \gamma >- \alpha$, we have for all $t >0$ and $x \in \R^d$ \begin{equation}\label{scalingdensityref}
			|x|^\gamma t^{-\frac{\gamma}{\alpha}} \rho^k(t,x) \leq \rho^{k-\gamma} (t,x).
		\end{equation}
		\item Let us fix $- \alpha < k_1 \leq k_2$. Then, for all function $y : (0,+ \infty) \to \R^d$ such that $t \in (0,+ \infty) \mapsto t^{-\frac{1}{\alpha}} y(t)$ is bounded, there exists $C>0$ such that for all $t>0$ and $x \in \R^d$ \begin{equation}\label{controldensityref}
			\rho^{k_2}(t,x+y(t)) \leq C \rho^{k_1}(t,x).
		\end{equation}
		
		\item For all $k > - \alpha$ and $R>0$, there exists $C$ such that for all $t >0$, $y \in \R^d$ and $x \in \R^d$ with $|x| \leq R$ \begin{equation}\label{uniformcontroldensityref}
		\rho^k(t,y+x) \leq (1+ct^{-\frac{1}{\alpha}}R)^{d+\alpha+k} \rho^k(t,y).
		\end{equation}
	
	\item For all $k_1,k_2 > - \alpha$, there exists $C>0$ such that for all $s\geq 0$, $t >s$ and $y \in \R^d$

\begin{equation}\label{convolineqdensityref}
\int_{\R^d}\rho^{k_1}(t-s,y-z)  \rho^{k_2}(s,z) \, dz \leq C\rho^{k_1 \wedge k_2}(t,y).
\end{equation} 

\end{itemize}
\end{Lemme}

The following lemma gathers the properties that we need on the proxy $\p$.

\begin{Lemme}\label{Lemmegradientestimatestable}
	
	For all $t \in (0,T]$, $y \in \R^d$, $\p(\cdot,t,\cdot,y)$ is of class $\CC^{1,\infty}$ on $[0,t)\times\R^d$. Moreover, it satisfies the following gradient estimates.
	
	\begin{itemize}
		\item For all $j \in \N$, there exists $C>0$ such that for all $t\in (0,T]$,  $s \in [0,t)$, $x,y \in \R^d$, we have \begin{equation}\label{gradientestimatestable}
		|\pa_x^j \p(s,t,x,y )| \leq C (t-s)^{-\frac{j}{\alpha}} \rho^j(t-s,y-x).
		\end{equation}
	
	\item  There exists a constant $C>0$ such that for all $j \in \{0,1\}$, $t\in (0,T]$, $s \in [0,t)$, $x,y \in \R^d$ \begin{equation}\label{eqproofproxy7}
		|\pa_s \pa^j_x \p(s,t,x,y)| \leq C(t-s)^{-1-\frac{j}{\alpha}} \rho^{j}(t-s,y-x).
	\end{equation}
		
		\item For all $j \in \N$, there exists $C>0$ such that for all $\gamma \in (0,1]$, $t\in (0,T]$,  $s \in [0,t)$, $ x_1,x_2,y \in \R^d$,  we have \begin{equation}\label{gradientHolderstable}
		|\pa_x^j \p(s,t,x_1,y) - \pa_x^j \p(s,t,x_2,y)| \leq C (t-s)^{-\frac{j+\gamma}{\alpha}} |x_1 - x_2|^\gamma \left[ \rho^{j}(t-s,y-x_1) + \rho^{j}(t-s,y-x_2)\right].
		\end{equation}
	
		\item For all $j \in \{0,1\}$, there exists $C>0$ such that for all $\gamma \in (0,1]$, $t\in (0,T]$, $s_1,s_2 \in [0,t)$, $ x,y\in \R^d$,  we have \begin{multline}\label{gradientHoldertimestable}
		|\pa_x^j \p(s_1,t,x,y) - \pa_x^j \p (s_2,t,x,y)| \\ \leq C \left[ \frac{|s_1 - s_2|^\gamma}{(t-s_1\wedge s_2)^{\gamma + \frac{j}{\alpha}}} \rho^{j}(t-s_1 \wedge s_2, y-x) + \frac{|s_1 - s_2|^\gamma}{(t-s_1\vee s_2)^{\gamma + \frac{j}{\alpha}}} \rho^{j}(t-s_1 \vee s_2, y-x) \right].
	\end{multline}

	\end{itemize}
	
\end{Lemme}

\begin{proof}
Notice that \eqref{gradientestimatestable} and \eqref{gradientHolderstable} are quite standard since $\p(s,t,x,y) = q(t-s,y-x),$ where $q(t-s,\cdot)$ is the density of the $\alpha$-stable random variable $Z_{t-s}$ (see \cite[Lemma $2.8$]{menozzi:hal-03088376}). We prove \eqref{eqproofproxy7}. To do this, we remark that the self-similarity of the $\alpha$-stable process $Z$ implies that \[\p(s,t,x,y) = q(t-s,y-x) = (t-s)^{-\frac{d}{\alpha}}q\left(1,\frac{y-x}{(t-s)^{\frac{1}{\alpha}}}\right) = (t-s)^{-\frac{d}{\alpha}} \p(0,1,(t-s)^{-\frac{1}{\alpha}}x,(t-s)^{-\frac{1}{\alpha}}y).\] 
This yields 
\begin{align}\label{time_der_stable_density}\pa_s \p(s,t,x,y) &= \frac{d}{\alpha}(t-s)^{-\frac{d}{\alpha}-1} \p(0,1,(t-s)^{-\frac{1}{\alpha}}x,(t-s)^{-\frac{1}{\alpha}}y)\\ \notag&\quad + (t-s)^{-\frac{d}{\alpha}} \pa_x\p(0,1,(t-s)^{-\frac{1}{\alpha}}x,(t-s)^{-\frac{1}{\alpha}}y)\cdot \left(\frac{1}{\alpha}(t-s)^{-\frac{1}{\alpha}-1}(x-y)\right).
\end{align}

Using \eqref{gradientestimatestable}, we obtain that 

\begin{align*}
	|\pa_s\p(,s,t,xy)| &\leq C(t-s)^{-1} \rho^0(t-s,y-x)+ C(t-s)^{-1}((t-s)^{-\frac{1}{\alpha}}|x-y|) \rho^1(1,(t-s)^{-\frac{1}{\alpha}}(y-x)).
	\end{align*}

The space-time inequality \eqref{scalingdensityref} finally yields \begin{align*}
	|\pa_s\p(,s,t,xy)| &\leq C(t-s)^{-1} \rho^0(t-s,y-x).
\end{align*}

By differentiating \eqref{time_der_stable_density} with respect to $x$, one has 
\begin{align}\pa_x \pa_s \p(s,t,x,y) &= \frac{d}{\alpha}(t-s)^{-\frac{1}{\alpha}} (t-s)^{-\frac{d}{\alpha}-1} \pa_x\p(0,1,(t-s)^{-\frac{1}{\alpha}}x,(t-s)^{-\frac{1}{\alpha}}y)\\ \notag&\quad + (t-s)^{-\frac{1}{\alpha}} (t-s)^{-\frac{d}{\alpha}} \pa^2_x\p(0,1,(t-s)^{-\frac{1}{\alpha}}x,(t-s)^{-\frac{1}{\alpha}}y) \left(\frac{1}{\alpha}(t-s)^{-\frac{1}{\alpha}-1}(x-y)\right) \\ \notag&\quad +\frac{1}{\alpha} (t-s)^{-\frac{1}{\alpha}-1} (t-s)^{-\frac{d}{\alpha}} \pa_x\p(0,1,(t-s)^{-\frac{1}{\alpha}}x,(t-s)^{-\frac{1}{\alpha}}y).
\end{align}

As previously, it follows from \eqref{gradientestimatestable} and the space-time inequality \eqref{scalingdensityref} that 

\begin{equation*}
	|\pa_s \pa_x \p(s,t,x,y)| \leq C(t-s)^{-1-\frac{j}{\alpha}} \rho^{1}(t-s,y-x).
\end{equation*}

We now use \eqref{eqproofproxy7} to prove \eqref{gradientHoldertimestable}. We fix $j \in \{0,1\}$, $\gamma \in (0,1]$ and we start with the case $|s_1 - s_2| > t - s_1 \vee s_2.$ In this case, using \eqref{gradientestimatestable}, we deduce that for some constant $C>0$, one has \begin{align*}
&|\pa_x^j \p (s_1,t,x,y) - \pa_x^j \p (s_2,t,x,y)| \\ &\leq C \left[(t-s_1\vee s_2)^{-\frac{j}{\alpha}}\rho^j(t-s_1\vee s_2,y-x) + (t-s_1\wedge s_2)^{-\frac{j}{\alpha}}\rho^j(t-s_1\wedge s_2,y-x)\right] \\ &\leq C \left[ \frac{|s_1-s_2|^\gamma}{|t-s_1\vee s_2|^{\gamma + \frac{j}{\alpha}}}\rho^j(t-s_1\vee s_2,y-x) + \frac{|s_1-s_2|^\gamma + |t-s_1 \vee s_2|^\gamma}{|t-s_1\wedge s_2|^{\gamma + \frac{j}{\alpha}}}\rho^j(t-s_1\wedge s_2,y-x)  \right] \\ &\leq C \left[ \frac{|s_1-s_2|^\gamma}{|t-s_1\vee s_2|^{\gamma + \frac{j}{\alpha}}}\rho^j(t-s_1\vee s_2,y-x) + \frac{|s_1-s_2|^\gamma }{|t-s_1\wedge s_2|^{\gamma + \frac{j}{\alpha}}}\rho^j(t-s_1\wedge s_2,y-x)  \right].
\end{align*}

We now focus on the case $|s_1-s_2| \leq t-s_1\vee s_2.$ For $\lambda \in [0,1]$, we set $ s_\lambda := \lambda s_1 + (1-\lambda)s_2$. We can thus write, thanks to \eqref{eqproofproxy7}, \begin{align*}
	&|\pa_x^j \p(s_1,t,x,y) - \pa_x^j \p (s_2,t,x,y)| \\ &\leq \int_0^1 |\pa_s \pa_x^j \p(s_\lambda,t,x,y)|\,|s_1-s_2| \, d\lambda \\ &\leq C|s_1-s_2| \int_0^1 (t-s_\lambda)^{-1-\frac{j}{\alpha}} \rho^{j}(t-s_\lambda,y-x) \, d\lambda \\ &\leq C|s_1-s_2| \int_0^1 (t-s_\lambda)^{-1-\frac{j+d}{\alpha}} (1+(t-s_\lambda)^{-\frac{1}{\alpha}} |y-x|)^{-d-\alpha -j}\, d\lambda \\ &\leq C |s_1-s_2|^\gamma (t-s_1\vee s_2)^{1-\gamma}(t-s_1\vee s_2)^{-1-\frac{j+d}{\alpha}} \left[(1+(t-s_1\vee s_2)^{-\frac{1}{\alpha}} |y-x|)^{-d-\alpha -j}\right. \\ &\left. \hspace{9cm} + (1+(t-s_1\wedge s_2)^{-\frac{1}{\alpha}} |y-x|)^{-d-\alpha -j}\right].
\end{align*}

Since $|s_1-s_2| \leq t-s_1\vee s_2$, we easily check that $(t-s_1\vee s_2)^{-1} \leq 2 (t-s_1\wedge s_2)^{-1}$. It follows that 

\begin{align*}
	&|\pa_x^j \p(s_1,t,x,y) - \pa_x^j \p (s_2,t,x,y)| \\ &\leq C |s_1-s_2|^\gamma\left[ (t-s_1\vee s_2)^{-\gamma-\frac{j+d}{\alpha}} (1+(t-s_1\vee s_2)^{-\frac{1}{\alpha}} |y-x|)^{-d-\alpha -j}\right. \\ &\left. \hspace{4cm} +  (t-s_1\wedge s_2)^{-\gamma-\frac{j+d}{\alpha}}(1+(t-s_1\wedge s_2)^{-\frac{1}{\alpha}} |y-x|)^{-d-\alpha -j}\right] \\ &  \leq C \left[ \frac{|s_1 - s_2|^\gamma}{(t-s_1\wedge s_2)^{\gamma + \frac{j}{\alpha}}} \rho^{j}(t-s_1 \wedge s_2, y-x) + \frac{|s_1 - s_2|^\gamma}{(t-s_1\vee s_2)^{\gamma + \frac{j}{\alpha}}} \rho^{j}(t-s_1 \vee s_2, y-x) \right].
\end{align*}
This concludes the proof.
\end{proof}

Recall that the Beta function $\BB$ is defined, for all $x,y >0$, by \[ \BB(x,y) := \int_0^1 (1-t)^{-1+x} t^{-1+y}  \, dt = \frac{\Gamma(x) \Gamma(y)}{\Gamma(x+y)},\] where $\Gamma$ is the Gamma function.\\

The next proposition gathers the controls that we need on the parametrix kernel $\H$ and the solution de the Volterra integral equation $\Phi$.

\begin{Prop}\label{Prop_H^k}
	
	The following estimates hold true. \begin{itemize}
		\item There exists $C>0$ such that for all $k \geq 1$, $0 \leq s < t\leq T$ and $x,y \in \R^d$
	
	\begin{equation}\label{H^k_bound}
		|\H^k(s,t,x,y)| \leq C^k (t-s)^{-\frac{1}{\alpha} + (k-1)\left( 1 - \frac{1}{\alpha} \right)} \prod_{j=1}^{k-1} \BB \left( j\left( 1 - \frac{1}{\alpha} \right),1 - \frac{1}{\alpha} \right) \rho^1(t-s,y-x).
	\end{equation}

\item For $\gamma \in (0, \eta]$ such that $\gamma < \alpha -1$, there exists $C>0$ depending on $\gamma$ such that for all $ k \geq 1$, $0 \leq s < t \leq T$ and $x_1,x_2,y \in \R^d$

\begin{multline}\label{H^k_Holder}
	| \H^k(s,t,x_1,y) - \H^k(s,t,x_2,y) | \leq C^k(t-s)^{-\frac{\gamma +1}{\alpha} + (k-1)\left( 1 - \frac{1}{\alpha} \right)} |x_1-x_2|^\gamma \prod_{j=1}^{k-1} \BB \left( -\frac{\gamma}{\alpha} + j\left( 1 - \frac{1}{\alpha} \right),1 - \frac{1}{\alpha} \right)\\  \left[\rho^{1 }(t-s,y-x_1) + \rho^{1 }(t-s,y-x_2) \right].
\end{multline}

\item There exists $C>0$ such that for all $k\geq 1$, $0 \leq s <t \leq T$ and $x,y \in \R^d$

\begin{equation}\label{potimesH^k_bound}
	|\p \otimes \H^k (s,t,x,y) | \leq C^{k+1} (t-s)^{k\left( 1 - \frac{1}{\alpha} \right)}\prod_{j=1}^{k} \BB \left( \frac{1}{\alpha} + j\left( 1 - \frac{1}{\alpha} \right),1 - \frac{1}{\alpha} \right) \rho^0(t-s,y-x).
\end{equation}

\item The series \eqref{defvolterra} defining $\Phi$ is absolutely convergent and there exists $C>0$ such that for all $0 \leq s <t\leq T$ and $x,y \in \R^d$

\begin{equation}\label{Phi_bound}
	|\Phi(s,t,x,y)| \leq C (t-s)^{-\frac{1}{\alpha}} \rho^1(t-s,y-x).
\end{equation}

\item For $\gamma \in (0, \eta]$ such that $\gamma < \alpha -1$, there exists $C>0$ depending on $\gamma$ such that for all $0\leq s < t\leq T$ and $x_1,x_2,y \in \R^d$

\begin{equation}\label{Phi_Holder}
	| \Phi(s,t,x_1,y) - \Phi(s,t,x_2,y) | \leq C(t-s)^{-\frac{\gamma +1}{\alpha}} |x_1-x_2|^\gamma \left[\rho^{1 }(t-s,y-x_1) + \rho^{1 }(t-s,y-x_2) \right].
\end{equation}

\end{itemize}

\end{Prop}

\begin{proof}
\noindent\textbf{Proof of \eqref{H^k_bound}.} We reason by induction on $k$. The base case $k=1$ is clear since $b$ is bounded and by Lemma \ref{Lemmegradientestimatestable}. We assume now that \eqref{H^k_bound} holds for $\H^k$ and we want to prove it for $\H^{k+1}$. We have thanks to Lemma \ref{Lemmedensityreferencecontrol}. \begin{align*}
	|\H^{k+1}(s,t,x,y)| &= \left\vert\int_s^t \int_{\R^d} \H(s,r,x,z) \H^k(r,t,z,y)  \, dz \, dr\right\vert\\ & \leq \int_s^t \int_{\R^d} C (r-s)^{-\frac{1}{\alpha}} \rho^1(r-s,z-x) C^k (t-r)^{-\frac{1}{\alpha} + (k-1)\left( 1 - \frac{1}{\alpha} \right)}\\ & \hspace{3cm} \prod_{j=1}^{k-1} \BB \left( j\left( 1 - \frac{1}{\alpha} \right),1 - \frac{1}{\alpha} \right) \rho^1(t-r,y-z) \, dz \, dr \\ &\leq C^{k+1} \left(\int_s^t (r-s)^{-\frac{1}{\alpha}} (t-r)^{-\frac{1}{\alpha} + (k-1)\left( 1 - \frac{1}{\alpha} \right)} \, dr \right) \prod_{j=1}^{k-1} \BB \left( j\left( 1 - \frac{1}{\alpha} \right),1 - \frac{1}{\alpha} \right) \rho^1(t-s,y-x).
\end{align*}

Changing variables in $r=s + \lambda (t-s)$ yields \begin{align*}	|\H^{k+1}(s,t,x,y)|& \leq C^{k+1}(t-s)^{-\frac{1}{\alpha} + k \left( 1 - \frac{1}{\alpha} \right)}\int_0^1 \lambda^{-\frac{1}{\alpha}} (1-\lambda)^{-\frac{1}{\alpha} + (k-1)\left( 1 - \frac{1}{\alpha} \right)} \, d\lambda \\ &\hspace{4cm} \prod_{j=1}^{k-1} \BB \left( j\left( 1 - \frac{1}{\alpha} \right),1 - \frac{1}{\alpha} \right) \rho^1(t-s,y-x)  \\ &\leq C^{k+1}  (t-s)^{-\frac{1}{\alpha} + k \left( 1 - \frac{1}{\alpha} \right)} \BB \left( k\left( 1 - \frac{1}{\alpha} \right),1 - \frac{1}{\alpha} \right)\prod_{j=1}^{k-1} \BB \left( j\left( 1 - \frac{1}{\alpha} \right),1 - \frac{1}{\alpha} \right)\rho^1(t-s,y-x).\end{align*}

\noindent\textbf{Proof of \eqref{H^k_Holder}.} We start with the case $k=1$.  We write \begin{align*} |\H(s,t,x_1,y) - \H(s,t,x_2,y)| &  = | b(s,x_1) \cdot \pa_x \p(s,t,x_1,y) -  b(s,x_2) \cdot \pa_x \p(s,t,x_2,y) | \\ &\leq |b(s,x_1)|\, |\pa_x \p(s,t,x_1,y) - \pa_x \p(s,t,x_2,y)| \\ &\quad + |\pa_x \p(s,t,x_2,y)|\, |b(s,x_1) - b(s,x_2)|\\ &=: I_1 + I_2. \end{align*}

Using that $b$ is bounded and Lemma \ref{Lemmegradientestimatestable}, we deduce that \[ I_1 \leq C (t-s)^{-\frac{1+\gamma}{\alpha}} |x_1 - x_2|^\gamma \left[ \rho^{1}(t-s,y-x_1) + \rho^{1}(t-s,y-x_2)\right]. \]
 Since $b$ is uniformly $\eta$-Hölder continuous and bounded, and thus uniformly $\gamma$-Hölder continuous because $\gamma \leq \eta$, one has by Lemma \ref{Lemmegradientestimatestable} \[ I_2  \leq C (t-s)^{-\frac{1}{\alpha}} |x_1-x_2|^\gamma\rho^1(t-s,y-x_2).\] We now prove that \eqref{H^k_Holder} holds for $\H^{k+1}$, with $ k \geq 1$. Using the case $k=1$, \eqref{H^k_bound} and Lemma \ref{Lemmedensityreferencecontrol}, we obtain that  \begin{align*}
	&|\H^{k+1}(s,t,x_1,y) - \H^{k+1}(s,t,x_2,y)| \\&= \left\vert\int_s^t \int_{\R^d} (\H(s,r,x_1,z) - \H(s,r,x_2,z)) \H^k(r,t,z,y) \,  dz\, dr \right\vert \\ & \leq \int_s^t \int_{\R^d} C (r-s)^{-\frac{1+\gamma}{\alpha}}|x_1-x_2|^\gamma \left[\rho^{1}(r-s,z-x_1) +\rho^{1}(r-s,z-x_2) \right] C^k (t-r)^{-\frac{1}{\alpha} + (k-1)\left( 1 - \frac{1}{\alpha} \right)}\\ & \hspace{3cm} \prod_{j=1}^{k-1} \BB \left( j\left( 1 - \frac{1}{\alpha} \right),1 - \frac{1}{\alpha} \right) \rho^1(t-r,y-z) \, dz \, dr \\ &\leq C^{k+1} (t-s)^{-\frac{1+\gamma}{\alpha} + k\left( 1 - \frac{1}{\alpha} \right)} \BB \left(k \left(1 - \frac{1}{\alpha}\right), 1 - \frac{1 + \gamma}{\alpha}\right) \prod_{j=1}^{k} \BB \left( j\left( 1 - \frac{1}{\alpha} \right),1 - \frac{1}{\alpha} \right) \\ & \hspace{8cm}\left[\rho^1(t-s,y-x_1)  + \rho^1(t-s,y-x_2) \right].
\end{align*}

We conclude noting that 

\begin{align*}
\BB \left(k \left(1 - \frac{1}{\alpha}\right), 1 - \frac{1 + \gamma}{\alpha}\right) \prod_{j=1}^{k} \BB \left( j\left( 1 - \frac{1}{\alpha} \right),1 - \frac{1}{\alpha} \right)  &= \frac{\Gamma \left(k \left(1 - \frac{1}{\alpha}\right)\right) \Gamma \left( 1 - \frac{1+\gamma}{\alpha}\right)}{\Gamma \left( 1 - \frac{1+\gamma}{\alpha} +k \left(1 - \frac{1}{\alpha}\right)\right)} \prod_{j=1}^{k-1} \frac{\Gamma \left(j \left(1 - \frac{1}{\alpha}\right)\right) \Gamma \left( 1 - \frac{1}{\alpha}\right)}{\Gamma \left((j+1) \left(1 - \frac{1}{\alpha}\right)\right)} \\ &= \frac{\Gamma \left(k \left(1 - \frac{1}{\alpha}\right)\right) \Gamma \left( 1 - \frac{1+\gamma}{\alpha}\right)}{\Gamma \left( 1 - \frac{1+\gamma}{\alpha} +k \left(1 - \frac{1}{\alpha}\right)\right)} \frac{ \left[\Gamma \left( 1 - \frac{1}{\alpha}\right)\right]^k}{\Gamma \left(k \left(1 - \frac{1}{\alpha}\right)\right)} \\ &= \frac{ \Gamma \left( 1 - \frac{1+\gamma}{\alpha}\right) \left[\Gamma \left( 1 - \frac{1}{\alpha}\right)\right]^k}{\Gamma \left( 1 - \frac{1+\gamma}{\alpha} +k \left(1 - \frac{1}{\alpha}\right)\right)} \\ & = \prod_{j=1}^k \frac{\Gamma \left( -\frac{\gamma}{\alpha } +j \left(1 - \frac{1}{\alpha}\right)\right) \Gamma \left( 1 - \frac{1}{\alpha}\right)}{\Gamma \left( -\frac{\gamma}{\alpha } +(j+1) \left(1 - \frac{1}{\alpha}\right)\right) } \\ &=\prod_{j=1}^{k} \BB \left( -\frac{\gamma}{\alpha} + j\left( 1 - \frac{1}{\alpha} \right),1 - \frac{1}{\alpha} \right).
\end{align*}

\noindent\textbf{Proof of \eqref{potimesH^k_bound}.} It follows from the control of $\p$ given by Lemma \ref{Lemmegradientestimatestable} and from \eqref{H^k_bound}. Indeed, one has thanks to Lemma \ref{Lemmedensityreferencecontrol}

\begin{align*}
	&|\p \otimes \H^k(s,t,x,y)| \\&\leq C\int_{s}^t \int_{\R^d} \rho^0(r-s,z-x)
	 C^k (t-r)^{-\frac{1}{\alpha} + (k-1)\left( 1 - \frac{1}{\alpha} \right)} \prod_{j=1}^{k-1} \BB \left( j\left( 1 - \frac{1}{\alpha} \right),1 - \frac{1}{\alpha} \right) \rho^1(t-r,y-z)\, dz \, dr \\ &\leq C^{k+1} \frac{(t-s)^{k\left( 1 - \frac{1}{\alpha}\right)}}{k\left(1 - \frac{1}{\alpha} \right)} \prod_{j=1}^{k-1} \BB \left( j\left( 1 - \frac{1}{\alpha} \right),1 - \frac{1}{\alpha} \right) \rho^0(t-s,y-x).
	 \end{align*}
We conclude noting that 

\begin{align*}
	\frac{1}{k\left(1 - \frac{1}{\alpha} \right)} \prod_{j=1}^{k-1} \BB \left( j\left( 1 - \frac{1}{\alpha} \right),1 - \frac{1}{\alpha} \right) &= \frac{1}{k\left(1 - \frac{1}{\alpha} \right)}\prod_{j=1}^{k-1} \frac{\Gamma \left(j \left(1 - \frac{1}{\alpha}\right)\right) \Gamma \left( 1 - \frac{1}{\alpha}\right)}{\Gamma \left((j+1) \left(1 - \frac{1}{\alpha}\right)\right)} \\ &= \frac{1}{k\left(1 - \frac{1}{\alpha} \right)} \frac{ \left[\Gamma \left( 1 - \frac{1}{\alpha}\right)\right]^k}{\Gamma \left(k \left(1 - \frac{1}{\alpha}\right)\right)} \\ &= \frac{ \left[\Gamma \left( 1 - \frac{1}{\alpha}\right)\right]^k}{\Gamma \left(1+k\left(1 - \frac{1}{\alpha}\right)\right)} \\ &=\prod_{j=1}^{k} \BB \left( \frac{1}{\alpha} + j\left( 1 - \frac{1}{\alpha} \right),1 - \frac{1}{\alpha} \right) 
\end{align*}

Finally, the proof of \eqref{Phi_bound} and \eqref{Phi_Holder} follows directly from \eqref{H^k_bound} and \eqref{H^k_Holder} using the asymptotic expansion of the Beta function.

\end{proof}

\begin{proof}[Proof of Theorem \ref{Thmdensityparametrix}]
	The existence of the density and its representation \eqref{representationdensityparametrix} is a consequence of \eqref{gradientestimatestable}, \eqref{potimesH^k_bound} and \eqref{Phi_bound}. Indeed, using the asymptotic expansion of the Beta function, we obtain that the series \eqref{representationdensityparametrix} is absolutely convergent, locally uniformly with respect to $(s,x) \in [0,t) \times \R^d$. We refer to \cite{Kulik_parametrix_stable} for a detailed presentation of the parametrix method for SDEs driven by stable processes. The permutation of the series and the convolution in the representation formula  \eqref{representationdensityparametrix} is clearly justified by the dominated convergence theorem. The regularity of the density $p$ with respect to $x$ and the controls \eqref{density_bound} and \eqref{gradientdensity_holder} follow from Lemma \ref{Lemmegradientestimatestable} for the proxy, i.e.\ the first term of the parametrix series \eqref{representationdensityparametrix}. We now prove that they also hold for the other term $\p \otimes \Phi$ of \eqref{representationdensityparametrix}. \\

	\noindent\textbf{Proof of \eqref{density_bound} for $j=1$.} We use \eqref{Phi_bound} and Lemma \ref{Lemmegradientestimatestable} to deduce that 
	
	\begin{align*}
		|\pa_x\p \otimes \Phi(s,t,x,y)|  & \leq C\int_s^t \int_{\R^d} (r-s)^{-\frac{1}{\alpha}} \rho^1(r-s,z-x)  (t-r)^{-\frac{1}{\alpha} } \rho^1(t-r,y-z) \, dz \, dr.
	\end{align*}
	
	Then, Lemma \ref{Lemmedensityreferencecontrol} yields
	
	\begin{align*}
		&|\pa_x\p \otimes \Phi(s,t,x,y)| \\ & \leq C\int_s^t  (r-s)^{-\frac{1}{\alpha}}  (t-r)^{-\frac{1}{\alpha}} \, dr   \rho^{1}(t-s,y-x) \\ & \leq C (t-s)^{-\frac{1}{\alpha}  1 - \frac{1}{\alpha} } \rho^{1}(t-s,y-x).
	\end{align*}
This concludes the proof of \eqref{density_bound} for $j=1$.

	\noindent\textbf{Proof of \eqref{gradientdensity_holder} for $j=0$.} We use again the control of $\p$ given by Lemma \ref{Lemmegradientestimatestable} and  \eqref{Phi_bound}. Thanks to Lemma \ref{Lemmedensityreferencecontrol}, we obtain 

\begin{align*}
	&|\p \otimes \Phi(s,t,x_1,y) - \p \otimes \Phi(s,t,x_2,y)| \\&= \left\vert\int_s^t \int_{\R^d} (\p(s,r,x_1,z) - \p(s,r,x_2,z)) \Phi(r,t,z,y) \,  dz\, dr \right\vert \\&\leq C\int_{s}^t \int_{\R^d} (r-s)^{- \frac{\gamma}{\alpha}} |x_1 - x_2|^\gamma \left[\rho^0(r-s,z-x_1) +  \rho^0(r-s,z-x_2)  \right] (t-r)^{-\frac{1}{\alpha}}  \rho^1(t-r,y-z) \, dz\,dr\\ &\leq C (t-s)^{- \frac{\gamma}{\alpha}+ 1 - \frac{1}{\alpha}} |x_1 - x_2|^\gamma  \left[ \rho^0(t-s,y-x_1)  +  \rho^0(t-s,y-x_2) \right].
\end{align*}

	\noindent\textbf{Proof of \eqref{gradientdensity_holder} for $j=1$.} We use the following decomposition \begin{align*}
		\pa_x\p \otimes \Phi(s,t,x_1,y) - \pa_x\p \otimes \Phi(s,t,x_2,y)&= \int_s^t \int_{\R^d} (\pa_x\p(s,r,x_1,z) - \pa_x\p(s,r,x_2,z)) \Phi(r,t,z,y) \,  dz\, dr   \\ &=  \int_{D_1} \int_{\R^d} (\pa_x\p(s,r,x_1,z) - \pa_x\p(s,r,x_2,z)) \Phi(r,t,z,y) \,  dz\, dr   \\ &\quad +  \int_{D_2} \int_{\R^d} (\pa_x\p(s,r,x_1,z) - \pa_x\p(s,r,x_2,z)) \Phi(r,t,z,y) \,  dz\, dr \\ &=: I_1 + I_2,
	\end{align*}
	where $D_1:= \{ r \in (s,t),\, |x_1 - x_2|> (r-s)^{\frac{1}{\alpha}}\}$ and $D_2:= \{ r \in (s,t),\, |x_1 - x_2|\leq (r-s)^{\frac{1}{\alpha}}\}$. For $I_1$, one can write \begin{align*}
		I_1 &=  \int_{D_1} \int_{\R^d} \pa_x\p(s,r,x_1,z) (\Phi(r,t,z,y) - \Phi(r,t,x_1,y)) \,  dz\, dr \\ &\quad- \int_{D_1} \int_{\R^d} \pa_x\p(s,r,x_2,z) (\Phi(r,t,z,y) - \Phi(r,t,x_2,y)) \,  dz\, dr,\\ &=: I_{1,1} + I_{1,2},
	\end{align*}
	since $\int_{\R^d} \pa_x\p(s,r,x,z) \, dz = 0$ for all $x \in \R^d$. Since $\alpha \in (1,2)$ and $\gamma < ( 2 \alpha - 2) \wedge (\eta + \alpha -1),$ we can pick $\delta \in (0,(\alpha-1)\wedge \eta)$ such that $\gamma < \delta + \alpha -1$. Then, we can use Lemma \ref{Lemmegradientestimatestable} and \eqref{Phi_Holder} with $\delta$, which yields
	
	\begin{multline*}
		|I_{1,1}|  \leq C \int_{D_1}\int_{\R^d} (r-s)^{-\frac{1}{\alpha}} \rho^1(r-s,z-x_1) (t-r)^{-\frac{\delta +1}{\alpha} } |z-x_1|^\delta \left[\rho^{1 }(t-r,y-z) + \rho^{1 }(t-r,y-x_1) \right] \, dz\,dr.
	\end{multline*}
	
	Using the space-time inequality \eqref{scalingdensityref} in Lemma \ref{Lemmedensityreferencecontrol}, we deduce that 
	
	\begin{multline*}
		|I_{1,1}|  \leq C \int_{D_1}\int_{\R^d} (r-s)^{\frac{ \delta-1}{\alpha}} \rho^{1-\delta}(r-s,z-x_1) (t-r)^{-\frac{\delta +1}{\alpha}}  \left[\rho^{1 }(t-r,y-z) + \rho^{1 }(t-r,y-x_1) \right] \, dz\,dr.
	\end{multline*}
	
	Let us note that $\int_{\R^d} \rho^{1-\delta}(r-s,z-x_1) \, dz$ is a constant independent of $r,s,x_1$, that if $r \in (s,t)$ \[\rho^{1 }(t-r,y-x_1) \leq (t-r)^{-\frac{d}{\alpha}} ( 1 + (t-s)^{-\frac{1}{\alpha}}|y-x_1|)^{-d-\alpha -1}\] and that if $r \in D_1$, $(r-s)^{\frac{\gamma}{\alpha}}<|x_1-x_2|^\gamma$. Applying Lemma \ref{Lemmedensityreferencecontrol}, we thus get that 
	
	\begin{align*}
		|I_{1,1}|  &\leq C\left[  \int_s^t (r-s)^{\frac{ \delta-1-\gamma}{\alpha}} (t-r)^{-\frac{\delta +1}{\alpha} } \,dr |x_1 -x_2|^\gamma \rho^{1-\delta }(t-s,y-x_1) \right. \\ &\left. \hspace{3cm}\int_s^t (r-s)^{\frac{ \delta-1-\gamma}{\alpha}} (t-r)^{-\frac{\delta +1 +d}{\alpha} } \,dr |x_1-x_2|^\gamma  ( 1 + (t-s)^{-\frac{1}{\alpha}}|y-x_1|)^{-d-\alpha -1} \right]\\ &\leq  C(t-s)^{-\frac{ 1+\gamma}{\alpha} + 1 - \frac{1}{\alpha}} |x_1-x_2|^\gamma \rho^{1-\delta }(t-s,y-x_1) . 
	\end{align*}

	Similarly, we obtain  \begin{equation*}
		|I_{1,2}| \leq C  (t-s)^{-\frac{ 1+\gamma}{\alpha} + 1 - \frac{1}{\alpha} } |x_1-x_2|^\gamma  \rho^{1- \delta}(t-s,y-x_2).
	\end{equation*}
	
	This proves \eqref{gradientdensity_holder} for $I_1$. We now focus on $I_2.$ One can write 
	
	\begin{align*}
		I_2 &=\int_{D_2} \int_{\R^d} (\pa_x\p(s,r,x_1,z) - \pa_x\p(s,r,x_2,z)) \Phi(r,t,z,y) \,  dz\, dr \\ &= \int_{D_2} \int_{\R^d} (\pa_x\p(s,r,x_1,z) - \pa_x\p(s,r,x_2,z)) (\Phi(r,t,z,y)- \Phi(r,t,x_2,y))\,  dz\, dr.
	\end{align*}
	Note that if $r \in D_2$, we have \begin{align*}|\pa_x\p(s,r,x_1,z) - \pa_x\p(s,r,x_2,z)| &\leq C(r-s)^{-\frac{2}{\alpha}} |x_1-x_2| \left[\rho^2(r-s,z-x_1) + \rho^2(r-s,z-x_1) \right] \\ &\leq C(r-s)^{-\frac{\gamma +1}{\alpha}} |x_1-x_2|^\gamma \rho^2(r-s,z-x_2),
	\end{align*}
	since $|x_1-x_2| \leq (r-s)^{\frac{1}{\alpha}}$ and by  \eqref{controldensityref}.
	Using \eqref{Phi_Holder} with $\delta \in (0,(\alpha -1) \wedge \eta)$ such that $\gamma < \delta + \alpha -1$ and the space-time inequality \eqref{scalingdensityref}, we get 
	
	\begin{align*}
		|I_2|&\leq C\int_{D_2}\int_{\R^d} (r-s)^{-\frac{ \gamma +1}{\alpha}}|x_1-x_2|^\gamma \rho^2(r-s,z-x_2) (t-r)^{-\frac{\delta +1}{\alpha}} |z-x_2|^\delta \\ &\hspace{7cm}\left[\rho^{1 }(t-r,y-z) + \rho^{1 }(t-r,y-x_2) \right] \, dz\,dr \\ &\leq C \int_{D_2}\int_{\R^d} (r-s)^{-\frac{ \gamma +1}{\alpha} + \frac{\delta}{\alpha}}|x_1-x_2|^\gamma \rho^{2- \delta}(r-s,z-x_2)  (t-r)^{-\frac{\delta +1}{\alpha}}  \\ &\hspace{7cm} \left[\rho^{1 }(t-r,y-z) + \rho^{1 }(t-r,y-x_2) \right] \, dz\,dr.
	\end{align*}
	
	As done previously to deal with $I_{1,1}$, Lemma \ref{Lemmedensityreferencecontrol} yields \begin{equation*}
		|I_2| \leq C (t-s)^{-\frac{ \gamma +1}{\alpha} +  1 - \frac{1}{\alpha}} \left[\rho^{1 }(t-r,y-x_1) + \rho^{1 }(t-r,y-x_2) \right].
	\end{equation*}
	
	It concludes the proof of \eqref{gradientdensity_holder} for $j=1$.\\

	\noindent\textbf{Proof of \eqref{EDPKolmogorovparametrix}.} Let us now prove that $p(\cdot,t,\cdot,y)$ is a fundamental solution to \eqref{EDPKolmogorovparametrix}. We fix $0 \leq s < t\leq T$ and $x,y \in \R^d$. From the Markov property satisfied by the SDE \eqref{SDEparametrix}, stemming from the well-posedness of the related martingale problem, one has for all $h >0$ such that $s-h \geq 0$ 
	
	\[ p(s-h,t,x,y) = \E (p(s,t,X^{s-h,x}_s,y)).\]
	
	Applying Itô's formula to the function $p(s,t,\cdot,y)$ which belongs to $\CC^{1+\gamma}_b(\R^d;\R)$ for $\gamma > \alpha -1$, we obtain that \begin{align*}
		p(s,t,X_s^{s-h,x},y) &= p(s,t,x,y) + \int_{s-h}^s b(r,X^{s-h,x}_r) \cdot \pa_x p(s,t,X^{s-h,x}_r,y) \, dr \\ &\quad + \int_{s-h}^s \int_{\R^d} p(s,t,X^{s-h,x}_{r^-} +z,y) - p(s,t,X^{s-h,x}_{r^-},y) \, \NNN(dr,dz) \\ &\quad +\int_{s-h}^s \int_{\R^d} p(s,t,X^{s-h,x}_{r^-} +z,y) - p(s,t,X^{s-h,x}_{r^-},y) - \pa_xp(s,t,X^{s-h,x}_{r^-},y) \cdot z  \, d\nu(z) \, dr
	\end{align*}
\end{proof}
We can take the expectation in the preceding formula using \eqref{density_bound}, \eqref{gradientdensity_holder}. It yields 

\begin{align*}
p(s-h,t,x,y) &= p(s,t,x,y) + \int_{s-h}^s \E (b(r,X^{s-h,x}_r) \cdot \pa_x p(s,t,X^{s-h,x}_r,y)) \, dr  \\ &\quad +\int_{s-h}^s \int_{\R^d} \E (p(s,t,X^{s-h,x}_{r} +z,y) - p(s,t,X^{s-h,x}_{r},y) - \pa_xp(s,t,X^{s-h,x}_{r},y) \cdot z ) \, d\nu(z) \, dr
\end{align*} 

Let us prove that \begin{equation}\label{eqproofparametrix1}
	 \frac1h  \int_{s-h}^s \E (b(r,X^{s-h,x}_r) \cdot \pa_x p(s,t,X^{s-h,x}_r,y)) \, dr  \underset{h \rightarrow 0}{\longrightarrow} b(s,x) \cdot \pa_x p(s,t,x,y).\end{equation}We can write \begin{align*}
&\left\vert\frac1h  \int_{s-h}^s \E (b(r,X^{s-h,x}_r) \cdot \pa_x p(s,t,X^{s-h,x}_r,y)) \, dr  - b(s,x) \cdot \pa_x p(s,t,x,y)\right\vert \\ & \leq  \left\vert\frac1h  \int_{s-h}^s \E (b(r,X^{s-h,x}_r) \cdot \pa_x p(s,t,X^{s-h,x}_r,y)) -  \E (b(r,X^{s-h,x}_s) \cdot \pa_x p(s,t,X^{s-h,x}_s,y)) \, dr \right\vert \\ &\quad +   \left\vert\frac1h  \int_{s-h}^s \E (b(r,X^{s-h,x}_s) \cdot \pa_x p(s,t,X^{s-h,x}_s,y)) -  (b(r,x) \cdot \pa_x p(s,t,x,y)) \, dr \right\vert \\ &\quad +  \left\vert\frac1h  \int_{s-h}^s   (b(r,x) \cdot \pa_x p(s,t,x,y)) - b(s,x)  \cdot \pa_x p(s,t,x,y)\, dr \right\vert \\ &=: I_1 + I_2 + I_3.
\end{align*}

It is clear that $I_3$ converges to $0$ as $h$ tends to $0$ since $b$ is continuous. Concerning $I_2$, one has \begin{align*}
	I_2 &\leq  \frac1h  \int_{s-h}^s \E |b(r,X^{s-h,x}_s) \cdot \pa_x p(s,t,X^{s-h,x}_s,y)) -  (b(r,x) \cdot \pa_x p(s,t,x,y)| \, dr \\ &\leq  \frac1h  \int_{s-h}^s \E| b(r,X^{s-h,x}_s)| \,| \pa_x p(s,t,X^{s-h,x}_s,y) - \pa_x p(s,t,x,y)| \, dr  \\ &\quad + \frac1h  \int_{s-h}^s \E  |b(r,X^{s-h,x}_s) - b(r,x)| \, |\pa_x p(s,t,x,y)| \, dr \\ &=: I_{2,1} + I_{2,2}.
\end{align*}

Since $b$ is bounded and by the Hölder control \eqref{gradientdensity_holder}, we obtain that for some constant $\gamma \in (0,1)$ \begin{align*}
	I_{2,1} &\leq C_{s,t,\gamma} \E |X^{s-h,x}_s - x|^\gamma \\&\leq C_{s,t,\gamma} h^{\frac{\gamma}{\alpha}}.
\end{align*}
The same reasoning based on the Hölder continuity of $b$ proves that $I_2 \underset{h \rightarrow 0}{\longrightarrow} 0.$

Finally, we decompose $I_1$ in the following way 

\begin{align*}
	I_1 &\leq \frac1h  \int_{s-h}^s \E |b(r,X^{s-h,x}_r)| \, |\pa_x p(s,t,X^{s-h,x}_r,y) - \pa_x p(s,t,X^{s-h,x}_s,y)| \, dr \\ &\quad+ \frac1h  \int_{s-h}^s \E |b(r,X^{s-h,x}_r) -b(r,X^{s-h,x}_s) | \, | \pa_x p(s,t,X^{s-h,x}_s,y)| \, dr \\ &=: I_{1,1} + I_{1,2}.
\end{align*} 
 Note that for all $\gamma \in (0,1]$, there exists $C_{s,t,\gamma} >0$ such that for all $r\in [s-h,s]$ \begin{equation*}
	\E  |X^{s-h,x}_r - X^{s-h,x}_s|^\gamma \leq C_{s,t,\gamma} (s-r)^{\frac{\gamma}{\alpha}},
\end{equation*}
The same reasoning as done for $I_2$ can be applied since $b$ and $\pa_x p(s,t,\cdot,y)$ are globally bounded and $b$ is uniformly Hölder continuous. It yields $I_1 \underset{h \rightarrow 0}{\longrightarrow} 0.$

Let us now show that \begin{multline}\label{eqproofparametrix2}
	\frac1h \int_{s-h}^s \int_{\R^d} \E(p(s,t,X^{s-h,x}_{r} +z,y) - p(s,t,X^{s-h,x}_{r},y) - \pa_xp(s,t,X^{s-h,x}_{r},y) \cdot z)  \, d\nu(z) \, dr\\ \underset{h \rightarrow 0}{\longrightarrow} \Delta^{\frac{\alpha}{2}} p(s,t,\cdot,y)(x).
\end{multline}

One can write \begin{align*}
	&\left\vert  \frac1h \int_{s-h}^s \int_{\R^d} \E(p(s,t,X^{s-h,x}_{r} +z,y) - p(s,t,X^{s-h,x}_{r},y) - \pa_xp(s,t,X^{s-h,x}_{r},y) \cdot z ) \, d\nu(z) \, dr - \Delta^{\frac{\alpha}{2}} p(s,t,\cdot,y)(x) \right\vert\\ &\leq  \left\vert   \int_{\R^d} \frac1h \int_{s-h}^s \E(p(s,t,X^{s-h,x}_{r} +z,y) - p(s,t,X^{s-h,x}_{r},y) - \pa_xp(s,t,X^{s-h,x}_{r},y) \cdot z)   \, dr\, d\nu(z)\right. \\ &\left.\hspace{4cm}-   \int_{\R^d} \E( p(s,t,X^{s-h,x}_{s} +z,y) - p(s,t,X^{s-h,x}_{s},y) - \pa_xp(s,t,X^{s-h,x}_{s},y) \cdot z)  \, d\nu(z) \right\vert \\ &\quad + \left\vert \int_{\R^d} \E( p(s,t,X^{s-h,x}_{s} +z,y) - p(s,t,X^{s-h,x}_{s},y) - \pa_xp(s,t,X^{s-h,x}_{s},y) \cdot z)  \, d\nu(z) - \Delta^{\frac{\alpha}{2}} p(s,t,\cdot,y)(x) \right\vert \\ &=: J_1 + J_2.
\end{align*}
For $J_1$, we obtain that \begin{align*}
	J_1 &\leq  \int_{\R^d} \frac1h \int_{s-h}^s\int_0^1 \E\left\vert(\pa_x p(s,t,X^{s-h,x}_{r} +\lambda z,y) -  \pa_xp(s,t,X^{s-h,x}_{r},y)) \right.   \\ &\left.\hspace{4cm}-   ( \pa_x p(s,t,X^{s-h,x}_{s} +\lambda z,y) - \pa_x p(s,t,X^{s-h,x}_{s},y) ) \right\vert \, |z| \, d\lambda \, dr \, d\nu(z) 
	\end{align*}
We are going to use the dominated convergence theorem in the integral with respect to $\nu$. By the Hölder control \eqref{gradientdensity_holder} on $\pa_x p$, we deduce that for some $\gamma > \alpha -1$, there exists a constant $C_{s,t,\gamma} >0$ such that for all $r \in [s-h,s]$, $ z \in \R^d$ \begin{align*}
&\frac1h \int_{s-h}^s\int_0^1 \E\left\vert(\pa_x p(s,t,X^{s-h,x}_{r} +\lambda z,y) -  \pa_xp(s,t,X^{s-h,x}_{r},y)) \right. \\ &\hspace{4cm} \left. -( \pa_x p(s,t,X^{s-h,x}_{s} +\lambda z,y) - \pa_x p(s,t,X^{s-h,x}_{s},y) ) \right\vert \, |z| \, d\lambda \, dr\\ &\leq C_{s,t,\gamma} \frac1h \int_{s-h}^s \E|X^{s-h,x}_r - X^{s-h,x}_s|^\gamma\, dr \, |z| \\ &\leq C_{s,t,\gamma} \frac1h \int_{s-h}^s (r-s)^{\frac{\gamma}{\alpha}} \, dr \, |z| \\ &\leq C_{s,t,\gamma}h^{\frac{\gamma}{\alpha}}|z|.
\end{align*}
The right-hand side term tends to $0$ when $h \to 0$. We start by justifying the domination on the ball $B_1$. In this case, we use again the Hölder continuity of $\pa_x p$ with respect to $x$, which yields

\begin{align*}
	&\frac1h \int_{s-h}^s\int_0^1 \E\left\vert(\pa_x p(s,t,X^{s-h,x}_{r} +\lambda z,y) -  \pa_xp(s,t,X^{s-h,x}_{r},y)) \right. \\ &\hspace{4cm} \left. -( \pa_x p(s,t,X^{s-h,x}_{s} +\lambda z,y) - \pa_x p(s,t,X^{s-h,x}_{s},y) ) \right\vert \, |z| \, d\lambda \, dr\\ &\leq  C_{s,t,\gamma} |z|^{1 + \gamma}.
\end{align*}

Since $\gamma > \alpha -1$, the map $z \in B_1 \mapsto |z|^{1 + \gamma}$ belongs to $L^1(B_1,\nu).$ The domination on $B_1^c$ is clear since $ \pa_x p(s,t,\cdot,y)$ is globally bounded by \eqref{density_bound} and $z \in B_1^c \mapsto |z| $ belongs to $L^1(B_1^c,\nu)$.\\

We now deal with $J_2.$ Thanks to the dominated convergence theorem, $ J_2$ converges to $0$ as $h$ tends to $0$. The Hölder control \eqref{gradientdensity_holder} of $\pa_x p$ ensures that for some constant $\gamma >\alpha -1$, there exists $C_{s,t,\gamma}>0$ such that \[ \E| p(s,t,X^{s-h,x}_{s} +z,y) - p(s,t,X^{s-h,x}_{s},y) - \pa_xp(s,t,X^{s-h,x}_{s},y) \cdot z| \leq C_{s,t,\gamma}|z|^{1 + \gamma}.\] This proves the domination on the ball $B_1$. The domination on $B_1^c$ is a consequence of \eqref{density_bound}, and the fact that $\alpha \in (1,2)$.\\

We have thus proved that the map $s \in [0,t) \mapsto p(s,t,x,y) $ is left-differentiable by \eqref{eqproofparametrix1} and \eqref{eqproofparametrix2}. Moreover, since the map $(s,x) \in [0,t) \times \R^d \mapsto b(s,x)\cdot \pa_x p(s,t,x,y) + \Delta^{\frac{\alpha}{2}} p(s,t,\cdot,y)(x)$ is continuous, we deduce that the map  $s \in [0,t) \mapsto p(s,t,x,y) $ is of class $\CC^1$ and that it solves \[\pa_s p(s,t,x,y) + b(s,x)\cdot \pa_x p(s,t,x,y) + \Delta^{\frac{\alpha}{2}} p(s,t,\cdot,y)(x) =0, \quad \forall (s,x) \in [0,t) \times \R^d.\]

Let us now fix $ f : \R^d \rightarrow \R$ a bounded and uniformly continuous function. We fix $\eps >0$. There exists $\delta >0$ such that for all $x,y \in \R^d$ with $|x-y|\leq \delta$, we have $ |f(x) - f(y)| \leq \eps$. Using \eqref{density_bound}, we obtain that \begin{align*}
	\sup_{x \in \R^d} \left\vert \int_{\R^d} f(y) p(s,t,x,y) \, dy - f(x) \right\vert& = 	\sup_{x \in \R^d} \left\vert \int_{\R^d} (f(y) - f(x)) p(s,t,x,y) \, dy \right\vert \\ &\leq  \eps + C\|f\|_{\infty}\int_{|y| >\delta} (t-s)^{-\frac{d}{\alpha}}( 1 + (t-s)^{-\frac{1}{\alpha}}|y|)^{-d-\alpha} \, dy \\ &= \eps + C\|f\|_{\infty}\int_{|z| >(t-s)^{-\frac{1}{\alpha}}\delta} ( 1 + |z|)^{-d-\alpha} \, dz.
\end{align*}
We conclude, taking the $\limsup$ when $s \rightarrow t$ in the preceding inequality, that $p(s,t,x,\cdot) \underset{s\rightarrow t^-}{\longrightarrow} \delta_x.$\\

\noindent\textbf{Proof of \eqref{timederivativedensity_bound}.}
We now prove that \[ |\Delta^{\frac{\alpha}{2}} p(s,t,\cdot,y)(x) | \leq C (t-s)^{-1} \rho^{0}(t-s,y-x),\] using the PDE \eqref{EDPKolmogorovparametrix} and the fact that $b$ is bounded and \eqref{density_bound}. Note that by symmetry, we have \[ \Delta^{\frac{\alpha}{2}} p(s,t,\cdot,y)(x) = \frac12 \int_{\R^d} p(s,t,x+z,y) + p(s,t,x-z,y) - 2 p(s,t,x,y) \, \frac{dz}{|z|^{d + \alpha}}. \] We decompose it in the following way \begin{align*}
	 |\Delta^{\frac{\alpha}{2}} p(s,t,\cdot,y)(x) |  &\leq \int_{|z| \leq (t-s)^{\frac{1}{\alpha}}} |p(s,t,x+z,y) + p(s,t,x-z,y) - 2 p(s,t,x,y) | \, \frac{dz}{|z|^{d + \alpha}} \\ &\quad + \int_{|z| > (t-s)^{\frac{1}{\alpha}}} | p(s,t,x+z,y) + p(s,t,x-z,y) - 2 p(s,t,x,y)| \, \frac{dz}{|z|^{d + \alpha}} \\ &=  \int_{|z| \leq (t-s)^{\frac{1}{\alpha}}} \left\vert \int_0^1 (\pa_x p(s,t,x +\lambda z,y) - \pa_xp(s,t,x-\lambda z,y) ) \cdot z \, d\lambda \right\vert \, \frac{dz}{|z|^{d + \alpha}} \\ &\quad + \int_{|z| > (t-s)^{\frac{1}{\alpha}}} \left\vert p(s,t,x+z,y) + p(s,t,x-z,y) - 2 p(s,t,x,y) \right\vert \, \frac{dz}{|z|^{d + \alpha}}
	 \\ &=: I_1 + I_2.
\end{align*}
We start with $I_1$. Using the Hölder control \eqref{gradientdensity_holder}, we obtain that for some $\gamma \in (0,(2\alpha -2)\wedge(\eta + \alpha -1))$ with $\gamma> \alpha -1$, there exists a constant $C>0$ such that \begin{align*}
	I_1 \leq C (t-s)^{-\frac{1 + \gamma}{\alpha}} \int_{|z| \leq (t-s)^{\frac{1}{\alpha}}} \int_0^1 (\rho^1 (t-s, y-x-\lambda z) + \rho^1(t-s,y-x + \lambda z)) \, d\lambda |z|^{1 + \gamma}  \, \frac{dz}{|z|^{d + \alpha}}.
\end{align*}
Since $|z|  \leq (t-s)^{\frac{1}{\alpha}}$, the space-time inequality \eqref{scalingdensityref} ensures that  \begin{align}\label{eqproofparametrix3}
	\notag I_1 &\leq C (t-s)^{-\frac{1 + \gamma}{\alpha}} \int_{|z| \leq (t-s)^{\frac{1}{\alpha}}} \rho^1(t-s,y-x) |z|^{1 + \gamma}  \, \frac{dz}{|z|^{d + \alpha}} \\ &\leq C (t-s)^{-1} \rho^1(t-s,y-x).
\end{align}
For $I_2$, we have with \eqref{density_bound}

\begin{align*}
	I_2 & \leq C\int_{ |z|> (t-s)^{\frac{1}{\alpha}}}  \left[\rho^0(t-s,y-x-z) + \rho^0(t-s,y-x+z)+ \rho^0(t-s,y-x) \right]  \, \frac{dz}{|z|^{d + \alpha}}.
\end{align*}

We distinguish two cases. Firstly, assume that $|y-x| \leq 2 (t-s)^{\frac{1}{\alpha}}$. Then, using that for all $s < t$ and $ y\in \R^d$ $\rho^0(t-s,x) \leq C(t-s)^{-\frac{d}{\alpha}},$ we get that 

\begin{align}\label{eqproofparametrix4}
\notag	I_2 & \leq C(t-s)^{-\frac{d}{\alpha}}\int_{ |z|> (t-s)^{\frac{1}{\alpha}}}  \frac{dz}{|z|^{d + \alpha}} \\ &\leq C (t-s)^{-1 - \frac{d}{\alpha}} \\ \notag&\leq C (t-s)^{-1} \rho^1(t-s,y-x),
\end{align} since $|y-x| \leq 2 (t-s)^{\frac{1}{\alpha}}.$ Finally, if $|y-x|> 2 (t-s)^{\frac{1}{\alpha}}$, we write 

\begin{align*}
	I_2 & \leq C\int_{  (t-s)^{\frac{1}{\alpha}}<|z| \leq \frac{|x-y|}{2}} \left[\rho^0(t-s,y-x-z)+ \rho^0(t-s,y-x+z) + \rho^0(t-s,y-x) \right]  \, \frac{dz}{|z|^{d + \alpha}} \\ &\quad +   C\int_{|z| > \frac{|x-y|}{2}}  \left[\rho^0(t-s,y-x-z) + \rho^0(t-s,y-x+z) + \rho^0(t-s,y-x) \right]  \, \frac{dz}{|z|^{d + \alpha}}\\ &=: I_{2,1} + I_{2,2}.
\end{align*}

Then, the reverse triangle inequality yields \begin{align}\label{eqproofparametrix5}
	\notag I_{2,1} &\leq C \int_{  (t-s)^{\frac{1}{\alpha}}<|z| \leq \frac{|x-y|}{2}} \left[ \rho^0(t-s,y-x)+  \rho^0(t-s,\frac{y-x}{2}) \right]  \, \frac{dz}{|z|^{d + \alpha}} \\ \notag&\leq  C  \rho^0(t-s,y-x)  \int_{  (t-s)^{\frac{1}{\alpha}}<|z|} \frac{dz}{|z|^{d + \alpha}} \\ &\leq C (t-s)^{-1}\rho^0(t-s,y-x).
\end{align}
For $I_{2,2}$, we have \begin{align*}
	 I_{2,2}  &\leq C|y-x|^{-d-\alpha} \int_{\R^d} \rho^0(t-s,y-x-z) + \rho^0(t-s,y-x+z)  \, dz \\& \quad + C \rho^0(t-s,y-x)\int_{|z| > \frac{|x-y|}{2}} \frac{dz}{|z|^{d + \alpha}}
	 \\ &\leq  C (t-s)^{\frac{-d-\alpha}{\alpha}} ((t-s)^{-\frac{1}{\alpha}}|y-x|)^{-d-\alpha }  + C \rho^0(t-s,y-x)|x-y|^{ - \alpha} \\ &\leq  C (t-s)^{\frac{-d-\alpha }{\alpha}} ((t-s)^{-\frac{1}{\alpha}}|y-x|)^{-d-\alpha }  + C(t-s)^{-1}  \rho^0(t-s,y-x).
\end{align*} 
Since $1 < \frac{1}{2}(t-s)^{-\frac{1}{\alpha}}|y-x|$, we deduce that for some constant $C>0$, one has \[((t-s)^{-\frac{1}{\alpha}}|y-x|)^{-d-\alpha } \leq C(1 + (t-s)^{-\frac{1}{\alpha}}|y-x|)^{-d-\alpha }. \]
It follows that 
\begin{align}\label{eqproofparametrix6}
	I_{2,2} &\leq C(t-s)^{-1}\rho^{0}(t-s,y-x).
\end{align} 

Combining \eqref{eqproofparametrix3}, \eqref{eqproofparametrix4}, \eqref{eqproofparametrix5} and \eqref{eqproofparametrix6}, we have proved that \eqref{timederivativedensity_bound} holds true.  

\subsection*{Acknowledgements}

I would like to thank Paul-Eric Chaudru de Raynal for his supervision, advices and for his careful reading of the paper. I am also grateful to Benjamin Jourdain for his careful reading of the paper and his relevant remarks.

\bibliographystyle{plain}
\bibliography{bibliography}

\end{document}